\documentclass[11pt]{article}

\usepackage{fancyhdr}
\usepackage[width=175mm,top=20mm,bottom=25mm,bindingoffset=0mm]{geometry}
\usepackage{graphicx}      
\usepackage{float}         
\usepackage{libertine}
\usepackage{amsmath}
\usepackage{mathrsfs}      
\usepackage{amsfonts}      
\usepackage{mathtools}     
\usepackage{amssymb}       
\usepackage{nccmath}       
\usepackage{bm}            
\usepackage{calligra}
\usepackage{longtable}
\usepackage{xcolor}        
\usepackage{indentfirst}   

\usepackage{hyperref}
\hypersetup{
  colorlinks = true,      
  urlcolor   = gray,      
  linkcolor  = blue,      
  citecolor  = red        
}

\usepackage[noprefix]{nomencl}

\renewcommand{\nomgroup}[1]{%
  \ifthenelse{\equal{#1}{A}}{
\item[\textbf{Symbols}]}{
    \ifthenelse{\equal{#1}{B}}{
\item[\textbf{Acronyms}]}{}}
}
\newcommand{\nomacro}[1][]{\nomenclature[B#1]}
\newcommand{\nomvar}[1][]{\nomenclature[A#1]}

\usepackage[numbers,square]{natbib}
\usepackage{filecontents}

\usepackage{mdframed}
\newmdenv[leftline=false,rightline=false]{topbot}
\newmdenv[topline=false,rightline=false]{leftbot}
\newmdenv[topline=false,bottomline=false,rightline=false]{leftonly}

\usepackage[skins]{tcolorbox}
\newtcolorbox{myframe}[2][]{%
  enhanced,colback=white,colframe=black,coltitle=black,
  sharp corners,boxrule=0.4pt,
  fonttitle=\itshape,
  attach boxed title to top left={yshift=-0.3\baselineskip-0.4pt,xshift=2mm},
  boxed title style={tile,size=minimal,left=0.5mm,right=0.5mm,
  colback=white,before upper=\strut},
  title=#2,#1
}

\usepackage{tabularray}    
\usepackage{amssymb}       
\usepackage{pifont}        
\newcommand{\cmark}{\ding{51}}%
\newcommand{\xmark}{\ding{55}}%
\usepackage{authblk}

\newcommand\numberthis{\addtocounter{equation}{1}\tag{\theequation}}
\newcommand{\vect}[1]{\bm{#1}}

\newcommand{\argmin}{\mathop{\mathrm{argmin}}\limits}
\newcommand*{\vertbar}{\rule[-1ex]{0.5pt}{2.5ex}}
\newcommand*{\horzbar}{\rule[.5ex]{2.5ex}{0.5pt}}

\renewcommand{\nompreamble}{We use bold symbols for vectors and tensors. Vectors are in bold mathtype, whereas matrices and tensors are in bold upright roman. We use uppercase Latin letters for vectors of dimension $p$ (the number of coefficients/weights in the discrete operator); Vectors of other dimensions are in lowercase. We use lowercase Latin letters for vector or tensor element indices (in row-major ordering), and lowercase Greek letters for coordinate directions. The most important symbols are defined as follows:\vspace{0.75em}}
\makenomenclature

\begin{document}

\fancyhead[L]{An Overview of Meshfree Collocation Methods}
\fancyhead[R]{ }
\pagestyle{fancy}

\begin{center}
  \vspace*{0.5cm}
  {\LARGE\bfseries An Overview of Meshfree Collocation Methods}\\[0.5cm]
  {\normalsize
    \begin{tabular}{c}
      Tomas Halada\textsuperscript{1,4}, Serhii Yaskovets\textsuperscript{1,2,3}, Abhinav Singh\textsuperscript{1,2,3}, Ludek Benes\textsuperscript{4}, \\
      Pratik Suchde\textsuperscript{5,6,\dag}, and Ivo~F.~Sbalzarini\textsuperscript{1,2,3,\dag,*}
    \end{tabular}
  }\\[0.3cm]
  {\small
    \begin{tabular}{c}
      \textsuperscript{1}Dresden University of Technology, Faculty of Computer Science, 01187 Dresden, Germany \\
      \textsuperscript{2}Max Planck Institute of Molecular Cell Biology and Genetics, 01307 Dresden, Germany \\
      \textsuperscript{3}Center for Systems Biology Dresden, 01307 Dresden, Germany \\
      \textsuperscript{4}Czech Technical University, Faculty of Mechanical Engineering, Department of Technical Mathematics, \\
      160\,00 Prague, Czech Republic \\
      \textsuperscript{5}University of Luxembourg, Faculty of Science, Technology and Medicine, 4365 Esch-sur-Alzette, Luxembourg \\
      \textsuperscript{6}Fraunhofer Institute for Industrial Mathematics, 67663 Kaiserslautern, Germany \\
      \\
      \textsuperscript{\dag}Joint senior authors \\
      \textsuperscript{*}Corresponding author {\footnotesize \url{sbalzarini@mpi-cbg.de}}
    \end{tabular}
  }
\end{center}

\begin{abstract}
  We provide a comprehensive overview of meshfree collocation methods for numerically approximating differential operators on continuously labeled unstructured point clouds. Meshfree collocation methods do not require a computational grid or mesh. Instead, they approximate smooth functions and their derivatives at potentially irregularly distributed collocation points, often called {\it particles}, to a desired order of consistency. We review several meshfree collocation methods from the literature, trace the historical development of key concepts, and propose a classification of methods according to their principle of derivation. Although some of the methods reviewed are similar or identical, there are subtle yet important differences between many, which we highlight and discuss. We present a unifying formulation of meshfree collocation methods that renders these differences apparent and show how each method can be derived from this formulation. Finally, we propose a generalized derivation for meshfree collocation methods going forward.
\end{abstract}

\tableofcontents

\begin{thenomenclature}
  \nomgroup{A}
\item [{ $\alpha, \beta, \gamma, \dots$ }]\begingroup indices used to index individual spatial components\nomeqref {0}\nompageref{10}
\item [{ $\bm{\alpha}$ }]\begingroup multi-index\nomeqref {2}\nompageref{11}
\item [{ $\bm{\bm{U}}_i$ }]\begingroup vector containing values of function $ u$ for all neighbors of particle $i$\nomeqref {19}\nompageref{13}
\item [{ $\bm{\lambda}$ }]\begingroup Lagrange multipliers [-]\nomeqref {144}\nompageref{42}
\item [{ $\bm{\lambda}$ }]\begingroup Lagrange multipliers [-]\nomeqref {97}\nompageref{28}
\item [{ $\bm{\Psi}_{i}^{\bm{\alpha}}$ }]\begingroup vector of weight coefficients used to construct weight functions\nomeqref {14}\nompageref{12}
\item [{ $\bm{B}_{i}^{\bm{\alpha}}$ }]\begingroup vector of moments of weight functions $w_{ji}^{d}$\nomeqref {11}\nompageref{11}
\item [{ $\bm{D}u$ }]\begingroup vector of all multi-index derivatives of variable $\psi$ up order of $m$\nomeqref {2}\nompageref{11}
\item [{ $\bm{E}_i$ }]\begingroup vector of reminders in Lagrange form for all neighbors of particle $i$\nomeqref {19}\nompageref{13}
\item [{ $\bm{H}_i$ }]\begingroup vector of inverse stencil size scales following Taylor monomial powers up to order of $m$\nomeqref {28}\nompageref{14}
\item [{ $\bm{P}_{ji}$ }]\begingroup vector of basis functions evaluated between particles $i$ and $j$\nomeqref {14}\nompageref{12}
\item [{ $\bm{W}$, $\bm{W}_{ji}$, $ W_{ji}^{q} $ }]\begingroup vector of Anisotropic Basis Functions (ABFs), single ABF\nomeqref {14}\nompageref{12}
\item [{ $\bm{w}_{\bullet i}^{\bm{\alpha}}$ }]\begingroup vector of weight functions for all neighbors of particle $i$\nomeqref {97}\nompageref{28}
\item [{ $\bm{X}$, $\bm{X}_{ji}$ }]\begingroup vector of Taylor monomials assumed up to order of $m$\nomeqref {6}\nompageref{11}
\item [{ $\delta_{\alpha\beta}$ }]\begingroup Kronecker delta\nomeqref {6}\nompageref{11}
\item [{ $\mathbb{B}_i$ }]\begingroup matrix with rows corresponding to weight function moments for all the derivatives up to order of $m$ at particle $i$\nomeqref {33}\nompageref{15}
\item [{ $\mathbb{M}_i$ }]\begingroup moment matrix related to the particle $ i $ [-]\nomeqref {27}\nompageref{14}
\item [{ $\mathbb{N}_{0}^{n}$ }]\begingroup $n$ dimensional space of natural numbers including zero\nomeqref {6}\nompageref{11}
\item [{ $\mathbb{P}_{m}$ }]\begingroup space of polynomials up to order of $m$\nomeqref {27}\nompageref{14}
\item [{ $\mathbb{R}^{n}$ }]\begingroup $n$-dimensional set of real numbers \nomeqref {0}\nompageref{10}
\item [{ $\mathbb{V}_{i}$ }]\begingroup diagonal matrix with diagonal entries  $ W_{ji}^{0} $ for all neighbors of particle $i$\nomeqref {51}\nompageref{19}
\item [{ $\mathbb{W}_i$ }]\begingroup matrix with lines corresponding to $ \bm{W}_{ji} $ evaluated for all neighbors of particle $i$\nomeqref {33}\nompageref{15}
\item [{ $\mathbb{Y}_i$ }]\begingroup matrix with rows corresponding to weight coefficients for all the derivatives up to order of $m$ at particle $i$\nomeqref {33}\nompageref{15}
\item [{ $\mathbf{C}_i$ }]\begingroup matrix of weight functions for all derivatives up to order of $m$ for all neighbors of particle $i$\nomeqref {19}\nompageref{13}
\item [{ $\mathbf{x}$, $\mathbf{x}_{i}$ }]\begingroup spatial position, spatial coordinates of discretization node $ i $ (or point/particles $i$)\nomeqref {0}\nompageref{10}
\item [{ $\mathbf{X}_i$ }]\begingroup matrix with lines corresponding to $ \bm{X}_{ji} $ evaluated for all neighbors of particle $i$\nomeqref {19}\nompageref{13}
\item [{ $\mathcal{C}^{d}$ }]\begingroup mapping vector which selects desired derivatives from vector of derivatives $\bm{D}\left.\left(\cdot\right)\right\rvert_{i}$\nomeqref {2}\nompageref{11}
\item [{ $\mathcal{D}^{\bm{\alpha}}u$ }]\begingroup multi-index derivative of variable $u$\nomeqref {2}\nompageref{11}
\item [{ $\mathcal{J}_{\mathcal{E},i}$ }]\begingroup functional for $\ell_2$ minimization of approximation error [-]\nomeqref {48}\nompageref{18}
\item [{ $\mathcal{J}_{\mathcal{P},i}$ }]\begingroup functional for $\ell_2$ minimization under the moment condition constraint [-]\nomeqref {144}\nompageref{42}
\item [{ $\mathcal{J}_{\mathcal{P},i}$ }]\begingroup functional for $\ell_2$ minimization under the moment condition constraint [-]\nomeqref {97}\nompageref{28}
\item [{ $\mathcal{J}_{\mathcal{P},i}$ }]\begingroup generalized functional for $\ell_2$ minimization under the moment condition constraint [-]\nomeqref {109}\nompageref{30}
\item [{ $\mathcal{J}_{\mathcal{P},i}$ }]\begingroup generalized functional for $\ell_2$ minimization under the moment condition constraint [-]\nomeqref {159}\nompageref{44}
\item [{ $\mathcal{N}_i$ }]\begingroup set of indices of nodes (or points/particles) neighboring with particle $ i $\nomeqref {0}\nompageref{10}
\item [{ $\mathcal{O}(h)$ }]\begingroup big O notation\nomeqref {27}\nompageref{14}
\item [{ $\mathcal{P}$ }]\begingroup set of indices of all discretization nodes (or point/particles)\nomeqref {0}\nompageref{10}
\item [{ $\mathcal{X}$ }]\begingroup set of spatial coordinates of all discretization nodes (or point/particles)\nomeqref {0}\nompageref{10}
\item [{ $\Omega$, $\overline{\Omega}$, $\partial\Omega$ }]\begingroup domain, closure of the domain, domain boundary\nomeqref {0}\nompageref{10}
\item [{ $\varphi_{ji}$ }]\begingroup radial basis function\nomeqref {14}\nompageref{12}
\item [{ $\widetilde{\bm{C}^{\bm{\alpha}}}$ }]\begingroup rescaled mapping vector\nomeqref {28}\nompageref{14}
\item [{ $\widetilde{\bm{X}}_{ji}$ }]\begingroup rescaled vector of Taylor monomials\nomeqref {28}\nompageref{14}
\item [{ $d$ }]\begingroup differential operator\nomeqref {2}\nompageref{11}
\item [{ $e_{ji}^{m}$ }]\begingroup Lagrange form of the remainder between particles $ i $ and $ j $ of order $ m $ \nomeqref {11}\nompageref{11}
\item [{ $h$ }]\begingroup characteristic size of computational stencil \nomeqref {0}\nompageref{10}
\item [{ $i, j, k,\dots$ }]\begingroup indices of individual discretization nodes (or point/particles)\nomeqref {0}\nompageref{10}
\item [{ $l(d)$ }]\begingroup degree the differential operator $ d $ \nomeqref {2}\nompageref{11}
\item [{ $L_{i}^{d}$ }]\begingroup general discrete differential operator (discretizing diff. op. $d$) \nomeqref {2}\nompageref{11}
\item [{ $m$ }]\begingroup  order of the approximation \nomeqref {0}\nompageref{10}
\item [{ $N$ }]\begingroup total number of discretization nodes (or points/particles)\nomeqref {0}\nompageref{10}
\item [{ $n$ }]\begingroup space dimension\nomeqref {0}\nompageref{10}
\item [{ $N_{i}$ }]\begingroup number of nodes (or points/particles) neighboring with particle $ i $\nomeqref {0}\nompageref{10}
\item [{ $p$, $p_\mathrm{2D} $, $p_\mathrm{3D} $ }]\begingroup number of elements in order and dimension dependent vectors, in 2D, in 3D \nomeqref {6}\nompageref{11}
\item [{ $s_{i}$ }]\begingroup average node spacing inside the neighborhood $ \mathcal{N}_{i} $ of node $ i $\nomeqref {0}\nompageref{10}
\item [{ $u$, $u(\vect{x})$, $ u_{i} $ }]\begingroup general variable, g. v. evaluated at $\mathbf{x}$, g. v. evaluated at $\mathbf{x}_i$ \nomeqref {0}\nompageref{10}
\item [{ $u^{*}$ }]\begingroup general function from space of polynomials\nomeqref {27}\nompageref{14}
\item [{ $w_{ji}^{d}$ }]\begingroup weight functions of general discrete differential operator $ d $\nomeqref {2}\nompageref{11}
  \nomgroup{B}
\item [{B-PIM}]\begingroup Boundary-Point Interpolation Method\nomeqref {0}\nompageref{9}
\item [{CMLS}]\begingroup Compact Moving Least Squares Method\nomeqref {0}\nompageref{5}
\item [{D-EM}]\begingroup Diffuse-Element Method\nomeqref {0}\nompageref{3}
\item [{DC-PSE}]\begingroup Discretization-Corrected Particle Strength Exchange\nomeqref {0}\nompageref{5}
\item [{EFG}]\begingroup Element Free Galerkin Method\nomeqref {0}\nompageref{3}
\item [{FPDM}]\begingroup Finite Difference Particle Method\nomeqref {0}\nompageref{5}
\item [{FPM}]\begingroup Finite Point Method\nomeqref {0}\nompageref{5}
\item [{FPsM}]\begingroup Finite Pointset Method\nomeqref {0}\nompageref{5}
\item [{GFDM}]\begingroup Generalized Finite-Difference Method\nomeqref {0}\nompageref{5}
\item [{GMLS}]\begingroup Generalized Moving Least Squares Method\nomeqref {0}\nompageref{4}
\item [{GRKCM}]\begingroup Gradient Reproducing Kernel Collocation Method\nomeqref {0}\nompageref{5}
\item [{HOCSPH}]\begingroup High-order Consistent SPH\nomeqref {0}\nompageref{6}
\item [{IMLS}]\begingroup Interpolating Moving Least Squares Method\nomeqref {0}\nompageref{5}
\item [{IMMLS}]\begingroup Interpolating Modified Moving Least Squares Method\nomeqref {0}\nompageref{5}
\item [{KE}]\begingroup Kernel Estimate\nomeqref {0}\nompageref{4}
\item [{KGF}]\begingroup Kernel Gradient Free\nomeqref {0}\nompageref{6}
\item [{KMM}]\begingroup Kinetic Meshfree Method\nomeqref {0}\nompageref{5}
\item [{LDD}]\begingroup Lagrangian Differencing Dynamics\nomeqref {0}\nompageref{6}
\item [{LFDM}]\begingroup Lagrangian Finite-Difference Method\nomeqref {0}\nompageref{6}
\item [{LSCM}]\begingroup Least Squares Collocation Meshless Method\nomeqref {0}\nompageref{5}
\item [{LSKUM}]\begingroup Least Squares Kinetic Upwind Method\nomeqref {0}\nompageref{5}
\item [{LSMFM}]\begingroup Least-Squares Meshfree Method\nomeqref {0}\nompageref{5}
\item [{MFDM}]\begingroup Meshfree Finite Difference Method\nomeqref {0}\nompageref{5}
\item [{MLPG}]\begingroup Local Petrov Galerkin Method\nomeqref {0}\nompageref{3}
\item [{MLS}]\begingroup Moving Least Squares\nomeqref {0}\nompageref{4}
\item [{MMLS}]\begingroup Modified Moving Least Squares Method\nomeqref {0}\nompageref{5}
\item [{MPS}]\begingroup Moving Particle Semi-Implicit\nomeqref {0}\nompageref{6}
\item [{PDA}]\begingroup Particle Derivative Approximation\nomeqref {0}\nompageref{5}
\item [{PDE}]\begingroup \textit{partial differential equations}\nomeqref {0}\nompageref{2}
\item [{PDM}]\begingroup Particle Difference Method\nomeqref {0}\nompageref{5}
\item [{PIC}]\begingroup Particle-In-Cell\nomeqref {0}\nompageref{2}
\item [{PIM}]\begingroup Point Interpolation Method\nomeqref {0}\nompageref{9}
\item [{PSE}]\begingroup Particle Strength Exchange\nomeqref {0}\nompageref{5}
\item [{PU}]\begingroup Partition of Unity\nomeqref {0}\nompageref{3}
\item [{RBCM}]\begingroup Radial Basis Collocation Method\nomeqref {0}\nompageref{4}
\item [{RKA}]\begingroup Reproducing Kernel Approximations\nomeqref {0}\nompageref{9}
\item [{RKCM}]\begingroup Reproducing Kernel Collocation Method\nomeqref {0}\nompageref{5}
\item [{RKPM}]\begingroup Reproducing Kernel Particle Method\nomeqref {0}\nompageref{4}
\item [{RMLS}]\begingroup Regularized Moving Least Squares Method\nomeqref {0}\nompageref{5}
\item [{SPH}]\begingroup Smoothed Particle Hydrodynamics\nomeqref {0}\nompageref{2}
\item [{WLS}]\begingroup Weighted Least Squares\nomeqref {0}\nompageref{5}

\end{thenomenclature}

\section{Introduction}\label{sec:introduction}
We present an overview and classification of meshfree collocation methods. Meshfree collocation methods are a class of numerical methods for approximating strong-form solutions to partial differential equations (PDEs). They do so without requiring a computational grid or mesh, by discretizing the unknown continuous functions on a set of irregularly placed collocation points. The collocation points therefore form a point cloud, on which the differential operators occurring in the PDE are approximated. In meshfree collocation methods, this is done in the strong sense by directly approximating the derivatives themselves from the discretized function values. Therefore, meshfree collocation methods can be seen as a generalization of finite differences to irregularly distributed collocation points, or nodes.
\nomacro{PDE}{\textit{partial differential equations}}

In addition to being a generalization of finite differences to scattered nodes, meshfree collocation methods are also a special sub-class of the larger algorithmic family of particle methods, although ``meshfree method'' and ``particle method'' are sometimes used synonymously in the literature. Particle methods are a family of numerical algorithms in which point-like objects---called particles---interact pairwise and then evolve as a result of these interactions, as has been rigorously defined~\cite{AUnifyingMathPahlke2023}. They are amongh the oldest and most computationally powerful numerical algorithms~\cite{Pahlke:2025}. In a particle method, particles can  represent discrete model entities, such as atoms in a molecular simulation~\cite{Alder1957, Rahman1964, Stillinger1974}, grains in simulations of granular materials~\cite{Cundall1971, Cundall1979, Cundall1987, Walter2009}, or molecules or material elements in more coarse-grained models~\cite{Hoogerbrugge1992, Espanol1995, Groot1997}. Alternatively, the particles in particle methods can also represent discretization elements or samples of a continuous field. This is the case, for example, for Smoothed Particle Hydrodynamics (SPH) simulations in astrophysics and fluid dynamics~\cite{monaghan1977,lucy1977,Monaghan1992,Monaghan2005}, Particle-In-Cell (PIC) simulations in plasma physics~\cite{Evans1957, Harlow1957, Birdsall1991, Arber2015}, and Vortex Methods for incompressible fluid flow~\cite{Chorin1978, Leonard1980, Gustafson1991, Cottet2000}. In all particle methods, the particles carry discrete or continuous properties, thus forming a labeled point cloud.
\nomacro{SPH}{Smoothed Particle Hydrodynamics}
\nomacro{PIC}{Particle-In-Cell}
The labeled point cloud in meshfree (sometimes also referred to as meshless) methods is one specific case~\cite{chen2017}. Meshfree methods numerically approximate the solution of a PDE on particles that represent the individual discretization points. Then, the particles have positions in the continuous domain of solution of the PDE and continuous properties discretizing the unknown functions.
Once the particle properties have been supplemented with appropriate boundary and/or initial conditions, the particles  evolve in time so as to approximate the solution of the PDE. During time evolution, the particles can change either their positions, their properties, or both. Methods in which particle properties are dynamic, but their positions are constant, are referred to as {\em Eulerian particle methods}. If the particle positions (also) change over time, for example to track a flow field or a deforming domain geometry, particle methods are called {\em Lagrangian}.

Both Eulerian and Lagrangian particle methods are meshfree numerical methods for continuous PDE problems. They can be formulated to approximate the solution of the PDE in either the strong or the weak sense. Meshfree Galerkin methods discretize the weak form of PDE. Even though they are not the topic of this article, we briefly mention them here for completeness. Meshfree Galerkin methods originated from the Diffuse-Element Method (D-EM)\footnote{
The commonly used abbreviation for the Diffuse-Element Method is DEM, but in the context of modern particle methods the abbreviation DEM usually means ``Discrete Element Method''.} as a generalization of the finite element method~\cite{nayroles1992}.
\nomacro{D-EM}{Diffuse-Element Method}
A meshfree Wavelet Galerkin method was proposed by Amaratunga et al.~\cite{amaratunga1994}.
Shortly thereafter Belytschko~\cite{belytschko1994efg, belytschko1995efg} and Lu~\cite{lu1994efg} introduced the Element Free Galerkin method (EFG), improving over the original D-EM.
In the late 1990s, meshfree Galerkin methods based on Partition of Unity (PU) were introduced by Babu\v{s}ka and Melek~\cite{babuska1997pum, melenk1996pum} and by Durate and Liszka~\cite{Duarte1996hp, Liszka1996hp}, known as \textit{hp}-clouds.
Most recently, the meshfree Local Petrov Galerkin Method (MLPG) has been introduced by Atluri et al.~\cite{atluri1998mlpg, atluri2002mlpg}.
In meshfree Galerkin methods, the particles carry (weighted) basis functions and integration volumes of the weak-form PDE solution.
\nomacro{EFG}{Element Free Galerkin Method}
\nomacro{MLPG}{Local Petrov Galerkin Method}
\nomacro{PU}{Partition of Unity}

In contrast to meshfree Galerkin methods, the particles in meshfree collocation methods are volume-less points~\cite{suchde2025} that directly sample the function values at the given locations, thus discretizing the strong form of the PDE.
In this text, we provide an overview of meshfree collocation methods.
For this class of numerical methods, we provide a unifying mathematical definition and show how all reviewed methods adhere to this definition for certain design choices. This makes explicit the differences between methods, but also shows under what conditions they are equivalent. In doing so, we bridge diverse nomenclatures, as different methods were developed by different scientific communities---such as solid mechanics, fluid mechanics, or applied mathematics---each using different terminology. This also led to several methods being re-invented under different names over time, while others were independently derived using novel approaches. We hope we have succeeded in collecting most of the existing methods. We classify all reviewed methods into four groups according to the mathematical principle underlying their derivation. Thereby, we hope to clarify equivalence, differences, and characteristics of meshfree collocation methods.

To better understand the mathematical basis and commonalities of meshfree collocation methods, we start by giving an overview of their historic development.

\subsection{A brief history of meshfree collocation methods\label{sec:brief_history_of_particle_methods}}

Meshfree collocation methods are historically rooted in works on function approximation and interpolation on scattered points from the early twentieth century~\cite{lanczos1938, frazer1937, slater1934}. The foundational reference specifically using these ideas for collocation methods is MacNeal, 1953~\cite{macneal1953}. Studying second-order boundary-value problems with curved boundary, MacNeal generalized finite differences to allow for determining weights for arbitrarily distributed nodes from a linear system of equations. This still is the basic concept underlying meshfree collocation methods today.

The concept of finite differences on scattered nodes was later discussed by Shepard in his work on interpolation of scattered data~\cite{shepard1968}. The function-approximation side of the problem was much later considered by Wendland in his seminal book on scattered data approximation~\cite{wendland2004}. Wendland's book addresses both the theoretical foundation and the practical applications of meshfree approximation, spanning from numerical analysis to computer graphics.

Focusing on numerical analysis and the approximation of differential operators in meshfree collocation methods, MacNeal's pioneering work~\cite{macneal1953} was  generalized and formalized by Forsyth and Wasow~\cite{forsythe1960book} and by Swartz and Wendroff~\cite{swartz1969generalized}. Meshfree collocation methods as we know them today were first formulated by Jensen~\cite{jensen1972}, where also the neighbor search procedure was  given attention.

Indeed, an important early difficulty for meshfree collocation methods was to efficiently find the neighbors of a point when constricting and evaluating the stencils~\cite{perrone1975}. The {\em stencil}, also sometimes referred to as ``star'' is the set of all neighboring points a given point needs to know the values of in order to build its local approximation. Constructing the stencil without having to search over {\em all} points, especially in Lagrangian discretizations where points move, became efficiently possible thanks to fast neighbor-search algorithms like cell lists~\cite{Hockney:1970} and Verlet lists~\cite{Verlet:1967}. The problem remained, though, to choose for each point the correct number of neighbors such that the linear system of equations for the stencil weights is determined.

The constraint of working with determined linear systems was relaxed in the early 1980s with the introduction of Moving Least Squares (MLS) approximation methods~\cite{lancaster1979,lancaster1981}. The concept of MLS had been inspired by the earlier work of Shepard~\cite{shepard1968} and, notably, MacLain~\cite{mclain1974, mclain1976, barnhill1977, gordon1978, franke1980} on irregular data interpolation. This extended the original method of Jensen~\cite{jensen1972} to considering more stencil points than the number of unknowns in the linear system, resulting in overdetermined stencil systems, as introduced by Liszka and Orkisz~\cite{liszka1980finite, liszka1984}.
The works by Liszka and Orkisz are therefore widely regarded as foundational for today's meshfree collocation methods.
\nomacro{MLS}{Moving Least Squares}

In parallel, and since the 1970s, alternative methods have been developed to approximate functions and their derivatives on scattered points.
Monaghan and Gingold~\cite{monaghan1977}, and independently Lucy~\cite{lucy1977}, introduced SPH, originally for applications in astrophysics, which brought the concept of a Kernel Estimate (KE).
A few years later, Chen and Liu elaborated on this concept, introducing the Reproducing Kernel Particle Method (RKPM)~\cite{Chen1996rkpm, Liu1995rkpm}. This was based on imposing a reproducing condition on the approximation kernels, with the aim of improving the accuracy of SPH. Another alternative approach is the Radial Basis Collocation method (RBCM) introduced by Kansa~\cite{kansa1990a, kansa1990b}, where radial basis functions are used to approximate the non-trivial domain integration problem.
\nomacro{KE}{Kernel Estimate}
\nomacro{RKPM}{Reproducing Kernel Particle Method}
\nomacro{RBCM}{Radial Basis Collocation Method}

Within the class of MLS-like methods, development continued with the Generalized Moving Least Squares Method (GMLS) introduced in the late 1990s~\cite{levin1998}. GMLS is not limited to approximating differential operators, but is able to handle more general operators on a Banach space~\cite{wendland2004,mirzaei2012,mirzaei2016,gross2020}.
It is often formulated in more rigorous mathematical language \cite{mirzaei2016, gross2020, mohhamadi2020}.
The connection between GMLS and RKPM was investigated in Ref.~\cite{behzadan2011}, and the connection between GMLS and RBF-FD was established in Ref.~\cite{jones2023}.
\nomacro{GMLS}{Generalized Moving Least Squares Method}

In summary, meshfree collocation methods have a long history. While mostly rooted in interpolation and approximation of scattered data, novel and original ideas have emerged over time. Some came from studying the numerical properties of the methods, some were inspired by particular applications. More recently, there is a clear trend to bridge across methods, study their mathematical connections, and generalize the theory. This is what we attempt in the remainder of this paper for meshfree collocation methods.

\subsection{Methods reviewed in this paper}

Meshfree collocation methods can roughly be divided into MLS-type methods and generalized finite-difference methods (GFDM). The former relate to the idea of operating on local least-squares regression approximations, whereas the latter take the more discrete viewpoint of weighted point-cloud stencils.

Modern MLS-type methods include Compact MLS (CMLS)~\cite{trask2016}, which introduced additional penalization terms to the MLS approximation that are specific to the PDE being solved.
Modified MLS (MMLS)~\cite{joldes2015mmls, joldes2019mmls} introduced additional terms based on the coefficients of the local approximation polynomials to automatically reduce the order of the approximation when not enough particles are available in the neighborhood.
The idea penalized MLS penalties was formalized in Regularized MLS (RMLS) \cite{wang2018rmls}. In the class of interpolating MLS methods, where the weight function is required to additionally satisfy a Kronecker delta property, Interpolating Modified MLS (IMMLS)~\cite{lohit2022immls} generalized Interpolating MLS (IMLS)~\cite{lancaster1981} in the same way as MMLS generalized the original MLS.
This was further extended GMLS to discretizing arbitrary linear bounded operators~\cite{levin1998}.
The numerical stability of these methods in time stepping was enhanced by the Adaptive Orthogonal Improved Interpolating MLS method~\cite{wang2019aoiimls}. Linking MLS and finite-difference approaches, the Least-Squares Meshfree Method (LSMFM)~\cite{park2001, zhang2005} was introduced, which uses a meshfree finite difference based on a MLS approximation.
\nomacro{CMLS}{Compact Moving Least Squares Method}
\nomacro{RMLS}{Regularized Moving Least Squares Method}
\nomacro{MMLS}{Modified Moving Least Squares Method}
\nomacro{IMLS}{Interpolating Moving Least Squares Method}
\nomacro{IMMLS}{Interpolating Modified Moving Least Squares Method}
\nomacro{LSMFM}{Least-Squares Meshfree Method}

GFDM trace back to the works of Jensen~\cite{jensen1972} and Perrone and Kao~\cite{perrone1975}, followed by significant advancements by Liszka and Orkisz~\cite{liszka1980finite, liszka1984, idelsohn1998computational}, eventually bridging them with MLS techniques~\cite{lancaster1981}.
More recent works~\cite{suchde2017, suchde2018, suchde2019} broadened the scope of GFDM and showed the formal connections between the methods in this class, including the Finite Point Method (FPM)~\cite{onate1996fpm}, Finite Pointset Method (FPsM)~\cite{tiwari2002fpsm}\footnote{We here use the acronym FPsM for the finite pointset method in order to avoid confusion with the finite point method, abbreviated as FPM.}, the Meshless Local Strong Form Method (MLSM) \cite{kosec2019mlsm},
Weighted Least Squares (WLS)~\cite{satyaprasad2025}, and Least Squares Collocation Meshless Method (LSCM)~\cite{zhang2001}.
Applications of GFDM include incompressible and free-surface flows~\cite{zhang2016gfdm, suchde2018incompressible}, phase transitions~\cite{suchde2024, lee2023}, geophysics~\cite{muelas2019}, soil mechanics~\cite{michel2017, kosec2019mlsm}, advection-diffusion problems~\cite{prieto2011gfdm}, shallow-water and thin-film flows~\cite{li2017swegfdm, li2021swegfdm, suchde2024b}, multiphase flow~\cite{zhan2023, bharadwaj2024eulerian, jiang2024weno}, sub-surface flow~\cite{michel2021, zhan2022, chen2021}, nonlinear obstacle problems~\cite{chan2013gfdm}, wave propagation~\cite{zhang2025, chang2025, xu2020}, inverse problems~\cite{fan2014gfdm, li2019gfdm}, and PDEs on curved surfaces~\cite{jiang2024manifolds, suchde2019, suchde2019b}. Comprehensive numerical analysis and error estimation is available for GFDMs~\cite{davydov2016error, davydov2018minimal, zheng2022gfdm}.
\nomacro{GFDM}{Generalized Finite-Difference Method}
\nomacro{FPM}{Finite Point Method}
\nomacro{FPsM}{Finite Pointset Method}
\nomacro{LSKUM}{Least Squares Kinetic Upwind Method}
\nomacro{KMM}{Kinetic Meshfree Method}
\nomacro{WLS}{Weighted Least Squares}
\nomacro{LSCM}{Least Squares Collocation Meshless Method}
The origins of GFDM for Lagrangian incompressible flow simulations are associated with FPsM~\cite{Kuhnert2001, Kuhnert2001_2}, whereas FPM~\cite{onate1996fpm, onate1998fpm} was mostly developed for compressible flows in an Eulerian discretization~\cite{onate1996finite, zhang2024}, hyperbolic conservation laws~\cite{Furst2001}, and free-surface flows~\cite{idelsohn2001}.
Despite this difference in their application, FPM and FSPM are mathematically identical in how they approximate derivatives.

A collocation version of RKPM, which avoids numerical quadrature, was introduced with the Reproducing Kernel Collocation Method (RKCM)~\cite{aluru2000}. Similar to MLS, RKCM first approximates the function itself, but the derivatives are obtained as full derivatives of the function approximation. A simplified derivative approximation is used in the Gradient Reproducing Kernel Collocation Method (GRKCM)~\cite{chi2013}, which corresponds to GFDM based on a MLS approximation.
\nomacro{RKCM}{Reproducing Kernel Collocation Method}
\nomacro{GRKCM}{Gradient Reproducing Kernel Collocation Method}

Another sub-class of GFDM originates from the Finite Difference Particle Method (FDPM) for wave problems~\cite{wang2018lagrangian}, which uses a variable smoothing length over Lagrangian particles.
This approach is similar to that of the Particle Difference Method (PDM)~\cite{yoon2021} and the Particle Derivative Approximation (PDA)~\cite{yoon2014extended}. The Meshless Local Strong Form Method (MLSM)~\cite{kosec2019mlsm} is based on similar principles but provides a link to meshfree Galerkin methods.
The same idea can also be found under the name Meshfree Finite Difference Method (MFDM) in solid mechanics~\cite{milewski2012meshless}. However, MFDM provides an original way to construct higher-order approximations without inverting complete moment matrices~\cite{milewski2021mfdm, milewski2022mfdm}.
\nomacro{FPDM}{Finite Difference Particle Method}
\nomacro{PDM}{Particle Difference Method}
\nomacro{PDA}{Particle Derivative Approximation}
\nomacro{MFDM}{Meshfree Finite Difference Method}

A meshfree collocation method that links both MLS-type and GFDM-type methods is Discretization-Corrected Particle Strength Exchange (DC-PSE)~\cite{schrader2010}.
Its derivation was motivated by the approximation errors of the classic Particle Strength Exchange (PSE) method~\cite{degondmasgallic1989, eldredge2002}, which prevented PSE from converging on asymmetric particle distributions.
Inspired by earlier works on discretization correction for particle methods~\cite{cottet1990, hieber2005, bergdorf2005}, DC-PSE results in a general formulation in which SPH, GFDM, MLS, RKPM, and grid-based finite differences are included as special cases~\cite{schrader2010, schrader2012}.
DC-PSE has been successfully used to solve incompressible flow problems~\cite{UsingDcPseOpBouran2016, EntropicallyDaSingh2023}, problems in solid mechanics~\cite{TaylorSeriesEJacque2020, RecoveryByDisZwick2023}, and biological active hydrodynamics~\cite{singh2023solv, singh2023inst, PhysRevX.14.041002}.
The method was also extended to solving PDEs on curved surfaces~\cite{AMeshfreeCollSingh2023}.
\nomacro{PSE}{Particle Strength Exchange}
\nomacro{DC-PSE}{Discretization-Corrected Particle Strength Exchange}

One of the most recent variants of GFDMs is the Local Anistropic Basis Function Method (LABFM)~\cite{HighOrderDiffKing2020, HighOrderSimuKing2022}.
Despite its name, it follows directly from the methods introduced so far, and it is not to be confused with particle methods employing anisotropic kernels for direction-dependent adaptive resolution~\cite{Hacki2015}.
Mathematically, LABFM is almost identical to DC-PSE but uses a different construction for the weight functions. In LABFM, weight functions are generated by orthogonal polynomials~\cite{HighOrderSimuKing2022}, which has clear advantages over the Taylor monomial basis used in DC-PSE.
LABFM so far has mostly been used for direct numerical simulation of combustion~\cite{king2023meshfreedns} and for viscoelastic flows~\cite{king2023meshfreeviscoelastic}, both in complex geometries.

In addition to MLS-type and GFDM-based methods, there are specialized meshfree collocation methods that were not originally designed to approximate derivatives. Instead, they are physics-inspired and usually targeted at a specific PDE or application. These methods take an intermediate position between a phenomenological simulation method and a numerical solver. This includes methods such as the Least Squares Kinetic Upwind Method (LSKUM)~\cite{ghosh1995least}, its variant q-LSKUM~\cite{ghosh2007qlskum}, and the Kinetic meshfree Method (KMM)~\cite{praveen2007kinetic},
all designed to solve inviscid compressible flows and supersonic flows~\cite{ghosh1995least}.
All three methods derive from the Boltzmann transport equation by taking statistical moments in an Eulerian coordinate system. Meshfree finite differences are then used to compute the fluxes from these moments.

Another group of physics-informed GFDM comprises the Lagrangian Finite-Difference Method (LFDM)~\cite{basic2018class, basic2020lagrangian, basic2022ldd} and Lagrangian Differencing Dynamics (LDD)~\cite{peng2021lagrangian}.
Motivated by discretizing the Laplace operator using SPH, and inspired by correction schemes based on Taylor expansion, these methods focus on first and second derivatives, introducing different formulations of correction matrices.
These correction matrices are usually smaller than the matrices obtained from  MLS-based schemes, which makes LFDM and LDD especially suitable for Lagrangian simulations, where the matrices have to be recomputed at each time step. The authors describe these methods as \textit{``a balance between oscillatory but faster SPH and Moving Particle Semi-Implicit (MPS) methods, and accurate but slower MLS–based methods''}~\cite{basic2022ldd}.
Applications of these methods focus on the incompressible Navier-Stokes equation in Lagrangian form. However, LDD was also used for mean-field granular flows~\cite{peng2021lagrangian}, viscoelastic flows~\cite{basic2024visco}, and fluid-structure interaction~\cite{paneer2024lddfsi}.
\nomacro{LDD}{Lagrangian Differencing Dynamics}
\nomacro{LFDM}{Lagrangian Finite-Difference Method}
\nomacro{MPS}{Moving Particle Semi-Implicit}

As is made evident in the context of LDD, meshfree collocation approaches are also commonly found within other particle methods, which are not strictly collocation methods themselves.
For SPH, this includes a large family of kernel correction methods that aim to make the derivative approximation in SPH higher-order consistent by introducing MLS-based correction matrices. This includes the Kernel Gradient Free (KGF) method~\cite{huang2015kgf}, the correction method by Randles and Libersky~\cite{randles1996sphmls} and Bonet and Look~\cite{bonet1999}, High-Order Consistent SPH (HOCSPH) \cite{nasar2021hoch}, and several others~\cite{ren2023correctedsph, zago2021sph, lei2016kgf, hardi2019}. Other examples of method that include meshfree collocation are remeshed particle methods~\cite{Koumoutsakos1995, Cottet2006, Chatelain2008}, where meshfree collocation approaches are often used to interpolate the function values from irregularly distributed particles to a regular Cartesian grid on which the derivatives are then approximated using finite differences~\cite{magni2012}.
This remeshing approach is frequently used in Vortex-Particle methods for incompressible flows~\cite{Koumoutsakos1995, Chatelain2008aircraft}, but also in SPH \cite{chaniotis2002,hieber2005} and in Wavelet Particle methods~\cite{bergdorf2006}.

\nomacro{KGF}{Kernel Gradient Free}
\nomacro{HOCSPH}{High-order Consistent SPH}

In addition to using meshfree collocation as a correction to ensure convergence of other particle methods, there is also a methodologically interesting connection with state-based peridynamics. State-based peridynamics presents an alternative formulation of continuum mechanics by introducing non-local forces in order to deal with large deformations and material separation~\cite{silling2000pd, madenci2013pd}.
There, the divergence of the stress tensor is replaced with an integral formulation of the action potential of the forces. Upon discretizing the integral by a numerical quadrature, this leads to a meshfree collocation method for the force calculation.
It is therefore not surprising that there have been attempts to connect peridynamics operators with other, more canonical, meshfree collocation formulations~\cite{bessa2014pd, hillman2020grkp}.

Taken together, there is a plethora of modern meshfree collocation methods in the literature today. While, as we have outlined in Section~\ref{sec:brief_history_of_particle_methods}, they all derive from a few common ancestors, they often diverged since then as they were developed by different scientific communities. For example, the methods used in computational fluid dyanmics (notably, GFDM, DC-PSE, LDD, LABFM) are disjoint from the methods used in solid mechanics (RKPM, GRKCM, MFDM).
While many of the methods are somewhat heuristic, rigorous numerical analysis is available for some of them, notably for GMLS, RKCM, and DC-PSE\footnote{For meshfree Galerkin methods, numerical analysis is more common in the literature than for meshfree collocation methods.}.
What stands out is that despite the diversity of methods presented in the literature, they are all based on common numerical principles and share many algorithmic properties. Understanding how all of these methods map onto one another could help translate algorithmic improvements and numerical analysis results from one method to the others. This is the main motivation for this paper. But before we present a unifying mathematical formulation for all meshfree collocation methods, we briefly highlight available software implementations and software libraries that can be used to implement meshfree collocation methods.

\subsection{Notable software implementations}
\label{sec:Software}

Many software libraries and codes implement meshfree collocation methods. While it is impossible to enumerate them all, we here provide a brief list of notable implementations with the purpose of highlighting the availability of the methods and the level of support for practical applications. Some of these software implementation specifically provide meshfree PDE solvers, whereas others are more general particle methods frameworks within which meshfree collocation have been implemented. Also, some of these codes are available as open-source free software, whereas others are commercially distributed or closed-source.

A generic open-source framework for particle methods is \textit{OpenFPM}~\cite{OpenfpmAScalIncard2019}, developed by the MOSAIC Group at Dresden University of Technology and the Max Planck Institute of Molecular Cell Biology and Genetics~\cite{openfpm-web}. The Open Framework for Particles and Meshes, \textit{OpenFPM}, provides a generic abstraction layer for implementing particles-only and hybrid particle--mesh codes on parallel computer systems. It uses compile-time code generation to support different hardware platforms, ranging from multi-core CPUs, over distributed-memory clusters, to Graphics Processing Units (GPUs) and distributed multi-GPU clusters~\cite{Incardona2023}. In its numerics module, among several other solvers, \textit{OpenFPM} provides implementations of DC-PSE~\cite{schrader2010} and Surface DC-PSE~\cite{AMeshfreeCollSingh2023} with a convenient template expression system for near-mathematical PDE notation~\cite{ACExpressioSingh2021}.

A comprehensive open-source toolkit for numerical computation on point clouds is provided by \textit{Compadre}~\cite{compadre_toolkit}, which is developed as part of the TRILINOS project~\cite{trilinos-website} for scientific computing. \textit{Compadre} focuses on assembly of the systems of equations for the stencil coefficients in meshfree collocation methods. It therefore excels at handling small dense matrices. For other parts of a meshfree collocation solver, such as for inverting the large sparse matrix of the global system, and for neighbor search, \textit{Compadre} relies on other components of the TRILINOS ecosystem. In the same spirit, \textit{Compadre} supports performance portability across different hardware platforms and GPUs through the \textit{Kokkos} core package of TRILINOS~\cite{Kokkos2022}. Besides providing a generic computational toolkit, \textit{Compadre} also implements the GMLS meshfree collocation method, both in bulk and on curved surfaces.

A more targeted open-source implementation of meshfree collocation methods for PDEs is \textit{Medusa}~\cite{slak2021medusa, Medusa}. \textit{Medusa} implements a range of methods, including FPM, GFDM, WLS, RBF-FD, and MLSM under a single code framework.
It provides wide range of predefined differential operators that can be used to straightforwardly discretize PDEs.

A compact stand-alone open-source implementation of the FPsM is provided by the \textit{gfdm-Toolkit}~\cite{gfdm-toolkit}. It supports two-dimensional simulations of incompressible hydrodynamics, linear elasticity, and heat conduction.
Finally, \textit{mFDlab}~\cite{D2020soft} is an open-source collection of MATLAB functions for developing and testing
meshfree collocation methods. While it supports implementing arbitrary meshfree collocation methods, RBF-FD and a range of MLS-type methods come predefined.

Complementing open-source implementations, there are commercial codes for meshfree collocation methods.
Among them is \textit{MESHFREE}~\cite{MESHFREE} by the Fraunhofer Institute for Industrial Mathematics, implementing GFDM and FPsM for problems from fluid and continuum mechanics. \textit{MESHFREE} enables simulations of free-surface flows, multiphase flows, non-Newtonian flows~\cite{veltmaat2022}, as well as problems with large deformations and moving boundaries. It supports incompressible, weakly compressible, and compressible flows regimes on parallel computing platforms. Coupling options include the Discrete Element Method (DEM), rigid-body solvers, and discrete phase models \cite{bharadwaj2022, bharadwaj2024lagrangian} within the software, as well as with several external structural mechanics software for fluid-structure interaction.

Another commercial code that includes meshfree collocation methods is~\textit{VPS} (formerly \textit{PAM-CRASH}) by the ESI Group \cite{VPS}. While \textit{VPS} is primarily a finite-element code for crash and impact simulations, it also implements the FPsM in a dedicated meshfree collocation module~\cite{tramecon2013, kuhnert2011}. This module originated as a branch of \textit{MESHFREE} and is mostly applied to fluid--structure interaction problems.

A third commercial code for meshfree collocation methods is~\textit{Points} by NOGRID GmbH~\cite{NOGRID}, which likewise originated as a branch of \textit{MESHFREE}. \textit{Points} implements the FPsM for problems across a broad spectrum of computational fluid mechanics. It is able to solve problems with coupled heat and mass transfer, as well as with plastic and elastic material properties.

Another commercial meshfree collocation methods software focusing on fluid mechanics applications is \textit{RHOXYZ}~\cite{RHOXYZ}, which is developed in collaboration with the University of Split, Croatia.
\textit{RHOXYZ} implements the LDD method with a focus on solving incompressible flow problems in a Lagrangian frame, also supporting GPU acceleration. The core solver is supplemented with fluid-structure interaction and with coupling to external Finite Element and Discrete Element solvers.

\subsection{Previous overview texts}

Before delving into the mathematical framework and classification of meshfree collocation methods, we provide a guide to previous review and overview articles as well as books on meshfree methods, in chronological order.

\begin{itemize}
  \item The early development of meshfree collocation methods has been comprehensively summarized by Franke, 1977~\cite{franke1977}, and compared with other early techniques for scattered data in 1982~\cite{franke1982}.
  \item Duarte, 1995~\cite{duarte1995overview}, presented an authoritative overview of meshfree Galerkin methods, covering MLS approaches, D-EM, EFG, Kernel Estimates (including SPH), and Wavelet Galerkin methods.
  \item Belytschko, 1996~\cite{belytschko1996overview}, provided a wide overview of methods ranging from SPH over MLS approaches to PU and Galerkin methods and GFDM. The methods are compared according to their main principles, and application examples are provided from solid mechanics.
  \item Li \& Liu, 2002~\cite{li2002}, gave an overview of meshfree and particle methods in general.
    This provided a cross-sectional view of particle methods with a focus on applications in structural mechanics. The same authors in 2004 then published a now-classic book on particle methods~\cite{li2004book}, particularly focusing on meshfree Galerkin and PU methods. However, the book is much broader and also touches on molecular dynamics methods and multi-scale particle methods. It includes theory, algorithms, and application cases.
  \item A reference book on meshfree methods was written by Liu \& Karamanlidis, 2003~\cite{liu2003}. It contains an introduction to several weak and strong form methods with more in-depth derivations of mostly weak-form methods, such as D-EM, EFG, MLPG, SPH, Point Interpolation Method (PIM), and the Boundary-Point Interpolation Method (B-PIM). The second part of the book focuses on applications in solid (beams, plates, shells) and fluid (convection, hydrodynamics) mechanics.
  \item An overview of meshfree Galerkin methods for PDEs was given by Fries \& Matthies, 2004~\cite{fries2004overview}.
    This contained PU methods and Galerkin approximations with and without extrinsic basis. The paper particularly discussed the choice of test functions. Meshfree collocation methods were not discussed, albeit some mathematical formulations of meshfree finite difference methods can be found in the paper, for example for GMLS, and their relation to meshfree Galerkin methods is discussed.
  \item Nguyen et al., 2008~\cite{nguyen2008overview}, reviewed meshfree methods for solid mechanics problems.
    This mainly included weak-form methods, covering a wide spectrum of aspects ranging from their mathematical formulation to software implementations.
  \item Liu, 2009~\cite{liu2009book} provided a large compendium of meshfree methods. It discusses fundamental concepts of function and derivative approximation on scattered nodes and then delves into several families of methods, including SPH, RKPM, MLS, WLS, PIM, EFG, and MLPG.
    This provides a broad overview of both strong- and weak-form methods, and it also explores their application in solid and fluid mechanics and discusses their integration with other methods in hybrid and multi-scale simulations.
  \item Chen, Hillman \& Chi, 2017~\cite{chen2017} presented a comprehensive overview of meshfree approximation techniques.
    They juxtaposed different ways of obtaining mehsfree approximations of a function, including approximations based on Least Squares (LS, WLS, MLS), KE, Reproducing Kernel Approximations (RKA), and PU methods.
    They categorize meshfree approximation of derivatives into five classes: Direct Derivatives, Diffuse Derivatives, Implicit Derivatives, Synchdonized Derivatives, and Generalized Finite Differences. We will refer to this classification further down when discussing methods.
    Although many meshfree methods are listed, this article did not highlight differences and commonalities across methods.
  \item Jacquemin et al., 2020~\cite{TaylorSeriesEJacque2020}, gave an overview of strong-form meshfree collocation methods derived from Taylor expansions. It is an excellent tutorial article for accessing the field of meshfree methods, as it also provides several examples and benchmark problems to compare methods. For many methods, detailed step-by-step derivations are given. This includes DC-PSE, GFDM, RBF-FD, MLS, and IMLS.
    Finally, this article unifies DC-PSE and GFDM into a single formulation.
  \item One of the most recent books on meshfree methods was authored by Belytschko, Chen \& Hillman in 2024~\cite{belytschko2024book}. This book focuses on meshfree Galerkin methods in MLS and RK formulations, which are introduced by considering problems from linear (diffusion, linear elasticity) and non-linear continuum mechanics. This book also follows the classification of derivative approximations given earlier by two of the authors~\cite{chen2017}, discussing methods using diffuse derivatives, implicit derivatives, synchronized derivatives, and direct derivatives. In this, the problem of quadrature is given particular attention, and various numerical integration methods are compared. Even though the book focuses on meshfree Galerkin methods, the comparison of nodally integrated meshfree Galerkin methods and nodally collocated meshfree methods also includes some material on meshfree collocation methods.
\end{itemize}

\nomacro{B-PIM}{Boundary-Point Interpolation Method}
\nomacro{PIM}{Point Interpolation Method}
\nomacro{RKA}{Reproducing Kernel Approximations}

\textbf{The remainder of this paper is structured as follows:}
We start in Section~\ref{sec:com_definition} by providing a unifying mathematical formulation of meshfree collocation methods for strong-form PDE solution. This allows us to define the notation and review the approximation principles that are common to all meshfree collocation methods. In a small subsection, we also introduce a compact matrix notation for this formulation.

All methods reviewed here follow this general formulation. They are, however, derived differently. We therefore propose a classification of meshfree collocation methods w.r.t.~their principle of derivation. This distinguishes methods derived by moment approximation (Section~\ref{sec:aom}), methods derived by $\ell_2$ minimization of the discretization error (Section~\ref{sec:l2e}), methods derived by $\ell_2$ minimization of the stencil weights (Section~\ref{sec:l2p}), and methods derived by generalized $\ell_2$ minimization (Section~\ref{sec:gl2p}).
In each of these sections, we provide the detailed derivation starting from the unifying definition of Section~\ref{sec:com_definition}, followed by a discussion of the methods that fall into the respective category. In this, we also mention further classifications of the methods, for example w.r.t.~the type of derivative approximation~\cite{chen2017}, and we cite examples of applications.

Finally, in Section~\ref{sec:summary}, we summarize the reviewed methods. This in particular entails a table of all reviewed methods. In this table, the general formulation of each method is given in the notation of Section~\ref{sec:com_definition}, along with its classification according to the optimization principle and the type of derivative approximation.

The main text is followed by four Appendices. Appendix~\ref{sec:appendix-oldAnd2D} recapitulates the general formulation from Section~\ref{sec:com_definition} in component-wise notation in two dimensions. This serves to further clarify notation and provde a didactic example of what the methods look like when all terms are explicitly written out. Appendix~\ref{sec:appendix-l2e} contains the formal derivation of the analytical solution for the $\ell_2$ error minimization from Section~\ref{sec:l2e}. Appendix~\ref{sec:appendix-l2p} does the same for the $\ell_2$ minimization problem over the weights from Section~\ref{sec:l2p}. Appendix~\ref{sec:appendix-gl2p} provides the formal solution for the generalized $\ell_2$ minimization problem from Section~\ref{sec:gl2p}. 
Finally, the Supplementary Material provides all details of every method discussed, serving as a reference.

\section{General formulation of meshfree collocation methods\label{sec:com_definition}}

All numerically consistent methods for approximating differential operators on  scattered nodes follow a common scheme. In this scheme, the function approximation at any given point can be interpreted as an interpolation or convolution-like operation over the values at neighboring particles within a certain finite neighborhood, multiplied with suitable weights. In all consistent methods, the moments of the weight function must fulfill certain conditions. Different meshfree collocation methods mainly differ in how these conditions are enforced (e.g., by minimization, constrained minimization, or exactly), and what ansatz is used for the weight function. Different design choices lead to approximations with varying accuracy and stability properties.

We formulate the general definition of meshfree collocation methods for an open continuous domain $\Omega \subset \mathbb{R}^n $ containing a set of $ N $ collocation points (nodes, particles) $\mathcal{X} \coloneq \{ \vect{x}_1,\, \ldots ,\, \vect{x}_N : \vect{x}_i \in \overline{\Omega} \}$, where $\overline{\Omega} = \Omega \cup \partial \Omega$, at locations $ \vect{x}_i \coloneq ( x_{i,1},\, \ldots ,\, x_{i,n} )^\top = ( x_{i,\beta})_{\beta=1}^{n} $.
Every point can be addressed by its index $ \mathcal{P} \coloneq \{i : \vect{x}_i \in \mathcal{X}\}$.

\nomvar{ $\Omega$, $\overline{\Omega}$, $\partial\Omega$ }{domain, closure of the domain, domain boundary}
\nomvar{ $\mathbb{R}^{n}$ }{$n$-dimensional set of real numbers }
\nomvar{ $n$ }{space dimension}
\nomvar{ $N$ }{total number of discretization nodes (or points/particles)}
\nomvar{ $\mathcal{X}$ }{set of spatial coordinates of all discretization nodes (or point/particles)}
\nomvar{ $\mathbf{x}$, $\mathbf{x}_{i}$ }{spatial position, spatial coordinates of discretization node $ i $ (or point/particles $i$)}
\nomvar{ $i, j, k,\dots$ }{indices of individual discretization nodes (or point/particles)}
\nomvar{ $\alpha, \beta, \gamma, \dots$ }{indices used to index individual spatial components}

Each particle $ i \in \mathcal{P} $ is surrounded by neighboring particles $ \mathcal{N}_i \coloneq \{ j \in \mathcal{P} :  \lVert \vect{x}_j - \vect{x}_i \rVert_2 < h_i \} $ within a disk defined by the \textit{stencil radius} $ h_i $, which is also sometimes called \textit{smoothing length}, \textit{cut-off radius} or \textit{interaction radius}\footnote{Note that $i \in \mathcal{N}_i$.}.
We similarly define neighboring particle sets $\mathcal{N}_{\bm{x}} $ for any point $\bm{x} \in \overline{\Omega}$, not necessarily coinciding with a particle location.
Different particles $i$ can have different $h_i$ to allow for adaptive-resolution methods. Therefore, and because particle distributions can be non-uniform, different points $i$ generally have different numbers of neighbors $ N_i \coloneq \lvert \mathcal{N}_i \rvert $.
A characteristic length scale of the neighborhood can be defined for arbitrary node distributions as $ \Delta x_i = \sqrt[n]{V_n(h_i)/N_i}$, where $ V_n(r) $ is $n$-dimensional volume of a Euclidean ball of radius $r$.

Together with the stencil radius, this defines the \textit{resolution ratio} $h_i/\Delta x_i$, allowing to estimate $ N_i \approx V_n(h_i / \Delta x_i) $ for locally homogeneous point distributions. 
\nomvar{ $\mathcal{P}$ }{set of indices of all discretization nodes (or point/particles)}
\nomvar{ $h$ }{characteristic size of computational stencil }
\nomvar{ $\mathcal{N}_i$ }{set of indices of nodes (or points/particles) neighboring with particle $ i $}
\nomvar{ $N_{i}$ }{number of nodes (or points/particles) neighboring with particle $ i $}
\nomvar{ $s_{i}$ }{average node spacing inside the neighborhood $ \mathcal{N}_{i} $ of node $ i $}

Now consider a sufficiently smooth continuous function $u(\vect{x}) : \Omega \mapsto \mathbb{R} $. For notational simplicity, we take $u(\vect{x})$ to be a scalar-valued function. For differential operators in flat spaces, all derivations apply element-wise for vector- and tensor-valued functions\footnote{For differential operators in curved spaces, such as Riemannian manifolds, additional metric terms appear~\cite{suchde2019b, AMeshfreeCollSingh2023}.}.
The function $u$ is sampled on the particles in $\mathcal{X}$ by assigning each particle $i$ a continuous label $ u_{i} \coloneq u(\vect{x}_i) $.
For convenience, we write the difference between $ u_j $ and $ u_i $ as $ u_{ji} \coloneq u_j - u_i $.
The goal of a meshfree collocation method is then twofold: (1) approximate the function $u(\vect{x})$ at off-particle locations $\vect{x}\neq \vect{x}_i\,\forall i$; (2) approximate derivatives of the function.

\nomvar{ $u$, $u(\vect{x})$, $ u_{i} $ }{general variable, g. v. evaluated at $\mathbf{x}$, g. v. evaluated at $\mathbf{x}_i$ }
\nomvar{ $m$ }{ order of the approximation }

For these approximations to exist, all methods assume that the neighborhood sets $ \mathcal{N}_i $ are $\mathbb{P}_{k}^{n}$-unisolvent for the linear space $\mathbb{P}_{k}^{n}$ of real-valued polynomials in $ n $ variables of degree at most $ k $.
Then, the local interpolant of $u(\vect{x})$ over $\mathcal{N}_i$ is unique, i.e., the zero polynomial is the only polynomial $\bm{p} \in \mathbb{P}_{k}^{n}$ that vanishes on $ \mathcal{N}_i $,
\begin{equation}
  \bm{p}(\vect{x}_j) = 0 \quad \forall \, \vect{x}_j \in \mathcal{N}_{i} \implies \bm{p} \equiv \bm{0}\, .
  \label{eq:unisolvency_property}
\end{equation}
This is sufficient to ensure that the points do not all lie in some lower-dimensional manifold~\cite{bos1991}.
Under this condition, the linear system of equations we define in the following has a unique solution. To ensure unisolvence in dimension $n=1$, $\vect{x}_i \neq \vect{x}_j\,\forall \,(i,j)$ is a sufficient condition. For $n>1$, distinct particle positions are necessary but not sufficient, and more complex algebraic conditions apply~\cite{bos1991}.
Indeed, unisolvence of the node set is equivalent with the Vandermonde matrix $ \mathbf{V} (\mathcal{N}_{i}) $ having full rank, such that $ \ker \mathbf{V} (\mathcal{N}_{i}) = \bm {0} $. This ensures that there exists no hypersurface $ H = Q^{-1}(\bm{0}) $ generated by a polynomial $ \bm{0} \ne Q \in \mathbb{P}_{k}^{n} $ such that $ \mathcal{N}_{i} \subseteq H $.
Unisolvence depends on the spatial distribution of points in the particle neighborhoods.
In particular, the points in $\mathcal{N}_i$ must not lie on some lower-dimensional hypersurface \cite{hecth2020}.

Assuming unisolvence in all neighborhoods $\mathcal{N}_i$, we define a discrete operator $ L^{\bm{\alpha}}_{i} $ at particle $ i $ as the weighted linear combination of function values within its neighborhood:
\begin{equation}
  L^{\bm{\alpha}}_{i}u \coloneq
  \displaystyle\sum_{j\in\mathcal{N}_i}
  u_{j}w^{\bm{\alpha}}_{ji} \, ,
  \label{eq:ba:ddo}
\end{equation}
where $ w^{\bm{\alpha}}_{ji} $ denote the particle weights. They are, in general, continuous functions of the distance between particles $i$ and $j$ and of the stencil radius, hence $ w^{\bm{\alpha}}_{ji} = w^{\bm{\alpha}}( \vect{x}_{j} - \vect{x}_i ,\, h_i) $. In some references, these weight functions are called ``shape functions''. The sum is only over the neighbors $\mathcal{N}_i$, which is important for the computational efficiency of the method. Therefore, the weight functions must be compactly supported with $\mathrm{supp~} w^{\bm{\alpha}}_{ji}(\cdot,\, h) = \mathcal{N}_i ,\, \forall i\in\mathcal{P}$.
\nomvar{ $L_{i}^{d}$ }{general discrete differential operator (discretizing diff. op. $d$) }
\nomvar{ $d$ }{differential operator}
\nomvar{ $l(d)$ }{degree the differential operator $ d $ }
\nomvar{ $w_{ji}^{d}$ }{weight functions of general discrete differential operator $ d $}
\nomvar{ $\bm{D}u$ }{vector of all multi-index derivatives of variable $\psi$ up order of $m$}
\nomvar{ $\mathcal{C}^{d}$ }{mapping vector which selects desired derivatives from vector of derivatives $\bm{D}\left.\left(\cdot\right)\right\rvert_{i}$}
\nomvar{ $\bm{\alpha}$ }{multi-index}
\nomvar{ $\mathcal{D}^{\bm{\alpha}}u$ }{multi-index derivative of variable $u$}

In a meshfree collocation method, we require the discrete operator $L_i^d$ to be a consistent approximation of some derivative of $u(\vect{x})$, with the zeroth derivative being the function itself. Let
\begin{equation}
  \bm{D}u
  \coloneq
  \left(
    \partial^{\bm{\alpha}} u
  \right)_{\bm{|\alpha|}=0}^m
  \text{~~~with~~~} \partial^{\bm{\alpha}}u \coloneq \partial^{\bm{|\alpha|}}u / \partial x_{1}^{\alpha_1} \partial x_{2}^{\alpha_2}\dots\partial x_{n}^{\alpha_n}\, ,
  \label{eq:derivatives_vect}
\end{equation}
be the vector of all derivatives up to and including the $m$-th derivative, indexed by the multi-index $ \bm{\alpha}\in \mathbb{N}_0^n = (\alpha_1,\, \ldots ,\, \alpha_n ) $ with notation $ |\bm{\alpha}| =  \alpha_1 + \ldots + \alpha_n $, $ \bm{\alpha}! = \alpha_1! \ldots \alpha_n! $, and $\vect{x}^{\bm{\alpha}} = x_1^{\alpha_1} \ldots x_{n}^{\alpha_n} $.
In order to select a specific derivative, or a linear differential operator, from this vector, we introduce the \textit{mapping vector} $ \bm{C}^{\bm{\alpha}} $
\begin{equation}
  \bm{C}^{\bm{\alpha}}
  \coloneq
  \left(
    \delta_{\bm{\alpha}\bm{\beta}}
  \right)_{|\bm{\beta}|=0}^m \, ,
  \label{eq:ba_def_mapping_vector}
\end{equation}
where $\delta_{\bm{\alpha}\bm{\beta}}$ is the Kronecker symbol.
For example, the derivatives $ \partial / \partial x $,\, $ \partial / \partial y $, and $ \nabla^2 $ in two dimensions correspond to:
\begin{equation}
  \bm{C}^{\bm{\alpha}}
  =
  \begin{cases}
    \left(
      0, 1, 0, 0, 0, 0, 0 ,0 ,0
    \dots\right)^\top
    &\text{for } \partial / \partial x\, , \\
    \left(
      0, 0, 1, 0, 0, 0, 0, 0 ,0
    \dots\right)^\top
    &\text{for } \partial / \partial y\, , \\
    \left(
      0, 0, 0, 1, 0, 1, 0, 0 ,0
    \dots\right)^\top
    &\text{for } \nabla^2 \,.
  \end{cases}
  \label{eq:selector_vect_example}
\end{equation}
We can then write:
\begin{equation}
  L^{\bm{\alpha}}_{i}u =
  \displaystyle\sum_{j\in\mathcal{N}_i}
  u_{j}w^{\bm{\alpha}}_{ji}
  \approx
  \bm{C}^{\bm{\alpha}} \cdot \left.\bm{D}u\right\rvert_{i}\, .
  \label{eq:general_discrete_diff_operator_zm}
\end{equation}

We require this approximation to have an order of consistency of $ m $, i.e., the numerical error when using the method in stable numerical scheme should converge as $O(h_i^{m})$ for all $h_i < 1$.
Approximating an $|\bm{\alpha}|$-th derivative then requires $\bm{D}u \in \mathbb{R}^p$ and $\bm{C}^{\bm{\alpha}} \in \mathbb{N}_0^p$ with
$ p = \sum_{k=0}^{k=m} (k + n - 1)!/( k! (n-1)!)$ in $n$ dimensions.

\nomvar{ $\delta_{\alpha\beta}$ }{Kronecker delta}
\nomvar{ $\mathbb{N}_{0}^{n}$ }{$n$ dimensional space of natural numbers including zero}
\nomvar{ $p$, $p_\mathrm{2D} $, $p_\mathrm{3D} $ }{number of elements in order and dimension dependent vectors, in 2D, in 3D }

\nomvar{ $\bm{X}$, $\bm{X}_{ji}$ }{vector of Taylor monomials assumed up to order of $m$}

To see the order of consistency of the approximation, we define the vector of scaled (by $1/\alpha_i!$) Taylor monomials
\begin{equation}
  \bm{X}_{ji}
  \coloneq
  \left(
    x_{ji,1}^{\alpha_1},\, \ldots ,\, x_{ji,n}^{\alpha_n} / \bm{\alpha}!
  \right)_{|\bm{\alpha}|=0}^m,
  \label{eq:taylor_monomials}
\end{equation}
with inter-particle distances $ x_{ji} \coloneq x_{j} - x_{i }$.
Thus, $ \bm{X}_{ji} \in \mathbb{R}^p $ and $ \bm{X}_{ji} \neq \bm{X}_{j} - \bm{X}_{i} $ since $\bm{X}(\bm{x})$ is nonlinear.
This allows us to write the Taylor expansion of $ u(\vect{x}) $ over $\mathcal{N}_i$ as:
\begin{equation}
  u_{j}
  =
  \bm{X}_{ji}\cdot\left.\bm{D}u\right\rvert_{i} + e_{ji}^{m+1},
  \label{eq:taylor_series_zm}
\end{equation}
where $e_{ji}^{m+1} $ denotes the approximation error \textit{(Lagrange form of the remainder)} of leading order \mbox{$m +1$}.

To obtain the approximation, we first substitute the Taylor expansion for $u_j$ from Eq.~\eqref{eq:taylor_series_zm} into Eq.~\eqref{eq:general_discrete_diff_operator_zm}.
Dropping the error term, this yields:
\begin{equation}
  L^{\bm{\alpha}}_{i}u
  =
  \displaystyle\sum_{j\in\mathcal{N}_i}
  \left.\bm{X}_{ji}\cdot\bm{D}u\right\rvert_{i}w^{\bm{\alpha}}_{ji} \, .
  \label{eq:error_short_zm}
\end{equation}
Introducing the vector of moments $ \bm{B}_{i}^{\bm{\alpha}} \in \mathbb{R}^{p} $ of the weights as
\begin{equation}
  \bm{B}^{\bm{\alpha}}_{i}
  \coloneq
  \displaystyle\sum_{j\in\mathcal{N}_i}
  \bm{X}_{ji}w^{\bm{\alpha}}_{ji}\, ,
  \label{eq:moments_vect_zm}
\end{equation}
this becomes
\begin{equation}
  L^{\bm{\alpha}}_{i}u
  =
  \left.\bm{D}u\right\rvert_{i} \cdot \bm{B}^{\bm{\alpha}}_{i} \, .
  \label{eq:error_moments_vect_zm}
\end{equation}
\nomvar{ $e_{ji}^{m}$ }{Lagrange form of the remainder between particles $ i $ and $ j $ of order $ m $ }
\nomvar{ $\bm{B}_{i}^{\bm{\alpha}}$ }{vector of moments of weight functions $w_{ji}^{d}$}

We assume the $ w_{ji}^{\bm{\alpha}} $ to be linear combinations of basis functions. In meshfree collocation methods, the basis can be general \textit{Anisotropic Basis Functions} (ABFs) $ W_{ji}^{q} = W^{q}(\vect{x}_j - \vect{x}_i, \, h_i) $. This is in contrast to \textit{Radial Basis Functions} (RBFs) $\varphi_{ji} = \varphi_{ji}(\| \vect{x}_{j} - \vect{x}_{i}\|_2,\, h_i)$ used for example in SPH. The ABFs $ W_{ji}^{q} $ depend on $ \vect{x}_{ji} $, whereas RBFs are radially symmetric and only depend on $ \| \vect{x}_{ji} \|_2 $.
We therefore have:
\begin{equation}
  w^{\bm{\alpha}}_{ji}
  \coloneq
  \bm{W}_{ji}\cdot\bm{\Psi}^{\bm{\alpha}}_{i}
  =
  W^{0}_{ji}\Psi^{\bm{\alpha}}_{i,0}
  +W^{1}_{ji}\Psi^{\bm{\alpha}}_{i,1}
  +W^{2}_{ji}\Psi^{\bm{\alpha}}_{i,2}
  +W^{3}_{ji}\Psi^{\bm{\alpha}}_{i,3}
  +W^{4}_{ji}\Psi^{\bm{\alpha}}_{i,4}
  +\dots+W^{n}_{ji}\Psi^{\bm{\alpha}}_{i,p}\, ,
  \label{eq:w_abf_zm}
\end{equation}
where $ \bm{W}_{ji} \in \mathbb{R}^p $ and $ \bm{\Psi}^{\bm{\alpha}}_{i} \in \mathbb{R}^p $ are the vector of basis functions and the vector of unknown coefficients, respectively:
\begin{equation}
  \bm{W}_{ji}
  \coloneq
  \left( W^{0}_{ji},\,W^{1}_{ji},\,W^{2}_{ji},\,W^{3}_{ji},\,\dots,\,W^{p}_{ji}\right)^\top
  \label{eq:w_vect_zm}
\end{equation}
\begin{equation}
  \bm{\Psi}_{i}^{\bm{\alpha}}
  \coloneq
  \left(\Psi^{\bm{\alpha}}_{i,0},\,\Psi^{\bm{\alpha}}_{i,1},\,\Psi^{\bm{\alpha}}_{i,2},\,\Psi^{\bm{\alpha}}_{i,3},\,\dots,\,\Psi^{\bm{\alpha}}_{i,p}\right)^\top.
  \label{eq:coef_vect_zm}
\end{equation}
The specific choice of ABFs is a pivotal design choice that distinguishes different meshfree collocation methods.
Most existing methods use a polynomial ABF basis $ \bm{P}_{ji} \in \mathbb{R}^p $ multiplied with a compactly supported (for computational efficiency) RBF $ \varphi_{ji} $, thus $ \bm{W}_{ji} \coloneq \bm{P}_{ji} \varphi_{ji} $. Some references refer to the RBFs $\varphi$ as ``basis functions'',  others confusingly call either the ABFs or the RBFs ``weights'', which we here reserve for the stencil coefficients $w^{\bm{\alpha}}_{ji}$.

\nomvar{ $\bm{W}$, $\bm{W}_{ji}$, $ W_{ji}^{q} $ }{vector of Anisotropic Basis Functions (ABFs), single ABF}
\nomvar{ $\bm{\Psi}_{i}^{\bm{\alpha}}$ }{vector of weight coefficients used to construct weight functions}
\nomvar{ $\bm{P}_{ji}$ }{vector of basis functions evaluated between particles $i$ and $j$}
\nomvar{ $\varphi_{ji}$ }{radial basis function}

The goal of every meshfree collocation method is to find weight functions $ w^{\bm{\alpha}}_{ji} $ in all $\mathcal{N}_i$ such that Eq.~\eqref{eq:general_discrete_diff_operator_zm} is fulfilled for the discrete operator $L^{\bm{\alpha}}_{i}$ from Eq.~\eqref{eq:ba:ddo}. This can be achieved in one of four ways.

The first way is based on assuming a particular parametric ansatz for the weight function, usually a linear combination of certain basis functions~\cite{degondmasgallic1989, eldredge2002, schrader2010}. Consistency is then ensured by determining the unknown coefficients of the weight function. Substituting the ansatz for the weight function into the Taylor expansion of Eq.~\eqref{eq:taylor_series_zm} leads to conditions on the moments $\sum_j \bm{X}_{ji} w_{ji}^{\bm{\alpha}}$ of the weight function given a certain derivative $\bm{C}^{\bm{\alpha}}\cdot\bm{D}u$. These moment conditions define a liner system of equations for the coefficients of the weight function. We refer to this approach as \textbf{approximation of moments} (AOM), and discuss it in more detail in Section~\ref{sec:aom}.

The second way to derive a consistent weight function starts with expressing the derivative values at the neighboring collocation points from the Taylor expansion in Eq.~\eqref{eq:taylor_series_zm}. This defines a linear system of equations where the system matrix depends on the positions of the neighboring particles~\cite{liszka1980finite}. The unknowns in this case are not the parameters of the weight function, like in AOM, but directly the values of the derivatives. There is one equation for each particle in $\mathcal{N}_i$. These systems are usually over-determined, and the unknown derivatives are approximated as a least-squares solution. We therefore refer to this approach as \textbf{$\ell_2$ minimization of the discretization error } ($\ell_2(\mathcal{E})$), and discuss it in more detail in Section~\ref{sec:l2e}.

The third way starts from Eq.~\eqref{eq:general_discrete_diff_operator_zm} and formulates an objective functional over a certain space of weight functions, for example polynomials. This functional is then optimized by minimizing some norm of the weights $w_{ji}^{\bm{\alpha}}(\bm{x}_j)$ at the neighboring points under the constraint of AOM-like moment conditions~\cite{suchde2017,suchde2018,suchde2019}.
We refer to this approach as \textbf{$\ell_2$ minimization of weights} ($\ell_2(\mathcal{P})$), and discuss it in more detail in Section~\ref{sec:l2p}.

The fourth approach generalizes $\ell_2$ minimization of weights by dropping the additional moment constraints. This is done by modifying the objective functional such that the resulting unconstrained solutions match the approximations obtained by AOM. While this approach is not found in the literature so far, and no methods exist that are explicitly based on it, we introduce it here to show the connection between AOM and minimization-based approaches. We call this the \textbf{generalized $\ell_2$ minimization} approach (g$\ell_2(\mathcal{P})$), and discuss it in more detail in Section~\ref{sec:gl2p}.

\subsection{Matrix notation}\label{sec:basic-approach-matrix-formulation}

Meshfree collocation methods can be written more compactly when using matrix notation, where summations over neighboring particles become matrix multiplications.
This is achieved by \textit{stacking} the individual particle pairs over rows, which requires defining an ordering over the neighbors $j\in\mathcal{N}_i$ of each particle $i$. Recall the definition of the general discrete operator from Eq.~\eqref{eq:general_discrete_diff_operator_zm} and write:
\begin{equation}
  L^{\bm{\alpha}}_{i} u
  =
  \sum_{j \in \mathcal{N}_{i}}
  u_{j}w^{\bm{\alpha}}_{ji}
  =
  \sum_{j=1}^{N_{i}}
  u_{j,i}w^{\bm{\alpha}}_{ji} \, ,
\end{equation}
where  $  u_{j, i} $ refers to the $ j $-th neighbor of particle $ i $.
Since the sum is commutative, the ordering introduced in the second step is arbitrary. While the method is mathematically invariant to the ordering, order may matter in practical software implementations, e.g., due to the non-commutativity of finite-precision arithmetics and for efficient memory access. With any ordering, though, we can write:
\begin{equation}
  \mathbf{X}_i
  =
  \left(
    \begin{array}{ccc}
      \horzbar & \bm{X}_{1i} & \horzbar \\
      \horzbar & \bm{X}_{2i} & \horzbar \\
      & \vdots    &          \\
      \horzbar & \bm{X}_{N_ii} & \horzbar
    \end{array}
  \right)
  \in \mathbb{R}^{N_{i}\times p}
  ,\qquad
  \bm{\bm{U}}_i
  =
  \left(
    \begin{array}{c}
      u_{1,i}   \\
      u_{2,i}   \\
      \vdots                    \\
      u_{N_i,i}
    \end{array}
  \right)
  \in \mathbb{R}^{N_{i}}
  ,\qquad
  \left.\bm{D} u\right\rvert_{i}
  \in \mathbb{R}^{p}\, ,
  \label{eq:appendix-introduction-of-matrix-formalism}
\end{equation}
where $ \mathbf{X}_i $ is the matrix of scaled Taylor monomials\footnote{
  The matrix $\mathbf{X}_{i}$ is related to the Vandermonde matrix $\mathbf{V}(\mathcal{N}_{i})$, but the elements of $\mathbf{X}_{i}$ are scaled by $1/\alpha_i!$, see Eq.~\eqref{eq:taylor_monomials}.
}.
Using this notation, Eq.~\eqref{eq:taylor_series_zm} becomes:
\begin{equation}
  \bm{\bm{U}}_i = \mathbf{X}_i \left.\bm{D} u\right\rvert_{i} + \bm{E}_i \, ,
\end{equation}
where $ \bm{E}_i \in \mathbb{R}^{N_{i}}$ is the vector of truncation errors of order \mbox{$m+1$} or higher.

Ordering the weights $ w_{ji}^{\bm{\alpha}} $ of the neighbors in the same way, $ \bm{w}_{\bullet i}^{\bm{\alpha}} \in \mathbb{R}^{N_{i}} $, leads to:
\begin{equation}
  \mathbf{C}_i
  =
  \left(
    \begin{array}{ccc}
      \horzbar & \bm{w}_{1i} & \horzbar \\
      \horzbar & \bm{w}_{2i} & \horzbar \\
      & \vdots      &          \\
      \horzbar & \bm{w}_{Ni} & \horzbar
    \end{array}
  \right)
  =
  \left(
    \begin{array}{cccc}
      \vertbar                & \vertbar                &        & \vertbar                   \\
      \bm{w}_{\bullet i}^{0}  & \bm{w}_{\bullet i}^{1}  & \ldots & \bm{w}_{\bullet i}^{p}  \\
      \vertbar                & \vertbar                &        & \vertbar
    \end{array}
  \right)
  \in \mathbb{R}^{N_{i}\times p} \, ,
  \label{eq:appendix_cmatrix}
\end{equation}
for all possible derivatives $\bm{\alpha}$. This brings Eq.~\eqref{eq:general_discrete_diff_operator_zm} into the final matrix form
\begin{equation}
  L^{\bm{\alpha}}_{i} u
  =
  \sum_{j\in\mathcal{N}_i}
  u_{j}w^{\bm{\alpha}}_{ji}
  =
  \bm{\bm{U}}_i \cdot \bm{w}_{\bullet i}^{\bm{\alpha}}
  =
  \bm{\bm{U}}_i \cdot
  \left(
    \mathbf{C}_i\bm{C}^{\bm{\alpha}}
  \right)\, ,
  \label{eq:general_discrete_diff_operator_zm_matrix_form}
\end{equation}
where $\bm{C}^{\bm{\alpha}} \in \mathbb{R}^{p} $ is the \textit{mapping vector} defined in Eq.~\eqref{eq:ba_def_mapping_vector}.

\nomvar{ $\mathbf{X}_i$ }{matrix with lines corresponding to $ \bm{X}_{ji} $ evaluated for all neighbors of particle $i$}
\nomvar{ $\bm{\bm{U}}_i$ }{vector containing values of function $ u$ for all neighbors of particle $i$}
\nomvar{ $\bm{E}_i$ }{vector of reminders in Lagrange form for all neighbors of particle $i$}
\nomvar{ $\mathbf{C}_i$ }{matrix of weight functions for all derivatives up to order of $m$ for all neighbors of particle $i$}

In matrix notation, the moments of the weight function are:
\begin{equation}
  \bm{B}_{\bullet i}^{\bm{\alpha}}
  =
  \displaystyle\sum_{j\in\mathcal{N}_i}
  \bm{X}_{ji}w^{\bm{\alpha}}_{ji}
  =
  \mathbf{X}_i \bm{w}_{\bullet i}^{\bm{\alpha}}
  =
  \mathbf{X}_i^{\top} \mathbf{C}_i\bm{C}^{\bm{\alpha}} \, .
  \label{eq:moments_vect_zm_matrix_form}
\end{equation}
We introduce the matrix $ \mathbf{W}_i$ of ABFs $ \bm{W}_{ji} \in \mathbb{R}^{p} $ and the matrix $\mathbf{Y}_i $ of the kernel coefficients $ \bm{\Psi}_{i}^{\bm{\alpha}} \in \mathbb{R}^{p} $, ordered as in  Eq.~\eqref{eq:appendix-introduction-of-matrix-formalism}:
\begin{equation}
  \mathbf{W}_i
  =
  \left[
    \begin{array}{ccc}
      \horzbar & \bm{W}_{1i} & \horzbar \\
      \horzbar & \bm{W}_{2i} & \horzbar \\
      & \vdots      &          \\
      \horzbar & \bm{W}_{Ni} & \horzbar
    \end{array}
  \right]
  \in \mathbb{R}^{N_i\times p}
  ,\qquad
  \mathbf{Y}_i
  =
  \left[
    \begin{array}{cccc}
      \vertbar      & \vertbar      &        & \vertbar       \\
      \bm{\Psi}^0_{i} & \bm{\Psi}^1_{i} & \ldots & \bm{\Psi}^{p}_{i}  \\
      \vertbar      & \vertbar      &        & \vertbar
    \end{array}
  \right]
  \in \mathbb{R}^{p\times p}
  \label{eq:appendix_wmatrix_ymatrix}
\end{equation}
and group the moments $ \bm{B}_{\bullet i}^{\bm{\alpha}} \in \mathbb{R}^N $ for all  derivatives
\begin{equation}
  \mathbf{B}_i
  =
  \left[
    \begin{array}{cccc}
      \vertbar               & \vertbar               &        & \vertbar                   \\
      \bm{B}_{\bullet i}^{0} & \bm{B}_{\bullet i}^{1} & \ldots & \bm{B}_{\bullet i}^{p}  \\
      \vertbar               & \vertbar               &        & \vertbar
    \end{array}
  \right]
  \in \mathbb{R}^{p\times p}\, ,
\end{equation}
such that element $ B_{ji}^{m} $ is the $m$-th moment related to $j$-th neighbor.

\section{Derivation by approximation of moments}\label{sec:aom}

Approximation of moments (AOM) methods find a weight function $ w^{\bm{\alpha}}_{ji} = w^{\bm{\alpha}}( \vect{x}_j - \vect{x}_i ,\, h_i)$ of the discrete operator in Eq.~\eqref{eq:general_discrete_diff_operator_zm} such that the moments of the desired derivatives in Eq.~\eqref{eq:error_moments_vect_zm} are equal to one, and all others (up to and including order $ m $) become zero.
For this, we substitute the ansatz for the weight function from Eq.~\eqref{eq:w_abf_zm} into Eq.~\eqref{eq:moments_vect_zm}.
The moments then become
\begin{equation}
  \bm{B}_{i}^{\bm{\alpha}}
  =
  \sum_{j\in\mathcal{N}_i}\bm{X}_{ji}\bm{W}_{ji}\cdot\bm{\Psi}_{i}^{\bm{\alpha}}\, ,
  \label{eq:moments_vect_subst_zm}
\end{equation}
which defines a linear system of equations for the unknown coefficients $ \bm{\Psi}_{i}^{\bm{\alpha}} $ with moment matrix $ \mathbf{M}_i \in \mathbb{R}^{p\times p}$:
\begin{equation}
  \mathbf{M}_{i}\bm{\Psi}^{\bm{\alpha}}_{i}
  =
  \bm{B}^{\bm{\alpha}}_{i}\, ,
  \quad
  \mathbf{M}_{i}
  \coloneq
  \sum_{j\in\mathcal{N}_i}
  \bm{X}_{ji}\bm{W}_{ji}^\top \, .
  \label{eq:linSystem_all_moments_zm}
\end{equation}
Imposing the values for the moments of a given target derivative $ \bm{C}^{\bm{\alpha}} $ amounts to replacing the vector $ \bm{B}^{\bm{\alpha}}_{i} $ with the appropriate $ \bm{C}^{\bm{\alpha}} $.
This then yields the weight-function coefficients
\begin{equation}
  \bm{\Psi}^{\bm{\alpha}}_{i}
  =
  \mathbf{M}_{i}^{-1} \bm{C}^{\bm{\alpha}}\, ,
  \label{eq:linSystem_zm}
\end{equation}
and thus the discrete operator over the neighborhood of each particle,
\begin{equation}
  L^{\bm{\alpha}}_{i}u
  =
  \sum_{j\in\mathcal{N}_i}
  u_{j}\bm{W}_{ji}\cdot \left(\mathbf{M}_{i}^{-1} \bm{C}^{\bm{\alpha}} \right).
  \label{eq:general_discrete_diff_operator_with_coefs_do}
\end{equation}

The operator in Eq.~\eqref{eq:general_discrete_diff_operator_with_coefs_do} is linear and amounts to a convolution operation with a matrix kernel. When fulfilling the moment conditions up to order $m$, it ensures polynomial consistency of order $ m $. Therefore, the error when approximating an $ l $-th derivative is $ o( h_{i}^{-m - 1 + l} ) $, where $o $ represents the standard Landau small-o notation.
For any resolution ratio $ h_i/\Delta x_i \in O(1) $, this implies an error of $ o( \Delta x_{i}^{-m - 1 + l} ) $, hence:
\begin{equation}
  L^{\bm{\alpha}}_{i}u
  =
  \displaystyle\sum_{j\in\mathcal{N}_i}
  u_{ji}w_{ji}^{d}
  =
  \bm{C}^{\bm{\alpha}} \cdot \left.\bm{D} u \right\rvert_{i} + o( \Delta x_{i}^{-m - 1 + l} )
  \iff
  L^{\bm{\alpha}}_{i}u^*
  = \bm{C}^{\bm{\alpha}} \cdot \left.\bm{D} u^* \right\rvert_{i} \quad \forall u^* \in \mathbb{P}_m \, ,
  \label{eq:general_discrete_diff_operator_theoret_orderd}
\end{equation}
where $ \mathbb{P}_m $ is the ring of polynomials of degree at most $ m $.
In other words, the discrete differential operator is exact for polynomials up to and including degree $ m $.
The condition expressed in Eq.~\eqref{eq:general_discrete_diff_operator_theoret_orderd}, more precisely the part $L^{\bm{\alpha}}_{i}u^* = \bm{C}^{\bm{\alpha}} \cdot \left.\bm{D}u^* \right\rvert_{i}$, is sometimes referred to as \textit{gradient reproducing condition}~\cite{chi2013}, \textit{discrete reproducing condition}~\cite{chen1998}, or \textit{discrete moment condition}~\cite{schrader2010}. In Appendix~\ref{sec:appendix-oldAnd2D}, the derivation is repeated in 2D with all components explicitly written out. This can help grasp the general concept of meshfree collocation methods.

\nomvar{ $\mathbb{M}_i$ }{moment matrix related to the particle $ i $ [-]}
\nomvar{ $\mathcal{O}(h)$ }{big O notation}
\nomvar{ $\mathbb{P}_{m}$ }{space of polynomials up to order of $m$}
\nomvar{ $u^{*}$ }{general function from space of polynomials}

The main difficulty in the AOM approach is the condition number of the moment matrix $ \mathbf{M}_i $ when numerically solving the linear system in Eq.~\eqref{eq:linSystem_all_moments_zm}. In general, this condition number depends on the ABF and on the distribution of points in $\mathcal{N}_i$. The closer the Vandermonde matrix from Section~\ref{sec:com_definition} is to singular, the worse the condition number of $\mathbf{M}_i$\footnote{Finding bounds on the condition number depending on the algebraic set of the neighbor particles and the ABF is, to our knowledge, an open problem.}. In practice, this often requires the use of pre-conditioning schemes. The simplest form of pre-conditioning scales each row with a power of $ h_i $, defining the diagonal pre-conditioning matrix
\begin{equation}
  \mathbf{H}_{i}
  \coloneq
  \mathrm{diag}\left(
    h_{i}^{-|\bm{\alpha}|}
  \right)_{|\bm{\alpha}|=0}^m \in \mathbb{R}^{p\times p} \, .
  \label{eq:h_vect_zm}
\end{equation}
We then define $ \widetilde{\bm{X}}_{ji} \coloneq \mathbf{H}_{i}\bm{X}_{ji} \in \mathbb{R}^p $ and $ \widetilde{\bm{C}^{\bm{\alpha}}} \coloneq \mathbf{H}_{i}\bm{C}^{\bm{\alpha}} \in \mathbb{R}^p $. Replacing $ \bm{X}_{ji} $ by $ \widetilde{\bm{X}}_{ji} $ and $ \bm{C}^{\bm{\alpha}} $ by $ \widetilde{\bm{C}^{\bm{\alpha}}} $ leaves the analytical solution of Eq.~\eqref{eq:linSystem_zm} unchanged, but improves the condition number of $ \mathbf{M}_i  $ in a numerical solver.

\nomvar{ $\bm{H}_i$ }{vector of inverse stencil size scales following Taylor monomial powers up to order of $m$}
\nomvar{ $\widetilde{\bm{C}^{\bm{\alpha}}}$ }{rescaled mapping vector}
\nomvar{ $\widetilde{\bm{X}}_{ji}$ }{rescaled vector of Taylor monomials}

\subsection{Formulation without the zeroth moment \label{sec:basic_approach_without_zeroth_moment}}

Following the original symmetric particle strength exchange (PSE) formulation~\cite{degond1989, degondmasgallic1989}, some works~\cite{schrader2010, milewski2012meshless, yoon2014extended, wang2018lagrangian}
write the operator $ L_i^{\bm{\alpha}} $ from Eq.~\eqref{eq:general_discrete_diff_operator_zm} as a function of the {\em differences of function values} $ u_{ji} = u_j - u_i $, rather than the values $u_j$ themselves, thus:
\begin{equation}
  L^{\bm{\alpha}}_{i}u
  \coloneq
  \displaystyle \sum_{j\in\mathcal{N}_i}
  u_{ji}w^{\bm{\alpha}}_{ji}
  \approx
  \bm{C}^{\bm{\alpha}}\cdot\left.\bm{D}u\right\rvert_{i}\, .
  \label{eq:general_discrete_diff_operator_without_zer_moment}
\end{equation}
This allows using weight functions with a singularity at $ w_{ii}^{\bm{\alpha}}$~\cite{HighOrderDiffKing2020}.
The indices in the vector of derivatives $ \bm{D}u $, the mapping vector $ \bm{C}^{\bm{\alpha}} $, and the vector of scaled monomials $ \bm{X}_{ji} $ then only start from 1:
\begin{equation}
  \bm{D}u
  \coloneq
  \left(
    \partial^{\bm{\alpha}} u
  \right)_{|\bm{\alpha}|=1}^m\, ,
  \quad
  \bm{C}^{(\bm{\alpha})}
  \coloneq
  \left(
    \delta_{\bm{\alpha}\bm{\beta}}
  \right)_{|\bm{\beta}|=1}^m \, ,
  \quad
  \bm{X}_{ji}
  \coloneq
  \left(
    x_{1,ji}^{\alpha_1},\, \ldots ,\, x_{n,ji}^{\alpha_n} / \bm{\alpha}!
  \right)_{|\bm{\alpha}|=1}^m \, .
  \label{eq:basic_vectors_without_zeroth_moment}
\end{equation}
This means that the zeroth moment of the weights is excluded and the length of the vectors is $ p = \sum_{k=1}^{k=m} (k + n - 1)!/( k! (n-1)!)$ for an $m$-th order approximation in $n$ dimensions.
This leads to an operator of the form:
\begin{equation}
  L^{\bm{\alpha}}_{i}u
  =
  \sum_{j\in\mathcal{N}_i}
  u_{ji}\bm{W}_{ji}\cdot \left(\mathbf{M}_{i}^{-1} \bm{C}^{\bm{\alpha}} \right)\, .
  \label{eq:general_discrete_diff_operator_with_coefs_withouth_zeroth_moment}
\end{equation}
The order of consistency remains as in Eq.~\eqref{eq:general_discrete_diff_operator_theoret_orderd}. Excluding the zeroth moment slightly reduces (by 1 row) the cost of solving the linear system in Eq.~\eqref{eq:linSystem_all_moments_zm}. More importantly, though, it improves the condition number of the moment matrix by eliminating the $O(1)$ entries~\cite{UsingDcPseOpBouran2016} and improves the stability of the discrete operator in explicit time stepping~\cite{schrader2010}.

\subsection{Matrix notation}
\label{sec:aom:matrix-formulation}

Following the notation introduced in Section~\ref{sec:basic-approach-matrix-formulation}, we see that $ \mathbf{C}_i = \mathbf{W}_i\mathbf{Y}_i $ and $ \mathbf{B}_i = \mathbf{X}_i^{\top}\mathbf{C}_i = \mathbf{X}_i^{\top} \mathbf{W}_i\mathbf{Y}_i $.
This defines a system of $p$ linear equations as in Eq.~\eqref{eq:linSystem_all_moments_zm} with moment matrix
\begin{equation}
  \mathbf{M}_{i}
  =
  \displaystyle\sum_{j\in\mathcal{N}_i}
  \bm{X}_{ji}\bm{W}_{ji}^{\top}
  =
  \mathbf{X}_i^{\top} \mathbf{W}_i
  \in \mathbb{R}^{p\times p}\, .
  \label{eq:momentum_matrix_matrix_form}
\end{equation}
The resulting form of the discrete differential operator then is:
\begin{equation}
  L^{\bm{\alpha}}_{i} u
  =
  \bm{\bm{U}}_i \cdot \bm{w}_{\bullet i}^{\bm{\alpha}}
  =
  \bm{\bm{U}}_i \cdot \left( \mathbf{C}_i\bm{C}^{\bm{\alpha}} \right)
  =
  \bm{\bm{U}}_i \cdot \left( \mathbf{W}_i \mathbf{Y}_{i}\bm{C}^{\bm{\alpha}} \right)
  =
  \bm{\bm{U}}_i \cdot \left( \mathbf{W}_i \mathbf{M}_{i}^{-1}\bm{C}^{\bm{\alpha}} \right)\, ,
  \label{eq:general_discrete_diff_operator_with_coefs_do_matrix_form}
\end{equation}
since $ \mathbf{Y}_{i} = \mathbf{M}_{i}^{-1} $.

\nomvar{ $\mathbb{W}_i$ }{matrix with lines corresponding to $ \bm{W}_{ji} $ evaluated for all neighbors of particle $i$}
\nomvar{ $\mathbb{Y}_i$ }{matrix with rows corresponding to weight coefficients for all the derivatives up to order of $m$ at particle $i$}
\nomvar{ $\mathbb{B}_i$ }{matrix with rows corresponding to weight function moments for all the derivatives up to order of $m$ at particle $i$}

\subsection{Methods derived by approximation of moments}

A typical moment-approximation method is \textbf{Discretization-Corrected Particle Strength Exchange (DC-PSE)}. This method was originally introduced by Schrader et al.~\cite{schrader2010} and later revisited in several other works~\cite{Schrader:2012, RecoveryByDisZwick2023, EntropicallyDaSingh2023} and  extended to differential operators on curved surfaces~\cite{AMeshfreeCollSingh2023}, spatially adaptive~\cite{Reboux:2012} and anisotropic~\cite{Hacki2015} resolution, and to PDEs with discontinuous coefficients~\cite{Fietier:2013}.
The original motivation for DC-PSE was to render the classic \textbf{Particle Strength Exchange (PSE)} operators~\cite{degond1989, eldredge2002} consistent on irregular particle distributions. This was achieved by replacing the continuous, integral moment conditions of PSE with discrete sums in the spirit of collocation methods as presented in Section~\ref{sec:com_definition}.
In DC-PSE, the discrete differential operator is written as:
\begin{equation}
  L^{\bm{\alpha}}_{i} u
  \coloneq
  \frac{1}{h_i^{|\bm{\alpha}|}}
  \sum_{j\in\mathcal{N}_i}
  \left( u_{j}\pm  u_{i}\right){w}^{\bm{\alpha}}_{ji}v_j
  \approx
  \bm{C}^{\bm{\alpha}} \cdot \left.\bm{D} u\right\rvert_{i}\, ,
  \label{eq:dcpse_discrete_diff_operator_with_volume}
\end{equation}
with weights $ {w}^{\bm{\alpha}}_{ji} \coloneq w_{ji}^{\bm{\alpha}}\left( \left( \bm{x}_{j} - \bm{x}_{i}  \right)/ h_i \right) $, and $ v_j $ the {\em volume} of particle $ j $.
The sign between the function values is positive for odd $|\bm{\alpha}|$ and negative for even $|\bm{\alpha}|$ and renders the method conservative on symmetric point distributions~\cite{eldredge2002}. The particle volume $v_j$ appears as the integration element when deriving the above equation by quadrature from the continuous moment integrals of PSE~\cite{eldredge2002}. As later noted~\cite{UsingDcPseOpBouran2016}, however, this is identical to Eq.~\eqref{eq:general_discrete_diff_operator_zm} when absorbing the unknown particle volumes into the ABF coefficients.
DC-PSE was derived for general ABFs~\cite{schrader2010}, in the same spirit as anisotropic-kernel PSE~\cite{Hacki2015}. However, all works on DC-PSE so far~\cite{schrader2010, UsingDcPseOpBouran2016, AMeshfreeCollSingh2023, RecoveryByDisZwick2023, EntropicallyDaSingh2023,TaylorSeriesEJacque2020} used ABFs that are $h$-scaled monomials multiplied with Gaussian RBFs:
\begin{equation}
  \bm{W}_{ji}
  =
  \widetilde{\bm{X}}_{ji} \varphi_{ji}^{\mathrm{Gauss}} \, .
  \label{eq:dcpse_w_taylor_monomials_with_rbf} 
\end{equation}
This choice allows writing the mapping vector as the differential operator $ \partial^{\bm{\alpha}} $
applied to the monomials with the result evaluated at zero, $ (\partial^{\bm{\alpha}}\bm{X}_{ji})|_0 $.
The scaling factor $ 1/h_{i}^{|\bm{\alpha}|} $ can be absorbed into the mapping vector, which amounts to scaling $ \bm{C}^{\bm{\alpha}} $ by $ \mathbf{H}_i $ according to Eq.~\eqref{eq:h_vect_zm}.
Thus expressing the mapping vector as
\begin{equation}
  \widetilde{\bm{C}^{\bm{\alpha}}}
  \coloneq
  (-1)^{\lvert\bm{\alpha}\rvert}\mathbf{H}_{i}(\partial^{\bm{\alpha}}\bm{X}_{ji})|_0 \, .
  \label{eq:dcpse_mapping_vector}
\end{equation}
Therefore, DC-PSE includes the pre-conditioning from Eq.~\eqref{eq:h_vect_zm}, rendering it better conditioned than MLS~\cite{UsingDcPseOpBouran2016}.
Due to its symmetric formulation w.r.t.~the function values $u_{i,j}$, DC-PSE can either include or exclude the zeroth moment (see Section~\ref{sec:basic_approach_without_zeroth_moment}). This can be exploited to improve its stability in explicit time stepping~\cite{schrader2010,schrader2012}. The final form of the DC-PSE operator then is:
\begin{equation}
  L^{\bm{\alpha}}_{i} u
  \coloneq
  \displaystyle\sum_{j\in\mathcal{N}_i}
  \left( u_{j}\pm  u_{i}\right) \bm{W}_{ji} \cdot \left(\mathbf{M}_{i}^{-1}\widetilde{\bm{C}^{\bm{\alpha}}}\right).
  \label{eq:aom:dcpse_discrete_diff_operator_hMod_recap}
\end{equation}

A similar approach, albeit formulated slightly differently, was used in deriving the \textbf{Reproducing Kernel Collocation Method (RKCM)}. Originally introduced by Aluru~\cite{aluru2000}, RKCM was inspired by the Reproducing Kernel Particle Method (RKPM)~\cite{Liu1995rkpm, Chen1996rkpm}. It aimed to eliminate the quadrature from RKPM in order to turn it into a collocation approach.
Following its RKPM lineage, RKCM is derived from \textit{reproducing conditions}, which is the RKPM name for moment conditions. Just like DC-PSE replaced the continuous moment integrals of PSE with discrete sums, RKCM replaces the continuous reproducing conditions of RKPM with discrete variants as in Eq.~\eqref{eq:general_discrete_diff_operator_theoret_orderd}.
This starts from approximating the function itself using Eq.~\eqref{eq:general_discrete_diff_operator_zm} for $\bm{\alpha}=\bm{0}$:
\begin{equation}
  L^{\bm{0}}_{i} u \coloneq
  \displaystyle\sum_{j\in\mathcal{N}_i}
  u_{j}w^{\bm{0}}_{ji}
  \approx
  \bm{C}^{\bm{0}} \cdot \left.\bm{D} u\right\rvert_{i}\, .
  \label{eq:aom:rkcm:function_interpolation}
\end{equation}
and determining the resulting weight function $w_{ji}^{\bm{0}}$ from the reproducing conditions.
The weight function for a differential operator $|\bm{\alpha}|>0$ is then obtained as the derivative $\partial^{\bm{\alpha}}\!\left( w_{ji}^{\bm{0}} \right)$.
This approach is referred to as \textit{direct derivatives} in the literature~\cite{chen2017}.
Then:
\begin{equation}
  L^{\bm{\alpha}}_{i} u
  =
  \partial^{\bm{\alpha}}
  L^{\bm{0}}_{i} u
  =
  \partial^{\bm{\alpha}}\!\!
  \left(
    \displaystyle\sum_{j\in\mathcal{N}_i}
    u_{j}w^{\bm{0}}_{ji}
  \right)
  =
  \displaystyle\sum_{j\in\mathcal{N}_i}
  u_{j}\partial^{\bm{\alpha}}w^{\bm{0}}_{ji}
  =
  \displaystyle\sum_{j\in\mathcal{N}_i}  u_{j}
  \partial^{\bm{\alpha}}\!
  \left(
    \bm{W}_{ji} \cdot (\mathbf{M}_i^{-1} \bm{C}^{\bm{0}})
  \right)
  =
  \displaystyle\sum_{j\in\mathcal{N}_i}
  u_{j}w^{\bm{\alpha}}_{ji}\, .
  \label{eq:aom:rkcm_discrete_diff_operator_derived}
\end{equation}
RKCM, similarly to RKPM, uses ABFs in form of Taylor monomials multiplied by suitable compactly supported RBFs, hence $\bm{W}_{ji} = \bm{X}_{ji} \varphi_{ji}$.
This choice turns Eq.~\eqref{eq:aom:rkcm_discrete_diff_operator_derived} into the RKCM approximation
\begin{equation}
  L^{\bm{\alpha}}_{i} u
  \coloneq
  \sum_{j\in\mathcal{N}_i}  u_{j}
  \left(
    \partial^{\bm{\alpha}}\! \bm{X}_{ji} \varphi_{ji} \cdot \mathbf{M}_{i}^{-1} \bm{C}^{\bm{0}} +
    \bm{X}_{ji} \partial^{\bm{\alpha}} \varphi_{ji} \cdot \mathbf{M}_{i}^{-1} \bm{C}^{\bm{0}} +
    \bm{X}_{ji} \varphi_{ji} \cdot \partial^{\bm{\alpha}} \mathbf{M}_{i}^{-1} \bm{C}^{\bm{0}}
  \right)\, .
  \label{eq:aom:rkcm:discrete_diff_operator}
\end{equation}
The order of consistency is the same as in Eq.~\eqref{eq:general_discrete_diff_operator_with_coefs_do}, and Eq.~\eqref{eq:general_discrete_diff_operator_theoret_orderd} holds (see Ref.~\cite{chen2017} for details).

The \textbf{Gradient Reproducing Kernel Collocation Method (GRKCM)} introduced by Chi et al.~\cite{chi2013} replaces the direct derivatives approximation in RKCM by  Eq.~\eqref{eq:general_discrete_diff_operator_with_coefs_do} with the zeroth moment included.
Still, the ABFs are Taylor monomials multiplied by compact RBFs, $\bm{W}_{ji} =\bm{X}_{ji} \varphi_{ji}$.

The \textbf{Local Anistropic Basis Function Method (LABFM)} by King et al.~\cite{HighOrderDiffKing2020, HighOrderSimuKing2022} is one of the most recent moment-approximation methods.
Originally, the method was formulated without the pre-conditionig of Eq.~\eqref{eq:h_vect_zm}~\cite{HighOrderDiffKing2020}, but this was added two years later~\cite{HighOrderSimuKing2022}.
LABFM excludes the zeroth moment of the weight function and emphasizes the choice of ABFs $ \bm{W}_{ji} $.
Since the elements of $ \bm{X}_{ji} $ scale as $ h_i,\, h_i^2,\, h_i^3,\,\ldots $, the ABFs should scale as $ h_i^{-1},\, h_i^{-2},\, h_i^{-3},\,\ldots $.
Originally, the ABFs in LABFM were taken to be the partial derivatives of some suitably scaled RBF $ \varphi_{ji} \coloneq \varphi\left( \| \vect{x}_{ji} \|_2,\, h_i \right) $ evaluated at $ \vect{x}_j $, thus~\cite{HighOrderDiffKing2020}
\begin{equation}
  \bm{W}_{ji}
  =
  \left.\bm{D}\varphi_{ji}\right\rvert_{\vect{x}_{j}} \, .
  \label{eq:aom:labfm_w_abf_from_rbf}
\end{equation}
While this works in two dimensions ($n=2$), it is unable to ensure linear independence of the rows of the moment matrix $ \mathbf{M}_i $ for higher dimensions and approximation orders. For example, for $n=2$ and $m=4$, the moment matrix is always singular.
Therefore, the more recent formulation of LABFM~\cite{HighOrderSimuKing2022} constructs ABFs from orthogonal polynomials, specifically Hermite or Legendre polynomials, multiplied by a compact RBF as:
\begin{equation}
  \bm{W}_{ji}
  =
  \bm{P}_{ji} \varphi_{ji}\, .
  \label{eq:aom:labfm_w_abf}
\end{equation}
Several choices of $ \bm{P}_{ji} $ and $ \varphi_{ji} $ have been studied~\cite{HighOrderDiffKing2020, HighOrderSimuKing2022}.

A similar approach using ABFs generated by RBF derivatives as in Eq.~\eqref{eq:aom:labfm_w_abf_from_rbf} is used in \textbf{High-Order Consistent SPH (HOCSPH)} as introduced by Nasar et al.~\cite{nasar2021hoch}.
Having roots in SPH, the method differs from LABFM in that it attributes finite volumes to the particles. However, this corresponds to direct node integration, rendering HOCSPH effectively a collocation scheme.

An independent and interesting view on moment approximation methods has resulted from studying the relationship between peridynamics~\cite{madenci2013pd, oterkus2021} and meshfree collocation. This resulted in the formulation of \textbf{Generalized Reproducing Kernel Peridynamics (GRKP)} by Hillman et al.~\cite{hillman2020grkp}.
Peridynamics is a non-local formulation of continuum mechanics based on the interaction potentials between material points.
The gradients of these potential functions then correspond to the forces acting in the material. They are approximated by weighed summation over neighboring points. This is akin to how meshfree collocation methods approximate derivatives.
GRKP assumes a differential operator of the form
\begin{equation}
  \partial^{\bm{\alpha}} u_i =  \int_{\mathcal{N}_i} u(\bm{x}') g_m^{\bm{\alpha}}(\bm{x}' - \bm{x}_i) \,\mathrm{d}\bm{x}' \, ,
  \label{eq:aom:peridinamic_derivative}
\end{equation}
where the $ g_m^{\bm{\alpha}}(\bm{x}' - \bm{x}_i)$ are orthogonal functions that need to be constructed specifically for the approximation order $m$.
The usual ansatz for the $g_m^{\bm{\alpha}}$ is a product of Taylor monomials $\bm{X}(\bm{x})$ and a (normalized and compactly supported) scalar weight function $\psi(\bm{x})$,
thus $\bm{W}(\bm{x}) = \bm{X}(\bm{x})\psi(\bm{x})$. With the vector of unknown coefficients
${\Psi}_{\bm{x}}^{\bm{\alpha}} $,
\begin{equation}
  g_m^{\bm{\alpha}}(\bm{x}) = \bm{W}(\bm{x})\cdot\bm{\Psi}_{\bm{x}}^{\bm{\alpha}}\, .
\end{equation}
This resembles the general ABF in meshfree collocation methods as given in  Eq.~\eqref{eq:w_abf_zm}.
The connection with meshfree collocation becomes obvious when approximating the integral in Eq.~\eqref{eq:aom:peridinamic_derivative} by direct nodal quadrature with an elementary volume associated with each point. Since the integration volumes (quadrature weights) can be absorbed into the weight functions, as is also done in DC-PSE, Eq.~\eqref{eq:aom:peridinamic_derivative} becomes identical to the definition of the meshfree collocation operator in Eq.~\eqref{eq:general_discrete_diff_operator_zm}, replacing $ \partial^{\bm{\alpha}} u_i$ with $ L^{\bm{\alpha}}_{i}u $.
Also, similarly to DC-PSE, the zeroth moment is optional, and both formulations are found in the literature~\cite{hillman2020grkp}.
This results in an approximation akin to that of DC-PSE in Eq.~\eqref{eq:aom:dcpse_discrete_diff_operator_hMod_recap}, albeit without the pre-conditioning:
\begin{equation}
  L^{\bm{\alpha}}_{i} u
  =
  \displaystyle\sum_{j\in\mathcal{N}_i}
  u_{ji}\bm{W}_{ji}\cdot \left(\mathbf{M}_{i}^{-1} \bm{C}^{\bm{\alpha}} \right).
  \label{eq:grkp_ld_approximation}
\end{equation}
The ABFs are usually chosen as Taylor monomials multiplied by a suitable RBF, $ \bm{W}_{ji} = \bm{X}_{ji} \varphi_{ji}$.

Taken together, the different moment approximation methods mainly differ in the form of the ABFs used and in whether the zeroth moment of the weights is included in the derivation.
Each method's derivation also brought its own mathematical contributions. For example, the DC-PSE paper~\cite{schrader2010} proved formal connections with (corrected) SPH, RKPM and, for the limit of particles lying on a Cartesian grid, compact finite differences. The LABFM papers~\cite{HighOrderDiffKing2020, HighOrderSimuKing2022} provided a detailed analysis of different ABF choices. RKCM introduced the concept of direct derivatives to strong-form collocation methods~\cite{aluru2000}, and GRKP provided a connection with peridynamics~\cite{hillman2020grkp}, showing a direct link between meshfree collocation methods and physical systems.

\section{Derivation by $\ell_2$ minimization of the discretization error}\label{sec:l2e}
An alternative way of obtaining the weights of the discrete operator in Eq.~\eqref{eq:general_discrete_diff_operator_zm} is to minimize the truncation error $ e_{ji}^{m+1} $ of the Taylor expansion in Eq.~\eqref{eq:taylor_series_zm} over $\mathcal{N}_i$ for the desired derivatives $ \left.\bm{D}u\right\rvert_{i} $.
In general, the error can be weighted by a scalar function $\varphi_{ji} $.
Therefore, we solve the following optimization problem:
\begin{align*}
  \mathcal{J}_{\mathcal{E},i}
  \coloneq
  \displaystyle\sum_{j\in\mathcal{N}_i} \varphi_{ji} \left( e_{ji}^{m+1} \right)^2
  =
  \displaystyle\sum_{j\in\mathcal{N}_i} \varphi_{ji} \left( \bm{X}_{ji} \cdot \left.\bm{D}u\right\rvert_{i} - u_{j} \right)^2 \\
  \mathrm{min} \left( \mathcal{J}_{\mathcal{E},i} \right)
  \quad
  \text{with respect to}
  \left.\bm{D}u\right\rvert_{i}.
  \numberthis
  \label{eq:l2e:emin-optimizationformulation}
\end{align*}
The formal solution to this problem as derived in Appendix~\ref{sec:appendix-l2e} is:
\begin{align*}
  \left.\bm{D}u\right\rvert_{i}
  &=
  \displaystyle\sum_{j\in\mathcal{N}_i} \varphi_{ji}  u_{j}
  \left(
    \displaystyle\sum_{k} \varphi_{ki} \bm{X}_{ki}  \bm{X}_{ki}^{\top}
  \right)^{-1} \!\!\!\bm{X}_{ji}
  \\
  &=
  \displaystyle\sum_{j\in\mathcal{N}_i} \varphi_{ji}  u_{j}
  \mathbf{M}_i^{-1}
  \bm{X}_{ji}\, .
  \numberthis
  \label{eq:l2e:emin-generaldiscreteoperator}
\end{align*}
The moment matrix $\mathbf{M}$ naturally appears in the solution of the optimization problem and does not depend on any choice of basis functions. Nevertheless, it is the same moment matrix as obtained in Eq.~\eqref{eq:momentum_matrix_matrix_form} for AOM methods:
\begin{equation}
  \mathbf{M}_i
  \coloneq
  \displaystyle\sum_{j\in\mathcal{N}_i} \varphi_{ji} \bm{X}_{ji}  \bm{X}_{ji}^{\top}
  =
  \displaystyle\sum_{j\in\mathcal{N}_i}\bm{X}_{ji}\bm{X}_{ji}^{\top} \varphi_{ji}
  =
  \displaystyle\sum_{j\in\mathcal{N}_i}\bm{X}_{ji}\bm{W}_{ji}^{\top}\, ,
  \label{eq:l2e:moment_matrix}
\end{equation}
with $ \bm{W}_{ji} \coloneq \bm{X}_{ji} \varphi_{ji}$. This reveals that the error weights $\varphi$ play the same role as the ABFs in AOM methods. In the present case, however, one is limited to radially symmetric functions, since the truncation error is a scalar.

\nomvar{ $\mathcal{J}_{\mathcal{E},i}$ }{functional for $\ell_2$ minimization of approximation error [-]}

The solution to the linear system of equations in Eq.~\eqref{eq:l2e:emin-generaldiscreteoperator} contains all derivatives up to and including the $m$-th derivative.
To obtain an approximation for a specific differential operator $\partial ^{\bm{\alpha}}$ in the form of Eq.~\eqref{eq:general_discrete_diff_operator_with_coefs_do}, we again use the mapping vector $\bm{C}^{\bm{\alpha}}$ to form the desired differential operator from $ \left.\bm{D}u\right\rvert_{i} $ as
\begin{equation}
  \bm{C}^{\bm{\alpha}} \cdot \left.\bm{D}u\right\rvert_{i}
  =
  \displaystyle\sum_{j\in\mathcal{N}_i} \varphi_{ji}  u_{j}
  \bm{C}^{\bm{\alpha}} \cdot
  \left(
    \mathbf{M}_i^{-1}
    \bm{X}_{ji}
  \right)\, .
  \label{eq:l2e:approx-error-l2-min-include-mapping-vect}
\end{equation}
Since in the present case the moment matrix $ \mathbf{M}_i $ (and hence its inverse $\mathbf{M}_i^{-1} $) is symmetric\footnote{$\mathbf{M} = \sum \bm{X}\bm{X}^{\top} \phi$, where $\bm{X}\bm{X}^{\top}$ is the outer product of two identical vectors, resulting in a symmetric matrix with each term scaled by $\phi$.}, this yields the final approximation:
\begin{equation}
  L^{\bm{\alpha}}_{i}u
  =
  \displaystyle\sum_{j\in\mathcal{N}_i} u_{j}
  \bm{W}_{ji} \cdot
  \left(
    \mathbf{M}_i^{-1}
    \bm{C}^{\bm{\alpha}}
  \right)\, .
  \label{eq:l2e:approx-error-l2-min_ld_result}
\end{equation}
Formally, this is the same discrete operator as obtained in Eq.~\eqref{eq:general_discrete_diff_operator_with_coefs_do} from approximation of moments. There are two important differences, though. First, the basis functions $\bm{W}_{ji} $ in Eq.~\eqref{eq:l2e:approx-error-l2-min_ld_result} are given by the Taylor monomials and the error weighting $ \bm{W}_{ji} = \bm{X}_{ji} \varphi_{ji}$. They cannot be chosen freely. Second, the moment matrix $\mathbf{M}_i$ in Eq.~\eqref{eq:l2e:approx-error-l2-min_ld_result} is symmetric, whereas it need not be symmetric in the AOM formulation (Section~\ref{sec:aom}).

The same derivation by error minimization can also be performed without the zeroth moment, following the procedure outlined in  Section~\ref{sec:basic_approach_without_zeroth_moment}.
In that case, the error $ e_{ji}^{m+1} $ in Eq.~\eqref{eq:l2e:emin-optimizationformulation} takes the form $ e_{ji}^{m+1} = \bm{X}_{ji} \cdot \left.\bm{D}u\right\rvert_{i} - u_{ji} $ with $ u_{ji} = u_{j} - u_{i} $. This replaces $ u_j $ with $ u_{ji} $ in Eq.~\eqref{eq:l2e:approx-error-l2-min_ld_result}, yielding the same result as in Eq.~\eqref{eq:general_discrete_diff_operator_with_coefs_withouth_zeroth_moment}.

\subsection{Matrix notation}

Following the notation introduced in Section \ref{sec:basic-approach-matrix-formulation}, the $\ell_2$-error minimization approach can also be written in matrix notation.
For this, we additionally define the diagonal matrix $\mathbf{V}_{i}$ of the error weights $\varphi_{ij}$ for individual particle pairs $(i,j)$:
\begin{equation}
  \mathbf{V}_{i} =
  \mathrm{diag}\left(
    \varphi_{ji}
  \right)_{j=1}^{N_i} \in \mathbb{R}^{N_i\times N_i} \, .
  \label{eq:l2e:appendix_matrixform_wmatrix}
\end{equation}

\nomvar{ $\mathbb{V}_{i}$ }{diagonal matrix with diagonal entries  $ W_{ji}^{0} $ for all neighbors of particle $i$}

With this matrix, the minimization problem in Eq.~\eqref{eq:l2e:emin-optimizationformulation} can be written as:
\begin{align*}
  \mathcal{J}_{\mathcal{E},i}
  =
  \displaystyle\sum_{j\in\mathcal{N}_i} \varphi_{ji} \left(e_{ji}^{m+1} \right)^2
  =
  \bm{E}_{i}^{\top}\mathbf{V}_{i} \bm{E}_{i}
  =
  \left( \mathbf{X}_i \left.\bm{D} u\right\rvert_{i} - \bm{\bm{U}}_{i} \right)^{\top}
  \mathbf{V}_{i}
  \left( \mathbf{X}_i \left.\bm{D} u\right\rvert_{i} - \bm{\bm{U}}_{i} \right)
  \\
  \mathrm{min} \left( \mathcal{J}_{\mathcal{E},i} \right)
  \quad
  \text{with respect to}
  \left.\bm{D} u\right\rvert_{i}.
  \numberthis
  \label{eq:l2e:emin-optimizationformulation-matrix}
\end{align*}
Formally solving this minimization problem (see Appendix~\ref{sec:AppB1}), we obtain:
\begin{equation}
  \left.\bm{D} u\right\rvert_{i}
  =
  \left[
    \left( \mathbf{X}_i^{\top} \mathbf{V}_{i} \mathbf{X}_i \right)^{-1} \mathbf{X}_i^{\top} \mathbf{V}_{i}
  \right]
  \bm{\bm{U}}_i \, .
\end{equation}
With the definition in Eq.~\eqref{eq:appendix_wmatrix_ymatrix} (here for $ \mathbf{W}_i = \mathbf{V}_i \mathbf{X}_i $) and the definition of the moment matrix from Eq.~\eqref{eq:momentum_matrix_matrix_form}, this yields the final operator:
\begin{equation}
  L^{\bm{\alpha}}_{i} u
  =
  \bm{\bm{U}}_i \cdot \bm{w}_{\bullet i}^{\bm{\alpha}}
  =
  \bm{\bm{U}}_i \cdot \left( \mathbf{C}_i\bm{C}^{\bm{\alpha}} \right)
  =
  \bm{\bm{U}}_i \cdot \left( \mathbf{W}_i \mathbf{Y}_{i}\bm{C}^{\bm{\alpha}} \right)
  =
  \bm{\bm{U}}_i \cdot \left( \mathbf{W}_i \mathbf{M}_{i}^{-1}\bm{C}^{\bm{\alpha}} \right),
  \label{eq:l2e:general_discrete_diff_operator_with_coefs_do_matrix_form_l2min}
\end{equation}
which is identical to Eq.~\eqref{eq:l2e:approx-error-l2-min_ld_result}.

\subsection{Methods derived by $\ell_2$ error minimization}
\label{sec:l2e methods}

The majority of meshfree collocation methods were derived by truncation error minimization. This notably includes the large class of methods based on moving least squares (MLS) minimization, and methods based on generalizing finite differences to scattered nodes. But it also encompasses more specialized methods, e.g., for advection problems or gradient and Laplacian approximation.

\subsubsection{Moving least squares methods}\label{sec:mlsmethods}

The largest and most popular group of methods in this category are the \textbf{Moving Least Squares (MLS)} methods. These start from approximating the function itself using Eq.~\eqref{eq:general_discrete_diff_operator_zm} for $\bm{\alpha}=\bm{0}$, as originally introduced by Lancester and Salkauskas~\cite{lancaster1981}:
\begin{equation}
  L_{i}^{\bm{0}} u
  =
  \bm{P}_i
  \cdot
  \displaystyle\sum_{j\in\mathcal{N}_i}
  \varphi_{ji}  u_j
  \mathbf{M}_{i}^{-1}
  \bm{P}_{j}\, .
  \label{eq:l2e:mls:function_approximation}
\end{equation}
Using relative coordinates with respect to $\vect{x}_i$ leads to a formulation with better condition number:
\begin{equation}
  L_{i}^{\bm{0}} u
  =
  \bm{P}(\bm{0})
  \cdot
  \displaystyle\sum_{j\in\mathcal{N}_i}
  \varphi_{ji}  u_j
  \mathbf{M}_{i}^{-1}
  \bm{P}_{ji}\, .
  \label{eq:fpm_approximation_of_function_itself}
\end{equation}
The differential operator $\partial^{\bm{\alpha}}$ is then approximated by taking the respective derivative of the basis,
\begin{equation}
  L^{\bm{\alpha}}_{i} u
  =
  \partial^{\bm{\alpha}} \bm{P}(\mathbf{0})
  \cdot
  \displaystyle\sum_{j\in\mathcal{N}_i}
  \varphi_{ji}  u_j
  \mathbf{M}_{i}^{-1}
  \bm{P}_{ji}
  =
  \displaystyle\sum_{j\in\mathcal{N}_i}
  u_j
  \bm{W}_{ji}
  \cdot
  \left(
    \mathbf{M}_{i}^{-1}
    \bm{C}^{\bm{\alpha}}
  \right)
  =
  \displaystyle\sum_{j\in\mathcal{N}_i}
  u_j
  w_{ji}^{\bm{\alpha}} \, ,
  \label{eq:fpm_derivative_of_basis}
\end{equation}
where we used $ \bm{C}^{\bm{\alpha}} = \partial^{\bm{\alpha}} \bm{P}(\bm{0}) $, the symmetry of $\mathbf{M}_{i}$, and $\bm{W}_{ji} = \bm{P}_{ji}\varphi_{ji}$.
The approximation in Eq.~\eqref{eq:fpm_derivative_of_basis} is sometimes called \textit{diffusive derivatives}~\cite{nayroles1992dem,Kim2003,Huerta2004,Yoon2006,chen2017}, since it uses the derivative of the basis instead of directly computing the derivative of the function $u$ itself (via the derivative of the $w_{ji}$ in the discrete approximation), as is for example done for RKCM in Eq.~\eqref{eq:aom:rkcm:discrete_diff_operator}.

The same derivative approximation is also used in the \textbf{Finite Point Method (FPM)} introduced by O\~{n}ate et al.~\cite{onate1995b} to describe a class of point-data interpolation schemes that includes Least Squares (LS), Weighted Least Squares (WLS), MLS, and RKPM.
The same authors later extended the MLS formalism by introducing \textbf{Multiple Fixed Least Squares (MFLS)}~\cite{onate1996fpm}, allowing the
$\varphi_{\bm{x}_i}$ to have different shapes at different particles $\bm{x}_i$.
When evaluating $\varphi(\bm{x})$ at an off-particle location, it is interpolated from the $\varphi_{\bm{x}_i}(\bm{x})$ at the neighboring particles.

To our knowledge, every available variant of MLS or FPM used Taylor monomials in the ABFs, i.e., $ \bm{P}_{ji} = \bm{X}_{ji} $ and thus $ \bm{W}_{ji} = \bm{X}_{ji} \varphi_{ji} $~\cite{onate1996finite, onate1996fpm, onate1998fpm, onate2000finite, onate2001fpm}.
While FPM is usually presented in matrix notation, MLS is more commonly found in vector form.

A variant of the classic MLS approach is the \textbf{Least Squares Meshfree Method (LSMFM)} introduced by Youn and co-workers~\cite{park2001, zhang2005}. It differs from classic MLS in that it replaces the diffuse derivative with a direct derivative, similar to RKCM.
Starting from Eq.~\eqref{eq:l2e:mls:function_approximation},
taking $\partial^{\bm{\alpha}}(w_{ji}^{\bm{0}})$ results in Eq.~\eqref{eq:aom:rkcm:discrete_diff_operator}.
While this is formally the same approximation as in RKCM, LSMFM derives it via approximation-error minimization instead of moment approximation. The result, however, is identical, and LSMFM also uses Taylor monomials as basis functions.

Levin~\cite{levin1998} generalized MLS approaches to arbitrary linear bounded operators in the \textbf{Generalized Moving Least Squares (GMLS)} method, which was further elaborated in several follow-up works~\cite{wendland2004, mirzaei2012, gross2020}. It is based on a thorough mathematical formalism. Assume the function $u(\vect{x})$ to be $m$ times continuously differentiable, $  u \in C^m(\Omega) $ for some $m \ge 0$\footnote{Here, $m=0$ is allowed because the method approximates the sampled values at the particles and then computes derivatives of the approximation. No derivatives of $u$ are computed.}.
In GMLS, the function $ u $ is approximated by a set of \textit{sampling functionals} $ \Lambda( u) = \left\{ \lambda_j( u)\right\}_{j=0}^{N} \in C^m(\Omega)^* $ where $ C^m(\Omega)^* $ is the dual space to $ C^m(\Omega) $.
In the simplest case, the sampling functionals are point-wise evaluations of $  u $, i.e., $ \lambda_i ( u) = \delta_{\vect{x}_i}( u ) =  u( \vect{x}_i) =  u_i $.
For a given $\lambda(u) \in C^m(\Omega)^*$, GMLS determines the approximation $\lambda_i(u)^h \approx \lambda(u)[\vect{x}_i]$ at resolution $h$ as a linear combination of the $ \left\{ \lambda_j( u)\right\}_{j\in\mathcal{N}_i}$:
\begin{equation}
  \lambda_i(u)^h \coloneq
  \sum_{j\in\mathcal{N}_i}
  \lambda_j( u)w_{j}(\lambda)\, ,
  \label{eq:l2e:gmls:general_approx_operator_gmls}
\end{equation}
which is analogous to $L_i^{\bm{0}}u$.
Following Eq.~\eqref{eq:l2e:gmls:general_approx_operator_gmls}, GMLS approximates any linear bounded operator $ \tau_i $ at $ \vect{x}_i $ as
\begin{equation}
  \tau_i^h u
  =
  \displaystyle\sum_{j\in\mathcal{N}_i}
  \lambda_{j}( u) \widetilde{w}_{j}(\lambda)\, .
  \label{eq:l2e:gmls:tau_approx_operator_gmls}
\end{equation}
In the GMLS literature, Eq.~\eqref{eq:l2e:gmls:tau_approx_operator_gmls} is mostly written in its functional matrix form as
\begin{equation}
  \lambda(u)^h = \lambda(u^*)
  =
  \bm{P}(\vect{x})^{\!\top} \widetilde{\bm{\Psi}}(u^*)\, ,
\end{equation}
where $\bm{P}(\bm{x})$ is the (stacked) vector of basis functions, $ \widetilde{\bm{\Psi}}(u^*) $ the (stacked) vector of coefficients, and $ u^* $ the solution to the local $\ell_2$-minimization problem
\begin{equation}
  u^*
  =
  \argmin_{q \in \mathcal{V}_h}
  \displaystyle\sum_{j\in\mathcal{N}_i}
  \left( \lambda_j( u) - \lambda_j(q) \right)^2 \varphi_{ji}\, ,
  \label{eq:l2e:cmls:gmls_l2_optim_problem}
\end{equation}
over the space $\mathcal{V}_h = \text{span}\{ P_1,\, \ldots ,\, P_p \} $ spanned by the basis $ P_i \in \bm{P}$ of dimension $ \dim({\mathcal{V}_h}) = p $ ($p$ depends on $m$ and $n$ as discussed after Eq.~\eqref{eq:general_discrete_diff_operator_zm}).
This optimization problem is solved as in Appendix~\ref{sec:appendix-l2e}, and the approximation $ \tau_i^h $ to the operator $ \tau_i $ is obtained by applying the operator to the basis functions while keeping the coefficients $ \widetilde{\bm{\Psi}}(u^*) $ constant.
In meshfree collocation methods, the operators of interest are differential operators at the locations of the collocation points $\vect{x}_i$, thus $\tau_i u = \left.\partial^{\bm{\alpha}} u\right\rvert_{i} = \bm{C}^{\bm{\alpha}} \cdot \left.\bm{D}  u\right\rvert_{i} $.
The GMLS procedure then yields the diffuse derivatives approximation~\cite{mirzaei2012}
\begin{equation}
  \tau_i^h u
  = (\tau_i \bm{P})^{\!\top} \widetilde{\bm{\Psi}}(u^*)
  \label{eq:gmls_diffuse_derivative}
  = \displaystyle\sum_{j\in\mathcal{N}_i}
  \lambda_{j}( u^*) \widetilde{w}_{j}(\lambda)\, ,
\end{equation}
which amounts to a functional form of Eq.~\eqref{eq:fpm_derivative_of_basis} as:
\begin{equation}
  L_{i}^{\bm{\alpha}}( u)
  =
  \partial^{\bm{\alpha}}\bm{P}(\mathbf{0})
  \cdot
  \displaystyle\sum_{j\in\mathcal{N}_i}
  \varphi_{ji} \lambda_j(u)
  \mathbf{M}_{i}^{-1}
  \lambda_j(\bm{P}_{i})\, ,
  \label{eq:gmls_ld_approximation_raw}
\end{equation}
with symmetric moment matrix
\begin{equation}
  \mathbf{M}_i
  =
  \displaystyle\sum_{j\in\mathcal{N}_i}
  \lambda_j(\bm{P}_i) \lambda_j(\bm{P}_i)^{\top} \varphi_{ji}\, .
  \label{eq:l2e:gmls:moment_matrix}
\end{equation}
For the simplest sampling functional $ \lambda_i = \delta_{\vect{x}_i} $, this simplifies to $ \lambda_j(\bm{P}_{i}) = \bm{P}(\vect{x}_j - \vect{x}_i) =\bm{P}_{ji} $ and $ \lambda_j( u) =  u_j $.
The mapping vector appears as $ \bm{C}^{\bm{\alpha}} = \partial^{\bm{\alpha}}\bm{P}(\mathbf{0}) $ with $ \bm{P}(\mathbf{0}) = \lambda_i(\bm{P}_i)  $.
Popular choice for the basis $ \bm{P}_{ji} $ are scaled Taylor monomials or Bernstein polynomials~\cite{compadre_toolkit}.

Trask et al.~\cite{trask2016} connected GMLS with compact finite differences in the \textbf{Compact Moving Least Squares (CMLS)} method. The goal was to derive more compactly supported meshfree stencils, taking inspiration from Compact Finite Difference Methods (CFDM)~\cite{lele1992cfdm, weinan1996cfdm}.
For CFDM, E proposed using the expression of the truncation error for the very equation to be solved to improve the approximation accuracy without increasing the stencil size~\cite{weinan1996cfdm}.
Another class of CFDM schemes, proposed by Lele~\cite{lele1992cfdm}, uses an implicit formulation in which the derivatives at grid points are not computed directly, but as part of a system of equations involving neighboring points.
Combining these two ideas, CMLS achieves compactly supported meshfree sencils by including an additional problem-dependent penalty into the optimization problem from Eq.~\eqref{eq:l2e:cmls:gmls_l2_optim_problem} while keeping the moment matrix invertible.
This has been reported to halve the support $|\mathcal{N}_i|$ of the MLS operators~\cite{trask2016}.
To illustrate the CMLS approach, consider the linear boundary value problem
\begin{align*}
  \mathcal{L}_{\Omega}  u &= f, \qquad \forall \vect{x} \in \Omega\, , \\
  \mathcal{L}_{\partial\Omega}  u &= g, \qquad \forall \vect{x} \in \partial\Omega\, .
  \numberthis \label{eq:cmls:example_problem}
\end{align*}
In CMLS, Eq.~\eqref{eq:l2e:cmls:gmls_l2_optim_problem} includes penalization terms reflecting the problem to be solved
\begin{equation}
  u^*
  =
  \argmin_{q \in \mathcal{V}_h}
  \sum_{j\in\mathcal{N}_i}
  \left[
    \left(
      \lambda_j( u) - \lambda_j(q)
    \right)^2
    +
    \varepsilon_{\Omega}
    \left(
      f_j - \mathcal{L}_{\Omega} q
    \right)^2
    +
    \chi_{j,\partial\Omega}
    \varepsilon_{\partial\Omega}
    \left(
      g_j - \mathcal{L}_{\partial\Omega} q
    \right)^2
  \right] \varphi_{ji}\, ,
  \label{eq:l2e:cmls:l2_optim_problem}
\end{equation}
where $ \varepsilon_{\Omega} $ and $ \varepsilon_{\partial\Omega} $ are interior and boundary regularization parameters, respectively, and $ \chi_{j,\partial\Omega} $ is the boundary indicator function $ \chi_{j,\partial\Omega} = 1 $ if $ \vect{x}_j \in \partial\Omega $ and $ 0 $ otherwise.
The first term is the standard GMLS term responsible for the approximation quality.
The second and third terms are the solution and boundary residuals as in CFDM.
The optimization problem in Eq.~\eqref{eq:l2e:cmls:l2_optim_problem} is solved as shown in Appendix~\ref{sec:appendix-l2e}, resulting in
\begin{equation}
  L^{\bm{\alpha}}_{i} u
  =
  \displaystyle\sum_{j\in\mathcal{N}_i}
  \left(
    u_{j} \bm{W}_{ji} +
    \epsilon_{\Omega} \mathcal{L}_{\Omega} \bm{P}_{ji} f_i \varphi_{ji} +
    \chi_{j,\partial\Omega} \epsilon_{\partial\Omega} \mathcal{L}_{\partial\Omega} \bm{P}_{ji} g_i \varphi_{ji}
  \right) \cdot (\mathbf{M}_i^{-1} \bm{C}^{\bm{\alpha}})\, ,
  \label{eq:l2e:cmls:discrete_diff_operator_no_enforced_bc}
\end{equation}
\begin{equation}
  \mathbf{M}_i =
  \displaystyle\sum_{j\in\mathcal{N}_i}
  \left(
    \bm{P}_{ji} \bm{W}_{ji}^{\top} +
    \epsilon_{\Omega} (\mathcal{L}_{\Omega} \bm{P}_{ji})(\mathcal{L}_{\Omega} \bm{P}_{ji}^{\top}) \varphi_{ji} +
    \chi_{j,\partial\Omega} \epsilon_{\partial\Omega} (\mathcal{L}_{\partial\Omega} \bm{P}_{ji})(\mathcal{L}_{\partial\Omega} \bm{P}_{ji}^{\top}) \varphi_{ji}
  \right).
\end{equation}
Both the resulting operator approximation and the moment matrix include contributions from the additional penalty terms.
The regularization parameters are used to weight the terms. Since the penalties act as soft constraints, however, the boundary condition may not be exactly fulfilled. Therefore, for points located on the boundary, an additional hard constraint
\begin{equation}
  \mathcal{L}_{\partial\Omega}q_i = g_i\, , \qquad \forall i:\; \bm{x}_i \in \partial\Omega\, ,
\end{equation}
is used to exactly enforce the boundary condition~\cite{trask2016}.
This ensures that the polynomial approximation constructed at the interior points, with the boundary condition imposed only as a soft constraint, fulfills the boundary condition exactly when evaluated at a boundary point.
This leads to an additional boundary constraint on the matrix $\mathbf{M}_i$.
Following the notation introduced in Section~\ref{sec:l2e} (in particular Eqs.~\eqref{eq:l2e:emin-generaldiscreteoperator} and \eqref{eq:l2e:approx-error-l2-min-include-mapping-vect}), we write Eq.~\eqref{eq:l2e:cmls:discrete_diff_operator_no_enforced_bc} as
\begin{equation}
  L_{i}^{\bm{\alpha}}u \approx \bm{C}^{\bm{\alpha}} \cdot \left.\bm{D}u\right\rvert_{i} \, ,
\end{equation}
where the derivatives $ \left.\bm{D}u\right\rvert_{i} $ are obtained by solving the linear system of equations $ \mathbf{M}_i \left.\bm{D}u\right\rvert_{i} = \bm{R}_{i} $ with right-hand side
\begin{equation}
  \bm{R}_{i} \coloneq
  \displaystyle\sum_{j\in\mathcal{N}_i}
  \left(
    u_{j} \bm{W}_{ji}^{\top} +
    \epsilon_{\Omega} \mathcal{L}_{\Omega} \bm{P}_{ji}^{\top} f_i \varphi_{ji} +
    \chi_{j,\partial\Omega} \epsilon_{\partial\Omega} \mathcal{L}_{\partial\Omega} \bm{P}_j^{\top} g_i \varphi_{ji}
  \right) .
\end{equation}
Including the boundary condition constraint, the problem takes the form
\begin{equation*}
  \widetilde{\mathbf{M}}_i \coloneq
  \begin{bmatrix}
    \mathbf{M}_i & \mathcal{L}_{\partial\Omega}\bm{P}_{ij} \\
    \mathcal{L}_{\partial\Omega}\bm{P}_{ij}^{\top} & 0
  \end{bmatrix},
  \quad
  \widetilde{\left.\bm{D}u\right\rvert}_{i} \coloneq
  \left[
    \left.\bm{D}u\right\rvert_{i},\, \lambda
  \right],
  \quad
  \widetilde{\bm{R}}_{i} \coloneq
  \left[
    \bm{R}_{i},\, g_i
  \right].
\end{equation*}
Solving the linear system of equations $ \widetilde{\mathbf{M}}_i \widetilde{\left.\bm{D}u\right\rvert}_{i} = \widetilde{\bm{R}}_{i} $, we obtain the approximation
\begin{equation}
  L_{i}^{\bm{\alpha}}u \approx \bm{C}^{\bm{\alpha}} \cdot \widetilde{\left.\bm{D}u\right\rvert}_{i}\, ,
\end{equation}
for which $ \mathcal{L}_{\Omega}^{h} u_i = L^{\bm{\alpha}}_{i} u = f_i $ for all $\{i:  \bm{x}_i \in \Omega\} $ and $ L_{i}^{\mathcal{L}_{\partial\Omega}}u = g_i $ for all $\{i:  \bm{x}_i \in  \partial \Omega\} $.
We can therefore write the resulting CMLS scheme as:
\begin{equation}
  L^{\bm{\alpha}}_{i} u
  =
  \displaystyle\sum_{j\in\mathcal{N}_i}
  \left( u_{j} \bm{W}_{ji} + \bm{Q}_{ji} \right) \cdot (\mathbf{M}_i^{-1} \bm{C}^{\bm{\alpha}})\, ,
  \quad
  \mathbf{M}_{i}
  =
  \displaystyle\sum_{j\in\mathcal{N}_i}\bm{P}_{ji}\bm{W}_{ji}^{\top} + \mathbf{Q}_i \, ,
  \label{eq:l2e:cmls:discrete_diff_operator}
\end{equation}
where $\bm{Q}_{ji} \in \mathbb{R}^{p} $
and $ \mathbf{Q}_i \in \mathbb{R}^{p \times p}$ represent the additional terms originating from the penalties.
Taylor monomials or Legendre polynomials are used as the basis functions \cite{trask2016}. For particles on a regular Cartesian grid, the CMLS approximation becomes identical to Lele CFDM~\cite{trask2016}.

Relaxing the conditions on the number and distribution of neighbors is the aim of \textbf{Modified Moving Least Squares (MMLS)} as introduced by Joldes et al.~\cite{joldes2015mmls, joldes2019mmls}.
Like CMLS, MMLS also includes additional penalties but uses a process akin to  \textit{Tikhonov--Miller} regularization to do so.
Unlike standard \textit{Tikhonov--Miller} regularization, though, MMLS can assign individual positive weights for the additional constraints.
This allows for higher-order approximations with less restrictive conditions on the size of the stencil.
MMLS also regularizes the moment matrix in cases where it becomes singular in classic MLS. For this, it uses sparsity to select the solution in which the most coefficients are zero.
This decreases the degree of the polynomial interpolation while preserving well-posedness of the problem.
In contrast to CMLS, the regularization in MMLS is problem-independent, as it directly acts on the coefficients of the polynomial basis functions.
This modifies Eq.~\eqref{eq:l2e:cmls:gmls_l2_optim_problem} to:
\begin{equation}
  u^*
  =
  \argmin_{q \in \mathcal{V}_h}
  \displaystyle\sum_{j\in\mathcal{N}_i}
  \left[
    \left(
      \lambda_j( u) - \lambda_j(q)
    \right)^2
    \varphi_{ji}
    +
    \bm{\mu}\cdot
    (
      \hat{\bm{\Psi}}_{i}
    )^{2}
  \right]\, ,
  \label{eq:l2e:mmls:l2_optim_problem}
\end{equation}
where $ (\mu_i )_{i=1}^s \coloneq \bm{\mu}\in\mathbb{R}^s $ is the vector of positive regularization weights for the subset of $s$ coefficients $ \hat{\bm{\Psi}}_{i} $ of the leading-order terms\footnote{
  We illustrate this in an example. Using a second-order basis ($ m = 2 $) in $n=2$ Cartesian dimensions, we have
  \begin{equation*}
    u^*
    =
    \argmin_{q \in \mathcal{V}_h}
    \displaystyle\sum_{j\in\mathcal{N}_i}
    \left(
      \lambda_j( u) - \lambda_j(q)
    \right)^2
    \varphi_{ji}
    +
    \mu_{x^2}
    \left(
      \Psi_{i,x^2}
    \right)^2
    +
    \mu_{xy}
    \left(
      \Psi_{i,xy}
    \right)^2
    +
    \mu_{y^2}
    \left(
      \Psi_{i,y^2}
    \right)^2
    \label{eq:mmls_l2_optim_problem_example_with_m2}
  \end{equation*}
  for the leading-order (here second-order) terms $x^2$, $xy$, and $y^2$.
}.
Using small positive $ 0 < \mu_i \ll 1 $ ensures that if the matrix $\mathbf{M}_i$ is close to singular, the coefficients for the leading-order terms become zero. Small regularization weights also ensure that the solution remains close to the standard MLS solution when the moment matrix is not singular. However, the important question of how to choose the regularization weights $\bm{\mu}$ is rarely discussed in the literature. In the original publication~\cite{joldes2019mmls}, the authors use a constant $ \mu = 10^{-7}$.
Solving then the optimization problem in Eq.~\eqref{eq:l2e:mmls:l2_optim_problem} as in Appendix~\ref{sec:appendix-l2e}, the MMLS moment matrix becomes:
\begin{equation}
  \widetilde{\mathbf{M}}_i
  =
  \displaystyle\sum_{j\in\mathcal{N}_i}
  \bm{X}_{ji} \bm{W}_{ji}^{\top} + \mathbf{H}_{\mu}
  =
  \mathbf{M}_i + \mathbf{H}_{\mu}\, ,
  \label{eq:l2e:mmls:momentum_matrix}
\end{equation}
with $ \mathbf{H}_{\mu} = \mathrm{diag}( \bm{0}_{p-s}, \bm{\mu} ) \in \mathbb{R}^{p\times p} $ and $\bm{0}_{p-s} \in \mathbb{R}^{p-s} $ a zero vector.
The resulting operator approximation is formally identical to Eq.~\eqref{eq:l2e:approx-error-l2-min_ld_result} but with the above moment matrix.
The MMLS moment matrix $\widetilde{\mathbf{M}}_i$ generally has the same dimensions
as the GMLS moment matrix $\mathbf{M}_i$ from Eq.~\eqref{eq:l2e:gmls:moment_matrix} but requires fewer neighboring particles.
As ABFs, MMLS also mostly uses polynomials $ \bm{P}_{ji} $ multiplied by a RBF $ \varphi_{ji} $.

Based on the general theme of MLS or regularized MLS, a large number of modified methods have been proposed.
The \textbf{Interpolating Moving Least Squares (IMLS}) method introduced by Lancester and Salkaukas~\cite{lancaster1981} adds the soft constraint $ w_{ji} \approx \delta_{ji} $ to favor interpolating weight functions.
This allows for direct node-value assignment when imposing Dirichlet boundary condition.
This idea was later also extended to MMLS, resulting in \textbf{Interpolating Modified Moving Least Squares (IMMLS)}~\cite{lohit2022immls}.
Another group of methods called \textbf{Improved Moving Least Squares (IMLS)}~\cite{cheng2005improvedmls, liew2006improvedmls} or \textbf{Improved Interpolating Moving Least Squares (IIMLS)}~\cite{wang2012iimls, li2015iimls} allows for the interpolating weight functions to be modified by any RBF $\varphi_{ji} $.
To avoid ill-conditioned moment matrices, \textbf{Adaptive Orthogonal IIMLS (AO-IIMLS)}~\cite{wang2019aoiimls} uses an adaptive orthogonalization algorithm to remove small entries.
Slightly different regularization techniques were presented as \textbf{Regularized Moving Least Squares (RMLS)} and \textbf{Regularized Improved Interpolating Moving Least Squares (RIIMLS)} by Wang et al.~\cite{wang2018rmls}. RMLS and RIIMLS are closely related to MMLS and IIMLS, but use
direct Tikhonov--Miller regularization.
The \textbf{Least Squares Collocation Meshless Method (LSCM)}~\cite{zhang2001} introduces auxiliary points where the approximation is not evaluated, but the solution must still be satisfied.
An attempt to unify the many MLS-type methods was made by introducing the \textbf{Meshless Local Strong Form Method (MLSM)}~\cite{kosec2019mlsm}, which also highlighted the wide range of basis functions $ \bm{P}_{ji} $ that can be used in MLS-type approaches.

\subsubsection{Generalized finite-difference methods}\label{sec:gfdm}

Besides MLS-type methods, another large class of meshfree collocation schemes derived by error minimization are \textbf{Meshless Finite Difference Methods (MFDM)}. Originally introduced by Milewsky~\cite{milewski2012meshless}, MFDM take the viewpoint of finite difference stencils, which they generalize to irregular collocation point distributions.
Recent work in this direction focused on deriving higher-order approximations that only require inverting lower-order moment matrices~\cite{milewski2021mfdm, milewski2022mfdm}.
Instead of approximating $ L^{\bm{\alpha}}_{i} $ with some order of consistency $m$, this is done by first obtaining a lower-order approximation $ L^{\bm{\alpha},L}_{i} $ of order $L<m$, which is then extrapolated to the higher-order approximation
\begin{equation}
  L^{\bm{\alpha}}_{i}u = L^{\bm{\alpha},L}_{i}u - \Delta_{i}^{\bm{\alpha}}\, .
\end{equation}
with a correction $ \Delta_{i}^{\bm{\alpha}} $.
The correction is computed using Richardson extrapolation. For this, the Taylor expansion in Eq.~\eqref{eq:taylor_series_zm} is split into terms of order $\leq L$ and terms of order $>L$:
\begin{equation}
  u_{j}
  =
  \bm{X}_{ji}^{\leq L}\cdot\left.\bm{D}^{\leq L} u\right\rvert_{i} +
  \bm{X}_{ji}^{>L}\cdot\left.\bm{D}^{>L} u\right\rvert_{i} +
  e_{ji}^{m+1},
  \label{eq:mfmd_taylor_series_split}
\end{equation}
where
\begin{equation}
  \bm{X}_{ji}^{\leq L}
  =
  \left\{
    x_{1,ji}^{\alpha_1}/ \bm{\alpha}! \ldots x_{n,ji}^{\alpha_n} / \bm{\alpha}!
  \right\}_{|\bm{\alpha}|=0}^L,
  \quad
  \bm{X}_{ji}^{>L}
  =
  \left\{
    x_{1,ji}^{\alpha_1}/ \bm{\alpha}! \ldots x_{n,ji}^{\alpha_n} / \bm{\alpha}!
  \right\}_{|\bm{\alpha}|>L}^m,
\end{equation}
\begin{equation}
  \bm{D}^{\leq L} u
  =
  \left\{
    \partial^{\bm{\alpha}}  u
  \right\}_{\bm{|\alpha|}=0}^L,
  \quad
  \bm{D}^{>L} u
  =
  \left\{
    \partial^{\bm{\alpha}}  u
  \right\}_{\bm{|\alpha|}>L}^m.
  \label{eq:mfdm_low_high_order_derivatives}
\end{equation}
A recommended choice for the split point is $L = |\bm{\alpha}|$ for a target order of $ m = 2L $~\cite{milewski2022mfdm}.
This then allows splitting the optimization problem in Eq.~\eqref{eq:l2e:emin-optimizationformulation} into lower- and higher-order terms:
\begin{align*}
  {\mathcal{J}}_{\mathcal{E},i}
  =
  \displaystyle\sum_{j\in\mathcal{N}_i}
  \varphi_{ji} \left( e_{ji}^{m+1} \right)^2
  =
  \displaystyle\sum_{j\in\mathcal{N}_i}
  \varphi_{ji}
  \left(
    \bm{X}_{ji}^{\leq L}\cdot\left.\bm{D}^{\leq L} u\right\rvert_{i}
    + \bm{X}_{ji}^{>L}\cdot\left.\bm{D}^{>L} u\right\rvert_{i}
    -  u_{ji}
  \right)^2
  \\
  \mathrm{min} \left( \mathcal{J}_{\mathcal{E},i} \right)
  \quad
  \text{with respect to}
  \left.\bm{D}^{\leq L} u\right\rvert_{i}.
  \numberthis
  \label{eq:mfdm_emin-optimizationformulation}
\end{align*}
Performing the minimization only with respect to the low-order terms results in
\begin{equation}
  \left.\bm{D}^{\leq L} u\right\rvert_{i}
  =
  \displaystyle\sum_{j\in\mathcal{N}_i}\mathbf{M}_{i}^{-1}\bm{W}_{ji} u_{j} - \bm{\Delta}_{i}\, ,
  \label{eq:l2e:mfdm:ho-lo derivatives}
\end{equation}
where the first term on the right-hand side represents the low-order derivatives, and $\bm{\Delta}_i$ is the higher-order correction
\begin{equation}
  \bm{\Delta}_i
  = \mathbf{M}_{i}^{-1}\bm{W}_{ji}
  \left(
    \bm{X}_{ji}^{>L}\cdot\left.\bm{D}^{>L} u\right\rvert_{i}
  \right).
\end{equation}
Following the general recipe of Richardson extrapolation, the unknown higher-order approximations $ \left.\bm{D}^{>L} u\right\rvert_{i} $ are expressed as nested low-order approximations:
\begin{align*}
  \left.\bm{D}^{>L} u\right\rvert_{i}
  &=
  \left(
    \bm{D}^{\leq L}\mathbf{P}_{t}
  \right)
  \left(
    \left.\bm{D}^{\leq L} u\right\rvert_{i}\mathbf{P}_{u}
  \right)
  \\
  &=
  \displaystyle\sum_{j\in\mathcal{N}_i}
  \left(
    \mathbf{M}_{i}^{-1}\bm{W}_{ji}\mathbf{P}_{t}
  \right)
  \left(
    \left.\bm{D}^{\leq L} u\right\rvert_{j}\mathbf{P}_{u}
  \right)
  \\
  &=
  \displaystyle\sum_{j\in\mathcal{N}_i}
  \left(
    \mathbf{M}_{i}^{-1}\bm{W}_{ji}\mathbf{P}_{t}
  \right)
  \left(
    \displaystyle\sum_{s}\mathbf{M}_{i}^{-1}\bm{W}_{sj} u_{s}\mathbf{P}_{u} - \bm{\Delta}_{j}\mathbf{P}_{u}
  \right)\, ,
\end{align*}
which we expanded in the second and third lines by substituting the discrete operator. The matrices $\mathbf{P}_{t}$ and $\mathbf{P}_{u}$ are non-unique permutation matrices.
The recursive dependence on the correction $\bm{\Delta}$ is approximated by performing a finite number of iterations $k=1,2,\ldots$ of
\begin{equation}
  \left( \left.\bm{D}^{>L} u\right\rvert_{i} \right)^{ (k + 1) }
  =
  \displaystyle\sum_{j\in\mathcal{N}_i}
  \left(
    \mathbf{M}_{i}^{-1}\bm{W}_{ji}\mathbf{P}_{t}
  \right)
  \left(
    \displaystyle\sum_{s}\mathbf{M}_{i}^{-1}\bm{W}_{sj} u_{s}\mathbf{P}_{u} - \bm{\Delta}_{j}^{(k)}\mathbf{P}_{u}
  \right)\, ,
\end{equation}
with the lower-order approximation as a starting point, $ \bm{\Delta}^{(0)} = \left.\bm{D}^{\leq L} u\right\rvert_{i} $.
This iteration converges for contractive maps. Whether the map above is contractive or not, howevever, depends on the distribution of particles in $\mathcal{N}_i$.
To select a particular derivative $ \bm{\alpha} $, one again uses the mapping vector $ \bm{C}^{\bm{\alpha}} $ as in Eq.~\eqref{eq:l2e:approx-error-l2-min-include-mapping-vect}.
For further details, as well as for higher-order corrections of boundary conditions, we refer to the original works~\cite{milewski2012meshless, milewski2022mfdm}.

Several related variants of meshfree finite-difference methods exist. Among them the \textbf{Finite Difference Particle Method (FDPM})~\cite{wang2018lagrangian}, which uses the same minimization procedure as MFDM but excludes the zeroth moment.
A variant that includes the zeroth moment is the \textbf{Particle Difference Method (PDM)}~\cite{yoon2014extended}, also sometimes called \textbf{Particle Derivative Approximation (PDA)}. While essentially identical to MFDM, an important contribution of PDM was the extension of the approach to problems with discontinuities by introducing non-smooth RBFs $\varphi_{ji}$; the method is then also sometimes referred to as \textbf{Extended Particle Difference Method (EPDM)}.

In addition, one branch of the \textbf{Finite Pointset Method (FPsM)} has been derived by minimization of the approximation error. The name {\em Finite Pointset Method} has been used in the literature to refer to two distinct approaches. The first branch uses error minimization~\cite{Kuhnert2001_2, Tiwari2002, Tiwari2003}, whereas the second branch follows from minimization of weights~\cite{drumm2008finite, uhlmann2013, lu2016, abdessalam2016} and is discussed in Section~\ref{sec:l2p}. Interestingly, both branches originate from the same authors.
The first branch of FPsM~\cite{Kuhnert2001}, which falls into the category of this chapter, accounts for the specific problem to be solved by augmenting the optimization functional with additional penalty terms~\cite{resendiz2018, saucedo2020, gedela2013}. This is similar to the approach in CMLS as presented in Section~\ref{sec:mlsmethods}, but it predates CMLS.
Therefore, the FPsM was probably the first to introduce problem-specific corrections to meshfree collocation methods. However, the motivation for doing so was different from that of CMLS.
In FPsM, the penalty terms were motivated by solving a second-order boundary problem~\cite{Kuhnert2001}, and not by reducing the stencil support.
FPsM therefore does not use the notion of stencil coefficients, but it directly assembles the global matrix of the system.
Considering again the model problem from Eq.~\eqref{eq:cmls:example_problem}, the minimization in FPsM is:
\begin{equation}
  u^*
  =
  \argmin_{q \in \mathcal{V}_h}
  \displaystyle\sum_{i\in\mathcal{N}_i}
  \left[
    \left(
      \lambda_j\left( u\right) - \lambda_j\left(q\right)
    \right)^2 \varphi_{ji}
    +
    \left(
      f_j - \mathcal{L}_{\Omega} q
    \right)^2
    +
    \chi_{j,\partial\Omega}
    \left(
      g_j - \mathcal{L}_{\partial\Omega} q
    \right)^2
  \right]\, ,
  \label{eq:l2p:FPsM:l2_optim_problem}
\end{equation}
where, in contrast to Eq.~\eqref{eq:l2e:cmls:l2_optim_problem}, the penalty terms do not have a relaxation prefactor and are independent of the weight function $\varphi_j$.
In some works, the problem-dependent penalty terms in Eq.~\eqref{eq:l2p:FPsM:l2_optim_problem} are understood to act globally~\cite{suchde2018}. This then results in a field-approximation method rather than a local differential operator, providing an interesting connection between meshfree collocation methods and global surrogate methods like {\em Optimizing a Discrete Loss (ODIL)}~\cite{karnakov2023} and global polynomial surrogates~\cite{suarezCardona2023,SuarezCardona2024}.
Ultimately, this links meshfree collocation methods to {\em Physics-Informed Neural Networks (PINN)}~\cite{raissi2019,suarezCardona2022}, a connection that might be worth exploring more deeply in the future.

Finally, we mention the \textbf{Generalized Finite Difference Method (GFDM)}~\cite{liszka1980finite, swartz1969generalized}, which can be derived by minimization of the approximation error~\cite{gavete2003, benito2001, zhang2016gfdm, li2017swegfdm}. However, more recent works derive it by minimization of weights~\cite{suchde2017,suchde2018,suchde2019},
which is why we discuss it in Section~\ref{sec:l2p}.

\subsubsection{One-sided and upwind methods}

A major limitation of collocation methods in Eulerian point clouds is that they are unconditionally unstable in explicit time stepping when discretizing the advection operator $\vect{v}\cdot\nabla$ for some flow velocity field $\vect{v}(\vect{x})$. The three classic remedies are to either: (1) use a Lagrangian frame of reference, where the collocation points move with $\vect{v}$, (2) use one-sided neighborhoods $\mathcal{N}_i$ containing only neighbors in the \textit{upwind} direction, i.e., the direction from which the flow comes, or (3) use staggered grids where $\vect{v}$ is discretized in-between grid points. Any meshfree collocation method can trivially be used in a Lagrangian setting. For upwinding and staggering, specific methods have been proposed.

An upwind meshfree collocation method is the \textbf{Least-Squares Kinetic Upwind Method (LSKUM)} introduced by Ghosh and Deshpande~\cite{ghosh1995least,ghosh1995robust} and later extended to higher order as q-LSKUM~\cite{deshpande1998, ghosh2007qlskum}.
The method was originally designed to numerically solve the Euler equation of gas dynamics and the Navier--Stokes equations of compressible viscous fluids.
LSKUM exploits the physical relation between Boltzmann transport and the Navier--Stokes equations to solve compressible flow problems by kinetic flux vector splitting~\cite{deshpande1986,mandal1994}. Hence, the method is also sometimes called q-KFVS referring to Kinetic Flux Vector Splitting~\cite{ghosh2007qlskum}. The fluxes in the fluid-flow equations are the derivatives to be approximated.
LSKUM uses one-sided neighborhoods to compute stable approximations of the flux vector.
Following classic grid-based upwind schemes, only half of the neighbors of any given particle $ i $ are used.
We denote this truncated neighborhood by $\mathcal{N}_i^{\beta\pm}$, where $ \beta $ stands for a direction in the coordinate system, e.g.,~$\beta=x,\, y,\,z$, and the $\pm$ denotes which half of the neighborhood is used.
A similar solution using one-sided particle neighborhoods was introduced by Eldredge et al.~\cite{eldredge2002} for PSE.
The coordinate directions are chosen from the signs of the flow velocity components $\vect{v} = (v_x,\, v_y,\, \dots)$ such that the flow propagates from the neighboring nodes $ j $ to the receiving node $ i $, but not vice versa.
Apart from this difference, LSKUM uses the standard truncation-error minimization derivation with ABFs $\bm{W}_{ji} = \bm{X}_{ji} \varphi_{ji}$ excluding the zeroth moment.

Largely based on LSKUM, the \textbf{Kinetic Meshless Method (KMM)} was introduced by Praveen and Deshpande~\cite{Praveen2003, praveen2007kinetic} to use staggering instead of one-sided kernels. This leads to smaller neighborhood radii and renders KMM rotationally invariant, simplifying boundary condition treatment. Since the neighborhoods are isotropic, and the condition number of the moment matrix is generally smaller than in LSKUM, the accuracy and robustness (wavenumber dissipation) of the method are also improved~\cite{praveen2007kinetic}.
KMM achieves this by replacing the Boltzmann-transport fluxes of LSKUM with a MLS-type approach to approximate derivatives over midpoint function values $ u_{ji/2} \coloneq u( ( \bm{x}_j - \bm{x}_i) / 2 ) $. This advection stabilization is akin to classic staggered-grid finite differences and is sometimes referred to as \textbf{Dual Least-Squares (DLS) method}~\cite{praveen2007kinetic}.
It is based on choosing the midpoint value from the appropriate side
\begin{equation}
  u_{ij/2} =
  \begin{cases}
    u_i, & \quad \bm{v} \cdot \hat{\bm{e}}_{ji} \ge 0 \\
    u_j, & \quad \bm{v} \cdot \hat{\bm{e}}_{ji} < 0\, ,
  \end{cases}
\end{equation}
where $ \vect{v} $ is the flow velocity and $ \hat{\bm{e}}_{ji} $ is the unit vector pointing from $ \bm{x}_i $ to $ \bm{x}_j $.
Evaluating the Taylor expansion in Eq.~\eqref{eq:taylor_series_zm} at $ u_{ji/2} $ instead of $ u_{j} $, the ABFs then become the Taylor monomials evaluated between $\vect{x}_i$ and $( \bm{x}_j - \bm{x}_i) / 2 $, multiplied by the RBFs $\varphi_{ji}$. Hence, with $ \bm{W}_{ji} = \bm{X}_{ji/2} \varphi_{ji} $, we obtain the final KMM approximation:
\begin{equation}
  L^{\bm{\alpha}}_{i} u
  =
  \displaystyle\sum_{j\in\mathcal{N}_i}
  (u_{ji/2} - u_{i})
  \bm{W}_{ji} \cdot (\mathbf{M}_i^{-1} \bm{C}^{\bm{\alpha}})\, ,
  \quad
  \mathbf{M}_{i}
  =
  \displaystyle\sum_{j\in\mathcal{N}_i}\bm{X}_{ji/2} \bm{W}_{ji}^{\top}\, .
  \label{eq:l2e:kmm_discrete_diff_operator_recap}
\end{equation}
Compared to the canonical approximation in Eq.~\eqref{eq:l2e:approx-error-l2-min_ld_result}, the values are supported at different locations, and the vector of Taylor monomials (and hence the moment matrix) is correspondingly modified.
Higher-order versions of KMM reconstruct $u_{ji/2}$ using a linear scheme from the seminal five-paper series {\it ``Towards the ultimate conservative difference scheme''} \cite{vanLeer1997series}, forming a one-parameter family known as the {\it k-schemes}.
For $ k = -1 $, this includes a one-sided second-order scheme to interpolate variables from neighboring particles to the midpoints $ ( \bm{x}_j - \bm{x}_i) / 2 $
before computing derivatives using Eq.~\eqref{eq:l2e:kmm_discrete_diff_operator_recap}~\cite{vanLeer2021}.

\subsubsection{Sepcialized gradient and Laplacian methods}

The methods described so far provide approximations to arbitrary linear differential operators $\partial ^{\bm{\alpha}}$. In addition to these general numerical methods, there are also several contributions in the literature that target only specific derivatives, mostly the gradient or Laplace operator.
These are usually subsumed under the name \textbf{Lagrangian Differencing Dynamics (LDD)}~\cite{basic2020lagrangian,peng2021lagrangian} and were originally proposed by Basic et al.~under the name \textbf{Renormalized Meshless Operator (RMO)}~\cite{basic2018class}.
LDD blends ideas from moment approximation (Section~\ref{sec:aom}) and MLS
with the aim of finding a compromise between SPH-type approximations, which do not provide sufficient accuracy, and accurate collocation schemes, which are computationally more expensive.
The central tool for this is a \textit{correction matrix}, also known from several corrected SPH methods~\cite{randles1996sphmls, bonet1999}.

We denote this correction matrix by $\mathbf{M}_{i}^{1}$, since it corresponds to the moment matrix $ \mathbf{M}_i $ with only linear ($ m = 1 $) terms present. The zeroth moment is excluded.
Together with the symmetry of the weight function $\varphi_{ji}$, this allows constructing compromise approximations. In the original work~\cite{basic2018class}, \textit{absolutely normalized} weight functions $ \psi_{ji} = \varphi_{ji} / \sum_{j\in\mathcal{N}_i} \varphi_{ji} $ were used, but in more recent works~\cite{basic2020lagrangian, peng2021lagrangian, basic2022ldd} the same authors use a canonical weight function $ \varphi_{ji} $ as introduced above in Section~\ref{sec:com_definition}.
While the derivation of LDD somewhat deviates from the derivation presented at the beginning of this section, it leads to the same approximation for the gradient $L_i^{\nabla}\approx \nabla$ as in Eq.~\eqref{eq:l2e:emin-generaldiscreteoperator} with $m=1$.
For the Laplace operator $\partial^{\bm{\alpha}}=\nabla ^2$, LDD presents three approximations.
The first, called \textit{na\"{i}ve approximation}, is:
\begin{equation}
  L_{i}^{\nabla^2\text{naive}} u
  =
  2n
  \displaystyle \sum_{j\in\mathcal{N}_i}
  \frac{\varphi_{ji}}{ \| \vect{x}_{ji} \|^2 }
  \left(
    u_{ji}
    -
    \vect{x}_{ji} \cdot L_i^{\nabla} u
  \right).
  \label{eq:l2e:ldd:naive}
\end{equation}
The approximation error of Eq.~\eqref{eq:l2e:ldd:naive} consists of the first-order error of the gradient approximation $ L_i^{\nabla} u $ and higher-order truncation contributions from the Taylor expansion~\cite{basic2018class}.
This is only accurate for large numbers of neighbors $|\mathcal{N}_i|$ and if the points in $\mathcal{N}_i$ are regularly distributed. Since these conditions are rather limiting in practice, the \textit{sum approximation} has been proposed:
\begin{equation}
  L_{i}^{\nabla^2\text{sum}} u
  =
  2n
  \frac{
    \displaystyle \sum_{j\in\mathcal{N}_i}
    \varphi_{ji}  u_{ji}
    \left(
      1 - \vect{x}_{ji} \cdot \left(\mathbf{M}_{i}^{1}\right)^{-1}\bm{O}_i
    \right)
  }
  {
    \displaystyle \sum_{j\in\mathcal{N}_i}
    \varphi_{ji}
    \lVert \vect{x}_{ji} \rVert^2
    \left(
      1 - \vect{x}_{ji} \cdot \left(\mathbf{M}_{i}^{1}\right)^{-1}\bm{O}_i
    \right)
  }\, ,
  \label{eq:l2e:ldd:sum}
\end{equation}
where $ \bm{O}_i $ is the \textit{offset vector} of particle $ i $
\begin{equation}
  \bm{O}_i =
  \displaystyle \sum_{j\in\mathcal{N}_i}
  \varphi_{ji}
  \vect{x}_{ji}\, ,
\end{equation}
pointing from $\vect{x}_i $ to the \textit{weighted center} of the neighborhood. For locally regular and symmetric (around $\vect{x}_i$) neighbor distributions, $\bm{O} = \mathbf{0}$. Then, and when using only the nearest neighbors, Eq.~\eqref{eq:l2e:ldd:sum} reduces to the compact second-order central finite difference stencil for the Laplacian.
Compared to Eq.~\eqref{eq:l2e:ldd:naive}, the sum approximation does not contain the first-order gradient error and does not require explicitly computing the gradient values. The resulting approximation error is an average of the error over the neighborhood. The third Laplacian approximation, called \textit{least-squares approximation}, further improves accuracy by employing a MLS procedure over the offset vector.
This defines a modified moment matrix $\widehat{\mathbf{M}}_i$, which in LDD is called \textit{renormalization matrix}
\begin{equation}
  \widehat{\mathbf{M}}_i
  =
  \displaystyle \sum_{j\in\mathcal{N}_i}
  \varphi_{ji}
  \bm{Q}_{ji}
  \bm{Q}_{ji}^{\top}\, ,
\end{equation}
where
\begin{equation}
  \bm{Q}_{ji}
  =
  \bm{X}_{ji}^{2}
  -
  \vect{x}_{ji}^{\top}
  \displaystyle \sum_{j\in\mathcal{N}_i}
  \varphi_{ji} \left(\mathbf{M}_{i}^{1}\right)^{-1}\vect{x}_{ji}
  \left( \bm{X}_{ji}^{2} \right) ^{\top}\, ,
  \label{eq:ldd-lsLaplace-q}
\end{equation}
and $ \bm{X}_{ji}^{2} $ are the square subcomponents of $\bm{X}_{ji}$
\begin{equation}
  \bm{X}_{ji}^{2}
  \coloneq
  \left\{
    x_{\beta,ji}^2
  \right\}_{\beta=1}^{n}\, .
\end{equation}
From this, one obtains a Laplacian approximation of the form~\cite{basic2018class}:
\begin{equation}
  L_{i}^{\nabla^2\text{LS}} u
  =
  \bm{I} \cdot \widehat{\mathbf{M}}_i^{-1}
  \displaystyle \sum_{j\in\mathcal{N}_i}
  2 \varphi_{ji} \vect{Q}_{ji}
  \left(
    u_{ji}
    -
    \vect{x}_{ji} \cdot L_i^{\nabla} u
  \right)\, ,
  \label{eq:l2e:ldd:ls_full}
\end{equation}
with $ \bm{I} \in \mathbb{R}^{n} $ a vector of ones.
In LDD, this is referred to as \textit{(full) inversion approximation} of the Laplace operator.
It is computationally inefficient, since several loops over the neighbors are required to compute all parts of the expression.
A computationally more efficient formulation of this approximation is obtained by modifying the moment matrix as
\begin{equation}
  \widetilde{\mathbf{M}}_i
  =
  \displaystyle \sum_{j\in\mathcal{N}_i}
  \varphi_{ji}
  \bm{X}_{ji}^{2}\left(\bm{X}_{ji}^{2}\right)^{\top}.
\end{equation}
This leads to:
\begin{equation}
  L_{i}^{\nabla^2\textrm{LS}} u
  =
  \bm{I} \cdot \widetilde{\mathbf{M}}_i^{-1}
  \displaystyle \sum_{j\in\mathcal{N}_i}
  2 \varphi_{ji} \vect{Q}_{ji}
  \left(
    u_{ji}
    -
    \vect{x}_{ji} \cdot L_i^{\nabla} u
  \right).
  \label{eq:l2e:ldd:ls_basic}
\end{equation}
Now, unlike in Eq.~\eqref{eq:l2e:ldd:ls_full}, both $ \widetilde{\mathbf{M}}_i $ and $ \mathbf{M}_{i}^{1} $ can be computed in the same loop, but the gradient $L_i^{\nabla} u$ is again required. In LDD, this is referred to as \textit{(basic) inversion approximation}.
Both least-squares approximations, Eqs.~\eqref{eq:l2e:ldd:ls_full} and \eqref{eq:l2e:ldd:ls_basic}, run into problems with the condition number of the matrix when some elements of the vectors $\bm{Q}_{ji}$ or $\bm{X}_{ji}^{2}$ are close to zero~\cite{basic2018class}.

\subsubsection{Summary}

Summarizing the many variants of meshfree collocation methods derived by $\ell_2$ minimization of the approximation error, we find that most of them build directly or indirectly on MLS~\cite{lancaster1981}.
With the exception of LSMFM~\cite{park2001, zhang2005}, which uses direct derivatives,
they all first approximate the function $u$ itself and then approximate derivatives of $u$ by computing derivatives of this approximation, an approach known as \textit{diffuse derivatives}~\cite{chen2017}.
Notable extensions to this recipe were provided by GMLS, CMLS, and MMLS.
GMLS~\cite{wendland2004, mirzaei2012}
extended the approach to general linear bounded operators and provided a more solid mathematical foundation.
CMLS~\cite{trask2016} and MMLS~\cite{joldes2015mmls, joldes2019mmls}
introduced additional penalty terms that favor more compactly supported stencils and improve the condition number of the MLS system.
MFDM~\cite{milewski2012meshless}
introduced the idea of iterative Richardson extrapolation to construct higher-order approximations by inverting only lower-order moment matrices.
One branch of FPsM~\cite{Kuhnert2001_2, Tiwari2002, Tiwari2003} introduced problem-specific penalties to solve second-order boundary problems, and some versions of GFDM were also derived by error minimization~\cite{gavete2003, benito2001, zhang2016gfdm, li2017swegfdm}. Upwind and staggering methods for Eulerian advection solvers were introduced by LSKUM~\cite{ghosh1995least,ghosh1995robust}
and KMM~\cite{Praveen2003, praveen2007kinetic},
respectively, and PDM~\cite{yoon2014extended}
showed an extension to discontinuous solutions with non-smooth weight functions.
LDD~\cite{basic2020lagrangian,peng2021lagrangian},
finally, presented a hybrid derivation for discretized gradient and Laplace operators following a conceptually different approach, which can avoid inversion of the moment matrix.

\section{Derivation by $\ell_2$ minimization of weights}\label{sec:l2p}
A third way to derive the weights for the discrete operator in Eq.~\eqref{eq:general_discrete_diff_operator_zm} is by applying  $\ell_2$ minimization directly to the weights.
This considers the optimization problem of finding the smallest weights such that the moment conditions from Section~\ref{sec:aom} are satisfied up to and including order $m$, hence:
\begin{align*}
  \mathcal{J}_{\mathcal{P},i}
  =
  \frac{1}{2}
  \displaystyle\sum_{j\in\mathcal{N}_i} \frac{ 1 }{ \varphi_{ji} }\left(w_{ji}^{\bm{\alpha}} \right)^{2} \\
  \mathrm{min} \left( \mathcal{J}_{\mathcal{P},i} \right)
  \quad
  \text{subject to}
  \quad
  \displaystyle\sum_{j\in\mathcal{N}_i}
  u_{ji}w_{ji}^{\bm{\alpha}}
  =
  \bm{C}^{\bm{\alpha}} \cdot \left.\bm{D} u^* \right\rvert_{i}
  \quad
  \forall u^* \in \mathbb{P}_m \, ,
  \numberthis
  \label{eq:l2p:mc-optimizationformulation}
\end{align*}
where $ \mathbb{P}_m $ is the linear space of polynomials of degree at most $ m $ with elements $u^*$.
The constraint enforcing derivative approximation can equivalently be written in terms of the scaled Taylor monomials, resulting in:
\begin{align*}
  \mathrm{min} \left( \mathcal{J}_{\mathcal{P},i} \right)
  \quad
  \text{subject to}
  \quad
  \displaystyle\sum_{j\in\mathcal{N}_i}
  \bm{X}_{ji}w_{ji}^{\bm{\alpha}}
  =
  \bm{B}_{i}^{\bm{\alpha}}
  =
  \bm{C}^{\bm{\alpha}}\, .
  \numberthis
  \label{eq:l2p:mc-optimizationformulation-withmoments}
\end{align*}
The reason why we divide by $\varphi_{ji}$ in Eq.~\eqref{eq:l2p:mc-optimizationformulation}, rather than multiplying, is purely notational, as this leads to the same form of the final approximation as in Section~\ref{sec:l2e}.

\nomvar{ $\mathcal{J}_{\mathcal{P},i}$ }{functional for $\ell_2$ minimization under the moment condition constraint [-]}
\nomvar{ $\bm{\lambda}$ }{Lagrange multipliers [-]}
\nomvar{ $\bm{w}_{\bullet i}^{\bm{\alpha}}$ }{vector of weight functions for all neighbors of particle $i$}

The formal solution as derived in Appendix~\ref{sec:appendix-l2p} is:
\begin{equation}
  w_{ji}^{\bm{\alpha}}
  =
  \varphi_{ji}
  \bm{X}_{ji}
  \cdot
  \mathbf{M}_i^{-1} \bm{C}^{\bm{\alpha}},
  \qquad
  \mathbf{M}_i
  =
  \displaystyle\sum_{j\in\mathcal{N}_i} \bm{X}_{ji} \bm{X}_{ji}^{\top} \varphi_{ji}\, .
\end{equation}
With this, we obtain the same discrete differential operator as in Section~\ref{sec:l2e}:
\begin{equation}
  L^{\bm{\alpha}}_{i}u
  =
  \displaystyle\sum_{j\in\mathcal{N}_i} u_{ji}
  \bm{W}_{ji} \cdot
  \left(
    \mathbf{M}_i^{-1}
    \bm{C}^{\bm{\alpha}}
  \right).
  \label{eq:l2p:approx-polynomial-method-l2-min}
\end{equation}
Also, as in Section~\ref{sec:l2e}, the basis of the weight functions is fixed to $ \bm{W}_{ji} = \bm{X}_{ji} \varphi_{ji} $, and the moment matrix $\mathbf{M}_i$ is symmetric. This makes the approach less general than the approximation of moments discussed in Section~\ref{sec:aom}.
A difference from Section~\ref{sec:l2e} is that in the present case the form of the weights is determined by the choice of the functional $\mathcal{J}_{\mathcal{P},i}$, and they amount to inverse MLS weights.
This approach is sometimes referred to as \textit{polynomial method}~\cite{suchde2018incompressible}, and its connection with approximation-error minimization has been explored~\cite{suchde2018}.

\subsection{Matrix notation}

In the matrix formulation introduced in Section \ref{sec:basic-approach-matrix-formulation}, the optimization problem from Eq.~\eqref{eq:l2p:mc-optimizationformulation} reads:
\begin{align*}
  \mathcal{J}_{\mathcal{P},i}
  =
  \left\lVert \left(\mathbf{V}_i \right)^{-1/2} \bm{w}_{\bullet i}^{\bm{\alpha}} \right\rVert^2  \\
  \mathrm{min} \left( \mathcal{J}_{\mathcal{P},i} \right)
  \quad
  \text{subject to}
  \quad
  \mathbf{X}^{\top}\bm{w}_{\bullet i}^{\bm{\alpha}} = \bm{C}^{\bm{\alpha}},
  \numberthis
  \label{eq:app_matrices_polynomial_opt}
\end{align*}
where $ \left(\mathbf{V}_i  \right)^{-1/2} $ is easy to compute, since $ \mathbf{V}_{i}  $ is a diagonal matrix.
The formal solution to this minimization problem, as derived in Appendix~\ref{sec:appendix-l2p}, is:
\begin{equation}
  \bm{w}_{\bullet i}^{\bm{\alpha}}
  =
  \mathbf{W}_i \mathbf{M}_i^{-1} \bm{C}^{\bm{\alpha}}.
  \label{eq:l2p_weights}
\end{equation}
With these weights, the discretized differential operator
\begin{equation}
  L^{\bm{\alpha}}_{i} u
  =
  \bm{\bm{U}}_i \cdot \bm{w}_{\bullet i}^{\bm{\alpha}}
  =
  \bm{\bm{U}}_i \cdot \left( \mathbf{C}_i\bm{C}^{\bm{\alpha}} \right)
  =
  \bm{\bm{U}}_i \cdot \left( \mathbf{W}_i \mathbf{Y}_{i}\bm{C}^{\bm{\alpha}} \right)
  =
  \bm{\bm{U}}_i \cdot \left( \mathbf{W}_i \mathbf{M}_{i}^{-1}\bm{C}^{\bm{\alpha}} \right)
  \label{eq:l2p:general_discrete_diff_operator_with_coefs_do_matrix_form_l2min}
\end{equation}
is exactly the same as in Eq.~\eqref{eq:l2p:approx-polynomial-method-l2-min}, albeit in matrix notation.

\subsection{Methods derived by $\ell_2$ minimization of weights}

The two methods discussed below were both originally derived by $\ell_2$ error minimization (see Section~\ref{sec:l2e methods}). However, they are nowadays more commonly presented by $\ell_2$ minimization of weights, which is why we discuss them here.

Reaching all the way back to the very roots of strong form collocation methods, the \textbf{Generalized Finite Difference Method (GFDM)} was introduced by Liszka and Orkisz \cite{liszka1980finite, liszka1984}. It still is a popular version of MFD as inspired by the work of Jensen~\cite{jensen1972} and of Swartz and Wendroff~\cite{swartz1969generalized}.
Since its inception, various improvements have been made to GFDM in different applicaition areas.
The method was originally formulated using approximation-error minimization (Section~\ref{sec:l2e}), but in recent years it was mostly presented according to the derivation in this section.
For recent formulations and applications of GFDM, we refer to Refs.~\cite{suchde2017,suchde2018,suchde2019}.
Therein, the method is derived exactly as presented above, with Taylor monomials multiplied by a suitable window function  as ABFs $ \bm{W}_{ji} = \bm{X}_{ji} \varphi_{ji} $.
The zeroth moment is typically excluded when deriving GFDM by approximation-error minimization but included when deriving the method by minimization of weights.

Since its introduction, GFDM has been extended in multiple directions. This included the possibility of using non-smooth basis functions in addition to polynomials (not RBFs)~\cite{suchde2019}, using {\em control cells}
to enforce local flux conservation ~\cite{suchde2017} or to address challenges arising from discontinuities in the coefficients of the PDE~\cite{kraus2023}. In addition, a central-weight control mechanism has been proposed to guarantee operator positivity \cite{suchde2018}, establishing a connection between GFDM and FPsM~\cite{Furst2001,jin2004positivity}.
Finally, positive GFDM operators have been derived with better computational efficiency and stability~\cite{seibold2008, bracco2025}, and higher-order GFDM operators have been derived~\cite{kraus2025}.

The second method in the present category is one branch of the \textbf{Finite Pointset Method (FPsM}). The FPsM was initially introduced by Kuhnert and Tiwari~\cite{Kuhnert2001} using error minimization (see Section~\ref{sec:l2e}) and has since undergone numerous improvements and extensions~\cite{Kuhnert2001_2, Tiwari2002, Tiwari2003, Tiwari2005, tiwari2007fpsm}.
There are two distinct branches of the FPsM method: the first follows an error-minimization approach, while the second uses minimizaton of weights. We discussed the first branch in Section~\ref{sec:gfdm}. Here, we focus on the second branch of FPsM, derived by minimization of weights~\cite{drumm2008finite, uhlmann2013, lu2016, abdessalam2016}.
This branch of FPsM, instead of including a problem-specific PDE penalty term, introduces a constraint requiring the weight function at the center particle $w^{\bm{\alpha}}_{ii}$ to be negative~\cite{drumm2008finite}. This improves the condition number of the linear equation system for the weight coefficients.
The additional constraint is added to Eq.~\eqref{eq:l2p:mc-optimizationformulation-withmoments} as:
\begin{align*}
  \mathcal{J}_{\mathcal{P},i}
  =
  \displaystyle\sum_{j\in\mathcal{N}_i} \frac{ 1 }{ \varphi_{ji} } \left(w_{ji}^{\bm{\alpha}} \right)^{2}\\
  \mathrm{min} \left( \mathcal{J}_{\mathcal{P},i} \right)
  \quad
  \text{subject to}
  \quad
  \displaystyle\sum_{j\in\mathcal{N}_i}
  \bm{X}_{ji}w_{ji}^{\bm{\alpha}}
  =
  \bm{B}_{i}^{\bm{\alpha}}
  =
  \bm{C}^{\bm{\alpha}},
  \;
  w_{ii}^{\bm{\alpha}} = - \frac{\epsilon}{h_i^2}\, ,
  \numberthis
  \label{eq:fpsm_mc-optimizationformulation-extended}
\end{align*}
where $\epsilon$ is a small non-negative number, and  $ h_i $ is the \textit{stencil size}, i.e., the radius of the support of $ \varphi_{ji} $.
To cover the additional constraint with the notation introduced above, we define an extended vector of Taylor monomials (which are the basis functions used in FPsM) $ \bm{X}_{ji}^{+} $ and an extended version of the mapping vector $ \bm{C}^{\bm{\alpha}} $, as:
\begin{equation}
  \bm{X}_{ji}^{+}
  =
  \left[
    \bm{X}_{ji},\, \delta_{ji}
  \right]
  \in \mathbb{R}^{p+1},
  \quad
  \bm{C}^{\bm{\alpha},+}
  =
  \left[
    \bm{C}^{\bm{\alpha}},\, \epsilon/ h_i^2
  \right]
  \in \mathbb{N}_0^{p+1},
\end{equation}
where $ \delta_{ji} $ is the Kronecker delta.
Using these \textit{extended} vectors, we rewrite the constraint in Eq.~\eqref{eq:fpsm_mc-optimizationformulation-extended} as:
\begin{equation}
  \sum_{j\in\mathcal{N}_i}
  \bm{X}_{ji}^{+}w_{ji}^{\bm{\alpha}}
  =
  \bm{C}^{d,+}.
\end{equation}
Denoting $ \bm{W}_{ji}^{+} = \bm{X}_{ji}^{+} \varphi_{ji} \in \mathbb{R}^{p+1} $, this leads to the \textit{extended} moment matrix
\begin{equation}
  \mathbf{M}_i^{+}
  = \sum_{j\in\mathcal{N}_i}\bm{X}_{ji}^{+}\left(\bm{W}_{ji}^{+}\right)^{\top} \in \mathbb{R}^{p+1,p+1}.
\end{equation}
The solution of the optimization problem then is:
\begin{equation}
  L^{\bm{\alpha}}_{i} u
  = \sum_{j\in\mathcal{N}_i}  u_{ji}
  \bm{W}_{ji}^{+} \cdot
  \left(
    \left(
      \mathbf{M}_i^{+}
    \right )^{-1}
    \bm{C}^{\bm{\alpha},+}
  \right)\, .
\end{equation}

A similar constraint on the central particle has also been used to ensure positivity of the discrete operator, in particular of a discrete Laplace operator~\cite{suchde2019}.
This uses two stencil coefficients with different values. One satisfies polynomial exactness, while the other resides in the null space of polynomial exactness. Therefore, any linear combination of these coefficients satisfies polynomial exactness.
The method finds a particular linear combination such that the resulting discrete operator is positive, yet accurately approximates the target differential operator~\cite{suchde2018}.

The differences between this branch of FPsM and the branch derived by error minimization (Section~\ref{sec:gfdm}) have been highlighted in Refs.~\cite{suchde2018incompressible, Iliev2003}.
The FPsM based on weight minimization is nowadays sometimes considered a special case of GFDM and referred to by this name~\cite{suchde2018incompressible, suchde2024}.

\section{Derivation by generalized $\ell_2$ minimization}\label{sec:gl2p}

The final way of deriving meshfree collocation methods is, at this point in time, purely theoretical. We are not aware of any methods in the literature that were derived or presented in this way. Nevertheless, it would be possible to do so, which is why we include this derivation here.

The derivation is motivated by the observation that the formulations obtained by $\ell_2$ minimization (Eqs.~\eqref{eq:l2e:approx-error-l2-min_ld_result} and \eqref{eq:l2p:approx-polynomial-method-l2-min} for error and weights minimization, respectively) are slightly different from what is obtained by approximating moment conditions (Eq.~\eqref{eq:general_discrete_diff_operator_with_coefs_do}).
Although the results are structurally alike, the minimization approaches led to a symmetric moment matrix composed of products of Taylor-monomial vectors $ \bm{X}_{ji}$. This corresponds to the choice of basis functions $ \bm{W}_{ji} = \bm{X}_{ji} \varphi_{ji} $.
The result is therefore less general than Eq.~\eqref{eq:general_discrete_diff_operator_with_coefs_do}, where the moment matrix need not be symmetric. This has been recognized by some authors, who posited that alternative choices of $\bm{W}_{ij} $ could improve the accuracy of the method~\cite{trask2016,HighOrderDiffKing2020, HighOrderSimuKing2022}.

In the case of $\ell_2$ error minimization, the symmetry of the moment matrix is a consequence of the procedure and cannot be changed. The structure of the matrix results from minimizing the square of the approximation error. In the case of minimizing for the \textit{best-fitting} weights with respect to the scaled moments, however, the structure of the moment matrix is a consequence of the formulation of the optimization problem, which can be generalized.
We generalize the formulation by allowing weight functions $w_{ji}^{\bm{\alpha}} = \bm{W}_{ji} \cdot \bm{\Psi}_i $ with arbitrary ABFs $ \bm{W}_{ji} $  scaled by the weight coefficients $ \bm{\Psi}_i $, as introduced in Eq.~\eqref{eq:w_abf_zm}.
Formulating the optimization problem from Section \ref{sec:l2p} with such weight functions yields:
\begin{align*}
  \widetilde{\mathcal{J}}_{\mathcal{P},i}
  =
  \frac{1}{2}
  \displaystyle\sum_{j\in\mathcal{N}_i}
  \left(
    w_{ji}^{\bm{\alpha}}
  \right)^2
  =
  \frac{1}{2}
  \displaystyle\sum_{j\in\mathcal{N}_i}
  \left(
    \bm{W}_{ji} \cdot \bm{\Psi}_{i}^{\bm{\alpha}}
  \right)^2
  \\
  \mathrm{min} \left( \widetilde{\mathcal{J}}_{\mathcal{P},i} \right)
  \quad
  \text{subject to}
  \quad
  \displaystyle\sum_{j\in\mathcal{N}_i}
  u_{ji}w_{ji}^{\bm{\alpha}}
  =
  \displaystyle\sum_{j\in\mathcal{N}_i}
  u_{ji} \bm{W}_{ji} \cdot \bm{\Psi}_{i}^{\bm{\alpha}}
  =
  \bm{C}^{\bm{\alpha}} \cdot \left.\bm{D} u^* \right\rvert_{i}
  \quad
  \forall u^* \in \mathbb{P}_m \, .
  \numberthis
  \label{eq:mcgeneralized-optimizationformulation}
\end{align*}
Due to the moment conditions in Eq.~\eqref{eq:general_discrete_diff_operator_theoret_orderd}, this can be rewritten as:
\begin{align*}
  \widetilde{\mathcal{J}}_{\mathcal{P},i}
  =
  \frac{1}{2}
  \displaystyle\sum_{j\in\mathcal{N}_i}
  \left(
    w_{ji}^{\bm{\alpha}}
  \right)^2
  =
  \frac{1}{2}
  \displaystyle\sum_{j\in\mathcal{N}_i}
  \left(
    \bm{W}_{ji} \cdot \bm{\Psi}_{i}^{\bm{\alpha}}
  \right)^2
  \\
  \mathrm{min} \left( \widetilde{\mathcal{J}}_{\mathcal{P},i} \right)
  \quad
  \text{subject to}
  \quad
  \displaystyle\sum_{j\in\mathcal{N}_i}
  \bm{X}_{ji}w_{ji}^{\bm{\alpha}}
  =
  \displaystyle\sum_{j\in\mathcal{N}_i}
  \bm{X}_{ji}\bm{W}_{ji} \cdot \bm{\Psi}_{i}^{\bm{\alpha}}
  =
  \mathbf{M}_i\bm{\Psi}_{i}^{\bm{\alpha}}
  =
  \bm{C}^{\bm{\alpha}}.
  \numberthis
  \label{eq:gl2p:mcgeneralized-optimizationformulation-mc}
\end{align*}
\nomvar{ $\mathcal{J}_{\mathcal{P},i}$ }{generalized functional for $\ell_2$ minimization under the moment condition constraint [-]}
The optimal coefficients $\bm{\Psi}_{i}^{\bm{\alpha}}$ are found by solving this minimization problem as shown in Appendix \ref{sec:appendix-gl2p}, leading to:
\begin{align*}
  w_{ji}^{\bm{\alpha}}
  =
  \bm{W}_{ji}
  \cdot
  \left(
    \mathbf{M}_i^{-1}
    \bm{C}^{\bm{\alpha}}
  \right),
  \quad
  \mathbf{M}_i
  =
  \displaystyle \sum_{j\in\mathcal{N}_i}
  \bm{X}_{ji} \bm{W}_{ji}^{\top}\, ,
\end{align*}
where $\mathbf{M}_i$ allows for arbitrary choices of $ \bm{W}_{ji} $, in contrast to the formulations from Sections~\ref{sec:l2e} and \ref{sec:l2p}.
The resulting discrete operator
\begin{equation}
  L^{\bm{\alpha}}_{i}u
  =
  \displaystyle\sum_{j\in\mathcal{N}_i}
  u_{ji}\bm{W}_{ji}\cdot \left(\mathbf{M}_{i}^{-1} \bm{C}^{\bm{\alpha}} \right)
  \label{eq:gl2p:general_discrete_diff_operator_with_coefs}
\end{equation}
is identical to that found by moment approximation without the zeroth moment, cf.~Eq.~\eqref{eq:general_discrete_diff_operator_with_coefs_withouth_zeroth_moment}. This not only suggests a generalization of miminization-based approaches, but it also links them to moment-based approximations, showing how and when the two are equivalent.

\subsection{Matrix notation}

In the matrix notation from Section \ref{sec:basic-approach-matrix-formulation}, the optimization problem in Eq.~\eqref{eq:mcgeneralized-optimizationformulation} reads:
\begin{align*}
  \widetilde{\mathcal{J}}_{\mathcal{P},i}
  =
  \frac{ 1 }{ 2 }
  \lVert \bm{w}_{i}^{\bm{\alpha}} \rVert^2
  =
  \frac{ 1 }{ 2 }
  \lVert \mathbf{W}_i \bm{\Psi}_{i}^{\bm{\alpha}} \rVert^2  \\
  \mathrm{min} \left( \widetilde{\mathcal{J}}_{\mathcal{P},i} \right)
  \quad
  \text{subject to}
  \quad
  \mathbf{M}_i\bm{\Psi}_i^{\bm{\alpha}} = \bm{C}^{\bm{\alpha}}\, .
  \numberthis
  \label{eq:app_matrices_polynomial_opt_generalized}
\end{align*}
Solving this as shown in Appendix \ref{sec:appendix-gl2p}, we find the optimal weight function
\begin{equation}
  \bm{w}_{\bullet i}^{\bm{\alpha}}
  =
  \mathbf{W}_i
  \mathbf{M}_i^{-1}
  \bm{C}^{\bm{\alpha}}\, .
  \label{eq:l2pg_weights}
\end{equation}
This has the same form as Eq.~\eqref{eq:l2p_weights}. But in the present form, the weight matrix
\begin{equation}
  \mathbf{M}_{i}
  =
  \displaystyle\sum_{j\in\mathcal{N}_i}
  \bm{X}_{ji}\bm{W}_{ji}^{\top}
  =
  \mathbf{X}_i^{\top} \mathbf{W}_i
\end{equation}
need not be symmetric and has the full flexibility of the definition in Eq.~\eqref{eq:momentum_matrix_matrix_form}, without any particular form of  ABFs $ \mathbf{W}_i $ (or $ \bm{W}_{ji} $).
This leads to the same discrete operator as in Eq.~\eqref{eq:gl2p:general_discrete_diff_operator_with_coefs} in matrix form:
\begin{equation}
  L^{\bm{\alpha}}_{i} u
  =
  \bm{\bm{U}}_i \cdot \bm{w}_{\bullet i}^{\bm{\alpha}}
  =
  \bm{\bm{U}}_i \cdot \left( \mathbf{C}_i\bm{C}^{\bm{\alpha}} \right)
  =
  \bm{\bm{U}}_i \cdot \left( \mathbf{W}_i \mathbf{Y}_{i}\bm{C}^{\bm{\alpha}} \right)
  =
  \bm{\bm{U}}_i \cdot \left( \mathbf{W}_i \mathbf{M}_{i}^{-1}\bm{C}^{\bm{\alpha}} \right)\, .
  \label{eq:general_discrete_diff_operator_with_coefs_do_matrix_form_l2mcg}
\end{equation}

\section{Summary\label{sec:summary}}

We comprehensively surveyed meshfree collocation methods for consistently (with some desired order of accuracy $m$) approximating differential operators on continuously labeled unstructured point clouds. We presented a notation that can accommodate all methods reviewed in a common general formulation. This suggested classifying the methods according to how they were derived, or are usually presented, into three groups: methods derived by approximation of moments, methods derived by minimization of the approximation error, and methods derived by minimization of weights. For each group, we presented the derivation in detail, showed how it fits the present formulation in both vector and matrix notation, and reviewed methods from the literature that fall within that group. This allowed us to highlight similarities and differences between the methods.

Comparing methods, we found that all use similar discrete approximation operators, which can be formulated within the same framework. There are, however, also important differences, for example in whether the zeroth moment is explicitly accounted for, the structure of the moment matrix (symmetric or not), or if additional penalty terms are included when deriving the approximation. These differences shape the numerical properties of the methods, such as their computational cost, stability in explicitly time stepping, wavenumber dissipation, boundary condition handling, and robustness to noise. Therefore, each method has its niche, and different applications may benefit from different methods. It would have been beyond the scope of this paper to discuss all numerical consequences of the differences between methods or to give advice on which methods are well suited for what type of applications. Instead, we focused on the mathematical formulation of the methods, their formal relationships, and the historical evolution of ideas and concepts.

An important difference between the three groups of methods defined here is whether the moment matrix is allowed to be asymmetric, which makes a wider choice of basis functions accessible. Minimizing the approximation error reaches the same order of accuracy as moment-approximation methods, but only the Taylor mononials can be used as basis functions. This observation led us to formulate a fourth, hitherto unexplored derivation of meshfree collocation. We named it {\em generalized $\ell _2$ minimization}, as it extends weight-minimization approaches to using arbitrary anisotropic basis functions. Interestingly, this yielded the exact same discrete operator as found by moment approximation. By showing that the same approximation can be derived in two conceptually different ways---by solving a generalized miminization problem or by canceling unwanted moments---we established a formal link between these groups of meshfree collocation methods. This hopefully allows theoretical results to be transfered between methods in the future.

We summarize all reviewed methods in Table~\ref{tab:summary}, together with their classification according to the presented method groups. Unsurprisingly, many of them are similar. Some of them are even identical and were re-discovered over the years and published under different names. In the Supplementary Material, we provide an encyclopedia of all methods for later reference.

\setcounter{table}{0}
\begin{longtblr}
  [ caption = {Summary of reviewed methods, listed alphabetically. For each method, we list the discrete operator used to approximate derivatives of a continuous function $u$ over unstructured point clouds. We then give the basis functions (ABFs) used in that method, and say if it includes (\cmark) the zeroth moment or not (\xmark). Both symbols simultaneously indicate that the method can do either, depending on parameter settings. Finally, we give the group each method belongs to (AOM: approximation of moments as in Section~\ref{sec:aom}; $\ell_2(\mathcal{E})$: $\ell_2$ minimization of approximation error as in Section~\ref{sec:l2e}; $\ell_2(\mathcal{P})$: $\ell_2$ minimization of weights as in Section~\ref{sec:l2p}), followed by comments. The acronyms of the methods are defined in the Nomenclature section at the beginning of the paper.}, label = {tab:summary} ]
  { colspec = {|  Q[0.2cm,valign=m] | Q[6.4cm,valign=m] | Q[3.3cm,valign=m] | Q[0.2cm,valign=m] | Q[0.7cm,valign=m] |Q[3.6cm,valign=m] | },
  rowhead = 1 }
  \hline
  \rotatebox[origin=c]{90}{\textbf{Method}} & \center{\textbf{Operator Approximation}} & \center{\textbf{ Basis }} & \rotatebox[origin=c]{90}{\textbf{0th Moment}} & \rotatebox[origin=c]{90}{\textbf{Group}} & \textbf{Comments} \\ \hline
  \hline
  \rotatebox[origin=c]{90}{CMLS}
  &
  \begin{fleqn}
    \begin{equation*}
      L^{\bm{\alpha}}_{i} u
      =
      \displaystyle\sum_{j\in\mathcal{N}_i}
      \left( u_{j} \bm{W}_{ji} + \bm{Q}_{ji} \right) \cdot (\mathbf{M}_i^{-1} \bm{C}^{\bm{\alpha}})
    \end{equation*}
    \begin{equation*}
      \mathbf{M}_{i}
      =
      \displaystyle\sum_{j\in\mathcal{N}_i}\bm{P}_{ji}\bm{W}_{ji}^{\top} + \mathbf{Q}_i
    \end{equation*}
  \end{fleqn}
  &
  \begin{equation*}
    \bm{W}_{ji}
    =
    \bm{P}_{ji} \varphi_{ji}
  \end{equation*}
  &
  \cmark
  &
  $\ell_2(\mathcal{E})$
  &
  Includes problem-dependent penalty terms for more compact operator support.
  \\ \hline
  \rotatebox[origin=c]{90}{DC-PSE}
  &
  \begin{fleqn}
    \begin{equation*}
      L^{\bm{\alpha}}_{i} u =
      \sum_{j\in\mathcal{N}_i} \left( u_{j}\pm  u_{i}\right) \bm{W}_{ji} \cdot \left(\mathbf{M}_{i}^{-1}\widetilde{\bm{C^{\bm{\alpha}}}}\right)
    \end{equation*}
    \begin{equation*}
      \mathbf{M}_{i}
      =
      \displaystyle\sum_{j\in\mathcal{N}_i}
      \widetilde{\bm{X}}_{ji}\bm{W}_{ji}^{\top}
    \end{equation*}
    with positive sign $+$ for $|\bm{\alpha}|$ odd and $-$ for $|\bm{\alpha}|$ even.
  \end{fleqn}
  &
  \begin{equation*}
    \bm{W}_{ji} =
    \widetilde{\bm{P}}_{ji} \varphi_{ji}
  \end{equation*}
  \begin{equation*}
    \text{usually } \bm{W}_{ji} = \widetilde{\bm{X}}_{ji} \varphi_{ji}
  \end{equation*}
  &
  \cmark \,\,\, \xmark
  &
  AOM
  &

  \\ \hline
  \rotatebox[origin=c]{90}{FDPM}
  &
  \begin{fleqn}
    \begin{equation*}
      L^{\bm{\alpha}}_{i} u
      =
      \displaystyle\sum_{j\in\mathcal{N}_i}
      u_{ji} \bm{W}_{ji} \cdot (\mathbf{M}_i^{-1} \bm{C}^{\bm{\alpha}}),
    \end{equation*}
    \begin{equation*}
      \mathbf{M}_{i}
      =
      \displaystyle\sum_{j\in\mathcal{N}_i}\bm{X}_{ji}\bm{W}_{ji}^{\top}
    \end{equation*}
  \end{fleqn}
  &
  \begin{equation*}
    \bm{W}_{ji}
    =
    \bm{X}_{ji} \varphi_{ji}
  \end{equation*}
  &
  \xmark
  &
  $\ell_2(\mathcal{E})$
  &
  \\ \hline
  \rotatebox[origin=c]{90}{FPM}
  &
  \begin{fleqn}
    \begin{equation*}
      L^{\bm{\alpha}}_{i} u
      =
      \displaystyle\sum_{j\in\mathcal{N}_i}
      u_{ji} \bm{W}_{ji} \cdot (\mathbf{M}_i^{-1} \bm{C}^{\bm{\alpha}}),
    \end{equation*}
    \begin{equation*}
      \mathbf{M}_{i}
      =
      \displaystyle\sum_{j\in\mathcal{N}_i}\bm{P}_{ji}\bm{W}_{ji}^{\top}
    \end{equation*}
  \end{fleqn}
  &
  \begin{equation*}
    \bm{W}_{ji}
    =
    \bm{P}_{ji} \varphi_{ji}
  \end{equation*}
  \begin{equation*}
    \text{usually } \bm{W}_{ji} = \bm{X}_{ji} \varphi_{ji}
  \end{equation*}
  &
  \xmark
  &
  $\ell_2(\mathcal{E})$
  &
  \\ \hline
  \rotatebox[origin=c]{90}{FPsM}
  &
  \begin{fleqn}
    \begin{equation*}
      L^{\bm{\alpha}}_{i} u
      =
      \displaystyle\sum_{j\in\mathcal{N}_i}  u_{ji}
      \bm{W}_{ji}^{+} \cdot
      \left(
        \left(
          \mathbf{M}_i^{+}
        \right )^{-1}
        \bm{C}^{\bm{\alpha},+} \cdot
      \right)
    \end{equation*}
    \begin{equation*}
      \mathbf{M}_i^{+}
      =
      \displaystyle\sum_{j\in\mathcal{N}_i}\bm{X}_{ji}^{+}\left(\bm{W}_{ji}^{+}\right)^{\top}
    \end{equation*}
  \end{fleqn}
  &
  \begin{equation*}
    \bm{W}_{ji}^{+}
    =
    \bm{X}_{ji}^{+} W^0_{ji}.
  \end{equation*}
  &
  \cmark
  &
  $\ell_2(\mathcal{E})$ or $\ell_2(\mathcal{P})$
  &
  Includes additional penalty to enforce negative center weight.
  In some works, the penalty terms are assumed globally, resulting in field approximation method.
  The $\ell_2(\mathcal{E})$ branch also includes problem dependent penalty terms.
  \\ \hline
  \rotatebox[origin=c]{90}{GFDM}
  &
  \begin{fleqn}
    \begin{equation*}
      L^{\bm{\alpha}}_{i} u
      =
      \displaystyle\sum_{j\in\mathcal{N}_i}
      u_{j} \bm{W}_{ji} \cdot (\mathbf{M}_i^{-1} \bm{C}^{\bm{\alpha}})
    \end{equation*}
    \begin{equation*}
      \mathbf{M}_{i}
      =
      \displaystyle\sum_{j\in\mathcal{N}_i}\bm{X}_{ji}\bm{W}_{ji}^{\top}
    \end{equation*}
  \end{fleqn}
  &
  \begin{equation*}
    \bm{W}_{ji}
    =
    \bm{X}_{ji} \varphi_{ji}
  \end{equation*}
  &
  \xmark

  or

  \cmark
  &
  $\ell_2(\mathcal{E})$ or $\ell_2(\mathcal{P})$
  &
  $\ell_2(\mathcal{E})$ usually exclude the 0th moment, $\ell_2(\mathcal{P})$ usually includes 0th moment
  \\ \hline
  \rotatebox[origin=c]{90}{GMLS}
  &
  \begin{fleqn}
    \begin{equation*}
      L^{\bm{\alpha}}_{i} u
      =
      \displaystyle\sum_{j\in\mathcal{N}_i}
      u_{j} \bm{W}_{ji} \cdot (\mathbf{M}_i^{-1} \bm{C}^{\bm{\alpha}}),
    \end{equation*}
    \begin{equation*}
      \mathbf{M}_{i}
      =
      \displaystyle\sum_{j\in\mathcal{N}_i}\bm{P}_{ji}\bm{W}_{ji}^{\top}
    \end{equation*}
  \end{fleqn}
  &
  \begin{equation*}
    \bm{W}_{ji}
    =
    \bm{P}_{ji} \varphi_{ji}
  \end{equation*}
  &
  \cmark
  &
  $\ell_2(\mathcal{E})$
  &
  Special case of approximating arbitrary linear bounded operators on Banach spaces.
  \\ \hline
  \rotatebox[origin=c]{90}{GRKCM}
  &
  \begin{fleqn}
    \begin{equation*}
      L^{\bm{\alpha}}_{i} u
      =
      \displaystyle\sum_{j\in\mathcal{N}_i}
      u_{j} \bm{W}_{ji} \cdot (\mathbf{M}_i^{-1} \bm{C}^{\bm{\alpha}})
    \end{equation*}
    \begin{equation*}
      \mathbf{M}_{i}
      =
      \displaystyle\sum_{j\in\mathcal{N}_i}\bm{X}_{ji}\bm{W}_{ji}^{\top}
    \end{equation*}
  \end{fleqn}
  &
  \begin{equation*}
    \bm{W}_{ji}
    =
    \bm{X}_{ji} \varphi_{ji}
  \end{equation*}
  &
  \cmark
  &
  AOM
  &
  \\ \hline
  \rotatebox[origin=c]{90}{GRKP}
  &
  \begin{fleqn}
    \begin{equation*}
      L^{\bm{\alpha}}_{i, [m]} u
      =
      \displaystyle\sum_{j\in\mathcal{N}_i}
      u_{ji}\bm{W}_{ji,[l,m]}\cdot \left(\mathbf{M}_{i,[l,m]}^{-1} \bm{C}^{\bm{\alpha}}_{[l,m]} \right)
    \end{equation*}
    \begin{equation*}
      \mathbf{M}_{i,[l,m]}
      =
      \displaystyle\sum_{j\in\mathcal{N}_i}
      \bm{X}_{ji,[l,m]} \bm{W}_{ji}^{\top}
    \end{equation*}
  \end{fleqn}
  &
  \begin{align*}
    &\bm{W}_{ji,[l,m]} = \\
    &\bm{X}_{ji,[l,m]} W^0_{ji,[l,m]}.
  \end{align*}
  &
  \cmark \,\,\, \xmark
  &
  AOM
  &

  \\ \hline
  \rotatebox[origin=c]{90}{HOCSPH}
  &
  \begin{fleqn}
    \begin{equation*}
      L^{\bm{\alpha}}_{i} u =
      \sum_{j\in\mathcal{N}_i} u_{ji}\bm{W}_{ji}\cdot \left(\mathbf{M}_{i}^{-1} \widetilde{\bm{C^{\bm{\alpha}}}}\right)v_j
    \end{equation*}
    \begin{equation*}
      \mathbf{M}_{i} = \sum_{j\in\mathcal{N}_i}\widetilde{\bm{X}}_{ji}\bm{W}_{ji}^{\top}v_j
    \end{equation*}
  \end{fleqn}
  &
  \begin{equation*}
    \bm{W}_{ji} =
    \left.\bm{D}\varphi\right\rvert_{\vect{x}_{ji}}
  \end{equation*}
  &
  \xmark
  &
  AOM
  &
  \\ \hline
  \rotatebox[origin=c]{90}{KKM}
  &
  \begin{fleqn}
    \begin{equation*}
      L^{\bm{\alpha}}_{i} u
      =
      \displaystyle\sum_{j\in\mathcal{N}_i}
      u_{ji} \bm{W}_{ji} \cdot (\mathbf{M}_i^{-1} \bm{C}^{\bm{\alpha}})
    \end{equation*}
    \begin{equation*}
      \mathbf{M}_{i}
      =
      \displaystyle\sum_{j\in\mathcal{N}_i}\bm{X}_{ji}\bm{W}_{ji}^{\top}
    \end{equation*}
  \end{fleqn}
  &
  \begin{equation*}
    \bm{W}_{ji}
    =
    \bm{X}_{ji} \varphi_{ji}
  \end{equation*}
  &
  \xmark
  &
  $\ell_2(\mathcal{E})$
  &
  \\ \hline
  \rotatebox[origin=c]{90}{LABFM}
  &
  \begin{fleqn}
    \begin{equation*}
      L^{\bm{\alpha}}_{i} u =
      \sum_{j\in\mathcal{N}_i} u_{ji}\bm{W}_{ji}\cdot \left(\mathbf{M}_{i}^{-1} \widetilde{\bm{C^{\bm{\alpha}}}}\right)
    \end{equation*}
    \begin{equation*}
      \mathbf{M}_{i} = \sum_{j\in\mathcal{N}_i}\widetilde{\bm{X}}_{ji}\bm{W}_{ji}^{\top}
    \end{equation*}
  \end{fleqn}
  &
  \begin{equation*}
    \bm{W}_{ji} =
    \left.\bm{D}\varphi\right\rvert_{\vect{x}_{ji}}
  \end{equation*}
  \begin{equation*}
    \bm{W}_{ji} =
    \bm{W}_{ji} = \widetilde{\bm{P}}_{ji} \varphi_{ji},
  \end{equation*}
  &
  \xmark
  &
  AOM
  &
  \\ \hline
  \rotatebox[origin=c]{90}{LDD}
  &
  \SetCell[c=2]{l}
  {
    $
    L_{i}^{\nabla^2,\textrm{naive}} u
    =
    2d
    \displaystyle \sum_{j\in\mathcal{N}_i}
    \frac{\varphi_{ji}}{ \lVert \vect{x}_{ji} \rVert^2 }
    \left(
      u_{ji}
      -
      \vect{x}_{ji} \cdot L_i^{\nabla} u
    \right)
    $ \\
    \vspace{0.2cm}
    $
    L_{i}^{\nabla^2,\textrm{sum}} u
    =
    2d
    \frac{
      \displaystyle \sum_{j\in\mathcal{N}_i}
      \varphi_{ji}  u_{ji}
      \left(
        1 - \vect{x}_{ji} \cdot \mathbf{M}_{1,i}^{-1}\bm{O}_i
      \right)
    }
    {
      \displaystyle \sum_{j\in\mathcal{N}_i}
      \varphi_{ji}
      \lVert \vect{x}_{ji} \rVert^2
      \left(
        1 - \vect{x}_{ji} \cdot \mathbf{M}_{1,i}^{-1}\bm{O}_i
      \right)
    },
    \quad
    $ \\
    $
    \qquad    \bm{O}_i =
    \displaystyle \sum_{j\in\mathcal{N}_i}
    \varphi_{ji}
    \vect{x}_{ji}.
    $ \\
    \vspace{0.2cm}
    $
    L_{i}^{\nabla^2,\textrm{LS}} u
    =
    \bm{I} \cdot \widetilde{\mathbf{M}}_i^{-1}
    \displaystyle \sum_{j\in\mathcal{N}_i}
    2 \varphi_{ji} \vect{Q}_{ji}
    \left(
      u_{ji}
      -
      \vect{x}_{ji} \cdot L_i^{\nabla} u
    \right),
    \quad
    $ \\
    $
    \qquad \bm{Q}_{ji}
    =
    \bm{X}_{ji}^{2}
    -
    \vect{x}_{ji}^{\top}
    \displaystyle \sum_{j\in\mathcal{N}_i}
    \varphi_{ji} \mathbf{M}_{1,i} \vect{x}_{ji}
    \left( \bm{X}_{ji}^{2} \right) ^{\top}
    $
  }
  &
  &
  \cmark
  &
  \textit{ren.}
  &
  Only first and second derivatives are considered.
  Derivation combines properties of AOM and $\ell_2(\mathcal{E})$ and it is inspired by \textit{renormalisation}
  \\ \hline
  \rotatebox[origin=c]{90}{LSKUM}
  &
  \begin{fleqn}
    \begin{equation*}
      L^{\bm{\alpha}}_{i} u
      =
      \displaystyle\sum_{j\in\mathcal{N}_i}
      u_{ji} \bm{W}_{ji} \cdot (\mathbf{M}_i^{-1} \bm{C}^{\bm{\alpha}})
    \end{equation*}
    \begin{equation*}
      \mathbf{M}_{i}
      =
      \displaystyle\sum_{j\in\mathcal{N}_i}\bm{X}_{ji}\bm{W}_{ji}^{\top}
    \end{equation*}
  \end{fleqn}
  &
  \begin{equation*}
    \bm{W}_{ji}
    =
    \bm{X}_{ji} \varphi_{ji}
  \end{equation*}
  &
  \xmark
  &
  $\ell_2(\mathcal{E})$
  &
  Exploits a physical connection with Boltzmann transport.
  \\ \hline
  \rotatebox[origin=c]{90}{LSMFM}
  &
  \begin{fleqn}
    \begin{align*}
      L^{\bm{\alpha}}_{i} u
      =
      \displaystyle\sum_{j\in\mathcal{N}_i}  u_{j}
      \bigg(
        & \partial^{\bm{\alpha}} \bm{X}_{ji} \varphi_{ji} \cdot \mathbf{M}_{i}^{-1} \bm{C}^{0} +\\
        & \bm{X}_{ji} \partial^{\bm{\alpha}} \varphi_{ji} \cdot \mathbf{M}_{i}^{-1} \bm{C}^{0} +\\
        & \bm{X}_{ji} \varphi_{ji} \cdot \partial^{\bm{\alpha}} \mathbf{M}_{i}^{-1} \bm{C}^{0}
      \bigg)
    \end{align*}
    \begin{equation*}
      \mathbf{M}_{i}
      =
      \displaystyle\sum_{j\in\mathcal{N}_i}
      \bm{X}_{ji}\bm{W}_{ji}^{\top}.
    \end{equation*}
  \end{fleqn}
  &
  \begin{equation*}
    \bm{W}_{ji} =
    \bm{P}_{ji} \varphi_{ji}
  \end{equation*}
  \begin{equation*}
    \text{usually } \bm{W}_{ji} = \bm{X}_{ji} \varphi_{ji}
  \end{equation*}
  &
  \cmark
  &
  $\ell_2(\mathcal{E})$
  &
  Referred to as \textit{direct derivatives}.
  \\ \hline
  \rotatebox[origin=c]{90}{MFDM}
  &
  \begin{fleqn}
    \begin{equation*}
      L^{d,(L)}_{i} u
      =
      \displaystyle\sum_{j\in\mathcal{N}_i}
      u_{ji} \bm{W}_{ji} \cdot (\mathbf{M}_i^{-1} \bm{C}^{\bm{\alpha}}),
    \end{equation*}
    \begin{equation*}
      L^{d,(H)}_{i} = L^{d,(L)}_{i} - \Delta_{i}^{\bm{\alpha}}.
    \end{equation*}
    \begin{equation*}
      \mathbf{M}_{i}
      =
      \displaystyle\sum_{j\in\mathcal{N}_i}\bm{X}_{ji}\bm{W}_{ji}^{\top}
    \end{equation*}
  \end{fleqn}
  &
  \begin{equation*}
    \bm{W}_{ji}
    =
    \bm{X}_{ji} \varphi_{ji}
  \end{equation*}
  &
  \cmark
  &
  $\ell_2(\mathcal{E})$
  &
  Provides an alternative treatment of high-order approximations.
  \\ \hline
  \rotatebox[origin=c]{90}{MLSM}
  &
  \begin{fleqn}
    \begin{equation*}
      L^{\bm{\alpha}}_{i} u
      =
      \displaystyle\sum_{j\in\mathcal{N}_i}
      u_{j} \bm{W}_{ji} \cdot (\mathbf{M}_i^{-1} \bm{C}^{\bm{\alpha}}),
    \end{equation*}
    \begin{equation*}
      \mathbf{M}_{i}
      =
      \displaystyle\sum_{j\in\mathcal{N}_i}\bm{P}_{ji} \bm{W}_{ji}^{\top}.
    \end{equation*}
  \end{fleqn}
  &
  \begin{equation*}
    \bm{W}_{ji}
    =
    \bm{P}_{ji} \varphi_{ji}.
  \end{equation*}
  &
  \cmark
  &
  $\ell_2(\mathcal{E})$
  &
  \\ \hline
  \rotatebox[origin=c]{90}{MMLS}
  &
  \begin{fleqn}
    \begin{equation*}
      L^{\bm{\alpha}}_{i} u
      =
      \displaystyle\sum_{j\in\mathcal{N}_i}
      u_{j} \bm{W}_{ji} \cdot (\mathbf{M}_i^{-1} \bm{C}^{\bm{\alpha}}),
    \end{equation*}
    \begin{equation*}
      \mathbf{M}_{i}
      =
      \displaystyle\sum_{j\in\mathcal{N}_i}\bm{P}_{ji}\bm{W}_{ji}^{\top} + \mathbf{H}_{\mu}
    \end{equation*}
  \end{fleqn}
  &
  \begin{equation*}
    \bm{W}_{ji}
    =
    \bm{P}_{ji} \varphi_{ji}
  \end{equation*}
  &
  \cmark
  &
  $\ell_2(\mathcal{E})$
  &
  Has reduced order of accuracy for singular penalty matrices.
  \\ \hline
  \rotatebox[origin=c]{90}{PDM/PDA}
  &
  \begin{fleqn}
    \begin{equation*}
      L^{\bm{\alpha}}_{i} u
      =
      \displaystyle\sum_{j\in\mathcal{N}_i}
      u_{ji} \bm{W}_{ji} \cdot (\mathbf{M}_i^{-1} \bm{C}^{\bm{\alpha}}),
    \end{equation*}
    \begin{equation*}
      \mathbf{M}_{i}
      =
      \displaystyle\sum_{j\in\mathcal{N}_i}\bm{X}_{ji}\bm{W}_{ji}^{\top}
    \end{equation*}
  \end{fleqn}
  &
  \begin{equation*}
    \bm{W}_{ji}
    =
    \bm{X}_{ji} \varphi_{ji}
  \end{equation*}
  &
  \xmark
  &
  $\ell_2(\mathcal{E})$
  &
  \\ \hline
  \rotatebox[origin=c]{90}{RKCM}
  &
  \begin{fleqn}
    \begin{align*}
      L^{\bm{\alpha}}_{i} u
      =
      \displaystyle\sum_{j\in\mathcal{N}_i}  u_{j}
      \bigg(
        & \partial^{\bm{\alpha}} \bm{X}_{ji} \varphi_{ji} \cdot \mathbf{M}_{i}^{-1} \bm{C}^{0} +\\
        & \bm{X}_{ji} \partial^{\bm{\alpha}} \varphi_{ji} \cdot \mathbf{M}_{i}^{-1} \bm{C}^{0} +\\
        & \bm{X}_{ji} \varphi_{ji} \cdot \partial^{\bm{\alpha}} \mathbf{M}_{i}^{-1} \bm{C}^{0}
      \bigg)
    \end{align*}
    \begin{equation*}
      \mathbf{M}_{i}
      =
      \displaystyle\sum_{j\in\mathcal{N}_i}
      \bm{X}_{ji}\bm{W}_{ji}^{\top}.
    \end{equation*}
  \end{fleqn}
  &
  \begin{equation*}
    \bm{W}_{ji} =
    \bm{P}_{ji} \varphi_{ji}
  \end{equation*}
  \begin{equation*}
    \text{usually } \bm{W}_{ji} = \bm{X}_{ji} \varphi_{ji}
  \end{equation*}
  &
  \cmark
  &
  AOM
  &
  Referred to as \textit{direct derivatives}.
  \\ \hline
\end{longtblr}

\appendix

\section{Component-wise example in 2D\label{sec:appendix-oldAnd2D}}
We find that to understand meshfree collocation methods, it often helps to go through the derivation from Section~\ref{sec:aom} in a low-dimensional case with all components explicitly written out. We therefore present the derivation in two dimensions (2D), separately for when the zeroth moment is included and when it is not.

Consider a 2D domain $\Omega \subset \mathbb{R}^2 $ with $ N $ arbitrarily placed points (nodes, particles) $\mathcal{X} = \{ \vect{x}_1, \dots, \vect{x}_N : \vect{x}_i \in \overline{\Omega} \}$, where $\overline{\Omega} = \Omega \cup \partial \Omega$, with positions $ \vect{x}_i = [ x_{i},\, y_{i} ]^\top $.
Every particle has a unique index $ \mathcal{P} = \{i : \vect{x}_i \in \mathcal{X}\}$.
Each particle interacts with other particles within a finite neighborhood $ \mathcal{N}_i $ around it, defined by a cutoff radius. The number of neighbors $ N_i $ in $ \mathcal{N}_i $ in 2D can be estimated as $ N_i \approx 4 \pi ( h_i / \Delta x_i )^2 $.

The particles collectively discretize a sufficiently smooth (at least $C^m$) function $  u(x,\, y)$. The value of this function at each particle $ i $ is $  u_{i} =  u(x_{i},\, y_{i}) $. We denote $  u_{ji} =  u_j -  u_i $. Moreover, by $ m $ we denote the order of consistency of the approximation.

\subsection{Derivation with the 0th moment\label{sec:appendix_2d_with_zm}}

We define a general discrete operator $ L^{\bm{\alpha}}_{i} $ at particle $ i $, which approximates a specific partial derivative $ \bm{\alpha} $ as
\begin{equation}
  L^{\bm{\alpha}}_{i} u=
  \displaystyle\sum_{j\in\mathcal{N}_i}
  u_{j}w^{\bm{\alpha}}_{ji}
  \approx
  \bm{C}^{\bm{\alpha}} \cdot \left.\bm{D} u\right\rvert_{i}\, ,
  \label{eq:appendix_2d_general_discrete_diff_operator_zm}
\end{equation}
where $ w^{\bm{\alpha}}_{ji} $ are the operator weights, and the right-hand side $ \bm{C}^{\bm{\alpha}} \cdot \left.\bm{D}  u\right\rvert_{i} $ is just a notation for the multi-index derivative $ \mathcal{D}^{|\bm{\alpha}|} u = \partial^{|\bm{\alpha}|} u / (\partial x^{\alpha_1} \partial y^{\alpha_2})$ with $ \bm{D} u $ denoting the vector of all possible derivatives up to and including order $m$, which in 2D is
\begin{equation}
  \bm{D} u
  =
  \left( u,
    \frac{\partial u}{\partial{x}},\,
    \frac{\partial u}{\partial{y}},\,
    \frac{\partial^{2} u}{\partial{x}^{2}},\,
    \frac{\partial^{2} u}{\partial{x}\partial{y}},\,
    \frac{\partial^{2} u}{\partial{y}^{2}},\,
    \frac{\partial^{3} u}{\partial{x}^{3}},\,
    \frac{\partial^{3} u}{\partial{x}^{2}\partial{y}},\,
    \frac{\partial^{3} u}{\partial{x}\partial{y}^{2}},\,
    \frac{\partial^{3} u}{\partial{y}^{3}},\,
    \frac{\partial^{4} u}{\partial{x}^{4}},\,
    \dots,\,
    \frac{\partial^{m} u}{\partial{y}^{m}}
  \right)^{\!\!\top}.
  \label{eq:appendix_2d_derivatives_vect_zm}
\end{equation}
$ \bm{C}^{\bm{\alpha}} $  is the\textit{mapping vector} that \textit{selects} the desired derivatives from the vector $ \bm{D}\left(\cdot\right) $, as introduced in Eq.~\ref{eq:selector_vect_example}.
The number of elements in $\bm{D} u$ is $ p = p_\mathrm{2D} = 1 + (m^2 + 3m)/2 $.
Thus, $ \bm{D} u \in \mathbb{R}^p $ and $ \bm{C}^d \in \mathbb{N}_0^p $.

The vector of 2D Taylor monomials is:
\begin{equation}
  \bm{X}_{ji}
  =
  \left(
    1,\,
    x_{ji},\,
    y_{ji},\,
    \frac{x_{ji}^{2}}{2},\,
    x_{ji}y_{ji},\,
    \frac{y_{ji}^{2}}{2},\,
    \frac{x_{ji}^{3}}{6},\,
    \frac{x_{ji}^{2}y_{ji}}{2},\,
    \frac{x_{ji}y_{ji}^{2}}{2},\,
    \frac{y_{ji}^{3}}{6},\,
    \frac{x_{ji}^{4}}{24},\,
    \dots \,
    \frac{y_{ji}^{m}}{m!}
  \right)^{\!\!\top},
  \label{eq:appendix_2d_taylor_monomials_vect_zm}
\end{equation}
where $ x_{ji} = x_{j} - x_{i}$ and $ y_{ji} = y_{j} - y_{i}$ are distances between particles. Thus, $ \bm{X}_{ji} \in \mathbb{R}^p $, $ \bm{X}_{ji} \neq \bm{X}_{j} - \bm{X}_{i} $.
This allows us to write the Taylor expansion of the function $  u $ in the neighborhood of particle $ i $ as
\begin{equation}
  u_{j}
  =
  \bm{X}_{ji}\cdot\left.\bm{D} u\right\rvert_{i} + e_{ji}^{m+1}.
  \label{eq:appendix_2d_taylor_series_zm}
\end{equation}
where $e_{ji}^{m+1} $ is the approximation error \textit{(Lagrange form of the remainder)} of orders $m +1$ and higher.
Dropping the error $e_{ji}^{m+1} $ and substituting the Taylor expansion into the operator approximation in Eq.~\eqref{eq:appendix_2d_general_discrete_diff_operator_zm}, we find:
\begin{equation}
  L^{\bm{\alpha}}_{i} u
  =
  \displaystyle\sum_{j\in\mathcal{N}_i}
  \left.\bm{X}_{ji}\cdot\bm{D} u\right\rvert_{i}w^{\bm{\alpha}}_{ji}\, ,
  \label{eq:appendix_2d_error_short_zm}
\end{equation}
which expands to
\begin{multline}
  L^{\bm{\alpha}}_{i} u
  =
  \left. u\right\rvert_{i}
  \displaystyle\sum_{j\in\mathcal{N}_i}
  w^{\bm{\alpha}}_{ji}
  +\left.\frac{\partial u}{\partial{x}}\right\rvert_{i}
  \displaystyle\sum_{j\in\mathcal{N}_i}
  x_{ji}w^{\bm{\alpha}}_{ji}
  +\left.\frac{\partial u}{\partial{y}}\right\rvert_{i}
  \displaystyle\sum_{j\in\mathcal{N}_i}
  y_{ji}w^{\bm{\alpha}}_{ji}
  +\left.\frac{\partial^{2} u}{\partial{x}^{2}}\right\rvert_{i}
  \displaystyle\sum_{j\in\mathcal{N}_i}
  \frac{x_{ji}^{2}}{2}w^{\bm{\alpha}}_{ji}
  \\
  +\left.\frac{\partial^{2} u}{\partial{x}\partial{y}}\right\rvert_{i}
  \displaystyle\sum_{j\in\mathcal{N}_i}
  x_{ji}y_{ji}w^{\bm{\alpha}}_{ji}
  +\left.\frac{\partial^{2} u}{\partial{y}^{2}}\right\rvert_{i}
  \displaystyle\sum_{j\in\mathcal{N}_i}
  \frac{y_{ji}^{2}}{2}w^{\bm{\alpha}}_{ji}
  +\left.\frac{\partial^{3} u}{\partial{x}^{3}}\right\rvert_{i}
  \displaystyle\sum_{j\in\mathcal{N}_i}
  \frac{x_{ji}^{3}}{6}w^{\bm{\alpha}}_{ji}
  +\left.\frac{\partial^{3} u}{\partial{x}^{2}\partial{y}}\right\rvert_{i}
  \displaystyle\sum_{j\in\mathcal{N}_i}
  \frac{x_{ji}^{2}y_{ji}}{2}w^{\bm{\alpha}}_{ji}
  \\
  +\left.\frac{\partial^{3} u}{\partial{x}\partial{y}^{2}}\right\rvert_{i}
  \displaystyle\sum_{j\in\mathcal{N}_i}
  \frac{x_{ji}y_{ji}^{2}}{2}w^{\bm{\alpha}}_{ji}
  +\left.\frac{\partial^{3} u}{\partial{y}^{3}}\right\rvert_{i}
  \displaystyle\sum_{j\in\mathcal{N}_i}
  \frac{y_{ji}^{3}}{6}w^{\bm{\alpha}}_{ji}
  +\left.\frac{\partial^{4} u}{\partial{x}^{4}}\right\rvert_{i}
  \displaystyle\sum_{j\in\mathcal{N}_i}
  \frac{x_{ji}^{4}}{24}w^{\bm{\alpha}}_{ji}
  +\dots.
  \label{eq:appendix_2d_error_long_zm}
\end{multline}
The right-hand side thus corresponds to a sum of moments of the weight coefficients $ w^{\bm{\alpha}} $.
We denote the vector of moments by
\begin{equation}
  \bm{B}^{\bm{\alpha}}_{i}
  =
  \displaystyle\sum_{j\in\mathcal{N}_i}
  \bm{X}_{ji}w^{\bm{\alpha}}_{ji}\, .
  \label{eq:appendix_2d_moments_vect_zm}
\end{equation}
With the vector $ \bm{B}^{\bm{\alpha}}_{i} $, Eq.~\eqref{eq:appendix_2d_error_long_zm} can be written compactly as:
\begin{equation}
  L^{\bm{\alpha}}_{i} u
  =
  \left.\bm{D} u\right\rvert_{i} \cdot \bm{B}^{\bm{\alpha}}_{i}\, ,
  \label{eq:appendix_2d_error_moments_vect_zm}
\end{equation}

The goal is to find a weight functions $ w^{\bm{\alpha}}_{ji} $ such that the moments associated with the desired derivatives are close to one, and all other moments (up to and including order $ m $) are zero.
We make an ansatz for $ w_{ji}^{\bm{\alpha}} $ to be a weighted sum of some anisotropic basis functions (ABFs) $ W_{ji} $ as introduced in Section~\ref{sec:com_definition}, thus:
\begin{equation}
  w^{\bm{\alpha}}_{ji}
  =
  \bm{W}_{ji}\cdot\bm{\Psi^{\bm{\alpha}}}_{i}
  =
  W^{0}_{ji}\Psi^{\bm{\alpha}}_{i,0}
  +W^{1}_{ji}\Psi^{\bm{\alpha}}_{i,1}
  +W^{2}_{ji}\Psi^{\bm{\alpha}}_{i,2}
  +W^{3}_{ji}\Psi^{\bm{\alpha}}_{i,3}
  +W^{4}_{ji}\Psi^{\bm{\alpha}}_{i,4}
  +\dots+W^{n}_{ji}\Psi^{\bm{\alpha}}_{i,p}
  \label{eq:appendix_2d_w_abf_zm}
\end{equation}
where $ \bm{W}_{ji} \in \mathbb{R}^p $ and $ \bm{\Psi^{\bm{\alpha}}}_{i} \in \mathbb{R}^p $ represent the vector of basis functions and the vector of weight coefficients, respectively.
\begin{equation}
  \bm{W}_{ji}
  =
  \left(W^{0}_{ji},W^{1}_{ji},W^{2}_{ji},W^{3}_{ji},\dots,W^{p}_{ji}\right)^{\top}
  \label{eq:appendix_2d_w_vect_zm}
\end{equation}
\begin{equation}
  \bm{\Psi}^{\bm{\alpha}}_{i}
  =
  \left(\Psi^{\bm{\alpha}}_{i,0},\Psi^{\bm{\alpha}}_{i,1},\Psi^{\bm{\alpha}}_{i,2},\Psi^{\bm{\alpha}}_{i,3},\dots,\Psi^{\bm{\alpha}}_{i,p}\right)^{\top}.
  \label{eq:appendix_2d_coef_vect_zm}
\end{equation}
To complete the 2D formulation, we substitute the weight functions into Eq.~\eqref{eq:appendix_2d_moments_vect_zm},
\begin{equation}
  \bm{B}^{\bm{\alpha}}_{i}
  =
  \displaystyle\sum_{j\in\mathcal{N}_i}\bm{X}_{ji}\bm{W}_{ji}\cdot\bm{\Psi}^{\bm{\alpha}}_{i}\, ,
  \label{eq:appendix_2d_moments_vect_subst_zm}
\end{equation}
which defines a system of linear equations for the weight coefficients $ \bm{\Psi}^{\bm{\alpha}}_{i} $ with system (moment) matrix $ \mathbf{M}_i \in \mathbb{R}^{p\times p}$
\begin{equation}
  \mathbf{M}_{i}\bm{\Psi}^{\bm{\alpha}}_{i}
  =
  \bm{B}^{\bm{\alpha}}_{i}\, ,
  \quad
  \mathbf{M}_{i}
  =
  \displaystyle\sum_{j\in\mathcal{N}_i}
  \bm{X}_{ji}\bm{W}_{ji}^{\top}\, .
  \label{eq:appendix_2d_linSystem_all_moments_zm}
\end{equation}
By specifying the values of the moments, that is, by replacing the vector $ \bm{B}^{\bm{\alpha}}_{i} $  with the vector $ \bm{C}^{\bm{\alpha}} $, we obtain the coefficients $ \bm{\Psi^{\bm{\alpha}}}_{i} $ for the desired derivative approximation by solving
\begin{equation}
  \mathbf{M}_{i}\bm{\Psi^{\bm{\alpha}}}_{i}
  =
  \bm{C}^{\bm{\alpha}}\, .
  \label{eq:appendix_2d_linSystem_zm}
\end{equation}
The discrete operator $ L^{\bm{\alpha}}_{i} $ over the neighborhood of each paricle then becomes
\begin{equation}
  L^{\bm{\alpha}}_{i} u
  =
  \displaystyle\sum_{j\in\mathcal{N}_i}
  u_{j}\bm{W}_{ji}\cdot \left(\mathbf{M}_{i}^{-1} \bm{C}^{\bm{\alpha}} \right)\, .
  \label{eq:appendix_2d_general_discrete_diff_operator_with_coefs_do}
\end{equation}

As a last step, we can optionally scale each row of the linear system with the corresponding power of the inter-particle spacing (spatial resolution) $ h_i $.
This improves the condition number of the system and, in 2D, amounts to preconditioning with the matrix
\begin{equation}
  \mathbf{H}_i
  =
  \text{diag} \left(
    1,\;
    h_{i}^{-1},
    h_{i}^{-1},
    h_{i}^{-2},
    h_{i}^{-2},
    h_{i}^{-2},
    h_{i}^{-3},
    h_{i}^{-3},
    h_{i}^{-3},
    h_{i}^{-3},
    \dots,
    h_{i}^{-m}
  \right)\in \mathbb{R}^{p\times p}
  \label{eq:appendix_2d_h_vect_zm}
\end{equation}
and define $ \widetilde{\bm{X}}_{ji} = \mathbf{H}_i \bm{X}_{ji} \in \mathbb{R}^p $ and $ \widetilde{\bm{C}^{\bm{\alpha}}} = \mathbf{H}_i \bm{C}^{\bm{\alpha}} \in \mathbb{R}^p $.

\subsection{Derivation without the 0th moment\label{sec:appendix_2d_without_zm}}

We repeat the above derivation with the zeroth moment excluded. It is implicitly assumed to be 1.
This leads to the {\em centered} (sometimes called {\em symmetric} or {\em conservative}\footnote{The term ``symmetric'' alludes to the fact that methods derived in this way exactly conserve the sum of $u$ over all particles if the particles are positioned on a regular grid.}) formulation in terms of $u_{ji}$ instead of $u_j$:
\begin{equation}
  L^{\bm{\alpha}}_{i} u
  =
  \displaystyle \sum_{j\in\mathcal{N}_i}
  u_{ji}w^{\bm{\alpha}}_{ji}
  \approx
  \bm{C}^{\bm{\alpha}}\cdot\left.\bm{D} u\right\rvert_{i}
  \label{eq:appendix_2d_general_discrete_diff_operator}
\end{equation}
In this case, the vectors $ \bm{D} u $ and $\bm{X}_{ji}$ omit the constant term and become:
\begin{equation}
  \bm{D} u
  =
  \left(
    \frac{\partial u}{\partial{x}}, \,
    \frac{\partial u}{\partial{y}}, \,
    \frac{\partial^{2} u}{\partial{x}^{2}}, \,
    \frac{\partial^{2} u}{\partial{x}\partial{y}}, \,
    \frac{\partial^{2} u}{\partial{y}^{2}}, \,
    \frac{\partial^{3} u}{\partial{x}^{3}}, \,
    \frac{\partial^{3} u}{\partial{x}^{2}\partial{y}}, \,
    \frac{\partial^{3} u}{\partial{x}\partial{y}^{2}}, \,
    \frac{\partial^{3} u}{\partial{y}^{3}}, \,
    \frac{\partial^{4} u}{\partial{x}^{4}}, \,
    \dots, \,
    \frac{\partial^{m} u}{\partial{y}^{m}}
  \right)^{\!\!\top},
  \label{eq:appendix_2d_derivatives_vect_expanded2D}
\end{equation}
\begin{equation}
  \bm{X}_{ji}
  =
  \left(
    x_{ji}, \,
    y_{ji}, \,
    \frac{x_{ji}^{2}}{2}, \,
    x_{ji}y_{ji}, \,
    \frac{y_{ji}^{2}}{2}, \,
    \frac{x_{ji}^{3}}{6}, \,
    \frac{x_{ji}^{2}y_{ji}}{2}, \,
    \frac{x_{ji}y_{ji}^{2}}{2}, \,
    \frac{y_{ji}^{3}}{6}, \,
    \frac{x_{ji}^{4}}{24}, \,
    \dots \,
    \frac{y_{ji}^{m}}{m!}
  \right)^{\top}.
  \label{eq:appendix_2d_taylor_monomials_vect}
\end{equation}
In this case, $ p =  p_\mathrm{2D} = (m^2 + 3m)/2 $ is one less.

Substituting the Taylor expansion
\begin{equation}
  u_{j}
  =
  u_{i}+\bm{X}_{ji}\cdot\left.\bm{D} u\right\rvert_{i} + e_{ji}^{m+1}.
  \label{eq:appendix_2d_taylor_series}
\end{equation}
without the remainder $e_{ji}^{m+1} $ into Eq.~\eqref{eq:appendix_2d_general_discrete_diff_operator} yields
\begin{equation}
  L^{\bm{\alpha}}_{i} u
  =
  \displaystyle\sum_{j\in\mathcal{N}_i}
  \left.\bm{X}_{ji}\cdot\bm{D} u\right\rvert_{i}w^{\bm{\alpha}}_{ji}\, ,
  \label{eq:appendix_2d_error_short}
\end{equation}
which is the same as in Eq.~\ref{eq:appendix_2d_error_short_zm}, but expands to the slightly different
\begin{multline}
  L^{\bm{\alpha}}_{i} u
  =
  \left.\frac{\partial u}{\partial{x}}\right\rvert_{i}
  \displaystyle\sum_{j\in\mathcal{N}_i}
  x_{ji}w^{\bm{\alpha}}_{ji}
  +\left.\frac{\partial u}{\partial{y}}\right\rvert_{i}
  \displaystyle\sum_{j\in\mathcal{N}_i}
  y_{ji}w^{\bm{\alpha}}_{ji}
  +\left.\frac{\partial^{2} u}{\partial{x}^{2}}\right\rvert_{i}
  \displaystyle\sum_{j\in\mathcal{N}_i}\frac{x_{ji}^{2}}{2}w^{\bm{\alpha}}_{ji}
  \\
  +\left.\frac{\partial^{2} u}{\partial{x}\partial{y}}\right\rvert_{i}
  \displaystyle\sum_{j\in\mathcal{N}_i}
  x_{ji}y_{ji}w^{\bm{\alpha}}_{ji}
  +\left.\frac{\partial^{2} u}{\partial{y}^{2}}\right\rvert_{i}
  \displaystyle\sum_{j\in\mathcal{N}_i}
  \frac{y_{ji}^{2}}{2}w^{\bm{\alpha}}_{ji}
  +\left.\frac{\partial^{3} u}{\partial{x}^{3}}\right\rvert_{i}
  \displaystyle\sum_{j\in\mathcal{N}_i}
  \frac{x_{ji}^{3}}{6}w^{\bm{\alpha}}_{ji}
  +\left.\frac{\partial^{3} u}{\partial{x}^{2}\partial{y}}\right\rvert_{i}
  \displaystyle\sum_{j\in\mathcal{N}_i}
  \frac{x_{ji}^{2}y_{ji}}{2}w^{\bm{\alpha}}_{ji}
  \\
  +\left.\frac{\partial^{3} u}{\partial{x}\partial{y}^{2}}\right\rvert_{i}
  \displaystyle\sum_{j\in\mathcal{N}_i}
  \frac{x_{ji}y_{ji}^{2}}{2}w^{\bm{\alpha}}_{ji}
  +\left.\frac{\partial^{3} u}{\partial{y}^{3}}\right\rvert_{i}
  \displaystyle\sum_{j\in\mathcal{N}_i}
  \frac{y_{ji}^{3}}{6}w^{\bm{\alpha}}_{ji}
  +\left.\frac{\partial^{4} u}{\partial{x}^{4}}\right\rvert_{i}
  \displaystyle\sum_{j\in\mathcal{N}_i}
  \frac{x_{ji}^{4}}{24}w^{\bm{\alpha}}_{ji}
  +\dots
  \label{eq:appendix_2d_error_long}
\end{multline}
As intended, the constant term is absent from the expansion, implicitly assuming a zeroth moment equal to 1.

The rest of the derivation remains the same as in the case with the zeroth moment above, but the discrete operator is now expressed in terms of $u_{ji}$ as:
The operator $ L^{\bm{\alpha}}_{i} $ can thus be expressed as
\begin{equation}
  L^{\bm{\alpha}}_{i} u
  =
  \displaystyle\sum_{j\in\mathcal{N}_i}
  u_{ji}\bm{W}_{ji}\cdot \left(\mathbf{M}_{i}^{-1} \bm{C}^{\bm{\alpha}} \right)\, .
  \label{eq:appendix_2d_general_discrete_diff_operator_with_coefs}
\end{equation}
The preconditioner, if used, shortens accordingly to
\begin{equation}
  \mathbf{H}_i
  =
  \text{diag} \left(
    h_{i}^{-1},
    h_{i}^{-1},
    h_{i}^{-2},
    h_{i}^{-2},
    h_{i}^{-2},
    h_{i}^{-3},
    h_{i}^{-3},
    h_{i}^{-3},
    h_{i}^{-3},
    \dots,
    h_{i}^{-m}
  \right)\, . 
  \label{eq:appendix_2d_h_vect}
\end{equation}
We thus have $ \mathbf{H}_i \in \mathbb{R}^{p\times p} $ and use $ \widetilde{\bm{X}}_{ji} = \mathbf{H}_i \bm{X}_{ji} \in \mathbb{R}^p $ and $ \widetilde{\bm{C}^{\bm{\alpha}}} = \mathbf{H}_i\bm{C}^{\bm{\alpha}} \in \mathbb{R}^p $.

\section{Formal $\ell_2$ error minimization\label{sec:appendix-l2e}}
We solve the optimization problem
\begin{align*}
  \mathcal{J}_{\mathcal{E},i}
  \coloneq
  \displaystyle\sum_{j\in\mathcal{N}_i} \varphi_{ji} \left( e_{ji}^{m+1} \right)^2
  =
  \displaystyle\sum_{j\in\mathcal{N}_i} \varphi_{ji} \left( \bm{X}_{ji} \cdot \left.\bm{D}u\right\rvert_{i} - u_{j} \right)^2 \\
  \mathrm{min} \left( \mathcal{J}_{\mathcal{E},i} \right)
  \quad
  \text{with respect to}
  \left.\bm{D}u\right\rvert_{i}\, .
  \numberthis
  \label{eq:emin-optimizationformulation}
\end{align*}
The partial derivative of the minimization functional is:
\begin{equation}
  \frac{\partial \mathcal{J}_{\mathcal{E},i}}{ \partial \left.\bm{D}u\right\rvert_{i} }
  =
  2 \displaystyle\sum_{j\in\mathcal{N}_i} \varphi_{ji} \left( \bm{X}_{ji} \cdot \left.\bm{D}u\right\rvert_{i} - u_{j} \right)\bm{X}_{ji}\, .
\end{equation}
Setting $\partial \mathcal{J}_i / \partial \left.\bm{D}u\right\rvert_{i} = 0 $ leads to:
\begin{align*}
  \displaystyle\sum_{j\in\mathcal{N}_i} \varphi_{ji} \left( \bm{X}_{ji} \cdot \left.\bm{D}u\right\rvert_{i} - u_{j} \right)\bm{X}_{ji}  &= 0 \\
  \displaystyle\sum_{j\in\mathcal{N}_i} \varphi_{ji} \left( \bm{X}_{ji} \cdot \left.\bm{D}u\right\rvert_{i} \right) \bm{X}_{ji}  &=  \displaystyle\sum_{j\in\mathcal{N}_i} \varphi_{ji}  u_{j} \bm{X}_{ji} \\
  \left(\displaystyle\sum_{j\in\mathcal{N}_i} \varphi_{ji} \bm{X}_{ji}  \bm{X}_{ji}^{\top} \right)\left.\bm{D}u\right\rvert_{i} &=  \displaystyle\sum_{j\in\mathcal{N}_i} \varphi_{ji}  u_{j} \bm{X}_{ji}\, ,
\end{align*}
from which we find the optimal
\begin{align*}
  \left.\bm{D}u\right\rvert_{i}
  &=
  \displaystyle\sum_{j\in\mathcal{N}_i} \varphi_{ji}  u_{j}
  \left(
    \displaystyle\sum_{k} \varphi_{ki} \bm{X}_{ki}  \bm{X}_{ki}^{\top}
  \right)^{-1}
  \bm{X}_{ji}
  \\
  &=
  \displaystyle\sum_{j\in\mathcal{N}_i} \varphi_{ji}  u_{j}
  \mathbf{M}_i^{-1}
  \bm{X}_{ji}.
  \numberthis
  \label{eq:emin-generaldiscreteoperator}
\end{align*}

\subsection{Matrix notation}\label{sec:AppB1}

The minimization problem in Eq.~\eqref{eq:emin-optimizationformulation} can be written in matrix form as:
\begin{align*}
  \mathcal{J}_{\mathcal{E},i}
  =
  \displaystyle\sum_{j\in\mathcal{N}_i} \varphi_{ji} \left(e_{ji}^{m+1} \right)^2
  =
  \bm{E}_{i}^{\top}\mathbf{V}_{i} \bm{E}_{i}
  =
  \left( \mathbf{X}_i \left.\bm{D} u\right\rvert_{i} - \bm{\bm{U}}_{i} \right)^{\top}
  \mathbf{V}_{i}
  \left( \mathbf{X}_i \left.\bm{D} u\right\rvert_{i} - \bm{\bm{U}}_{i} \right)
  \\
  \mathrm{min} \left( \mathcal{J}_{\mathcal{E},i} \right)
  \quad
  \text{with respect to}
  \left.\bm{D} u\right\rvert_{i}\, .
  \numberthis
  \label{eq:appendix:emin-optimizationformulation-matrix}
\end{align*}
Following the same steps as above, we find:
\begin{align*}
  \frac{\partial \mathcal{J}_{\mathcal{E},i}}{ \partial \left.\bm{D} u\right\rvert_{i} }
  \left(
    \left( \mathbf{X}_i \left.\bm{D} u\right\rvert_{i} - \bm{\bm{U}}_{i} \right)^{\top}
    \mathbf{V}_{i}
    \left( \mathbf{X}_i \left.\bm{D} u\right\rvert_{i} - \bm{\bm{U}}_{i} \right)
  \right)
  &= \mathbf{0}
  \\
  \frac{\partial \mathcal{J}_{\mathcal{E},i}}{ \partial \left.\bm{D} u\right\rvert_{i} }
  \left(
    \left.\bm{D} u\right\rvert_{i}^{\top}
    \left( \mathbf{X}_i^{\top} \mathbf{V}_{i} \mathbf{X}_i \right)
    \left.\bm{D} u\right\rvert_{i}
    -
    \left.\bm{D} u\right\rvert_{i}^{\top}
    \mathbf{X}_i^{\top}
    \left( \mathbf{V}_{i} \bm{\bm{U}}_i \right)
    -
    \bm{\bm{U}}_i^{\top}
    \left( \mathbf{V}_{i} \mathbf{X}_i  \left.\bm{D} u\right\rvert_{i} \right)
    +
    \bm{\bm{U}}_i^{\top} \mathbf{V}_{i} \bm{\bm{U}}_i
  \right)
  &= \mathbf{0}
  \\
  2 \left( \mathbf{X}_i^{\top} \mathbf{V}_{i} \mathbf{X}_i \right)
  \left.\bm{D} u\right\rvert_{i}
  -
  \mathbf{X}_i^{\top} \mathbf{V}_{i} \bm{\bm{U}}_i
  -
  \left(
    \bm{\bm{U}}_i^{\top}
    \mathbf{V}_{i} \mathbf{X}_i
  \right)^{\top}
  &= \mathbf{0}
  \\
  2 \left( \mathbf{X}_i^{\top} \mathbf{V}_{i} \mathbf{X}_i \right)
  \left.\bm{D} u\right\rvert_{i}
  -
  2\mathbf{X}_i^{\top} \mathbf{V}_{i} \bm{\bm{U}}_i
  &= \mathbf{0}\, ,
\end{align*}
from which the optimal solution follows as:
\begin{equation}
  \left.\bm{D} u\right\rvert_{i}
  =
  \left[
    \left( \mathbf{X}_i^{\top} \mathbf{V}_{i} \mathbf{X}_i \right)^{-1} \mathbf{X}_i^{\top} \mathbf{V}_{i}
  \right]
  \bm{\bm{U}}_i \, .
\end{equation}

\section{Formal $\ell_2$ minimization of weights\label{sec:appendix-l2p}}
We solve the optimization problem
\begin{align*}
  \mathcal{J}_{\mathcal{P},i}
  =
  \frac{1}{2}
  \displaystyle\sum_{j\in\mathcal{N}_i} \frac{ 1 }{ \varphi_{ji} }\left(w_{ji}^{\bm{\alpha}} \right)^{2} \\
  \mathrm{min} \left( \mathcal{J}_{\mathcal{P},i} \right)
  \quad
  \text{subject to}
  \quad
  \displaystyle\sum_{j\in\mathcal{N}_i}
  \bm{X}_{ji}w_{ji}^{\bm{\alpha}}
  =
  \bm{B}_{i}^{\bm{\alpha}}
  =
  \bm{C}^{\bm{\alpha}}\, .
  \numberthis
  \label{eq:mc-optimizationformulation-withmoments}
\end{align*}

\nomvar{ $\mathcal{J}_{\mathcal{P},i}$ }{functional for $\ell_2$ minimization under the moment condition constraint [-]}
\nomvar{ $\bm{\lambda}$ }{Lagrange multipliers [-]}

The Lagrange function for this problem is
\begin{equation}
  \mathcal{L}\left( \bm{w}_{\bullet i}^{\bm{\alpha}}, \bm{\lambda} \right)
  =
  \frac{1}{2}
  \displaystyle\sum_{j\in\mathcal{N}_i}
  \frac{1}{\varphi_{ji}}
  \left(
    w_{ji}^{\bm{\alpha}}
  \right)^2
  +
  \bm{\lambda}^{\!\top}
  \left(
    \displaystyle\sum_{j\in\mathcal{N}_i}
    \bm{X}_{ji}w_{ji}^{\bm{\alpha}}
    - \bm{C}^{\bm{\alpha}}
  \right).
  \label{eq:mc_lagrangeFunction}
\end{equation}
for which we find the derivative
\begin{equation}
  \frac{\partial \mathcal{L}\left( \bm{w}_{\bullet i}^{\bm{\alpha}}, \bm{\lambda} \right)}{\partial w_{ki}^{\bm{\alpha}}}
  =
  \displaystyle\sum_{j\in\mathcal{N}_i}
  \left(
    \frac{1}{\varphi_{ji}}
  \right)
  w_{ji}^{\bm{\alpha}}
  \delta_{jk}
  +
  \bm{\lambda}^{\!\top}
  \left(
    \displaystyle\sum_{j\in\mathcal{N}_i}
    \bm{X}_{ji}
    \delta_{jk}
  \right)
  \quad
  \forall k \in \mathcal{N}_i \, .
\end{equation}
Setting $\partial \mathcal{L}\left( \bm{w}_{\bullet i}^{\bm{\alpha}}, \bm{\lambda} \right)/ \partial  w_{ki}^{\bm{\alpha}} = 0, \ \forall k \in \mathcal{N}_i $, leads to:
\begin{align*}
  \left(
    \frac{1}{\varphi_{ji}}
  \right)
  w_{ji}^{\bm{\alpha}}
  +
  \bm{\lambda}^{\!\top}
  \bm{X}_{ji}
  &= 0
  \quad
  \forall j \in \mathcal{N}_i
  \\
  w_{ji}^{\bm{\alpha}}
  &=
  -
  \varphi_{ji}
  \bm{\lambda}^{\!\top}
  \bm{X}_{ji}\, .
  \numberthis
  \label{eq:momentumMethod_weightFunctionsTemp}
\end{align*}
Substituting this into the constraint, we find:
\begin{equation}
  \displaystyle\sum_{j\in\mathcal{N}_i} \bm{X}_{ji}  w_{ji}^{\bm{\alpha}}
  =
  - \displaystyle\sum_{j\in\mathcal{N}_i} \bm{X}_{ji}
  \varphi_{ji}
  \bm{\lambda}^{\!\top}
  \bm{X}_{ji}
  =
  - \displaystyle\sum_{j\in\mathcal{N}_i} \bm{X}_{ji}
  \varphi_{ji}
  \bm{X}_{ji}^{\top}
  \bm{\lambda}
  =
  \bm{C}^{\bm{\alpha}}.
\end{equation}
Using $ \bm{W}_{ji} = \bm{X}_{ji} \varphi_{ji} $, we recognize the moment matrix $ \mathbf{M}_i $ defined in Eq.~\eqref{eq:linSystem_all_moments_zm} and hence write:
\begin{equation}
  - \displaystyle\sum_{j\in\mathcal{N}_i} \bm{X}_{ji}
  \varphi_{ji}
  \bm{X}_{ji}^{\top}
  \bm{\lambda}
  =
  - \displaystyle\sum_{j\in\mathcal{N}_i} \bm{W}_{ji}
  \bm{X}_{ji}^{\top}
  \bm{\lambda}
  =
  \mathbf{M}_i
  \bm{\lambda}
  =
  \bm{C}^{\bm{\alpha}}.
  \label{eq:momentummethod_temp_constrains}
\end{equation}
Solving for $ \bm{\lambda}$ yields:
\begin{equation}
  \bm{\lambda} = - \mathbf{M}_i^{-1} \bm{C}^{\bm{\alpha}}\, ,
\end{equation}
which we substitute back into Eq.~\eqref{eq:momentumMethod_weightFunctionsTemp} to find the optimal solution:
\begin{equation}
  w_{ji}^{\bm{\alpha}}
  =
  \varphi_{ji}
  \bm{X}_{ji}
  \cdot
  \mathbf{M}_i^{-1} \bm{C}^{\bm{\alpha}},
  \qquad
  \mathbf{M}_i
  =
  \displaystyle\sum_{j\in\mathcal{N}_i} \bm{X}_{ji} \bm{X}_{ji}^{\top} \varphi_{ji}\, .
\end{equation}

\subsection{Matrix notation}

In matrix notation, the optimization problem from Eq.~\eqref{eq:mc-optimizationformulation-withmoments} reads:
\begin{align*}
  \mathcal{J}_{\mathcal{P},i}
  =
  \left\lVert \left(\mathbf{V}_i \right)^{-1/2} \bm{w}_{\bullet i}^{\bm{\alpha}} \right\rVert^2  \\
  \mathrm{min} \left( \mathcal{J}_{\mathcal{P},i} \right)
  \quad
  \text{subject to}
  \quad
  \mathbf{X}^{\top}\bm{w}_{\bullet i}^{\bm{\alpha}} = \bm{C}^{\bm{\alpha}},
  \numberthis
  \label{eq:appendix:app_matrices_polynomial_opt}
\end{align*}
where $ \left(\mathbf{V}_i  \right)^{-1/2} $ is easy to compute, since $ \mathbf{V}_{ji}$ is a diagonal matrix.
In this notation, the Lagrange function from Eq.~\eqref{eq:mc_lagrangeFunction} becomes:
\begin{align*}
  \mathcal{L}\left( \bm{w}_{\bullet i}^{\bm{\alpha}}, \bm{\lambda} \right)
  &=
  \frac{1}{2}
  \left( \bm{w}_{\bullet i}^{\bm{\alpha}} \right)^{\top}
  \left(\left( \mathbf{V}_{i} \right)^{-1/2}\right)^{\top}
  \left( \mathbf{V}_{i} \right)^{-1/2}
  \bm{w}_{\bullet i}^{\bm{\alpha}}
  +
  \bm{\lambda}^{\top}
  \left(
    \mathbf{X}_i^{\top}  \bm{w}_{\bullet i}^{\bm{\alpha}}  - \bm{C}^{\bm{\alpha}}
  \right)
\end{align*}
With this, we formally solve the minimization problem as:
\begin{equation}
  \frac{\partial \mathcal{L}\left( \bm{w}_{\bullet i}^{\bm{\alpha}}, \bm{\lambda} \right) }{ \partial \bm{w}_{\bullet i}^{\bm{\alpha}} }
  =
  (\mathbf{V}_i)^{-1} \bm{w}_{\bullet i}^{\bm{\alpha}} + \mathbf{X}_i \bm{\lambda}
  =
  \bm{0}\, ,
\end{equation}
\begin{equation}
  \bm{w}_{\bullet i}^{\bm{\alpha}} = -\mathbf{V}_{i} \mathbf{X}_i \bm{\lambda}\, .
\end{equation}
This vector of weight functions $\bm{\omega}_{\bullet i}^{\bm{\alpha}} $ is then  substituted into the constraint in Eq.~\eqref{eq:appendix:app_matrices_polynomial_opt} to find:
\begin{equation}
  - \mathbf{X}_i^{\top} \mathbf{V}_{i} \mathbf{X}_i \bm{\lambda} = \bm{C}^{\bm{\alpha}}.
  \label{eq:mcmatrix-wdvectorwithlambda}
\end{equation}
This allows us to obtain
\begin{equation}
  \bm{\lambda}
  =
  - \left(
    \mathbf{X}_i^{\top}
    \mathbf{V}_{i}
    \mathbf{X}_i
  \right)^{-1} \bm{C}^{\bm{\alpha}}
  =
  - \left(
    \mathbf{X}_i^{\top}
    \mathbf{W}_i
  \right)^{-1} \bm{C}^{\bm{\alpha}}
  =
  - \mathbf{M}_{i}^{-1}\bm{C}^{\bm{\alpha}}\, ,
  \label{eq:mcmatrix-lambdas}
\end{equation}
where we have used that $ \mathbf{W}_i = \mathbf{V}_i \mathbf{X}_i = \mathbf{X}_i \mathbf{V}_i $ for the moment matrix from Eq.~\eqref{eq:momentum_matrix_matrix_form}.
Substituting Eq.~\eqref{eq:mcmatrix-lambdas} back into Eq.~\eqref{eq:mcmatrix-wdvectorwithlambda} results in a vector containing weight for all neighboring particles:
\begin{equation}
  \bm{w}_{\bullet i}^{\bm{\alpha}}
  =
  \mathbf{V}_i  \mathbf{X}_i \mathbf{M}_i^{-1} \bm{C}^{\bm{\alpha}}
  =
  \mathbf{W}_i \mathbf{M}_i^{-1} \bm{C}^{\bm{\alpha}}.
\end{equation}
Expanding the matrix multiplication corresponding to the summation of the neighboring particle, we can see that this is the same as in the original derivation in Section~\ref{sec:l2p},
\begin{equation}
  L^{\bm{\alpha}}_{i} u
  =
  \bm{\bm{U}}_i \cdot \bm{w}_{\bullet i}^{\bm{\alpha}}
  =
  \bm{\bm{U}}_i \cdot \left( \mathbf{C}_i\bm{C}^{\bm{\alpha}} \right)
  =
  \bm{\bm{U}}_i \cdot \left( \mathbf{W}_i \mathbf{Y}_{i}\bm{C}^{\bm{\alpha}} \right)
  =
  \bm{\bm{U}}_i \cdot \left( \mathbf{W}_i \mathbf{M}_{i}^{-1}\bm{C}^{\bm{\alpha}} \right).
  \label{eq:general_discrete_diff_operator_with_coefs_do_matrix_form_l2min}
\end{equation}

\section{Formal generalized $\ell_2$ minimization\label{sec:appendix-gl2p}}
We solve optimization problem
\begin{align*}
  \widetilde{\mathcal{J}}_{\mathcal{P},i}
  =
  \frac{1}{2}
  \displaystyle\sum_{j\in\mathcal{N}_i}
  \left(
    w_{ji}^{\bm{\alpha}}
  \right)^2
  =
  \frac{1}{2}
  \displaystyle\sum_{j\in\mathcal{N}_i}
  \left(
    \bm{W}_{ji} \cdot \bm{\Psi}_{i}^{\bm{\alpha}}
  \right)^2
  \\
  \mathrm{min} \left( \widetilde{\mathcal{J}}_{\mathcal{P},i} \right)
  \quad
  \text{subject to}
  \quad
  \displaystyle\sum_{j\in\mathcal{N}_i}
  \bm{X}_{ji}w_{ji}^{\bm{\alpha}}
  =
  \displaystyle\sum_{j\in\mathcal{N}_i}
  \bm{X}_{ji}\bm{W}_{ji} \cdot \bm{\Psi}_{i}^{\bm{\alpha}}
  =
  \mathbf{M}_i\bm{\Psi}_{i}^{\bm{\alpha}}
  =
  \bm{C}^{\bm{\alpha}}.
  \numberthis
  \label{eq:mcgeneralized-optimizationformulation-mc2}
\end{align*}

\nomvar{ $\mathcal{J}_{\mathcal{P},i}$ }{generalized functional for $\ell_2$ minimization under the moment condition constraint [-]}

The Lagrange function for this problem is:
\begin{align*}
  \mathcal{L}\left( \bm{\Psi}_{i}^{\bm{\alpha}}, \bm{\lambda} \right)
  &=
  \frac{1}{2}
  \displaystyle\sum_{j\in\mathcal{N}_i}
  \left(
    \bm{W}_{ji} \cdot \bm{\Psi}_{i}^{\bm{\alpha}}
  \right)^2
  +
  \bm{\lambda}^{\top}
  \left(
    \displaystyle\sum_{j\in\mathcal{N}_i}
    \bm{X}_{ji}\bm{W}_{ji} \cdot \bm{\Psi}_{i}^{\bm{\alpha}}
    - \bm{C}^{\bm{\alpha}}
  \right)
  \\
  &=
  \frac{1}{2}
  \displaystyle\sum_{j\in\mathcal{N}_i}
  \left(
    \bm{W}_{ji} \cdot \bm{\Psi}_{i}^{\bm{\alpha}}
  \right)^2
  +
  \bm{\lambda}^{\top}
  \left(
    \mathbf{M}_i \bm{\Psi}_{i}^{\bm{\alpha}}
    - \bm{C}^{\bm{\alpha}}
  \right)\, ,
  \numberthis
  \label{eq:l2mcg_lagrangeFunction}
\end{align*}
for which we find the derivative
\begin{equation}
  \frac{\partial \mathcal{L}\left( \bm{\Psi}_{i}^{\bm{\alpha}}, \bm{\lambda} \right)}{\partial \bm{\Psi}_{i}^{\bm{\alpha}}}
  =
  \displaystyle\sum_{j\in\mathcal{N}_i}
  \left(
    \bm{W}_{ji} \cdot \bm{\Psi}_{i}^{\bm{\alpha}}
  \right)
  \bm{W}_{ji}
  +
  \mathbf{M}_{i}^{\top}\bm{\lambda}\, .
  \label{eq:l2mcg_derivativeOfLagrangeFunction}
\end{equation}
Setting $\partial \mathcal{L}\left( \bm{\Psi}_{i}^{\bm{\alpha}}, \bm{\lambda} \right)/ \partial  \bm{\Psi}_{i}^{\bm{\alpha}} = 0 $ leads to
\begin{align*}
  \displaystyle\sum_{j\in\mathcal{N}_i}
  \left(
    \bm{W}_{ji} \cdot \bm{\Psi}_{i}^{\bm{\alpha}}
  \right)
  \bm{W}_{ji}
  &=
  - \mathbf{M}_{i}^{\top}\bm{\lambda} \\
  \displaystyle\sum_{j\in\mathcal{N}_i}
  \bm{W}_{ji}
  \bm{W}_{ji}^{\top} \bm{\Psi}_{i}^{\bm{\alpha}}
  &=   - \mathbf{M}_{i}^{\top}\bm{\lambda}
  \\
  \bm{\Psi}_{i}^{\bm{\alpha}}
  &=
  -
  \left(
    \displaystyle\sum_{j\in\mathcal{N}_i}
    \bm{W}_{ji}
    \bm{W}_{ji}^{\top}
  \right)^{-1}
  \mathbf{M}_{i}^{\top}\bm{\lambda}\, .
\end{align*}

Using the constraint from Eq.~\eqref{eq:mcgeneralized-optimizationformulation-mc2}, we substitute for $w_{ji}^{\bm{\alpha}}$,
\begin{equation}
  - \displaystyle\sum_{j\in\mathcal{N}_i}
  \bm{X}_{ji}\bm{W}_{ji}^{\top}
  \left(
    \displaystyle\sum_{j\in\mathcal{N}_i}
    \bm{W}_{ji}
    \bm{W}_{ji}^{\top}
  \right)^{-1}
  \mathbf{M}_{i}^{\top}\bm{\lambda}
  =
  \bm{C}^{\bm{\alpha}}\, .
\end{equation}
This allows expressing the Lagrange multipliers $ \bm{\lambda} $ as:
\begin{equation}
  \bm{\lambda}
  =
  - \mathbf{M}_i^{-T}
  \left(
    \displaystyle\sum_{j\in\mathcal{N}_i}
    \bm{W}_{ji}
    \bm{W}_{ji}^{\top}
  \right)
  \mathbf{M}_i^{-1}
  \bm{C}^{\bm{\alpha}}\, .
\end{equation}
Substituting this back into the expression for the weights, we find:
\begin{align*}
  w_{ji}^{\bm{\alpha}}
  &=
  \bm{W}_{ji} \cdot \bm{\Psi}_{i}^{\bm{\alpha}}
  \\
  &=
  -
  \bm{W}_{ji}^{\top}
  \left(
    \displaystyle\sum_{j\in\mathcal{N}_i}
    \bm{W}_{ji}
    \bm{W}_{ji}^{\top}
  \right)^{-1}
  \mathbf{M}_{i}^{\top}\bm{\lambda}
  \\
  &=
  \bm{W}_{ji}^{\top}
  \left(
    \displaystyle\sum_{j\in\mathcal{N}_i}
    \bm{W}_{ji}
    \bm{W}_{ji}^{\top}
  \right)^{-1}
  \mathbf{M}_{i}^{\top}
  \mathbf{M}_i^{-T}
  \left(
    \displaystyle\sum_{j\in\mathcal{N}_i}
    \bm{W}_{ji}
    \bm{W}_{ji}^{\top}
  \right)
  \mathbf{M}_i^{-1}
  \bm{C}^{\bm{\alpha}}
  \\
  &=
  \bm{W}_{ji}
  \cdot
  \left(
    \mathbf{M}_i^{-1}
    \bm{C}^{\bm{\alpha}}
  \right)\, .
\end{align*}
This leads to the optimal solution
\begin{equation}
  L^{\bm{\alpha}}_{i}u
  =
  \displaystyle\sum_{j\in\mathcal{N}_i}
  u_{ji}\bm{W}_{ji}\cdot \left(\mathbf{M}_{i}^{-1} \bm{C}^{\bm{\alpha}} \right)\, ,
  \qquad
  \mathbf{M}_i
  =
  \displaystyle \sum_{j\in\mathcal{N}_i}
  \bm{X}_{ji} \bm{W}_{ji}^{\top} \, .
  \label{eq:appendix:l2mcg-general_discrete_diff_operator_with_coefs}
\end{equation}

\subsection{Matrix notation}

In the matrix notation introduced in Section~\ref{sec:basic-approach-matrix-formulation}, the optimization problem from Eq.~\eqref{eq:mcgeneralized-optimizationformulation-mc2} reads:
\begin{align*}
  \widetilde{\mathcal{J}}_{\mathcal{P},i}
  =
  \frac{ 1 }{ 2 }
  \lVert \bm{w}_{i}^{\bm{\alpha}} \rVert^2
  =
  \frac{ 1 }{ 2 }
  \lVert \mathbf{W}_i \bm{\Psi}_{i}^{\bm{\alpha}} \rVert^2  \\
  \mathrm{min} \left( \widetilde{\mathcal{J}}_{\mathcal{P},i} \right)
  \quad
  \text{subject to}
  \quad
  \mathbf{M}_i\bm{\Psi}_i^{\bm{\alpha}} = \bm{C}^{\bm{\alpha}}
  \numberthis
  \label{eq:appendix:app_matrices_polynomial_opt_generalized}
\end{align*}
The matrix form of the Lagrange function and its derivative are:
\begin{align*}
  \mathcal{L}\left( \bm{\Psi}_{i}^{\bm{\alpha}}, \bm{\lambda} \right)
  &=
  \frac{1}{2}
  \left( \bm{\Psi}_{i}^{\bm{\alpha}} \right)^{\top}
  \left( \mathbf{W}_{i} \right)^{\top}
  \left( \mathbf{W}_{i} \right)
  \bm{\Psi}_{i}^{\bm{\alpha}}
  +
  \bm{\lambda}^{\top}
  \left(
    \mathbf{M}_i  \bm{\Psi}_{i}^{\bm{\alpha}}  - \bm{C}^{\bm{\alpha}}
  \right) \, ,
  \\
  \frac{\partial \mathcal{L}\left( \bm{\Psi}_{i}^{\bm{\alpha}}, \bm{\lambda} \right) }{ \partial \bm{\Psi}_{i}^{\bm{\alpha}} }
  &=
  (\mathbf{W}_i)^{\top}\mathbf{W}_i \bm{\Psi}_{i}^{\bm{\alpha}} + \mathbf{M}_i^{\top} \bm{\lambda}
  =
  \bm{0} \, ,\numberthis
  \\
  \bm{\Psi}_i^{\bm{\alpha}}
  &=
  - \left(
    \mathbf{W}_{i}^{\top} \mathbf{W}_i
  \right)^{-1}
  \mathbf{M}_i^{\top} \bm{\lambda}\, . \numberthis
\end{align*}
Substituting the weight coefficients $ \bm{\Psi}_{i}^{\bm{\alpha}} $ in the constraint in Eq.~\eqref{eq:appendix:app_matrices_polynomial_opt_generalized} yields:
\begin{equation}
  - \mathbf{M}_i \left( \mathbf{W}_{i}^{\top} \mathbf{W}_i \right)^{-1} \mathbf{M}_i^{\top} \bm{\lambda}
  =
  \bm{C}^{\bm{\alpha}}\, .
  \numberthis \label{eq:l2mcg_matrix_weight_temp_relation}
\end{equation}
From this, we find
\begin{equation}
  \bm{\lambda}
  =
  -
  \mathbf{M}_i^{-T} \left( \mathbf{W}_{i}^{\top} \mathbf{W}_i \right) \mathbf{M}_i^{-1}
  \bm{C}^{\bm{\alpha}}\, .
\end{equation}
Subsituting this back into the expression for the weights,
\begin{align*}
  \bm{w}_{\bullet i}^{\bm{\alpha}}
  &=
  \mathbf{W}_i \bm{\Psi}_{i}^{\bm{\alpha}}
  \\
  &=
  \mathbf{W}_i
  \left(
    \mathbf{W}_{i}^{\top} \mathbf{W}_i
  \right)^{-1}
  \mathbf{M}_i^{\top} \bm{\lambda}
  \\
  &=
  \mathbf{W}_i
  \left(
    \mathbf{W}_{i}^{\top} \mathbf{W}_i
  \right)^{-1}
  \mathbf{M}_i^{\top}
  \mathbf{M}_i^{-T} \left( \mathbf{W}_{i}^{\top} \mathbf{W}_i \right) \mathbf{M}_i^{-1}
  \bm{C}^{\bm{\alpha}}\, .
  \\
\end{align*}
In contrast to the error or weights minimization approaches, the moment matrix here has a structure as in the definition in Eq.~\eqref{eq:momentum_matrix_matrix_form}, without any particular form of anisotropic basis functions $ \mathbf{W}_i $ (or $ \bm{W}_{ji} $ ). Therefore, the resulting discrete operator becomes:
\begin{equation}
  \mathbf{M}_{i}
  =
  \displaystyle\sum_{j\in\mathcal{N}_i}
  \bm{X}_{ji}\bm{W}_{ji}^{\top}
  =
  \mathbf{X}_i^{\top} \mathbf{W}_i\, ,
\end{equation}
\begin{equation}
  L^{\bm{\alpha}}_{i} u
  =
  \bm{\bm{U}}_i \cdot \bm{w}_{\bullet i}^{\bm{\alpha}}
  =
  \bm{\bm{U}}_i \cdot \left( \mathbf{C}_i\bm{C}^{\bm{\alpha}} \right)
  =
  \bm{\bm{U}}_i \cdot \left( \mathbf{W}_i \mathbf{Y}_{i}\bm{C}^{\bm{\alpha}} \right)
  =
  \bm{\bm{U}}_i \cdot \left( \mathbf{W}_i \mathbf{M}_{i}^{-1}\bm{C}^{\bm{\alpha}} \right)\, .
\end{equation}

\end{document}


\fancyhead[L]{An Overview of Meshfree Collocation Methods}
\fancyhead[R]{Supplementary Material}
\pagestyle{fancy}

\begin{center}
  \vspace*{0.5cm}
  {\LARGE\bfseries An Overview of Meshfree Collocation Methods\\ --- \\Supplementary Material}\\[0.5cm]
  {\normalsize
    \begin{tabular}{c}
      Tomas Halada\textsuperscript{1,4}, Serhii Yaskovets\textsuperscript{1,2,3}, Abhinav Singh\textsuperscript{1,2,3}, Ludek Benes\textsuperscript{4}, \\
      Pratik Suchde\textsuperscript{5,6,\dag}, and Ivo~F.~Sbalzarini\textsuperscript{1,2,3,\dag,*}
    \end{tabular}
  }\\[0.3cm]
  {\small
    \begin{tabular}{c}
      \textsuperscript{1}Dresden University of Technology, Faculty of Computer Science, 01187 Dresden, Germany \\
      \textsuperscript{2}Max Planck Institute of Molecular Cell Biology and Genetics, 01307 Dresden, Germany \\
      \textsuperscript{3}Center for Systems Biology Dresden, 01307 Dresden, Germany \\
      \textsuperscript{4}Czech Technical University, Faculty of Mechanical Engineering, Department of Technical Mathematics, \\
      160\,00 Prague, Czech Republic \\
      \textsuperscript{5}University of Luxembourg, Faculty of Science, Technology and Medicine, 4365 Esch-sur-Alzette, Luxembourg \\
      \textsuperscript{6}Fraunhofer Institute for Industrial Mathematics, 67663 Kaiserslautern, Germany \\
      \\
      \textsuperscript{\dag}Joint senior authors \\
      \textsuperscript{*}Corresponding author {\footnotesize \url{sbalzarini@mpi-cbg.de}}
    \end{tabular}
  }
\end{center}

\tableofcontents

\section{Encyclopedia of meshfree collocation methods\label{sec:com}}

In this Supplement, we provide an encyclopedia of all reviewed meshfree collocation methods for convenient reference.
In this, we also present the specific features and approaches of each method and point toward applications where suitable. This Supplement is intended for information look-up, and not for continuous reading.

\subsection{Local Anistropic Basis Function Method (LABFM)\label{sec:com_labfm}}

The Local Anistropic Basis Function Method (LABFM) introduced by King in \cite{HighOrderDiffKing2020, HighOrderSimuKing2022} is one of the newest additions to the family of finite-difference methods in scattered nodes.
A discrete approximation of the differential operator in \textbf{LABFM originates from the direct approximation of moments 2} with the excluded 0th moment and is given as follows.
\begin{equation}
  L^{\bm{\alpha}}_{i} u
  \coloneq
  \displaystyle\sum_{j\in\mathcal{N}_i}
  u_{ji}\bm{W}_{ji}\cdot \left(\mathbf{M}_{i}^{-1} \widetilde{\bm{C}^{\bm{\alpha}}} \right),
  \quad
  \mathbf{M}_{i}
  \coloneq
  \displaystyle\sum_{j\in\mathcal{N}_i}
  \widetilde{\bm{X}}_{ji}\bm{W}_{ji}^{\top}.
  \label{eq:labfm_discrete_diff_operator}
\end{equation}
Originally, in \cite{HighOrderDiffKing2020}, the method was formulated without scaling $\bm{X}_{ji} $ and $ \bm{C}^{\bm{\alpha}} $ with stencil sizes $h$, which was included two years later in \cite{HighOrderSimuKing2022}.

Compared to similar methods, the introduction to LABFM emphasizes the choice of anisotropic basis functions $ \bm{W}_{ji} $.
Since the elements of $ \bm{X}_{ji} $ scale with $ h_i, h_i^2, h_i^3,... $, we need to ensure that corresponding ABFs scales with $ h_i^{-1}, h_i^{-2}, h_i^{-3},... $ also.
Originally, in \cite{HighOrderDiffKing2020} the basis functions $ \bm{W}_{ji} $ are taken as partial derivatives of some suitably scaled RBF $ \varphi_{ji} \coloneq \varphi\left( | \vect{x}_{ji} |, h \right) $ with respect to $ \vect{x}_j $
\begin{equation}
  \bm{W}_{ji}
  =
  \left.\bm{D}\varphi_{ji}\right\rvert_{\vect{x}_{j}}.  
  \label{eq:labfm_w_abf_from_rbf}
\end{equation}
However, two years later, the authors of the original work \cite{HighOrderSimuKing2022} published improvements over the original choice \eqref{eq:labfm_w_abf_from_rbf}.
Although the original approach works well in 2D, it is unable to ensure linear independence of the rows of the moment matrix $ \mathbf{M}_i $ in general.
Based on \cite{HighOrderSimuKing2022}, this eventually breaks down into three dimensions, of order $ m = 4 $, where regardless of the node distribution, $ \mathbf{M}_i $ has no full rank.
In \cite{HighOrderSimuKing2022} they propose to construct ABFs using orthogonal polynomials (Hermite or Legendre polynomials) multiplied by a suitable RBF $ \varphi_{ji} $.
Thus, ABFs can be written as
\begin{equation}
  \bm{W}_{ji}
  =
  \bm{P}_{ji} \varphi_{ji}
  \label{eq:labfm_w_abf_from_base}
\end{equation}
where $ \bm{P}_{ji} $ are suitable basis functions (preferably orthogonal).
Several choices of $ \bm{P}_{ji} $ and $ \varphi_{ji} $ were investigated at \cite{HighOrderDiffKing2020, HighOrderSimuKing2022}.

\subsection{High-order Consistent SPH (HOCSPH)\label{sec:com_hocsph}}

Around the same time as LABFM, \textbf{High-order Consistent SPH (HOCSPH)} introduced by Nasar et al. \cite{nasar2021hoch}, which is very similar to LABFM, sharing even some of the authors.
Having roots in SPH, the method in contrast to LABFM assumes a volume defined for each particle, which corresponds to direct node integration, and also makes the HOCSPH effectively a collocation scheme.
Besides that, the strategy is the same.
\textbf{HOCSPH originates from the direct approximation of moments 2} with the excluded 0th moment.
In the notation introduced above, ABFs in HOCSPH are obtained as derivatives of SPH weight function which corresponds to \eqref{eq:labfm_w_abf_from_rbf}.

Considering the particle volume $v_i$ associated with particle $i$, with our notation, the general discrete operator \eqref{eq:ba:ddo} of HOCSPH takes from
\begin{equation}
  L^{\bm{\alpha}}_{i}u
  \coloneq
  \displaystyle \sum_{j\in\mathcal{N}_i}
  u_{ji}w^{\bm{\alpha}}_{ji}v_j
  \label{eq:appendix:variations:hochs ddo}.
\end{equation}
and -- following approximation of moments 2 -- the resulting approximation is
\begin{equation}
  L^{\bm{\alpha}}_{i} u
  \coloneq
  \displaystyle\sum_{j\in\mathcal{N}_i}
  u_{ji}\bm{W}_{ji}\cdot \left(\mathbf{M}_{i}^{-1} \bm{C}^{\bm{\alpha}} \right)v_j,
  \quad
  \mathbf{M}_{i}
  \coloneq
  \displaystyle\sum_{j\in\mathcal{N}_i}
  \bm{X}_{ji}\bm{W}_{ji}^{\top}v_j.
  \label{eq:appendix:variations:hochs ddo approx}
\end{equation}
As we show the in section~\ref{sec:com_dcpse} describing DC-PSE method, considering direct nodal integration, the volume can be absorbed by the weight functions.

\subsection{Discretization-Corrected Particle Strength Exchange Method (DC-PSE)\label{sec:com_dcpse}}

Discretization-Corrected Particle Strength Exchange Method (DC-PSE) is a generalization of the integral particle strength exchange (PSE) operator approach to approximate spatial derivatives of any degree \cite{degond1989, eldredge2002}. \textbf{It originates from the direct approximation of moments 2.}
The method was introduced by Schrader et al. \cite{schrader2010} and later was presented in \cite{UsingDcPseOpBouran2016, AMeshfreeCollSingh2023, RecoveryByDisZwick2023, EntropicallyDaSingh2023}.
In its original formulation, the discrete approximation of differential operator is introduced as
\begin{equation}
  L^{\bm{\alpha}}_{i} u
  \coloneq
  \frac{1}{h^{|\bm{\alpha}|}}
  \displaystyle\sum_{j\in\mathcal{N}_i}
  \left( u_{j}\pm  u_{i}\right)w^{\bm{\alpha}}_{ji}v_j
  \approx
  \bm{C}^{\bm{\alpha}} \cdot \left.\bm{D} u\right\rvert_{i}
  \label{eq:dcpse_discrete_diff_operator_with_volume}
\end{equation}
with weights $ \widetilde{w}^{\bm{\alpha}}_{ji} \coloneq w_{ji}^{\bm{\alpha}}\left( \left( \mathbf{x}_{j} - \mathbf{x}_{i}  \right)/ h \right) $ and $ v_j $ representing volume of particle $ j $ from the quadrature of PSE.
The sign is positive for odd derivatives (in case $ |\bm{\alpha}| $ is odd) and negative for even derivatives (in case $ |\bm{\alpha}| $ is even), however as mentioned in \cite{UsingDcPseOpBouran2016}, form exactly corresponding to \eqref{eq:general_discrete_diff_operator_zm} can also be taken as the default.
Later the volume $ v_j $ was omitted since it was relic of PSE and the volume and the normalization factor, which is usually hidden in the weight functions, can be absorbed into the normalization term (and weight coefficients).
As it is clear from derivation 2, volume is not part of the formulas.
The common form to present DC-PSE is thus
\begin{equation}
  L^{\bm{\alpha}}_{i} u
  =
  \frac{1}{h_{i}^{|\bm{\alpha}|}}
  \displaystyle\sum_{j\in\mathcal{N}_i}
  \left( u_{j}\pm  u_{i}\right)w^{\bm{\alpha}}_{ji}
  \approx
  \bm{C}^{\bm{\alpha}} \cdot \left.\bm{D} u\right\rvert_{i}
  \label{eq:dcpse_discrete_diff_operator}
\end{equation}
and it is clear that the use of factor $1/h^{|\bm{\alpha}|}$ and normalized weight function arguments is equivalent with stencil-size normalisation introduced at the end of section 2.

The feature highlighted in DC-PSE compared to other methods of this kind is the work with zero moment.
Using the Taylor expansion with normalized Taylor monomials $ \widetilde{\bm{X}}_{ji} \coloneq \bm{X}(\mathbf{x}_{ji}/h) = \mathbf{H}_{i}\bm{X}_{ji} $
\begin{equation}
  L^{\bm{\alpha}}_{i} u
  =
  \frac{1}{h_{i}^{|\bm{\alpha}|}}
  \displaystyle\sum_{j\in\mathcal{N}_i}
  \left(
    \left.\widetilde{\bm{X}}_{ji}\cdot\bm{D} u\right\rvert_{i}
    \pm u_{i}
  \right) w^{\bm{\alpha}}_{ji},
  \label{eq:dcpse_error_short}
\end{equation}
we can see, that the zeroth moment is automatically canceling out for $|\bm{\alpha}|$ odd, while stays present for $|\bm{\alpha}|$ even.
A slightly different approach is taken in \cite{TaylorSeriesEJacque2020} where the DC-PSE approximation is assumed in form
\begin{equation}
  L^{\bm{\alpha}}_{i} u
  =
  \displaystyle\sum_{j\in\mathcal{N}_i}
  u_{j}w^{\bm{\alpha}}_{ji}
  \approx
  \bm{C}^{\bm{\alpha}} \cdot \left.\bm{D} u\right\rvert_{i}.
  \label{eq:dcpse_discrete_diff_operator_no_substr_form}
\end{equation}
This however omits the nuances associated with zeroth moment and it is completely similar to the derivation carried out in section presenting basic approach 2.

Following \cite{schrader2010, UsingDcPseOpBouran2016}, the weight function are composed by vector of monomials divided by stencil sizes $ \widetilde{\bm{X}}_{ji} $ multiplied by a suitable radial base function $ \varphi_{ji} $ (for example, in \cite{schrader2010,UsingDcPseOpBouran2016} an exponential window function is used) and unknown weight coefficients $ \bm{\Psi^{\bm{\alpha}}}_{i} $
\begin{equation}
  L^{\bm{\alpha}}_{i} u
  =
  \frac{1}{h_{i}^{|\bm{\alpha}|}}
  \displaystyle\sum_{j\in\mathcal{N}_i}
  \left( u_{j}\pm  u_{i}\right) \widetilde{\bm{X}}_{ji} \cdot \bm{\Psi^{\bm{\alpha}}}_{i} \varphi_{ji}
  =
  \displaystyle\sum_{j\in\mathcal{N}_i}
  \left( u_{j}\pm  u_{i}\right) \bm{W}_{ji} \cdot \bm{\Psi^{\bm{\alpha}}}_{i}.
  \label{eq:dcpse_discrete_diff_operator_hMod}
\end{equation}
Taking into account the notation of general form of MFD introduced above, both \cite{schrader2010, UsingDcPseOpBouran2016, AMeshfreeCollSingh2023, RecoveryByDisZwick2023, EntropicallyDaSingh2023} and \cite{TaylorSeriesEJacque2020} assume the ABFs as vector of monomials $ \widetilde{\bm{X}}_{ji} $ multiplied by a radial base function $ \varphi_{ji} $
\begin{equation}
  \bm{W}_{ji}
  =
  \widetilde{\bm{X}}_{ji} \varphi_{ji}^{\text{Gauss}}
  \label{eq:dcpse_w_taylor_monomials_with_rbf} 
\end{equation}
In the cited works, the right-hand side for the system of coefficients in the DC-PSE method is usually written as the differential operator $ \mathcal{D}^{\bm{\alpha}} $ acting on vector with monomials $ \bm{X}(0) $, which however is not anything else then different expression for our mapping vector $ \bm{C}^{\bm{\alpha}} $.
Moreover, the scaling factor $ 1/h_{i}^{|\bm{\alpha}|} $ can by absorbed by the mapping vector $ \bm{C}^{\bm{\alpha}} $, which is equivalent to scaling $ \bm{C}^{\bm{\alpha}} $ by $ \mathbf{H}_i $.
Since the $ \bm{C}^{\bm{\alpha}} $ selects the approximated derivative, the normalisation coefficient $ 1/h_{i}^{|\bm{\alpha}|} $ is automatically part of the mapping vector.
Finally DC-PSE can be written as
\begin{equation}
  \widetilde{\bm{C}^{\bm{\alpha}}}
  \coloneq
  (-1)^{\lvert\bm{\alpha}\rvert}\mathbf{H}_{i}(\partial^{\bm{\alpha}}\bm{X}_{ji})|_0 \, .
  \label{eq:dcpse_linSystem_rhs} 
\end{equation}

\begin{equation}
  L^{\bm{\alpha}}_{i} u
  \coloneq
  \displaystyle\sum_{j\in\mathcal{N}_i}
  \left( u_{j}\pm  u_{i}\right) \bm{W}_{ji} \cdot \left(\mathbf{M}_{i}^{-1}\widetilde{\bm{C}^{\bm{\alpha}}}\right),
  \quad
  \mathbf{M}_{i}
  \coloneq
  \displaystyle\sum_{j\in\mathcal{N}_i}\widetilde{\bm{X}}_{ji}\bm{W}_{ji}^{\top}.
  \label{eq:dcpse_discrete_diff_operator_hMod_recap}
\end{equation}

\subsubsection*{Notation remark}

To connect our notation with notation usually used by DC-PSE community, we also present small \textit{dictionary} to make the relation clear.
The cited papers use $\bm{Q}^\alpha  u (\vect{x}_i) $ instead of our $ L^{\bm{\alpha}}_{i} u $ to express the differential operator, $\eta^{\bm{\alpha}}( \vect{x}_i - \vect{x}_j )$ is used instead of our $ w_{ji}^{\bm{\alpha}} $ to express the weight functions (which is called a \textit{kernel function} in DC-PSE).
DC-PSE uses the kernel or stencil size $ \epsilon(\vect{x}_i) $ which we replaced with our $ h_i $ and in the scaling procedure with the vector of characteristic stencil sizes $ \mathbf{H}_i $.
In fact, all the arguments $ \vect{x}_i - \vect{x}_j $ are divided by $ \epsilon(\vect{x}_i) $ which is equivalent to the scaling with $ \mathbf{H}_i $. The searched coefficients $ \bm{a}(\vect{x}_i) $ corresponds to our $ \bm{\Psi}_i^{\bm{\alpha}} $ and the base (although only Taylor monomials are mentioned) $ \bm{p}(\mathbf{x}_j - \mathbf{x}_i) $ can be viewed as our $ \bm{P}_{ji} $ (or $ \bm{X}_{ji} $ in this context).
DC-PSE is often presented also in matrix form.

\subsection{Reproducing Kernel Collocation Method (RKCM)\label{sec:com_rkcm}}

Reproducing Kernel Collocation Method (RKCM) is the next member of the family of meshfree finite difference methods introduced by Aluru \cite{aluru2000}.
RKCM is inspired by Reproducing Kernel Particle method (RKPM) \cite{Liu1995rkpm, Chen1996rkpm}.
The motivation was to get rid of numerical integration (or quadrature) required in RKPM and utilize fully the collocation approach.
RKCM thus follows so called \textit{reproducing conditions} which are same with RKPM, but formulates them directly in discrete form.
RKCM (RKPM), like other methods, relies on the canceling of unwanted moments.
However, in contrast to the methods presented so far, RKCM (RKPM) first approximates the function itself (which is equivalent with equation (2) for $\bm{\alpha}=\bm{0}$ resulting in search for weight functions $w_{ji}^{\bm{0}}$).
The resulting weight functions $w_{ji}^{\bm{\alpha}}$ for a given differential operator $\bm{\alpha}$ are obtained as derivatives of the weight functions which approximates the function itself $w_{ji}^{\bm{0}}$, i.e. $\partial^{\bm{\alpha}}\left( w_{ji}^{\bm{0}} \right)$.
This approach is usually referred to as direct derivatives (see for example \cite{chen2017}).

The idea behind reproducing kernel approximation is to obtain weight functions $ w^{\bm{\alpha}}_{ji} $ in such a form that we correct the RBF weight function $ \varphi_{ji} $ (i.e. the kernel approximation) to force consistency to the desired order.
The correction is assumed in form $ \bm{X}_{ji}\cdot\bm{\Psi}_{i}^{\bm{\alpha}} $ where $ \bm{\Psi}_{i}^{\bm{\alpha}} $ are the unknown coefficients, that we search for by canceling selected moments.
The way the resulting function approximation is obtained \textbf{corresponds to the principle of direct approximation of moments described in section 2} which enforces the reproducing conditions for $\bm{\alpha}=\bm{0}$.

Following the procedure outlined above, to obtained RKCM approximation of derivatives we start with approximation of function itself
\begin{equation}
  L^{\bm{0}}_{i} u \coloneq
  \displaystyle\sum_{j\in\mathcal{N}_i}^{N}
  u_{j}w^{\bm{0}}_{ji}
  \approx
  \bm{C}^{\bm{0}} \cdot \left.\bm{D} u\right\rvert_{i}
  \label{eq:rkcm_function_interpolation}
\end{equation}
and following the procedure described step by step in section 2 we end up with relation
\begin{equation}
  \bm{W}_{ji}
  \coloneq
  \bm{X}_{ji} \varphi_{ji}
\end{equation}
\begin{equation}
  L^{\bm{0}}_{i} u
  =
  \displaystyle\sum_{j\in\mathcal{N}_i}  u_{j}
  \bm{W}_{ji} \cdot (\mathbf{M}_i^{-1} \bm{C}^{\bm{0}}),
  \quad
  \mathbf{M}_{i}
  \coloneq
  \displaystyle\sum_{j\in\mathcal{N}_i}
  \bm{X}_{ji}\bm{W}_{ji}^{\top}
  \label{eq:rkcm_function_interpolation_substituted weights}
\end{equation}
In RKCM, ABFs are usually assumed as Taylor monomials $\bm{X}_{ji}$ multiplied by some radial base function $ \varphi_{ji} $ and by $ \bm{C}^{\bm{0}} $ we mean mapping vector with a one in the first position and zeros in the rest.
(Note that in case of formulation with $0$th moment mentioned in section 2, $\bm{C} = \mathbf{0}$).

Derivatives of the function are then obtained by taking derivatives of the approximation
\begin{equation}
  L^{\bm{\alpha}}_{i} u
  =
  \partial^{\bm{\alpha}}
  L^{\bm{0}}_{i} u
  =
  \partial^{\bm{\alpha}}
  \left(
    \displaystyle\sum_{j\in\mathcal{N}_i}^{N}
    u_{j}w^{\bm{0}}_{ji}
  \right)
  =
  \displaystyle\sum_{j\in\mathcal{N}_i}^{N}
  u_{j}\partial^{\bm{\alpha}}w^{\bm{0}}_{ji}
  =
  \displaystyle\sum_{j\in\mathcal{N}_i}^{N}
  u_{j}w^{\bm{\alpha}}_{ji}.
  \label{eq:rkcm_discrete_diff_operator_not_substituted}
\end{equation}
By substituting for $w_{ji}^{\bm{0}}$ from \eqref{eq:rkcm_function_interpolation_substituted weights} to \eqref{eq:rkcm_discrete_diff_operator_not_substituted}, we get
\begin{equation}
  L^{\bm{\alpha}}_{i} u
  =
  \displaystyle\sum_{j\in\mathcal{N}_i}  u_{j}
  \partial^{\bm{\alpha}}
  \left(
    \bm{W}_{ji} \cdot (\mathbf{M}_i^{-1} \bm{C}^{\bm{0}})
  \right)
\end{equation}
and considering $\bm{W}_{ji} = \bm{X}_{ji} \varphi_{ji}$, this expands to resulting differential operator approximation of derivative used in RKCM
\begin{equation}
  \bm{W}_{ji}
  \coloneq
  \bm{X}_{ji} \varphi_{ji}
\end{equation}
\begin{equation}
  L^{\bm{\alpha}}_{i} u
  \coloneq
  \displaystyle\sum_{j\in\mathcal{N}_i}  u_{j}
  \left(
    \partial^{\bm{\alpha}} \bm{X}_{ji} \varphi_{ji} \cdot \mathbf{M}_{i}^{-1} \bm{C}^{\bm{0}} +
    \bm{X}_{ji} \partial^{\bm{\alpha}} \varphi_{ji} \cdot \mathbf{M}_{i}^{-1} \bm{C}^{\bm{0}} +
    \bm{X}_{ji} \varphi_{ji} \cdot \partial^{\bm{\alpha}} \mathbf{M}_{i}^{-1} \bm{C}^{\bm{0}}
  \right),
  \quad
  \mathbf{M}_{i}
  \coloneq
  \displaystyle\sum_{j\in\mathcal{N}_i}
  \bm{X}_{ji}\bm{W}_{ji}^{\top}.
  \label{eq:rkcm_discrete_diff_operator}
\end{equation}
Note also that the reproducing condition is formulated in RKCM using Taylor monomials.

\subsubsection*{Notation remark}

As in the previous case, we would like to link our notation with notation established in RKCM (which is common also to RKPM).
In RKCM/RKPM, the vector of Taylor monomials is usually denoted by $ \bm{h}(\vect{x}_j - \vect{x}_i) $ (or sometimes $ \bm{H}(\vect{x}_j - \vect{x}_i) $) instead of our $ \bm{X}_{ji} $.
The correction functions are usually denoted by $ C(\vect{x}_j, \vect{x}_j -\vect{x}_i) $ and consist of Taylor monomials $ \bm{h}(\vect{x}_j - \vect{x}_i) $ and unknown coefficients $ c(\vect{x}_j, \vect{x}_i) $.
In our notation, this corresponds to the product $ \bm{X}_{ji}\cdot\bm{\Psi}_{i}^{\bm{\alpha}} $.

\subsection{Gradient Reproducing Kernel Collocation Method (GRKCM)\label{sec:com_grkcm}}

Gradient Kernel Collocation Method (GRKCM) introduced by Chi \cite{chi2013} aims to simplify the derivative approximation used in RKCM introduced in previous part \ref{sec:com_rkcm}.
Like the RKCM, the GRKCM is based on the ideas of the RKPM, but does not use direct derivatives.
Instead, derivatives are approximated directly, in the same way we presented in the introduction (which is sometimes referred to as synchronous or implicit derivation, see \cite{chen2017}).
The method \textbf{is derived using the direct approximation of moments} exactly as presented above 2.
The 0th moment is included.

We start with discrete approximation in form
\begin{equation}
  L^{\bm{\alpha}}_{i} u
  \coloneq
  \displaystyle\sum_{j\in\mathcal{N}_i}
  u_{j}w^{\bm{\alpha}}_{ji}
  \approx
  \bm{C}^{\bm{\alpha}} \cdot \left.\bm{D} u\right\rvert_{i}.
  \label{eq:grkcm_discrete_diff_operator}
\end{equation}
Similarly to the previous section, we aim to construct correction which ensures desired order of the approximation and we assume the correction in form  $ \bm{X}_{ji}\cdot\bm{\Psi}_{i}^{\bm{\alpha}} $.
Canceling unwanted moments (or in the words of RKPM, by enforcing the condition of reproduction), we end up with
\begin{equation}
  L^{\bm{\alpha}}_{i} u
  =
  \displaystyle\sum_{j\in\mathcal{N}_i}
  u_{j}\bm{X}_{ji} \cdot \bm{\Psi^{\bm{\alpha}}}_{i} \varphi_{ji}
  \label{eq:grkcm_discrete_diff_operator_subs_weights}
\end{equation}
where the coefficients $ \bm{\Psi}_{i}^{\bm{\alpha}} $ are obtained solving linear system with symmetric matrix
\begin{equation}
  \mathbf{M}_i\bm{\Psi}_{i}^{\bm{\alpha}}
  =
  \bm{C}^{\bm{\alpha}}.
  \label{eq:grkcm_linSystem}
\end{equation}
Using the notation introduced above, the RKCM is expressed as
\begin{equation}
  \bm{W}_{ji}
  \coloneq
  \bm{X}_{ji} \varphi_{ji}
  \label{eq:grkcm_w_taylor_monomials_with_rbf}
\end{equation}

\begin{equation}
  L^{\bm{\alpha}}_{i} u
  \coloneq
  \displaystyle\sum_{j\in\mathcal{N}_i}  u_{j} \bm{W}_{ji} \cdot (\mathbf{M}_i^{-1} \bm{C}^{\bm{\alpha}}),
  \quad
  \mathbf{M}_{i}
  \coloneq
  \displaystyle\sum_{j\in\mathcal{N}_i}
  \bm{X}_{ji}\bm{W}_{ji}^{\top}.
  \label{eq:grkcm_discrete_diff_operator_recap}
\end{equation}

\subsubsection*{Notation remark}

Same as above, in GRKCM the vector of Taylor monomials is usually denoted by $ \bm{h}(\vect{x}_j - \vect{x}_i) $ (or sometimes $ \bm{H}(\vect{x}_j - \vect{x}_i) $) instead of our $ \bm{X}_{ji} $.
The mapping vector is usually denoted using the vector of Taylor monomials $ \bm{C}^{\bm{\alpha}} = (-1)^{\lvert \bm{\alpha} \rvert}\partial^{\bm{\alpha}}\bm{H}(\mathbf{0}) $.
The correction functions are usually denoted by $ C(\vect{x}_j, \vect{x}_j - \vect{x}_i) $ and consist of Taylor monomials $ \bm{h}(\vect{x}_j - \vect{x}_i) $ and unknown coefficients $ c(\vect{x}_j, \vect{x}_i) $.
In our notation, this corresponds to the product $ \bm{X}_{ji}\cdot\bm{\Psi}_{i}^{\bm{\alpha}} $.

\subsection{Generalized Finite Difference Method (GFDM)\label{sec:com_gfdm}}

Generalized Finite Difference Method (GFDM) is popular versions of MFD introduced by Liszka and Orkisz \cite{liszka1980finite, liszka1984} inspired by work by Jensen \cite{jensen1972} or paper by Swartz and Wendroff \cite{swartz1969generalized}.
Since then, various improvements have been made across different areas and a number of additional techniques have been developed.
Here we focus on the current concept of GFDM as it can be found under this particular name.
For this purpose, we refer to the work \cite{suchde2018}.
The mechanism of GFDM is exactly the same as for the previous methods but \textbf{it was build on the principle of $\ell_2$ minimization of approximation error} which is described in section 4 and later on \textbf{using the polynomial method $\ell_2$ minimization} as we presented in section 5.

In the case of approximation error $ \ell_2 $ we start with the discrete version of differential operator
\begin{equation}
  L^{\bm{\alpha}}_{i} u  \coloneq
  \displaystyle\sum_{j\in\mathcal{N}_i}
  u_{j}w^{\bm{\alpha}}_{ji}
  \label{eq:gfdm_general_discrete_diff_operator_zm}
\end{equation}
and using Taylor expansion of field $  u $ in the vicinity of point $ i $
\begin{equation}
  u_{j}
  =
  \bm{X}_{ji}\cdot\left.\bm{D} u\right\rvert_{i} + \bm{E}_{ji},
  \label{eq:gfdm_taylor_expansion}
\end{equation}
with aim to minimize the square of the error of the approximation $\lVert \bm{E}_{ji} \rVert_{2}^{2}$ weighted by some suitable weight function $ \varphi_{ji} $.
This leads to a optimization problem which is step by step solved in section 4 which results in the known system of equations for the vector of derivatives $ \left.\bm{D} u\right\rvert_{i} $
\begin{equation}
  \mathbf{M}_{i}\left.\bm{D} u\right\rvert_{i}
  =
  \displaystyle\sum_{j\in\mathcal{N}_i}
  u_{ji}\bm{X}_{ji}\varphi_{ji},
  \quad
  \mathbf{M}_{i}
  \coloneq
  \displaystyle\sum_{j\in\mathcal{N}_i}\bm{X}_{ji}\bm{X}_{ji}^{\top} \varphi_{ji}.
  \label{eq:gfdm_linSystem}
\end{equation}
Expressing this by the vector of basis functions $ \bm{W}_{ji} $, which is in this case again composed from the vector of monomials $ \bm{X}_{ji} $ multiplied by a suitable RBF $ \varphi_{ji} $
\begin{equation}
  \bm{W}_{ji}
  \coloneq
  \bm{X}_{ji} \varphi_{ji}
  \label{eq:gfdm_w_taylor_monomials_with_rbf}
\end{equation}
we can write \eqref{eq:gfdm_linSystem} as
\begin{equation}
  \mathbf{M}_{i}\left.\bm{D} u\right\rvert_{i}
  =
  \displaystyle\sum_{j\in\mathcal{N}_i} u_{ji}\bm{W}_{ji},
  \quad
  \mathbf{M}_{i}
  =
  \displaystyle\sum_{j\in\mathcal{N}_i}\bm{X}_{ji}\bm{W}_{ji}^{\top}.
  \label{eq:gfdm_linSystem_with_w}
\end{equation}

In GFDM, all first $m$ derivatives are computed together at once. To construct particular differential operator $ d $, we include $ \bm{C}^{\bm{\alpha}} $ to select desired derivatives.
The result is
\begin{equation}
  L^{\bm{\alpha}}_{i} u
  =
  \displaystyle\sum_{j\in\mathcal{N}_i} \bm{C}^{\bm{\alpha}} \cdot (\mathbf{M}_i^{-1} \bm{W}_{ji})  u_{ji}
  \label{eq:gfdm_discrete_diff_operator}
\end{equation}
which we can do due to the symmetry of the matrix $ \mathbf{M} $. Reshuffle to the form \eqref{eq:general_discrete_diff_operator_with_coefs_do}
\begin{equation}
  \bm{W}_{ji}
  \coloneq
  \bm{X}_{ji} \varphi_{ji}.
  \label{eq:gfdm_w_taylor_monomials_with_rbf_recap}
\end{equation}
\begin{equation}
  L^{\bm{\alpha}}_{i} u
  \coloneq
  \displaystyle\sum_{j\in\mathcal{N}_i}
  u_{ji} \bm{W}_{ji} \cdot (\mathbf{M}_i^{-1} \bm{C}^{\bm{\alpha}}),
  \quad
  \mathbf{M}_{i}
  \coloneq
  \displaystyle\sum_{j\in\mathcal{N}_i}\bm{X}_{ji}\bm{W}_{ji}^{\top}.
  \label{eq:gfdm_discrete_diff_operator_recap}
\end{equation}

In the case of polynomial method $\ell_2$ minimization, we start from the same assumed form of discrete differential operator \eqref{eq:gfdm_general_discrete_diff_operator_zm}, but instead of minimizing the approximation error, we determine the weight functions $ w_{ji}^{\bm{\alpha}} $ by enforcing consistency up to order of $ m $ of gradient reproducing conditions. The minimization problem given by equation (97) and the section 5 presents a step-by-step procedure leading to the same equation \eqref{eq:gfdm_discrete_diff_operator_recap}.
Following the derivation, we obtained the same result (as for the approximation error $ \ell_2 $) that the basis function are given as $ \bm{W}_{ji} = \bm{X}_{ji} \varphi_{ji} $.

As we can observe, deriving the method using minimization of approximation error leads in fact to less general results, since the choice of the basis is less flexible (the vector of ABFs is given) than for the methods introduced so far.
Although the polynomial minimization produces the same result, we showed in section 5 that this approach can be generalized such that arbitrary basis functions can be chosen in the same way as in 2.

\subsection{Moving least squares collocation methods (MLS) \label{sec:com_mls} }

One particular class of meshfree finite difference methods is based directly on the Moving Least Squares (MLS) method introduced by Lancester \cite{lancaster1981}.
This approach was quickly adopted for approximation of differential operators \textbf{following the derivations using optimization procedure from sections 4 and 5}.
Over the years, a number of variants, modifications and improvements have been introduced to improve the properties of the moment matrix and reduce the number of points needed for approximation.
MLS is also widely used in meshfree methods employed in with weak form, however, in the following section we present purely collocation methods based on MLS for approximating differential operators in strong form.
Namely we introduce Finite Point Method (FPM) \ref{sec:com_fpm}, Generalized Moving Least Squares (GMLS) \ref{sec:com_gmls}, Compact Moving Least Squares (CMLS) \ref{sec:com_cmls}, Modified Moving Least squares (MMLS) \ref{sec:com_mmls}, Interpolating Moving Least Squares (IMLS) \ref{sec:com_moreOnMMLS} and, used with the same abbreviation as the previous one Improved Moving Least Squares (IMLS) \ref{sec:com_moreOnMMLS}.

\subsubsection{Finite Point Method (FPM)\label{sec:com_fpm}}

The term \textit{Finite Point Methods} is used originally by the authors \cite{onate1995b} to describe a class of point data interpolation-based procedures that encompass the methods like Least Squares Methods (LSQ), Weighted Least Squares (WLS), Moving Least Squares (MLS) and Reproducing Kernel Particle (RKP) Methods. In the later works the authors use MLS as the main method under the title \textit{Finite Point Method} (FPM) to simulate the equations of self-adjoint and non-self-adjoint nature like advective-diffusive transport problems, compressible fluid flow problems and elasticity problems in structural mechanics\cite{onate1996finite, onate1996fpm, onate1998fpm, onate2000finite, onate2001fpm}. In addition, the authors thoroughly elaborate on stabilization techniques for the non-self-adjoint problems applicable to the meshfree formulation.

Given MLS is assumed to be the method of choice for FPM, the derivation is identical to the definition given in the section 4 for $\bm{D}u = L^{\bm{0}}u$, i.e. it first builds the approximation of the function $u_i$ which is then used to compute its derivatives.

Additionally, the authors have extended \cite{onate1996fpm} MLS formalism by introducing Multiple Fixed Least Squares (MFLS) approximation by replacing the weighting function $ \varphi_k(\mathbf{x}_j) $ centered at an arbitrary point $\mathbf{x}_k$ (not necessarily one of the nodes $\mathbf{x}_i$) where the function approximation $u_k$ or the derivative approximation $L_k^{\bm{\alpha}}u$ is evaluated, to $ \varphi_j(\mathbf{x}_k) $, i.e. defining the weighting function at each node location $ \varphi_j $. This formalism is equivalent to the one used by us (see section 2).

Following the formal minimization presented in section 4 we end up with the expression
\begin{equation}
  u_{i}^{*}
  =
  L_{i}^{\bm{0}} u
  =
  \bm{P}_i
  \cdot
  \displaystyle\sum_{j\in\mathcal{N}_i}
  \left[
    \left(
      \displaystyle\sum_{k\in\mathcal{N}_i}
      \bm{P}_{ki} \bm{P}_{ki}^{\top} W_{ki}^{\bm{0}}
    \right)^{-1}
    \bm{P}_{ji} \varphi_{ji}  u_j
  \right].
  \label{eq:fpm_l0_approximation_raw}
\end{equation}

The authors consider two approaches to compute derivatives of the shape functions: \textit{full derivative} (also known as \textit{direct derivative}) and \textit{diffuse derivative} \cite{nayroles1992dem,onate1995b} (for further details see \cite{chen2017}). \textit{Direct derivative} of the expression \eqref{eq:fpm_l0_approximation_raw} leads to a complicated and computationally expensive formula since both $ \bm{P}_{ji}$ and $ \varphi_{ji} $ depend on the position $ \mathbf{x} $.
Instead with the method of \textit{diffuse derivative} the coefficients of the basis $ \varphi_{ji} $ are assumed to be constant and their derivatives are therefore zero.
\footnote{
  Another view - not given in the FPM derivation, however - of diffusion derivations may be that instead of primary function, we approximate its derivatives directly. For further details see \cite{nayroles1992dem, chen2017}.
}
The approximation of the corresponding differential operator $ L_{i}^{\bm{\alpha}} u $ is thus obtained by replacing the basis functions $ \bm{P}_i $ in equation \eqref{eq:fpm_l0_approximation_raw} by their respective derivatives $ \left.\partial^{\bm{\alpha}}\bm{P}\right\rvert_{i} $ which gives us
\begin{equation}
  L_{i}^{\bm{\alpha}} u
  =
  \left.\partial^{\bm{\alpha}}\bm{P}\right\rvert_{i}
  \cdot
  \displaystyle\sum_{j\in\mathcal{N}_i}
  \left[
    \left(
      \displaystyle\sum_{k\in\mathcal{N}_i}
      \bm{P}_{ki} \bm{P}_{ki}^{\top} W_{ki}^{\bm{0}}
    \right)^{-1}
    \bm{P}_{ji} \varphi_{ji}  u_j
  \right].
  \label{eq:fpm_ld_approximation_raw}
\end{equation}
Taking the ABFs in form $ \bm{W}_{ji} = \bm{P}_{ji} \varphi_{ji} $, we recognize the moment matrix inside the parentheses
\begin{equation}
  \mathbf{M}_i
  =
  \displaystyle\sum_{j\in\mathcal{N}_i}
  \bm{P}_{ji}\bm{W}_{ji}^{\top}.
  \label{eq:fpm_momentum_matrix}
\end{equation}
Moreover $ \left.\partial^{\bm{\alpha}}\bm{P}_i\right\rvert_{i} $ can be interpreted as the \textit{mapping vector} $ \bm{C}^{\bm{\alpha}} $.
This results in final form of the discrete differential operator
\begin{equation}
  \bm{W}_{ji}
  =
  \bm{P}_{ji} \varphi_{ji}.
\end{equation}
\begin{equation}
  L^{\bm{\alpha}}_{i} u
  =
  \displaystyle\sum_{j\in\mathcal{N}_i}
  u_{j} \bm{W}_{ji} \cdot (\mathbf{M}_i^{-1} \bm{C}^{\bm{\alpha}}),
  \quad
  \mathbf{M}_{i}
  =
  \displaystyle\sum_{j\in\mathcal{N}_i}\bm{P}_{ji}\bm{W}_{ji}^{\top}.
  \label{eq:fpm_ld_approximation}
\end{equation}
which is in agreement with the results presented above.

To conclude the FPM section, we would like to point out, that to the best of our knowledge, every available application has used Taylor monomials as basis for ABFs, i.e. $ \bm{P}_{ji} = \bm{X}_{ji} $ and thus $ \bm{W}_{ji} = \bm{X}_{ji} \varphi_{ji} $.
FPM is also usually presented in the matrix notation.

\paragraph{Notation remark}

Last, we provide link between our notation and notation consistently used in FPM.
Papers cited above usually use $ \bm{p}(\mathbf{x}) $ to denote the basis (or vector of Taylor monomials in particular cases) which is denoted by our $\bm{P}_{ji} $ or by $\bm{X}_{ji} $ in case of Taylor monomials.
Weight coefficients in FPM are usually denoted by $ \bm{\alpha} $ while we use $ \Psi^{\bm{\alpha}}_{i} $.
The weight functions are usually represented by vector $ \bm{N}_i $ corresponds to our $ w_{ji}^{\bm{\alpha}} $ or in matrix notation $ \bm{w}_{i}^{\bm{\alpha}} $.
The RBF weight function is usually denoted as $ \varphi_i(\mathbf{x}_j) $ which we replace with $ \varphi_{ji} $.

\subsubsection{Meshless Local Strong Form Method (LSMFM)\label{sec:com_lsmfm}}

Meshless Local Strong Form Method (LSMFM) introduced by Park and Youn \cite{park2001} was introduced with the aim to remove the numerical integration required in meshfree methods using Galerkin formulation.
LSMFM based on MLS \ref{sec:com_mls}, so similarly as other methods of this kind, it \textbf{uses minimization of approximation error 4, to approximate the function itself, which is then used to compute desired derivatives}.
LSMFM is thus very similar to RKCM introduced in section \ref{sec:com_rkcm}, with that difference it follows derivation based on MLS instead of approximation of moments.
\textbf{The zero-th moment is included.}

To obtain the approximation, first, function itself is approximated which corresponds to equation (2) for $\bm{\alpha}=\bm{0}$
\begin{equation}
  L^{\bm{0}}_{i} u \coloneq
  \displaystyle\sum_{j\in\mathcal{N}_i}^{N}
  u_{j}w^{\bm{0}}_{ji}
  \approx
  \bm{C}^{\bm{0}} \cdot \left.\bm{D} u\right\rvert_{i}
  \label{eq:app:variations:lsmfm:funcion approx general}
\end{equation}
Using ABFs  given by Taylor monomials multiplied by eight order RBF $\bm{W}_{ji} = \bm{P}_{ji} \varphi_{ji}$ and momentum matrix defined, the approximation of function results
\begin{equation}
  L^{\bm{0}}_{i} u
  =
  \displaystyle\sum_{j\in\mathcal{N}_i}  u_{j}
  \bm{W}_{ji} \cdot (\mathbf{M}_i^{-1} \bm{C}^{\bm{0}}),
  \label{eq:app:variations:lsmfm:funcion approx}
\end{equation}
and taking the derivatives of the function approximation
\begin{equation}
  L^{\bm{\alpha}}_{i} u
  =
  \displaystyle\sum_{j\in\mathcal{N}_i}  u_{j}
  \partial^{\bm{\alpha}}
  \left(
    \bm{W}_{ji} \cdot (\mathbf{M}_i^{-1} \bm{C}^{\bm{0}})
  \right)
\end{equation}
gives the final approximation of derivatives
\begin{equation}
  \bm{W}_{ji}
  \coloneq
  \bm{X}_{ji} \varphi_{ji}
\end{equation}
\begin{equation}
  L^{\bm{\alpha}}_{i} u
  \coloneq
  \displaystyle\sum_{j\in\mathcal{N}_i}  u_{j}
  \left(
    \partial^{\bm{\alpha}} \bm{X}_{ji} \varphi_{ji} \cdot \mathbf{M}_{i}^{-1} \bm{C}^{\bm{0}} +
    \bm{X}_{ji} \partial^{\bm{\alpha}} \varphi_{ji} \cdot \mathbf{M}_{i}^{-1} \bm{C}^{\bm{0}} +
    \bm{X}_{ji} \varphi_{ji} \cdot \partial^{\bm{\alpha}} \mathbf{M}_{i}^{-1} \bm{C}^{\bm{0}}
  \right),
  \quad
  \mathbf{M}_{i}
  \coloneq
  \displaystyle\sum_{j\in\mathcal{N}_i}
  \bm{X}_{ji}\bm{W}_{ji}^{\top}.
  \label{eq:app:variations:lsmfm:ddo approx}
\end{equation}
which corresponds to the forumla \eqref{eq:rkcm_discrete_diff_operator} of the RKCM method introduced above \ref{sec:com_rkcm}.

\subsubsection{Generalized Moving Least Squares (GMLS) \label{sec:com_gmls} }

Taking into account all these strong form methods constructed as some variation on MLS, Generalized Moving Least Square method proposed in \cite{ levin1998, wendland2004, mirzaei2012, gross2020} is very popular one.
The approach used in \textbf{GMLS is based on the minimization of approximation error, but in contrast to GFDM and section 4, it focuses on the approximation of reconstruction of target function, which is then used to compute desired derivatives}.

A more thorough mathematical formalism and rigor can be found in GMLS.
The methods is also presented as a more general tool since it is focus on arbitrary linear bounded operator in general.
In the following lines we present a derivation of the GMLS method, emphasizing the precise classification of the individual terms, in contrast to the derivations and methods presented above.
The notation introduced above is used for the derivation, while the connection to the classical notation usually used in GMLS is given at the end of the section.

Assume function $u(\vect{x})$ to be $m$ times differentiable function $  u \in C^m(\Omega) $ for some $m \ge 0$. Moreover, assume, that $ u $ can be characterized by a collection of scattered so called \textit{sampling} functionals $ \Lambda( u) = \left\{ \lambda_j( u)\right\}_{j=0}^{N} \in C^m(\Omega)^* $ where $ C^m(\Omega)^* $ is the dual space to $ C^m(\Omega) $.
Usually, these sampling functionals are chosen as simple evaluation of $  u $ in given point i.e. $ \lambda_i ( u) = \delta_{\vect{x}_i}[  u ] =  u( \vect{x}_i) $ or simply $ \lambda_i ( u) =  u_i $.

For a given evaluation functional $\lambda(u) \in C^m(\Omega)^*$ the problem of approximating $\lambda(u)^h$ of the function value $\lambda(u)$ in the context of GMLS consits of searching for a linear function of the data $ \left\{ \lambda_j( u)\right\}_{j=0}^{N}$:

\begin{equation}
  \lambda(u)^h \coloneq
  \displaystyle\sum_{j\in\mathcal{N}_i}^{N}
  \lambda_j( u)w_{j}(\lambda).
  \label{eq:general_approx_operator_gmls}
\end{equation}

Now, we consider linear operator  $ \tau_{\vect{x}} $ at location $ \vect{x} $.
In our case, we are interested in differential operators $ \tau_i = \left.\partial^{\bm{\alpha}} u\right\rvert_{i} = \bm{C}^{\bm{\alpha}} \cdot \left.\bm{D}  u\right\rvert_{i} $ but we keep in mind that this can be in practice arbitrary linear operator.
We aim for approximation $ \tau_i^h $ of the operator $ \tau_i $ which we search in the form
\begin{equation}
  \tau_i^h u
  =
  \displaystyle\sum_{j\in\mathcal{N}_i}
  \lambda_{j}( u) \widetilde{w}_{j}(\lambda).
  \label{eq:tau_approx_operator_gmls}
\end{equation}
Alternatively, the solution to the function approximation problems \eqref{eq:general_approx_operator_gmls} can be expressed in the matrix form (see Appendix \ref{sec:appendix-matrixForm}) using a vector of coefficients $ \widetilde{\bm{\Psi}}(u^*) $ as

\begin{equation}
  \lambda(u)^h = \lambda(u^*)
  =
  \bm{P}(\vect{x})^{\top} \widetilde{\bm{\Psi}}(u^*)
\end{equation}
where $\bm{P}(\mathbf{x})$ denotes the vector of basis functions and $ u^* $ is the solutions to the local weight $\ell_2$-optimization problem over finite dimensional subspace $ \mathcal{V}_h \in \mathcal{V} $.
For $  u^* \in \mathcal{V}_h $ we need to solve
\begin{equation}
  u^*
  =
  \argmin_{q \in \mathcal{V}_h}
  \displaystyle\sum_{j\in\mathcal{N}_i}
  \left( \lambda_j\left( u\right) - \lambda_j\left(q\right) \right)^2 \varphi_{ji}.
  \label{eq:gmls_l2_optim_problem}
\end{equation}
where $\mathcal{V}_h = \text{span}\{ P_1, ..., P_{\text{dim}(\mathcal{V}_h)} \} $, we denote $ P_i $ denotes individual basis functions of basis $\bm{P}$.
The dimension $ \dim({\mathcal{V}_h}) = p $ depends on the desired order of approximation $ m $.

To obtain an approximation $ \tau_i^h $ of the operator $ \tau_i $ we compute derivatives of the vector of basis functions $ \bm{P} $ while keeping the vector of coefficients $ \widetilde{\bm{\Psi}}(u^*) $ constant.
This approach has been connected to the notion of \textit{diffuse derivatives} (see \cite{mirzaei2012}). For further details on \textit{diffuse derivatives} and their convergence to exact derivatives see \cite{nayroles1992dem,Kim2003,Huerta2004,Yoon2006,chen2017}

\begin{equation}
  \tau_i^h u
  = \tau_i (\bm{P})^{\top} \widetilde{\bm{\Psi}}(u).
  \label{eq:gmls_diffuse_derivative}
\end{equation}

Few notes about the formulas above.
Fist, in the minimization problem  \eqref{eq:gmls_l2_optim_problem}, some of the well known radial functions $ \varphi_{ji} $ are used to weight the square of an error.
Second, in contrast to optimization problem presented in 4, here we focus on minimizing the error of the approximation of the target function instead of minimizing the error of a particular approximation of the differential operator.

The two operator approximations \eqref{eq:tau_approx_operator_gmls} and \eqref{eq:gmls_diffuse_derivative} are connected by the relation
\begin{equation}
  \tau_i^h u
  = \tau_i (\bm{P})^{\top} \widetilde{\bm{\Psi}}(u^*) = \displaystyle\sum_{j\in\mathcal{N}_i}
  \lambda_{j}( u) \widetilde{w}_{j}(\lambda).
\end{equation}
Since we focus on differential operators and $ \tau_i = \left.\partial^{\bm{\alpha}} u\right\rvert_{i} \approx \bm{C}^{\bm{\alpha}} \cdot \left.\bm{D}  u\right\rvert_{i} $, the resulting approximation corresponds to our $L_{i}^{\bm{\alpha}} u $ and we would like to express the result in form
\begin{equation}
  L^{\bm{\alpha}}_{i}u \coloneq
  \displaystyle\sum_{j\in\mathcal{N}_i}^{N}
  u_{j}w^{\bm{\alpha}}_{ji}
  \label{eq:gmls_gddo}
\end{equation}
The procedure to obtain the coefficient vector is the same as presented in the section 4.
Following the formal minimization, we end up with approximation of the differential operator in form
\begin{equation}
  L_{i}^{\bm{\alpha}}( u)
  =
  \partial^{\bm{\alpha}}\bm{P}(\mathbf{0})
  \cdot
  \displaystyle\sum_{j\in\mathcal{N}_i}
  \varphi_{ji} \lambda_j(u)
  \left(
    \displaystyle\sum_{k\in\mathcal{N}_i}
    \lambda_k(\bm{P}_i) \lambda_k(\bm{P}_i)^{\top} \varphi_{ki}
  \right)^{-1}
  \lambda_j(\bm{P}_{i}).
  \label{eq:gmls_ld_approximation_raw}
\end{equation}
Since we assume that the functional $ \lambda_j = \delta_{\vect{x}_i} $ quantifies the point value of the variable at point $ j $ (however, one should beware that it may in principle be a more complex functional), following notation introduced above, we can simply write $ \lambda_j(\bm{P}_{i}) = \bm{P}(\vect{x}_j - \vect{x}_i) = \bm{P}_{ji} $ (again $\bm{P}_{ji} \neq \bm{P}_j - \bm{P}_i $) and $ \lambda_j( u) =  u_j $.
Next, we cam see that the first term of the dot product is the \textit{mapping vector} $ \partial^{\bm{\alpha}}\bm{P}(\mathbf{0}) = \bm{C}^{\bm{\alpha}} $  where $ \partial^{\bm{\alpha}}\bm{P}(\mathbf{0}) $ originates from applying the differential operator $ \partial^{\bm{\alpha}} $ to $\bm{P}_i $ evaluated at point $ \mathbf{x}_i $ which can be written as $ \lambda_i(\bm{P}_i) = \bm{P}(\mathbf{0}) $.
The bracket is nothing else then moment matrix
\begin{equation}
  \mathbf{M}_i
  =
  \displaystyle\sum_{k\in\mathcal{N}_i}
  \lambda_k(\bm{P}_i) \lambda_k(\bm{P}_i)^{\top} \varphi_{ki}.
  \label{eq:gmls_momentum_matrix}
\end{equation}
In the context of vector of ABFs, the basis $\bm{P}$ is the \textit{building block} for the vector $ \bm{W}_{ji} $.
So far we assumed mostly Taylor monomials (i.e. $\bm{P}_{ji} = \bm{X}_{ji} $) and for example in case of LABFM in \ref{sec:com_labfm} $ \bm{P} $ was represented by Hermite or Legendre polynonomials.
As in the previous cases, the vector of ABFs is given as
\begin{equation}
  \bm{W}_{ji}
  =
  \bm{P}_{ji} \varphi_{ji}.
  \label{eq:gmls_w_basis}
\end{equation}
Using \eqref{eq:gmls_w_basis} the moment matrix can be rewritten as
\begin{equation}
  \mathbf{M}_i
  =
  \displaystyle\sum_{j\in\mathcal{N}_i}\bm{P}_{ji} \bm{W}_{ji}^{\top}.
  \label{eq:gmls_m_matrix}
\end{equation}

Using the symmetry of the moment matrix, we can pull the vector $\bm{C}^{\bm{\alpha}}$  inside the summation and swapping the matrix multiplication results in the approximation of differential operator given by
\begin{equation}
  \bm{W}_{ji}
  =
  \bm{P}_{ji} \varphi_{ji}.
\end{equation}
\begin{equation}
  L^{\bm{\alpha}}_{i} u
  =
  \displaystyle\sum_{j\in\mathcal{N}_i}
  u_{j} \bm{W}_{ji} \cdot (\mathbf{M}_i^{-1} \bm{C}^{\bm{\alpha}}),
  \quad
  \mathbf{M}_{i}
  =
  \displaystyle\sum_{j\in\mathcal{N}_i}\bm{P}_{ji} \bm{W}_{ji}^{\top}.
  \label{eq:gmls_ld_approximation}
\end{equation}
which is almost the same result, as for the methods above.
There are however slight nuances.
As in the case of the GFDM \ref{sec:com_gfdm} and RKCM \ref{sec:com_rkcm}, the result is a matrix that is always symmetric since it is formed by the product of the base $ \mathcal{P_{ji}} $.
This is a consequence of the derivation via $\ell_2$ error minimization.
In contrast, DC-PSE \ref{sec:com_dcpse} and LABFM \ref{sec:com_labfm} allow the matrix to be chosen in other ways.
Moreover, in contrast to the GFDM, where we start from a Taylor expansion and the ABFs are formed by Taylor monomials, here we have the choice of $ \bm{W}_{ji} $ functions, hence the basis $ \bm{P}_{ji} $ functions.

For the basis $ \bm{P}_{ji} $ bivariate polynomials are among the popular choices.
In the literature, orthogonal polynomials such as Legendre, Laguerre or Hermite polynomials are also mentioned.
Taylor monomials are as well an option.

When solving for the operator approximation $ \tau_i^h $ in the matrix formulation \eqref{eq:gmls_diffuse_derivative}, resulting equation \eqref{eq:gmls_ld_approximation_raw} can be expressed as
\begin{equation}
  \tau_i^h u
  =
  \tau_i (\bm{P})^{\top} \left( \Lambda_i(\bm{P})^{\top} \mathbf{V}_i \Lambda(\bm{P}) \right)^{-1} \Lambda_i(\bm{P})^{\top} \mathbf{V}_i \Lambda_i( u)
\end{equation}
where $ \Lambda( u) = \left\{ \lambda_j( u)\right\}_{j=0}^{N} \in \mathbb{R}^{N} $ \textit{stacks} all the values from the neighboring points $ \mathcal{N}_i $ similarly to the equation \eqref{eq:appendix-introduction-of-matrix-formalism} and $ \mathbf{V}_i $ is the diagonal matrix with stores RBF values for all neighbors defined by equation \eqref{eq:l2e:appendix_matrixform_wmatrix}.

\subsubsection{Compact Moving Least Squares (CMLS)\label{sec:com_cmls}}

Compact Moving Least Squares (CMLS) introduced by Trask \cite{trask2016} builds on the standard MLS approach and enriches it with procedures for constructing compact finite difference schemes \cite{lele1992cfdm, weinan1996cfdm}.
The strategy is to consider the formulation of problem the when calculating the coefficients of the weight functions.
In practice, this is achieved by augmenting the functional for minimization from which we derive the reconstruction with additional conditions that take into account the initial problem.
This aims to reduce the size of the stencils and keep moment matrices invertible while using smaller number of neighbors, then is required for the standard case.

First of all, assume the standard MLS approach \cite{chen2017} or simply the same method as described in previous section \ref{sec:com_gmls}.
Assume for example following linear boundary value problem
\begin{align*}
  \mathcal{L}_{\Omega}  u &= f, \qquad \forall \bm{x} \in \Omega \\
  \mathcal{L}_{\partial\Omega}  u &= g, \qquad \forall \bm{x} \in \partial\Omega
\end{align*}
In case of CMLS, minimization problem for optimal polynomial reconstruction \eqref{eq:gmls_l2_optim_problem} is enhanced by additional penalization terms which reflects the solved problem.
For $  u^* \in \mathcal{V}_h $ we need to solve
\begin{equation}
  u^*
  =
  \argmin_{q \in \mathcal{V}_h}
  \displaystyle\sum_{j\in\mathcal{N}_i}
  \left[
    \left(
      \lambda_j\left( u\right) - \lambda_j\left(q\right)
    \right)^2
    +
    \varepsilon_{\Omega}
    \left(
      f_j - \mathcal{L}_{\Omega} q
    \right)^2
    +
    \chi_{j,\partial\Omega}
    \varepsilon_{\partial\Omega}
    \left(
      g_j - \mathcal{L}_{\partial\Omega} q
    \right)^2
  \right] \varphi_{ji}.
  \label{eq:cmls_l2_optim_problem}
\end{equation}
where $ \varepsilon_{\Omega} $ and $ \varepsilon_{\partial\Omega} $ stand for the regularization parameters and $ \chi_{j,\partial\Omega} $ is the characteristic function which $ \chi_{j,\partial\Omega} = 1 $ if $ \vect{x}_j \in \partial\Omega $ and $ \chi_{j,\partial\Omega} = 0 $ otherwise.
The first penalty term is same as in the previous case (or as in the standard MLS) and it is responsible for the interpolation quality of the reconstruction.
The second penalty term takes into account the correspondence between reconstruction and the governing equation and the third penalty term reflects how faithful is the reconstruction to boundary conditions.
The regularization parameters are used to tune the significance of individual terms.
To enforce boundary conditions, an additional constraint is necessary which requires that reconstruction satisfy the boundary conditions exactly
\begin{equation}
  \mathcal{L}_{\partial\Omega}q_i = g_i, \qquad \forall i:\; \bm{x}_i \in \partial\Omega.
\end{equation}

Optimization problem \eqref{eq:cmls_l2_optim_problem} is solved in the same fashion as it is done in detail in sections 4, 5 and \ref{sec:gl2p}.
The structure of the moment matrix is more complex in this case, as it includes contributions from additional penalty terms.
\begin{equation}
  \bm{W}_{ji}
  =
  \bm{P}_{ji} \varphi_{ji}
\end{equation}
\begin{equation}
  L^{\bm{\alpha}}_{i} u
  =
  \displaystyle\sum_{j\in\mathcal{N}_i}
  \left(
    u_{j} \bm{W}_{ji} +
    \epsilon_{\Omega} \mathcal{L}_{\Omega} \bm{P}_{ji} f_i \varphi_{ji} +
    \chi_{j,\partial\Omega} \epsilon_{\partial\Omega} \mathcal{L}_{\partial\Omega} \bm{P}_{ji} g_i \varphi_{ji}
  \right) \cdot (\mathbf{M}_i^{-1} \bm{C}^{\bm{\alpha}}),
  \label{eq:cmls_discrete_diff_operator_no_enforced_bc}
\end{equation}
\begin{equation}
  \mathbf{M}_i =
  \displaystyle\sum_{j\in\mathcal{N}_i}
  \left(
    \bm{P}_{ji} \bm{W}_{ji}^{\top} +
    \epsilon_{\Omega} (\mathcal{L}_{\Omega} \bm{P}_{ji})(\mathcal{L}_{\Omega} \bm{P}_{ji}^{\top}) \varphi_{ji} +
    \chi_{j,\partial\Omega} \epsilon_{\partial\Omega} (\mathcal{L}_{\partial\Omega} \bm{P}_{ji})(\mathcal{L}_{\partial\Omega} \bm{P}_{ji}^{\top}) \varphi_{ji}
  \right)
\end{equation}

To enforce exactly boundary conditions, following \cite{trask2016}, we need to include boundary condition constraint to the matrix $M_i$ and right hand side vector.
Following notation introduced in section 4, we write equation \eqref{eq:cmls_discrete_diff_operator_no_enforced_bc} as
\begin{equation}
  L_{i}^{\bm{\alpha}}u = \bm{C}^{\bm{\alpha}} \cdot \left.\bm{D}u\right\rvert_{i}
\end{equation}
where the derivatives coefficient $ \left.\bm{D}u\right\rvert_{i} $ are obtained as solution of system of linear equations $ \mathbf{M}_i \left.\bm{D}u\right\rvert_{i} = \bm{R}_{i} $, where
\begin{equation}
  \bm{R}_{i} \coloneq
  \left[
    \displaystyle\sum_{j\in\mathcal{N}_i}
    \left(
      u_{j} \bm{W}_{ji}^{\top} +
      \epsilon_{\Omega} \mathcal{L}_{\Omega} \bm{P}_{ji}^{\top} f_i \varphi_{ji} +
      \chi_{j,\partial\Omega} \epsilon_{\partial\Omega} \mathcal{L}_{\partial\Omega} \bm{P}_j^{\top} g_i \varphi_{ji}
    \right)
  \right].
\end{equation}.
Including the boundary condition constraint, the problem takes form
\begin{equation*}
  \widetilde{\mathbf{M}_i} \coloneq
  \begin{bmatrix}
    \mathbf{M}_i & \mathcal{L}_{\partial\Omega}\bm{P}_{ij} \\
    \mathcal{L}_{\partial\Omega}\bm{P}_{ij}^{\top} & 0
  \end{bmatrix},
  \quad
  \widetilde{\left.\bm{D}u\right\rvert_{i}} \coloneq
  \left[
    \left.\bm{D}u\right\rvert_{i}, \lambda
  \right],
  \quad
  \widetilde{\bm{R}_{i}} \coloneq
  \left[
    \bm{R}_{i}, g_i
  \right].
\end{equation*}
Solving system of linear equations $ \widetilde{\mathbf{M}_i} \widetilde{\left.\bm{D}u\right\rvert_{i}} = \widetilde{\bm{R}_{i}} $, we obtain the approximation
\begin{equation}
  L_{i}^{\bm{\alpha}}u = \bm{C}^{\bm{\alpha}} \cdot \widetilde{\left.\bm{D}u\right\rvert_{i}}
\end{equation}
for which $ L_{i}^{\mathcal{L}_{\Omega}}u = f_i $ for all $i: \; \mathbf{x}_i \in \Omega \cup \partial \Omega $.

Formally, we can write the resulting scheme in general form as
\begin{equation}
  \bm{W}_{ji}
  =
  \bm{P}_{ji} \varphi_{ji}
\end{equation}
\begin{equation}
  L^{\bm{\alpha}}_{i} u
  =
  \displaystyle\sum_{j\in\mathcal{N}_i}
  \left( u_{j} \bm{W}_{ji} + \bm{Q}_{ji} \right) \cdot (\mathbf{M}_i^{-1} \bm{C}^{\bm{\alpha}}),
  \quad
  \mathbf{M}_{i}
  =
  \displaystyle\sum_{j\in\mathcal{N}_i}\bm{P}_{ji}\bm{W}_{ji}^{\top} + \mathbf{Q}_i
  \label{eq:cmls_discrete_diff_operator}
\end{equation}
where $\bm{Q}_i \in \mathbb{R}^{p} $ and $ \mathbf{Q}_i \in \mathbb{R}^{p \times p}$ represent the additional terms originating in solved problem.

In \cite{trask2016}, the authors report that this can reduce the number of required neighbors by up to half compared to the standard approach for particles inside the domain.

\subsubsection{Modified Moving Least Squares (MMLS)\label{sec:com_mmls}}

Modified Moving Least Squares (MMLS) introduced by Joldes \cite{joldes2015mmls, joldes2019mmls} as well as CMLS \ref{sec:com_cmls} brings the idea of applying \textit{Tikhonov–Miller} alike regularization to allow the use of higher-order polynomial bases while having less restrictive conditions on the size of the \textit{stencil}, or the number of neighbors needed for interpolation and their distribution compared to the standard MLS method. Additionally, the method introduces constraints for the situation when the classical MLS moment matrix is singular (i.e. there exist infinitely many solutions), that ensure the solution with zero coefficients for the higher order monomials in the basis is favoured. This has the practical effect of decreasing the order of polynomial interpolation, but preserves well-posedness of the problem.

As in the previous case with CMLS, these properties are achieved by modifying the approximation error functional.
However, in contrast to the previous case, in MMLS this modification does not take into account the system of governing equations, but instead the method adds additional terms based on the coefficients of the polynomial base functions.

Assume a modification of the minimization problem \eqref{eq:gmls_l2_optim_problem} to the form
\begin{equation}
  u^*
  =
  \argmin_{q \in \mathcal{V}_h}
  \displaystyle\sum_{j\in\mathcal{N}_i}
  \left[
    \left(
      \lambda_j\left( u\right) - \lambda_j\left(q\right)
    \right)^2
  \right] \varphi_{ji}
  +
  \bm{\mu}_s\cdot
  \left(
    \hat{\bm{\Psi}}_{s,i}
  \right)^2.
  \label{eq:mmls_l2_optim_problem}
\end{equation}
where $ \bm{\mu} $ is a vector containing positive regularization parameters for the additional constraints and $ \hat{\bm{\Psi}}_{s,i} $ is a subset of the coefficients $ \bm{\Phi}_i^{\bm{\alpha}} $ connected with the higher orders terms.
To make this clear, for example using second order basis $ m = 2 $ in 2D, we have
\begin{equation}
  u^*
  =
  \argmin_{q \in \mathcal{V}_h}
  \displaystyle\sum_{j\in\mathcal{N}_i}
  \left[
    \left(
      \lambda_j\left( u\right) - \lambda_j\left(q\right)
    \right)^2
  \right] \varphi_{ji}
  +
  \mu_{x^2}
  \left(
    \Psi_{i,x^2}
  \right)^2
  +
  \mu_{xy}
  \left(
    \Psi_{i,xy}
  \right)^2
  +
  \mu_{y^2}
  \left(
    \Psi_{i,y^2}
  \right)^2
  .
  \label{eq:mmls_l2_optim_problem_example_with_m2}
\end{equation}

By multiplying small positive ($ \ll 1 $) coefficients $ \bm{\mu} $ to the higher order terms $ \hat{\bm{\Psi}}_{s,i} $, in the case the matrix $\mathbf{M}_i$ tends to be singular, the strategy is to obtain a solution with zero coefficients for these higher order terms. By choosing the coefficients small enough it is ensured that the solution remains nearly identical to the standard MLS solution when the moment matrix is not singular. Let us add that there are only a few comments on the specific choice of regularization parameters $ \bm{\mu} $ to be found.
In \cite{joldes2019mmls} the authors use a constant value $ \mu = 10^{-7}$.

The optimization problem \eqref{eq:mmls_l2_optim_problem} is solved in the same way as presented in detail in sections 4, 5 and 6.
For MMLS the moment matrix is obtained in the form
\begin{equation}
  \widetilde{\mathbf{M}}_i
  =
  \displaystyle\sum_{j\in\mathcal{N}_i}
  \bm{X}_{ji} \bm{W}_{ji}^{\top} + \mathbf{H}
  =
  \mathbf{M}_i + \mathbf{H}_{\mu}
  \label{eq:mmls_momentum_matrix}
\end{equation}
where $ \mathbf{H}_{\mu} = \mathrm{diag}( \bm{0}_{p-s}, \bm{\mu}_s ) \in \mathbb{R}^{p\times p} $, with $\bm{0}_{p-s} \in \mathbb{R}^{p-s} $ representing a zero-vector and $ \bm{\mu}_s \in \mathbb{R}^s $ representing the regularization coefficients.
For the ABFs, same as above, arbitrary polynomials $ \bm{P}_{ji} $ multiplied by a suitable radial weight function $ \varphi_{ji} $ can be used.
We end up with the resulting scheme
\begin{equation}
  \bm{W}_{ji}
  =
  \bm{P}_{ji} \varphi_{ji}.
  \label{eq:mmls_general_abfs}
\end{equation}
\begin{equation}
  L^{\bm{\alpha}}_{i} u
  =
  \displaystyle\sum_{j\in\mathcal{N}_i}
  u_{ji} \bm{W}_{ji} \cdot (\widetilde{\mathbf{M}}_i^{-1} \bm{C}^{\bm{\alpha}}),
  \quad
  \widetilde{\mathbf{M}}_{i}
  =
  \displaystyle\sum_{j\in\mathcal{N}_i}\bm{P}_{ji}\bm{W}_{ji}^{\top} + \mathbf{H}_{\mu}.
  \label{eq:mmls_discrete_diff_operator}
\end{equation}
with the symmetric \textit{modified} moment matrix $\widetilde{\mathbf{M}_i}$ built using arbitrary basis functions.
The matrix $\widetilde{\mathbf{M}_i}$ generally has the same rank as the matrix $\mathbf{M}_i$ from the equation \eqref{eq:gmls_m_matrix}, but requires for it fewer neighboring nodes. \\

\subsubsection{Other variants based on MLS: Regularized Moving Least Squares (RMLS), Interpolating Moving Least Squares (IMLS), Improved Interpolating Moving Least Squares (IIMLS), Interpolating Modified Moving Least Squares (IMMLS) and more\label{sec:com_moreOnMMLS}}

Following the introduced methods based on MLS, large number of modifications and improvements exist and some of them are as well adopted as a strong form collocation schemes.
We do not discuss these methods in detail, but at least briefly mention their differences and refer to a more detailed description.

Another regularization techniques were presented as Regularized Moving Least Squares (RMLS) in \cite{wang2018rmls}.
These are close to MMLS but are based directly on the Tikhonov-Miller regularization.
Group of methods with name Interpolating Moving Least Squares (IMLS) introduced by Lancester \cite{lancaster1981} adds the property that the weight functions have to satisfy property $ w_{ji} \approx \delta_{ji} $.
The \textit{interpolating} property was also extended to MMLS, which is named Interpolating Modified Moving Least Squares (IMMLS) \cite{lohit2022immls}.
Another group of methods is called Improved Moving Least Squares (IMLS) \cite{zhang2015immls} or Improved Interpolating Moving Least Squares (IIMLS) \cite{wang2012iimls, li2015iimls}, where the weight functions have the property of delta distributions but any standard radial functions $\varphi_{ji} $ can be used.
In order to avoid ill-conditioned moment matrices, Adaptive Orthogonal IIMLS (AO-IIMLS) with additional extensions was proposed \cite{wang2019aoiimls}.

\subsection{Finite Pointset Method (FPsM)\label{sec:com_FPsM}}

Finite Pointset Method (FPsM) is based on the method of weighted least squares and was first  introduced in \cite{Kuhnert2001} with further works on application of the method to flow problems \cite{Kuhnert2001_2, Tiwari2002, Tiwari2003, Tiwari2005, tiwari2007fpsm, drumm2008finite}.

To construct the approximation of differential operators, similarly to GFDM, both \textbf{polynomial method $\ell_2$ minimization which corresponds to the derivation introduced in section 5} is used. Given a second order differential model equation in 3D:

\begin{equation}
  Au + \bm{B}\cdot \nabla u + C \Delta \psi = f(\bm{x})
  \label{eq:general_second_order_eq}
\end{equation}

where the coeﬃcients $A$, $B$, $C$ and $f(\bm{x}$ are given real valued function.
  To solve this equation with Dirichlet $u = \phi$ or Neumann boundary conditions

  \begin{equation}
    \frac{\partial u}{\partial \bm{n}} = \phi
    \label{eq:general_second_order_eq_Neumann_BC}
  \end{equation}

  We extend the definition of the matrix $\bm{X}_{ji}$ of Taylor monomials defined in the section \ref{sec::com_definition} to add equations \eqref{eq:general_second_order_eq} and \eqref{eq:general_second_order_eq_Neumann_BC} as constraints in the form

  \begin{equation}
    Aa_0 + B_1 a_1 + B_2 a_2 + B_3 a_3 + C(a_4 + a_7 + a_9 ) = f
  \end{equation}

  \begin{equation}
    n_1 a_1 + n_2 a_2 + n_3 a_3 = \phi
  \end{equation}

  where $B = (B_1 , B_2 , B_3 )$, $n = (n_1 , n_2 , n_3 )$ and vector $a$ contains the elements of vector $ \bm{D}u $ defined in \eqref{eq:derivatives_vect}, i.e.

  \begin{equation}
    \begin{aligned}
      a_1 &= \frac{\partial u}{\partial x}\\
      a_2 &= \frac{\partial u}{\partial y}\\
      ...\\
      a_9 &= \frac{\partial^2 u}{\partial^2 z}
    \end{aligned}
  \end{equation}

  The matrix $\tilde{\mathbb{X}}_{i}$ in matrix notation \eqref{eq:appendix-introduction-of-matrix-formalism} becomes (see more on matrix notation \ref{sec:appendix-matrixForm})

  \begin{equation}
    \tilde{\mathbb{X}}_{i} =
    \begin{bmatrix}
      1 & dx_1 & dy_1 & dz_1 & \frac{1}{2}dx_1^2 & dx_1 dy_1 & dx_1 dz_1 & \frac{1}{2} dy_1^2 & dy_1 dz_1 & \frac{1}{2} dz_1^2 \\
      \vdots  & \vdots & \vdots & \vdots & \vdots & \vdots & \vdots & \vdots & \vdots & \vdots \\
      1 & dx_{N_i} & dy_{N_i} & dz_{N_i} & \frac{1}{2}dx_{N_i}^2 & dx_{N_i} dy_{N_i} & dx_{N_i} dz_{N_i} & \frac{1}{2} dy_{N_i}^2 & dy_{N_i} dz_{N_i} & \frac{1}{2} dz_{N_i}^2\\
      A & B_1 & B_2 & B_3 & C & 0 & 0 & C & 0 & C\\
      0 & n_1 & n_2 & n_3 & 0 & 0 & 0 & 0 & 0 & 0
    \end{bmatrix}
    \label{eq:FPsM_matrix_X}
  \end{equation}

  To construct the approximation of differential operators, similarly to GFDM, both \textbf{polynomial method $\ell_2$ minimization which corresponds to the derivation introduced in section 5} was used \cite{drumm2008finite} and also \textbf{derivation based on $\ell_2$ minimization of approximation error} given in section 4 can be found in \cite{Kuhnert2001}. Analytical and numerical comparison of numerical stability of FPsM and GFDM methods is given in \cite{liszka1980finite}.

  Taking \cite{Tiwari2003} as the referential work (although only Laplace operator is considered, it is clear, that this approach can be used for any particular differential operator $d$ ), we assume approximate differential operator
  \begin{equation}
    L^{\bm{\alpha}}_{i} u=
    \displaystyle\sum_{j\in\mathcal{N}_i}
    u_{j}w^{\bm{\alpha}}_{ji}
    \approx
    \bm{C}^{\bm{\alpha}} \cdot \left.\bm{D} u\right\rvert_{i}.
    \label{eq:FPsM_ddo}
  \end{equation}
  The weights are obtained similarly to 4 from the following minimization problem
  \begin{align*}
    \tilde{\mathcal{J}}_{\mathcal{E},i}
    \coloneq
    \displaystyle\sum_{j\in\mathcal{N}_i} \tilde{\varphi}_{ji} \left( \tilde{e}_{ji}^{m+1} \right)^2
    \quad
    \text{with respect to}
    \left.\bm{D}u\right\rvert_{i}.
    \numberthis
    \label{eq:FPsM_mc-optimizationformulation-extended2}
  \end{align*}

  where
  \begin{equation}
    \begin{aligned}
      \tilde{e}_{N_i+1, i}^{m+1} &= Au + \bm{B}\cdot \nabla u + C \Delta \psi - f(\bm{x}_i)\\
      \tilde{e}_{N_i+2, i}^{m+1} &= \frac{\partial u}{\partial \bm{n}} - \phi\\
      \tilde{\varphi}_{N_i+1, i} &= 1\\
      \tilde{\varphi}_{N_i+2, i} &= 1
    \end{aligned}
  \end{equation}

  The final scheme reads as:

  \begin{equation}
    \tilde{\mathbf{M}}_i
    =
    \displaystyle\sum_{j\in\mathcal{N}_i}\tilde{\bm{X}}_{ji}\left(\tilde{\bm{W}}_{ji}\right)^{\top},
    \quad
    \tilde{\bm{W}}_{ji} = \tilde{\bm{X}}_{ji} \tilde{\varphi_{ji}}.
  \end{equation}
  \begin{equation}
    L^{\bm{\alpha}}_{i} u
    =
    \displaystyle\sum_{j\in\mathcal{N}_i}  u_{ji}
    \tilde{\bm{W}}_{ji} \cdot
    \left(
      \left(
        \tilde{\mathbf{M}_i}
      \right )^{-1}
      \bm{C}^{\bm{\alpha}}
    \right)
  \end{equation}
  which is consistent with the methods and procedures presented above.

  \subsection{Least Squares Kinetic Upwind Method (LSKUM)\label{sec:com_lskum}}

  Another methods that employs finite differences on scattered nodes is Least-squares Kinetic Upwind Method (LSKUM) introduced by Ghosh and Deshpande \cite{ghosh1995least,ghosh1995robust} and extended later on for high order schemes as q-LSKUM \cite{deshpande1998, ghosh2007qlskum} (denoted as well as q-KVFS referring to Kinetic Flux Vector Splitting).
  This method itself is not just a tool for approximating derivatives.
  The method was designed originally for obtaining numerical solution of Euler equations of gas dynamics and of Navier–Stokes equations of compressible viscous heat conducting fluid.
  LSKUM uses a relation between the Boltzmann transport equation and the Euler equations for gas dynamics to solve inviscid compressible flow and works with the so-called kinetic flux vector splitting \cite{deshpande1986,mandal1994} to compute the fluxes.
  It comes with a technique for discretizing derivatives at arbitrarily spaced points.
  Here, we focus only on the discretization of the spatial derivatives, and we refer to \cite{ghosh1995least,deshpande2002} for details about how the discrete derivatives are employed.

  LSKUM uses approach \textbf{which $ \ell_2 $ minimization of approximation error} and thus the presented technique to obtain approximation of derivatives is equivalent with approach presented in section 4.
  \textbf{Zeroth moment is not included}.

  Since the method relays on minimization of approximation error of Taylor expansion, introducing LSKUM within the presented frameworks, the ABFs are given by Taylor monomials multiplied by suitable weight function  $ \bm{W}_{ji} = \bm{X}_{ji} \varphi_{ji} $.
  Thus, the method for approxmating the derivatives can be summarized as follows
  \begin{equation}
    \bm{W}_{ji}
    =
    \bm{X}_{ji} \varphi_{ji}.
    \label{eq:lskum_w_taylor_monomials_with_rbf_recap}
  \end{equation}
  \begin{equation}
    L^{\bm{\alpha}}_{i} u
    =
    \displaystyle\sum_{j\in\mathcal{N}_i}
    u_{ji} \bm{W}_{ji} \cdot (\mathbf{M}_i^{-1} \bm{C}^{\bm{\alpha}}),
    \quad
    \mathbf{M}_{i}
    =
    \displaystyle\sum_{j\in\mathcal{N}_i}\bm{X}_{ji} \bm{W}_{ji}^{\top}.
    \label{eq:lskum_discrete_diff_operator_recap}
  \end{equation}
  To recover the expressions for first and second order from original paper \cite{ghosh1995least} (where the authors symbolically expressed the solution), equation \eqref{eq:lskum_discrete_diff_operator_recap} can be used with $ m = 1 $ and $ m=2 $ for the first and second order respectively (see section \ref{sec:com_definition}).

  For the needs of LSKUM, not all the particles in $ \mathcal{N}_i $ are used to obtain the derivatives.
  \footnote{
    Derivatives of fluxes in the context of the method.
  }
  Since the method builds on the upwind scheme, only half of the neighboring nodes of the considered node $ i $ are used along each of the coordinate axes $ k $.
  We denote this truncated neighborhood as $\mathcal{N}_i^{\beta\pm}$ where $ \beta $ stands for a direction along the coordinate axes, i.e. $\beta=x, y, \dots n$ and the plus-minus sign denotes which half of the neighborhood is used.
  This indicates that the propagation of information from the neighboring nodes $ j $ to the receiving node $ i $ depends on the location $ j $ in relation to $ i $ and the signs of the imposed velocity $v_x, v_y \dots v_z$
  \footnote{
    Molecular velocity in the context of the method
  }
  , thus ensuring the inclusion in the computation of the derivative operator of only those neighbors that influence the solution of the advection equation at $ i $.

  \subsection{Kinetic Meshless Method (KMM)\label{sec:com_kkm}}

  Kinetic Meshless Method is largely based on the LSKUM. It utilizes finite differences on scattered nodes in an effort to improve accuracy and robustness of the original method.
  KMM was introduced by Praveen and Deshpande \cite{Praveen2003, praveen2007kinetic} and similarly to LSKUM, the method was designed to solve compressible inviscid flows and discretization of derivatives is only one particular step in the procedure (see LSKUM section \ref{sec:com_lskum} and the original papers \cite{Praveen2003, praveen2007kinetic}.
    Approximation of derivatives \textbf{is based on $ \ell_2 $ minimization of approximation error} presented in section $ 4 $ and \textbf{zeroth moment is not included}.

    The approach to approximate derivatives is the same as we derived in section 4 or described in LSKUM  section \ref{sec:com_lskum}.
    However, the procedure presented in the KMM paper is called dual least-squares (DLS) method. The idea is to use a modified least squares scheme in which mid-point function values $ u_{j/2} := u( ( \bm{x}_j - \bm{x}_i) / 2 ) $ (i.e. fluxes) are used to introduce upwind convective stabilization. These mid-point values are used instead of the values $  u_j $ to approximate derivatives at the position $ \bm{x}_i $ by replacing $ u_{ji} =  u_{j} -  u_{i} $ in the Taylor expansion with $ u_{ji/i} =  u_{j/2} -  u_{i} $, thus avoiding the flux splitting employed before in the LSKUM. Hence, similarly to all methods introduced before it uses all nodes from the surroundings $ \mathcal{N}_i $.

    Additionally, in comparison to the (q-)LSKUM, rotational invariance due to the DLS approximation is proven that simplifies setting up boundary conditions \cite{praveen2007kinetic}.

    To introduce the upwind convective stabilization, the values of $ u_{j/2} $ take the following approximation:

    \begin{equation}
      u_{ij/2} =
      \begin{cases}
        u_i, & \quad \bm{v} \cdot \hat{\bm{e}}_{ji} \ge 0 \\
        u_j, & \quad \bm{v} \cdot \hat{\bm{e}}_{ji} < 0\, ,
      \end{cases}
    \end{equation}

    where $ \mathbf{v} $ is the imposed velocity, $ \hat{\bm{e}} $ is a unit vector directed from $ \bm{x}_i $ to $ \bm{x}_j $.

    Due to the minimization of approximation error, the structure of ABFs is given by Taylor monomials multiplied by suitable weight function  $ \bm{W}_{ji} = \bm{X}_{ji/2} \varphi_{ji} $ which leads to approximation formula.
    \begin{equation}
      \bm{W}_{ji}
      =
      \bm{X}_{j/2i} \varphi_{ji}.
      \label{eq:kmm_w_taylor_monomials_with_rbf_recap}
    \end{equation}
    \begin{equation}
      L^{\bm{\alpha}}_{i} u
      =
      \displaystyle\sum_{j\in\mathcal{N}_i}
      (u_{ji/2} - u_{i})
      \bm{W}_{ji} \cdot (\mathbf{M}_i^{-1} \bm{C}^{\bm{\alpha}})\, ,
      \quad
      \mathbf{M}_{i}
      =
      \displaystyle\sum_{j\in\mathcal{N}_i}\bm{X}_{ji/2} \bm{W}_{ji}^{\top}\, .
      \label{eq:kmm_discrete_diff_operator_recap}
    \end{equation}

    Compared to the general approach taken in section 4, the method modifies the matrix $ \mathbf{M} $ and the vector $ \bm{X} $ to account for the halved node distance.

    The method described so far corresponds in the original work \cite{praveen2007kinetic} to a kinetic least-squares formulation as opposed to a least-squares formulation at the macroscopic level. The latter is also derived and shown to have a higher order of approximation, but more prone to ill-conditioning.

    For the higher order approximation at the kinetic level, KMM utilizes the scheme from the seminal series {\it Towards the ultimate conservative difference scheme.}
    \cite{vanleer1977,vanLeer1979muscl,vanleer1977_2} that form a one-parameter family known as the {\it k-schemes}, after the parameter $ k $, with $ k = -1 $ generating the fully one-sided second-order scheme to interpolate variables from neighboring particles to the middle points $  u_{ji/2} $ where the flux is evaluated.
    This interpolation uses derivatives computed by the standard approximation given by equation \eqref{eq:kmm_discrete_diff_operator_recap}.

    \subsection{Particle Difference Method (PDM)\label{sec:com_pdm}}

    Particle Difference Method (PDM) \cite{yoon2014extended} (sometimes also called Particle Derivative Approximation - PDA) is derived \textbf{using the $\ell_2$ minimization of approximation error and exactly corresponds to the derivation introduced in section 4.} In this case, zeroth moment is included and the resulting approximation is obtained by minimizing functional $ \mathcal{J}_{\mathcal{E},i} $ and solving the optimization problem. Using this approach, all derivatives are obtained simultaneously

    \begin{equation}
      \left.\bm{D} u\right\rvert_{i}
      =
      \displaystyle\sum_{j\in\mathcal{N}_i} \varphi_{ji}   u_{ji}
      \left(
        \mathbf{M}_i^{-1}
        \bm{X}_{ji}
      \right)
      \label{eq:pdm-approximation-all-derivatives}
    \end{equation}
    \begin{equation}
      \mathbf{M}_i
      =
      \displaystyle\sum_{j\in\mathcal{N}_i}\bm{X}_{ji}\bm{W}_{ji}^{\top},
      \quad
      \bm{W}_{ji} = \bm{X}_{ji} \varphi_{ji}
      .
    \end{equation}

    To obtain approximation of particular derivative in the form of (2), the\textit{mapping vector} $\bm{C}^{\bm{\alpha}}$ is employed
    \begin{equation}
      L^{\bm{\alpha}}_{i} u
      =
      \displaystyle\sum_{j\in\mathcal{N}_i}  u_{ji}
      \bm{W}_{ji} \cdot
      \left(
        \mathbf{M}_i^{-1}
        \bm{C}^{\bm{\alpha}} \cdot
      \right).
    \end{equation}

    An important contribution of \cite{yoon2014extended} is the extension of this approach to model weak and strong discontinuities - so called EPDM or EPDA. Note also that the paper considers non-smooth $\varphi_{ji}$ functions as opposed to the standard \textit{bell shape} functions.

    \subsection{Finite Difference Particle Method (FDPM)\label{sec:com_fdpm}}

    Finite Difference Particle Method (FDMP) \cite{wang2018lagrangian} is derived \textbf{using the $\ell_2$ minimization of approximation error and exactly corresponds to the derivation introduced in section 4.} The authors essentially present this method as a GFDM-based scheme. To keep our review consistent, zeroth moment is not included, which means that the resulting scheme is

    \begin{equation}
      \left.\bm{D} u\right\rvert_{i}
      =
      \displaystyle\sum_{j\in\mathcal{N}_i} \varphi_{ji}   u_{ji}
      \left(
        \mathbf{M}_i^{-1}
        \bm{X}_{ji}
      \right)
      \label{eq:fdpm-approximation-all-derivatives}
    \end{equation}
    \begin{equation}
      \mathbf{M}_i
      =
      \displaystyle\sum_{j\in\mathcal{N}_i}\bm{X}_{ji}\bm{W}_{ji}^{\top},
      \quad
      \bm{W}_{ji} = \bm{X}_{ji} \varphi_{ji}.
    \end{equation}

    To obtain approximation of particular derivative in the form of \eqref{eq:general_discrete_diff_operator_zm}, the\textit{mapping vector} $\bm{C}^{\bm{\alpha}}$ is employed
    \begin{equation}
      L^{\bm{\alpha}}_{i} u
      =
      \displaystyle\sum_{j\in\mathcal{N}_i}  u_{ji}
      \bm{W}_{ji} \cdot
      \left(
        \mathbf{M}_i^{-1}
        \bm{C}^{\bm{\alpha}} \cdot
      \right).
    \end{equation}

    The work of \cite{wang2018lagrangian} extends the standard scheme by introducing a dynamic method to update variable smoothing length so that the number of the neighbor particles remains relatively constant. In turn it makes the parameter $ h_i $ to the \textit{radial basis function} $ \varphi_{ij} $ vary spatially, and the function becomes anisotropic.

    \subsection{Meshless Finite Difference Method (MFDM)\label{sec:com_mfdm}}

    Meshless Finite Difference Method MFDM is not a new variant of inter-network MFD methods but rather an umbrella term for a class of methods using the same principles of use by Milewsky in \cite{milewski2012meshless}.
    To cover all the various occurrences of the same or similar methods under different names, we decided to include MFDM aswell, although it is not presented as o novel method.
    The article summarize discretization of differential operators and their use in methods with strong and weak formulation.

    Considering the strong form method, MFDM is similar to GFDM.
    The method is derived \textbf{using $\ell_2$ minimization of approximation error} which we presented in detail in section 4.
    \textbf{The zeroth moment is included}.
    Moreover, high order extensions are presented by MFDM, which are done by canceling higher order terms up to certain order $ m $.
    However, MFDM brings different approach to high order corrections \cite{milewski2021mfdm, milewski2022mfdm}, which are not done through simply zeroing out of higher order terms as we presented above.

    MFDM is presented in vector formalism which is similar to what we introduced in opening section 2.
    Starting from the formulation of the minimization problem, similarly to the summary presented in GFDM section \ref{sec:com_gfdm}, we obtain ABFs given by Taylor monomials $ \bm{X}_{ji} $ multiplied appropriate RBF weight function $ \varphi_{ji} $, which originates in minimization of Taylor expansion approximation error.
    Following derivation from section 4, we end up with approximation of all derivatives up to order of $ m $
    \begin{equation}
      \bm{W}_{ji}
      =
      \bm{X}_{ji} \varphi_{ji}\, ,
      \label{eq:mfdm_w_taylor_monomials_with_rbf}
    \end{equation}
    we can write
    \begin{equation}
      \left.\bm{D} u\right\rvert_{i}
      =
      \displaystyle\sum_{j\in\mathcal{N}_i}\mathbf{M}_{i}^{-1}\bm{W}_{ji} u_{ji},
      \quad
      \mathbf{M}_{i}
      =
      \displaystyle\sum_{j\in\mathcal{N}_i}\bm{X}_{ji}\bm{W}_{ji}^{\top}.
      \label{eq:mfdm_linSystem_with_w}
    \end{equation}
    To express approximation of particular differential operator $ d $, we include the \textit{mapping vector} to select desired derivatives
    \begin{equation}
      \bm{C}^{\bm{\alpha}} \cdot \left.\bm{D} u\right\rvert_{i}
      =
      \displaystyle\sum_{j\in\mathcal{N}_i}
      \bm{C}^{\bm{\alpha}}
      \cdot
      \left(
        \mathbf{M}_{i}^{-1}\bm{W}_{ji}
      \right)
      u_{ji},
      \label{eq:mfdm_selec_derivative}
    \end{equation}
    and thanks to symmetry of moment matrix $ \mathbf{M}_i $, we can summarize the method as
    \begin{equation}
      \bm{W}_{ji}
      =
      \bm{X}_{ji} \varphi_{ji}
      \label{eq:mfdm_w_taylor_monomials_with_rbf_recap}
    \end{equation}
    \begin{equation}
      L^{\bm{\alpha}}_{i} u
      =
      \displaystyle\sum_{j\in\mathcal{N}_i}
      u_{ji} \bm{W}_{ji} \cdot (\mathbf{M}_i^{-1} \bm{C}^{\bm{\alpha}}),
      \quad
      \mathbf{M}_{i}
      =
      \displaystyle\sum_{j\in\mathcal{N}_i}\bm{X}_{ji}\bm{W}_{ji}^{\top}.
      \label{eq:mfdm_discrete_diff_operator_recap}
    \end{equation}

    Although this approach is common to all methods presented so far, in order to obtain high order approximation, MFDM takes different strategy then just simple cancel higher order terms by setting approximation order $ m $ higher then is the order of approximated derivative $L = |\bm{\alpha}| $.
    Approach to treat high order approximations is introduced \cite{milewski2012meshless} and more recently in \cite{milewski2021mfdm, milewski2022mfdm}.
    High order approximation $ L^{\bm{\alpha},>L}_{i} $ denoted by superscript $>L$ is obtained by introducing correction $ \Delta^{\bm{\alpha}} $ to lower order approximation $ L^{\bm{\alpha},{\leq L}}_{i} $ denoted by superscript ${\leq L}$
    \begin{equation}
      L^{\bm{\alpha},>L}_{i} = L^{\bm{\alpha},{\leq L}}_{i} - \Delta_{i}^{\bm{\alpha}}.
    \end{equation}
    The $ \Delta_{i}^{\bm{\alpha}} $ correction factor is obtained as follows.
    Assuming Taylor expansion \eqref{eq:taylor_series_zm}, we can split the terms to lower order and higher order
    \begin{equation}
      u_{j}
      =
      \bm{X}_{ji}^{{\leq L}}\cdot\left.\bm{D}^{\leq L} u\right\rvert_{i} +
      \bm{X}_{ji}^{>L}\cdot\left.\bm{D}^{>L} u\right\rvert_{i} +
      e_{ji}^{m+1},
      \label{eq:mfmd_taylor_series_split}
    \end{equation}
    where
    \begin{equation}
      \bm{X}_{ji}^{{\leq L}}
      =
      \left\{
        x_{1,ji}^{\alpha_1} \ldots x_{n,ji}^{\alpha_n} / \bm{\alpha}!
      \right\}_{|\bm{\alpha}|=0}^L,
      \quad
      \bm{X}_{ji}^{>L}
      =
      \left\{
        x_{1,ji}^{\alpha_1} \ldots x_{n,ji}^{\alpha_n} / \bm{\alpha}!
      \right\}_{|\bm{\alpha}|>L}^m,
    \end{equation}
    \begin{equation}
      \bm{D}^{\leq L} u
      =
      \left\{
        \partial^{\bm{\alpha}}  u
      \right\}_{\bm{|\alpha|}=0}^L,
      \quad
      \bm{D}^{>L} u
      =
      \left\{
        \partial^{\bm{\alpha}}  u
      \right\}_{\bm{|\alpha|}>L}^m.
      \label{eq:mfdm_low_high_order_derivatives}
    \end{equation}
    In the work \cite{milewski2022mfdm} author states that usually $L = |\bm{\alpha}|$ and in most cases $ m = 2L $.

    Next, partitioning to lower and higher order terms is introduced as well in optimization problem \eqref{eq:l2e:emin-optimizationformulation} to minimize approximation error
    \begin{align*}
      \widetilde{\mathcal{J}}_{\mathcal{E},i}
      =
      \displaystyle\sum_{j\in\mathcal{N}_i}
      \varphi_{ji} \left( e_{ji}^{m+1} \right)^2
      =
      \displaystyle\sum_{j\in\mathcal{N}_i}
      \varphi_{ji}
      \left(
        \bm{X}_{ji}^{{\leq L}}\cdot\left.\bm{D}^{\leq L} u\right\rvert_{i}
        + \bm{X}_{ji}^{>L}\cdot\left.\bm{D}^{>L} u\right\rvert_{i}
        -  u_{ji}
      \right)^2
      \\
      \mathrm{min} \left( \mathcal{J}_{\mathcal{E},i} \right)
      \quad
      \text{with respect to}
      \left.\bm{D}^{\leq L} u\right\rvert_{i}.
      \numberthis
      \label{eq:mfdm_emin-optimizationformulation}
    \end{align*}
    Following formal minimization $\partial \mathcal{J}_i / \partial \left.\bm{D}^{\leq L} u\right\rvert_{i} = 0 $  performed in section 4 we obtain following result
    \begin{equation}
      \left.\bm{D}^{\leq L} u\right\rvert_{i}
      =
      \displaystyle\sum_{j\in\mathcal{N}_i}\mathbf{M}_{i}^{-1}\bm{W}_{ji} u_{j} - \bm{\Delta}_{i}\, ,
      \label{eq:}
    \end{equation}
    where the first term on the right hand side corresponds to the expression for derivatives above \eqref{eq:mfdm_linSystem_with_w} and
    \begin{equation}
      \bm{\Delta}_i
      = \mathbf{M}_{i}^{-1}\bm{W}_{ji}
      \left(
        \bm{X}_{ji}^{>L}\cdot\left.\bm{D}^{>L} u\right\rvert_{i}
      \right)\, ,
    \end{equation}
    groups all the additional terms arising from the presence of unknown higher order derivatives $ \left.\bm{D}^{>L} u\right\rvert_{i} $.
    $ \Delta_{i} $ represents a vector of correction terms, each element corrects the corresponding derivative in $ \bm{\Delta}_{i} $.
    In case we are interested only in particular derivative $ d $, we proceed in the same way as in equation \eqref{eq:mfdm_selec_derivative} using \text{mapping vector} $ \bm{C}^{\bm{\alpha}} $.
    The problem of correction terms $ \bm{\Delta}_i $ is that it contains unknown vector of higher order derivatives $ \left.\bm{D}^{>L} u\right\rvert_{i} $.
    High order derivatives $ \left.\bm{D}^{>L} u\right\rvert_{i} $ can be expressed by appropriate composition law of low order derivatives
    \begin{align*}
      \left.\bm{D}^{>L} u\right\rvert_{i}
      &=
      \left(
        \bm{D}^{\leq L}\mathbf{P}_{t}
      \right)
      \left(
        \left.\bm{D}^{\leq L} u\right\rvert_{i}\mathbf{P}_{u}
      \right)
      \\
      &=
      \displaystyle\sum_{j\in\mathcal{N}_i}
      \left(
        \mathbf{M}_{i}^{-1}\bm{W}_{ji}\mathbf{P}_{t}
      \right)
      \left(
        \left.\bm{D}^{\leq L} u\right\rvert_{j}\mathbf{P}_{u}
      \right)
      \\
      &=
      \displaystyle\sum_{j\in\mathcal{N}_i}
      \left(
        \mathbf{M}_{i}^{-1}\bm{W}_{ji}\mathbf{P}_{t}
      \right)
      \left(
        \displaystyle\sum_{s}\mathbf{M}_{i}^{-1}\bm{W}_{sj} u_{s}\mathbf{P}_{u} - \bm{\Delta}_{j}\mathbf{P}_{u}
      \right)\, ,
    \end{align*}
    where the matrices $\mathbf{P}_{t}$ and $\mathbf{P}_{u}$ are non-unique permutation matrices that fulfill the property $p - t - u = 0$ for every higher-order approximation $p = | \bm{\alpha} |$ and the lower-order derivatives $t = | \bm{\alpha}^{'} |$ and $u = | \bm{\alpha}^{''} |$. The permutation matrices select $\bm{\alpha}^{'}$ and $\bm{\alpha}^{''}$ from the  corresponding rows of the matrix $\mathbf{M}_{i}$.

    However, this resulting formula again depends on the correction terms $ \bm{\Delta}_{j} $.
    The problem is solved by employing iterative procedure
    \begin{equation}
      \left( \left.\bm{D}^{>L} u\right\rvert_{i} \right)^{ (k + 1) }
      =
      \displaystyle\sum_{j\in\mathcal{N}_i}
      \left(
        \mathbf{M}_{i}^{-1}\bm{W}_{ji}\mathbf{P}_{t}
      \right)
      \left(
        \displaystyle\sum_{s}\mathbf{M}_{i}^{-1}\bm{W}_{sj} u_{s}\mathbf{P}_{u} - \bm{\Delta}_{j}^{(k)}\mathbf{P}_{u}
      \right)\, ,
    \end{equation}
    with the initial solution corresponding to low order approximation $ \bm{\Delta}^{(0)} = \left.\bm{D}^{\leq L} u\right\rvert_{i} $.

    For further details as well as for treatment of boundary conditions and handling high order correction for the boundary conditions, we refer to the original works \cite{milewski2012meshless, milewski2022mfdm}.
    Last, we should mention, that MFDM is presented mostly in matrix formalism.

    \subsubsection*{Notation remark}

    As we already mentioned, the way how is MFDM written is close to the notation used through the paper.
    We would like to only point out that original works use $ \mathbf{P} $ to denote vector of Taylor monomials which we denote by $ \bm{X}_{ji} $.
    Moreover, original paper use matrix $ \mathbf{M} $ to denote matrix of coefficients of all derivatives instead of moment matrix (there is no any particular symbol for a moment matrix on its own).
    Thus, $ \mathbf{M} $ from the original papers corresponds to our $ \mathbb{C}_{i}^{\top} $ introduced in the section 2.

    \subsection{Lagrangian Differencing Dynamics (LDD)\label{sec:com_ldd}}

    In contrast to all methods presented so far, Lagrangian Differencing Dynamics (Originally Proposed as Renormalised meshless operator \cite{basic2018class}, then renamed to Lagrangian finite difference method \cite{basic2020lagrangian} and finally renamed to Lagrangian differencing dynamics \cite{peng2021lagrangian}) does not consider general formula to approximate general differential operator, but it focuses only on gradient and Laplace operator approximations providing second order formulas.
    The motivation behind LDD is to find suitable compromise between SPH/MPS-like approximations, which do not provide sufficient accuracy and accurate collocation schemes, which are computationally more expensive.
    This is achieved only with the help of a correction matrix, known for example from SPH.
    In our notation, this is $ \mathbf{M}_i $ for $ m = 1 $, where zeroth terms are not considered.
    For the purpose of this section, we denote this matrix as $\mathbf{M}_{i}^{1}$, to highlight that only linear terms are present.
    Together with the symmetry of the weight function $\varphi_{ji}$, this allows us to construct the desired compromise approximation.

    At this point, short comment about the weight function $ \varphi $ has to be made.
    In the original work \cite{basic2018class} \textit{absolutely normalized} weight function $ \psi_{ji} = \varphi_{ji} / \sum_{j\in\mathcal{N}_i} \varphi_{ji} $ is used.
    This function acts the same was as the $ \varphi_{ji} $ introduced above.
    However, in more recent works \cite{basic2020lagrangian, peng2021lagrangian, basic2022ldd} authors operates with the standard normalized weight function $ \varphi_{ji} $.
    To keep consistency with the notation introduced above, we adopt the approach from recent works and weight function $ \varphi_{ji} $ through the expressions.

    For the purpose of derivation LDD, we write Taylor expansion with split second derivatives
    \begin{equation}
      u_{i}
      =
      u_{j} +
      \vect{x}_{ji} \cdot \left.\nabla  u\right\rvert_{i} +
      \frac{1}{2} \left.\bm{D}^{2} u\right\rvert_{i} \cdot  \bm{X}_{ji}^{2} + \left.\bm{D}^{2,\textrm{mixed}} u\right\rvert_{i} \cdot  \bm{X}_{ji}^{2,\textrm{mixed}} + e_{ji}^{3}
      \label{eq:ldd-taylorExpansionSplitedDerivatives}
    \end{equation}
    where $ e_{ji}^{3} $ is third order remainder term and where we divided $ \bm{D} u $ and $ \bm{X}_{ji} $ to sub-components
    \footnote{
      To make the procedure clear, assuming two dimensions, the sub-components of $ \bm{D} u $ and $ \bm{X}_{ji} $ have form
      \begin{equation}
        \bm{D}^{2} u
        =
        \begin{bmatrix}
          \frac{\partial^{2} u}{\partial{x}^{2}},&
          \frac{\partial^{2} u}{\partial{y}^{2}}
        \end{bmatrix}^{\top},
        \quad
        \bm{D}^{2,\textrm{mixed}} u
        =
        \begin{bmatrix}
          \frac{\partial^{2} u}{\partial{x}\partial{y}},&
          \frac{\partial^{2} u}{\partial{y}\partial{x}}
        \end{bmatrix}^{\top}
      \end{equation}
      \begin{equation}
        \bm{X}_{ji}^{2}
        =
        \begin{bmatrix}
          x_{ji}^2,&
          y_{ji}^2
        \end{bmatrix}^{\top},
        \quad
        \bm{X}_{ji}^{2,\textrm{mixed}}
        =
        \begin{bmatrix}
          x_{ji} y_{ji},&
          y_{ji} x_{ji}
        \end{bmatrix}^{\top}.
      \end{equation}
    }
    \begin{equation}
      \bm{D}^{2} u
      =
      \left\{
        \frac{\partial^{2} u}{\partial{x}^{2}_{\beta}}
      \right\}_{\beta=1}^{n}
      \quad
      \bm{D}^{2,\textrm{mixed}} u
      =
      \left\{
        \frac{\partial^{2} u}{\partial{x}_{\alpha}\partial{x}_{\beta}}
      \right\}_{\alpha, \beta =1, \alpha\neq\beta}^{n}
      \label{eq:ldd-selectedDerivativesAndMonomials}
    \end{equation}
    \begin{equation}
      \bm{X}_{ji}^{2}
      =
      \left\{
        x_{\beta,ji}^2
      \right\}_{\beta=1}^{n}
      \quad
      \bm{X}_{ji}^{2,\textrm{mixed}}
      =
      \left\{
        x_{\alpha,ji} x_{\beta,ji}
      \right\}_{\alpha, \beta =1, \alpha\neq\beta}^{n}
      \label{eq:ldd-selectedDerivativesAndMonomials}
    \end{equation}
    where $ \bm{D}^{2} u \in \mathbb{R}^{n} $, $ \bm{D}^{2,\textrm{mixed}} u \in \mathbb{R}^{\widetilde{p}} $ and $ \bm{X}_{ji}^{2} \in \mathbb{R}^{n} $,  $ \bm{X}_{ji}^{2,\textrm{mixed}} \in \mathbb{R}^{\widetilde{p}} $ with $ \widetilde{p} = \sum_{k=1}^{n-1} k $.

    Considering identity vector $ \bm{I} \in \mathbb{R}^{n} $, the Laplace operator can be written as
    \begin{equation}
      \left.\nabla^2  u\right\rvert_{i} = \left.\bm{D}^{2} u\right\rvert_{i} \cdot \bm{I}.
      \label{eq:ldd-taylorExpansionSelectedLaplace}
    \end{equation}
    Isolating the second derivatives $ \left.\bm{D}^{2} u\right\rvert_{i} $ from Taylor expansion \eqref{eq:ldd-taylorExpansionSplitedDerivatives}
    \begin{equation}
      \left.\bm{D}^{2} u\right\rvert_{i} \cdot \bm{X}_{ji}^{2}
      =
      2
      \left(
        u_{ji}
        -
        \vect{x}_{ji} \cdot \left.\nabla u\right\rvert_{i}
        -
        \left.\bm{D}^{2,\textrm{mixed}} u\right\rvert_{i}
        \cdot
        \bm{X}_{ji}^{2,\textrm{mixed}}
        -
        e_{ji}^{3}
      \right)
      \label{eq:ldd-taylorExpansionExpressedSecondDerivatives}
    \end{equation}
    there are three ways how we can construct the approximation.

    \subsubsection{Laplacian approximation - naive version}

    In order to obtain Laplacian, equation \eqref{eq:ldd-taylorExpansionExpressedSecondDerivatives} is multiplied by $ \bm{I} \cdot \bm{I} / \lVert \vect{x}_{ji} \rVert^2$, which turns the left hand side into \eqref{eq:ldd-taylorExpansionSelectedLaplace}
    \begin{equation}
      \left.\nabla^2  u\right\rvert_{i}
      =
      \frac{2n}{ \lVert \vect{x}_{ji} \rVert^2 }
      \left(
        u_{ji}
        -
        \vect{x}_{ji} \cdot \left.\nabla u\right\rvert_{i}
        -
        \left.\bm{D}^{2,\textrm{mixed}} u\right\rvert_{i}
        \cdot
        \bm{X}_{ji}^{2,\textrm{mixed}}
        -
        e_{ji}^{3}
      \right),
      \label{eq:ldd-naiveLaplacianPairwise}
    \end{equation}
    where the dimension factor $ n $ is product of $ n = \bm{I} \cdot \bm{I} $.
    Multiplying equation \eqref{eq:ldd-naiveLaplacianPairwise} with the weight function $ \varphi_{ji} $ and performing summation over all neighboring particles, due to symmetry of $ \varphi_{ji} $, the term with mixed derivatives is canceled and dropping the Taylor expansion remainder, we obtain formula
    \begin{equation}
      L_{i}^{\nabla^2,\textrm{naive}} u
      =
      2n
      \displaystyle \sum_{j\in\mathcal{N}_i}
      \frac{\varphi_{ji}}{ \lVert \vect{x}_{ji} \rVert^2 }
      \left(
        u_{ji}
        -
        \vect{x}_{ji} \cdot L_i^{\nabla} u
      \right).
      \label{eq:ldd_naive}
    \end{equation}
    As show in \cite{basic2018class}, the error of approximation \eqref{eq:ldd_naive} contains error of gradient approximation $ L_i^{\nabla} u $ and truncated higher order terms from the Taylor series.
    Moreover as the authors highlights, to obtain sufficient accuracy, high number of neighbors of utmost regular distribution is required.

    \subsubsection{Laplacian approximation - sum version}
    The naive approximation is strongly affected by the accuracy of approximation of gradient $ L_i^{\nabla} u $.
    Substituting the first order approximation of gradient obtained from (2) (which corresponds to LDD approximation of gradient) directly to the Taylor expansion \eqref{eq:ldd-taylorExpansionSelectedLaplace}, we obtain following expression
    \begin{equation}
      u_j
      =
      u_i +
      \vect{x}_{ji} \cdot \left(\mathbf{M}_{i}^{1}\right)^{-1}
      \displaystyle \sum_{j\in\mathcal{N}_i} \varphi_{ji} \vect{x}_{ji}
      \left(
        u_{ji}
        -
        \frac{1}{2} \left.\bm{D}^{2} u\right\rvert_{i} \cdot  \bm{X}_{ji}^{2}
        - e_{ji}^{3}
      \right)
      +
      \frac{1}{2} \left.\bm{D}^{2} u\right\rvert_{i} \cdot  \bm{X}_{ji}^{2}
      +
      e_{ji}^{3}
    \end{equation}
    Reshuffling the terms in order to group the terms with second derivatives $ \left.\bm{D}^{2} u\right\rvert_{i}$ on one side, we obtain
    \begin{align*}
      \left.\bm{D}^{2} u\right\rvert_{i} \cdot  \bm{X}_{ji}^{2}
      -
      \vect{x}_{ji} \cdot \left(\mathbf{M}_{i}^{1} \right)^{-1}
      \displaystyle \sum_{j\in\mathcal{N}_i} \varphi_{ji} \vect{x}_{ji}
      \left(
        \left.\bm{D}^{2} u\right\rvert_{i} \cdot  \bm{X}_{ji}^{2}
      \right)
      = \\
      2
      \left(
        u_{ji}
        -
        \vect{x}_{ji} \cdot \left(\mathbf{M}_{i}^{1} \right)^{-1}
        \displaystyle \sum_{j\in\mathcal{N}_i}
        \varphi_{ji} \vect{x}_{ji}
        \left(
          u_{ji} - e_{ji}^{3}
        \right)
        -
        e_{ji}^{3}
      \right).
      \numberthis \label{eq:ldd-sumlaplacian-expressedLaplacianTerms}
    \end{align*}
    In order to obtain the Laplacian, we multiply the equation \eqref{eq:ldd-sumlaplacian-expressedLaplacianTerms} by $ \bm{I} \cdot \bm{I} $.
    Since $ \bm{X}_{ji}^{2} \cdot{I}  = \lVert \vect{x}_{ji} \rVert ^2 $, omitting the Taylor expansion reminders we have
    \begin{equation}
      \nabla^2  u_i
      \left(
        \lVert \vect{x}_{ji} \rVert ^2
        -
        \vect{x}_{ji} \cdot \left(\mathbf{M}_{i}^{1}\right)^{-1}
        \displaystyle \sum_{j\in\mathcal{N}_i}
        \varphi_{ji} \vect{x}_{ji} \lVert \vect{x}_{ji} \rVert ^2
      \right)
      \approx
      2n
      \left(
        u_{ji}
        -
        \vect{x}_{ji} \cdot \left(\mathbf{M}_{i}^{1}\right)^{-1}
        \displaystyle \sum_{j\in\mathcal{N}_i}
        \varphi_{ji}
        \vect{x}_{ji}
        u_{ji}
      \right).
    \end{equation}
    Next, we multiply both sides by the weight function $ \varphi_{ji} $ and sum over the neighbors.
    Due to normalization of the weight functions (summing weight function over all particles is approximately one), the value of Laplacian remains unchanged.
    Reordering the terms once again results in
    \begin{equation}
      \nabla^2  u_i
      \left(
        \displaystyle \sum_{j\in\mathcal{N}_i}
        \varphi_{ji} \lVert \vect{x}_{ji} \rVert ^2
        -
        \left(
          \displaystyle \sum_{j\in\mathcal{N}_i}
          \varphi_{ji} \lVert \vect{x}_{ji} \rVert ^2
        \right)
        \cdot \left(\mathbf{M}_{i}^{1}\right)^{-1}
        \displaystyle \sum_{j\in\mathcal{N}_i}
        \psi_{ki} \vect{x}_{ki} \lVert \vect{x}_{ki} \rVert ^2
      \right)
      \approx
      2n
      \displaystyle \sum_{j\in\mathcal{N}_i}
      \varphi_{ji} u_{ji}
      \left(
        1
        -
        \vect{x}_{ji} \cdot \left(\mathbf{M}_{i}^{1}\right)^{-1}
        \bm{O}_i
      \right)
      \label{eq:ldd-sumLaplacian-multipliedAndSummed}
    \end{equation}
    where $ \bm{O}_i $ is so called $\textit{offset} $ vector related with particle $ i $
    \begin{equation}
      \bm{O}_i =
      \displaystyle \sum_{j\in\mathcal{N}_i}
      \varphi_{ji}
      \vect{x}_{ji},
    \end{equation}
    which points from the position of particle $\vect{x}_i $ to the \textit{weighted center} of the neighborhood where the distribution of neighbor points dominates and which is zero for regular distribution of neighbors.

    It can be seen from the equation \eqref{eq:ldd-sumLaplacian-multipliedAndSummed} that the resulting shape of the Laplace approximation is is
    \begin{equation}
      L_{i}^{\nabla^2,\textrm{sum}} u
      =
      2n
      \frac{
        \displaystyle \sum_{j\in\mathcal{N}_i}
        \varphi_{ji}  u_{ji}
        \left(
          1 - \vect{x}_{ji} \cdot \left(\mathbf{M}_{i}^{1}\right)^{-1}\bm{O}_i
        \right)
      }
      {
        \displaystyle \sum_{j\in\mathcal{N}_i}
        \varphi_{ji}
        \lVert \vect{x}_{ji} \rVert^2
        \left(
          1 - \vect{x}_{ji} \cdot \left(\mathbf{M}_{i}^{1}\right)^{-1}\bm{O}_i
        \right)
      }.
      \label{eq:ldd_sum}
    \end{equation}

    Assuming regular distribution ($\bm{O} = \mathbf{0}$) and only closest neighbors are taken into account, approximation \eqref{eq:ldd_sum} reduces to central difference.
    The offset vector $\bm{O}$ inside the sums acts as correction, which within the scalar products \textit{virtually corrects} relative position of each neighborhood to \textit{compensate} iregularity of the distribution \cite{basic2018class}.
    As shown in \cite{basic2018class}, in contrast to \eqref{eq:ldd_naive} sum version \eqref{eq:ldd_sum} does not contain the first order error but the average error collected from the neighbor points and it does not require to precompute the gradient values.

    \subsubsection{Laplacian approximation - LS version}
    An alternative way to get the Laplace operator from equation \eqref{eq:ldd-sumlaplacian-expressedLaplacianTerms} is to use the moving least square procedure, which stated by the authors \cite{basic2018class} results in more accurate approximation.
    This is closer to the approach used in LDD to approximate the first derivative or to the approach introduced in the Basic approach section 2.
    From the left side of the equation \eqref{eq:ldd-sumlaplacian-expressedLaplacianTerms}, we can express $ \left.\bm{D}^{2} u\right\rvert_{i} $ and group the remaining terms left side as
    \begin{equation}
      \left.\bm{D}^{2} u\right\rvert_{i} \cdot \bm{Q}_{ji}
      =
      2
      \left(
        u_{ji}
        -
        \vect{x}_{ji} \cdot \left(\mathbf{M}_{i}^{1}\right)^{-1}
        \displaystyle \sum_{j\in\mathcal{N}_i}
        \varphi_{ji} \vect{x}_{ji}
        \left(
          u_{ji} - e_{ji}^{3}
        \right)
        -
        e_{ji}^{3}
      \right).
      \label{eq:ldd-lsLaplace-expressedWithVectorQ}
    \end{equation}
    where
    \begin{equation}
      \bm{Q}_{ji}
      =
      \bm{X}_{ji}^{2}
      -
      \vect{x}_{ji}^{\top}
      \displaystyle \sum_{j\in\mathcal{N}_i}
      \varphi_{ji} \left(\mathbf{M}_{i}^{1}\right)^{-1}\vect{x}_{ji}
      \left( \bm{X}_{ji}^{2} \right) ^{\top}.
      \label{eq:ldd-lsLaplace-q}
    \end{equation}
    Multiplying equation \eqref{eq:ldd-lsLaplace-expressedWithVectorQ} by $\varphi_{ji} \bm{Q}_{ji}^{\top} $ and summing over the particles, we end up with
    \begin{equation}
      \displaystyle \sum_{j\in\mathcal{N}_i}
      \left.\bm{D}^{2} u\right\rvert_{i}
      \varphi_{ji} \cdot \bm{Q}_{ji} \bm{Q}_{ji}^{\top}
      =
      \displaystyle \sum_{j\in\mathcal{N}_i}
      2 \varphi_{ji}
      \left(
        u_{ji}
        -
        \vect{x}_{ji} \cdot \left(\mathbf{M}_{i}^{1}\right)^{-1}
        \displaystyle \sum_{j\in\mathcal{N}_i}
        \varphi_{ji} \vect{x}_{ji}
        \left(
          u_{ji} - e_{ji}^{3}
        \right)
        -
        e_{ji}^{3}
      \right) \bm{Q}_{ji}^{\top}.
      \label{eq:ldd-lsLaplace-expressedWithVectorQMultipliedAndSummed}
    \end{equation}
    Defining modified correction matrix $\widehat{\mathbf{M}}_i$ (in LDD called renormalisation matrix)
    \begin{equation}
      \widehat{\mathbf{M}}_i
      =
      \displaystyle \sum_{j\in\mathcal{N}_i}
      \varphi_{ji}
      \bm{Q}_{ji}
      \bm{Q}_{ji}^{\top}
    \end{equation}
    and multiplying both sides of equation \eqref{eq:ldd-lsLaplace-expressedWithVectorQMultipliedAndSummed} with identity vector $\bm{I}$, we obtain discretization of Laplacian in form
    \begin{equation}
      L_{i}^{\nabla^2,\textrm{LS}} u
      =
      \bm{I} \cdot \widehat{\mathbf{M}}_i^{-1}
      \displaystyle \sum_{j\in\mathcal{N}_i}
      2 \varphi_{ji} \vect{Q}_{ji}
      \left(
        u_{ji}
        -
        \vect{x}_{ji} \cdot L_i^{\nabla} u
      \right).
      \label{eq:ldd_ls_full}
    \end{equation}
    where we recognized the approximation for gradient
    \begin{equation}
      L_i^{\nabla} u
      =
      \left(\mathbf{M}_{i}^{1} \right)^{-1}
      \displaystyle \sum_{j\in\mathcal{N}_i}
      u_{ji}
      \bm{W}_{ji}.
    \end{equation}
    The equation \eqref{eq:ldd_ls_full} is referred to as \textit{full inversion version} of Laplace operator.
    The full inversion version is computationally not very efficient, since several loops over the neighbors need to be done to compute all the parts of the expression.
    To provide computationally more advantageous discretization, authors in \cite{basic2018class} proposed different approach.
    Instead of the steps performed above, equation \eqref{eq:ldd-sumlaplacian-expressedLaplacianTerms} is multiplied by $\bm{X}_{ji}^{2}$, weight function and summed over all neighbors

    \begin{align*}
      \left.\bm{D}^{2} u\right\rvert_{i}
      \left[
        \displaystyle \sum_{j\in\mathcal{N}_i}
        \varphi_{ji}\bm{X}_{ji}^{2}\left(\bm{X}_{ji}^{2}\right)^{\top}
        +
        \displaystyle \sum_{j\in\mathcal{N}_i}
        \varphi_{ji}\mathbf{x}_{ji}\left(\bm{X}_{ji}^{2}\right)^{\top}
        \displaystyle \sum_{k}
        \varphi_{ji}\left( \mathbf{M}_{i}^{1} \right)^{-1}\mathbf{x}_{ji}\left(\bm{X}_{ji}^{2}\right)^{\top}
      \right]
      \approx \\
      \displaystyle \sum_{j\in\mathcal{N}_i}
      2 \varphi_{ji} \bm{X}_{ji}^{2}
      \left(
        u_{ji}
        -
        \vect{x}_{ji} \cdot L_i^{\nabla} u
      \right)
      \numberthis \label{eq:ldd_ls_efficent_strategy}
    \end{align*}
    where we similarly to \eqref{eq:ldd_ls_full} used the notation for the first derivatives.
    In the second term on the left side of \eqref{eq:ldd_ls_efficent_strategy}, the product of two similar tensors is recovered.
    Following \cite{basic2018class} and taking the cubic order elements as zero (due to symmetry of weight function $\varphi$) we can define the modified correction matrix $ \widetilde{\mathbf{M}}_i $ simply as
    \begin{equation}
      \widetilde{\mathbf{M}}_i
      =
      \displaystyle \sum_{j\in\mathcal{N}_i}
      \varphi_{ji}
      \bm{X}_{ji}^{2}\left(\bm{X}_{ji}^{2}\right)^{\top}.
    \end{equation}
    This leads to approximation of Laplacian in form
    \begin{equation}
      L_{i}^{\nabla^2,\textrm{LS}} u
      =
      \bm{I} \cdot \widetilde{\mathbf{M}}_i^{-1}
      \displaystyle \sum_{j\in\mathcal{N}_i}
      2 \varphi_{ji} \vect{Q}_{ji}
      \left(
        u_{ji}
        -
        \vect{x}_{ji} \cdot L_i^{\nabla} u
      \right).
      \label{eq:ldd_ls_basic}
    \end{equation}
    which is denoted as \textit{basic inversion version} of the Laplace operator. In contrast to \eqref{eq:ldd_ls_full} both $ \widetilde{\mathbf{M}}_i $ and $ \mathbf{M}_{i}^{1} $ can be computed in the same loop.
    The gradient $L_i^{\nabla} u$ is in contrast to sum version \eqref{eq:ldd_sum} still required.
    Potential drawback reported by the authors \cite{basic2018class} with regard to formula \eqref{eq:ldd_ls_full} and \eqref{eq:ldd_ls_basic} is problem with conditioning of the matrix in case that elements of vectors $\bm{Q}_{ji}$ or $\bm{X}_{ji}^{2}$ are close to zero.

    \subsection{Generalized reproducing kernel peridynamics (GRKP)\label{sec:com_grkp}}

    Generalized reproducing kernel peridynamics (GRKP) \cite{hillman2020grkp} was proposed based on the relationship between meshfree pridynamic method \cite{madenci2013pd, oterkus2021}, used to solve new class of problems from solid mechanics, and classical meshfree methods discussed in this paper.
    Peridynamics is a mathematical framework that describes the behavior of materials and provides non-local formulation of continuum mechanics.
    Peridynamics describes the material as set of material points and focuses on interactions between them in a similar way how particle methods approximate function values and derivatives.
    The gradients in peridynamics are obtained as weighed summation over surrounding points, originating from Taylor expansion, strikingly resembling the approximation of derivatives in collocation methods, which motivate the authors to explore the connection.
    On one hand, GRKP bridges two completely different fields.
    On the other hand, peridynamics is based on particles and the need to represent derivatives on scattered data hints the intuitive link.
    GRKP which unifies peridynamics differential operators with RKPM methods using implicit gradients (see GRKCM section \ref{sec:com_grkcm}) \textbf{is based on approximation of moments 2}. Including the zeroth moment is optional.

    The theory of peridynamics assumes that variation of  function $ u(\mathbf{x}_i) $ at point $ \mathbf{x}_i $ (or simply $u_i$) is influenced by interactions with another points in its neighborhood $N_i$.
    As stated in \cite{Madenci2016}, in the original formulations of peridynamics, the neighborhoods were symmetric and uniform with constant cut-off distance $r_c$.
    In recent variants, the sizes and shapes of the neighborhoods of different points can vary.
    Similarly to particle methods, the interaction between peridynamics points are local thanks to the weight function $\varphi_{ji} = \varphi(\lvert x_j - x_i \rvert, h )$ which weights the interactions.
  For detail description of derivatives in the peridynamics context we refer to \cite{Madenci2016} and for more information about the peridynamics, we refer to original GRKP paper \cite{hillman2020grkp} or specialized peridynamics literature \cite{madenci2013pd, oterkus2021}.)

  The peridynamics assumes the differential operators in form
  \begin{equation}
    \partial^{\bm{\alpha}} u_i =  \int_{\mathcal{N}_i} u(\mathbf{x}') g_m^{\bm{\alpha}}(\mathbf{x}' - \mathbf{x}_i) \mathrm{d}\mathbf{x}'
    \label{eq:peridinamic_derivative}
  \end{equation}
  where $\bm{\alpha} $ is the multi-index and $ g_m^{\bm{\alpha}}(\mathbf{x}' - \mathbf{x}_i) = g_m^{\alpha_1, \alpha_2,\dots\alpha_n}(\mathbf{x}' - \mathbf{x}_i)$ are orthogonal functions which needs to be constructed in order to represent the operator.
  The $g_m^{\bm{\alpha}}$ functions are usually assumed as product of Taylor monomials $\mathbf{X}(\mathbf{x})$ multiplied by weight functions $\psi(\mathbf{x})$, which we denote $\mathbf{W}(\mathbf{x}) = \mathbf{X}(\mathbf{x})\psi(\mathbf{x})$, and vector of unknown coefficients ${\Psi}_{\mathbf{x}}^{\bm{\alpha}} $
  \begin{equation}
    g_m^{\bm{\alpha}}(\mathbf{x}) = \bm{W}(\mathbf{x})\cdot\bm{\Psi}_{\mathbf{x}}^{\bm{\alpha}}.
  \end{equation}
  This is similar to our weight functions defined by equation \eqref{eq:w_abf_zm}.

  Formulating equation \label{eq:peridinamic_derivative} in discrete collocation sense or assuming direct nodal integration, i.e. that there is an elementary volume associated with each point, which can be incorporated to the weight functions, we can spot that equation \eqref{eq:peridinamic_derivative} become equivalent with our definition of general discrete operator \eqref{eq:general_discrete_diff_operator_zm}.
  Replacing $ \partial^{\bm{\alpha}} u_i$ with $ L^{\bm{\alpha}}_{i}u $, we have
  \begin{equation}
    L^{\bm{\alpha}}_{i}u =
    \displaystyle\sum_{j\in\mathcal{N}_i}^{N}
    u_{j}w^{\bm{\alpha}}_{ji}
    \approx
    \bm{C}^{\bm{\alpha}} \cdot \left.\bm{D}u\right\rvert_{i}
    \label{eq:grkp_differential_operator}
  \end{equation}
  where we assume that the weights have the form $ w_{ji}^{\bm{\alpha}}  = \bm{W}_{ji} \cdot \bm{\Psi}_i^{\bm{\alpha}} $ of anisotropic basis functions $ \bm{W}_{ji} $ multiplied by unknown coefficients $ \bm{\Psi}_i^{\bm{\alpha}} $.

  As shown in detail in \cite{hillman2020grkp}, the idea of peridynamic differential operators is similar to the approximation of a differential operator using reproducing conditions.
  The authors moreover introduce Generalized Reproducing Kernel Peridynamics (GRKP) which combines the peridynamics differential operators and implicit gradient-based RKPM (for RKPM and implicit gradients see GRKCM \ref{sec:com_grkcm}).
  GRKP follows directly the derivation presented in Section 2.
  The method is presented in both continuous and discrete collocation versions, showing that the resulting formula corresponds to a scheme which is equivalent with GRKCM and which can also be interpreted as a peridynamics differential operator employed in the peridynamics theory.

  Including the zeroth moment in the GRKP is optional.
  In the referential paper \cite{hillman2020grkp}, both situations are considered, and this feature is directly built into the GRKP notation.
  If the moment is used, it is denoted by a square bracket subscript containing two indices $ [l,m] $ - $l$ for the \textit{starting order of the basis} i.e., $l=0$ using the unit element or $l=1$ to exclude the unit element, and $m$ which denotes the order of the basis used and hence the approximation.

  Summarizing the results in our notation, the operator $ L^{\bm{\alpha}}_{i} $ can thus be expressed as
  \begin{equation}
    \bm{W}_{ji,[l,m]}
    =
    \bm{X}_{ji,[l,m]} \varphi_{ji,[l,m]}
    \label{eq:grkp_abfs}
  \end{equation}
  \begin{equation}
    L^{\bm{\alpha}}_{i, [m]} u
    =
    \displaystyle\sum_{j\in\mathcal{N}_i}
    u_{ji}\bm{W}_{ji,[l,m]}\cdot \left(\mathbf{M}_{i,[l,m]}^{-1} \bm{C}^{\bm{\alpha}}_{[l,m]} \right),
    \quad
    \mathbf{M}_{i,[l,m]}
    =
    \displaystyle\sum_{j\in\mathcal{N}_i}
    \bm{X}_{ji,[l,m]} \bm{W}_{ji}^{\top}
    \quad
    \label{eq:grkp_ld_approximation}
  \end{equation}
  where Taylor monoials $\bm{X}_{ji}$ are used to construct anizotropic basis functions.

  \subsection{Meshless Local Strong Form Method (MLSM)\label{sec:com_mlsm}}

  The name Meshless Local Strong Form Method is introduced in \cite{kosec2019mlsm} as a generalization of several strong form methods presented in the literature.
  The approach presented in the paper is based on MLS and corresponds exactly to the methods presented in section \ref{sec:com_mls}.
  MLSM is thus \textbf{ based on minimization of approximation error, but focuses on the reconstruction of the target function, which is then used to compute its derivatives}.
  As presented in the referred paper, \textbf{zeroth moment is included.}
  The paper also highlights that a wide range of basis function $ \bm{P}_{ji} $ can be used, which allows one to write the method as
  \begin{equation}
    \bm{W}_{ji}
    =
    \bm{P}_{ji} \varphi_{ji}.
  \end{equation}
  \begin{equation}
    L^{\bm{\alpha}}_{i} u
    =
    \displaystyle\sum_{j\in\mathcal{N}_i}
    u_{j} \bm{W}_{ji} \cdot (\mathbf{M}_i^{-1} \bm{C}^{\bm{\alpha}}),
    \quad
    \mathbf{M}_{i}
    =
    \displaystyle\sum_{j\in\mathcal{N}_i}\bm{P}_{ji} \bm{W}_{ji}^{\top}.
    \label{eq:mlsm_ld_approximation}
  \end{equation}
